%% file: TRUnderlying.tex
\documentclass[10pt]{article}

\usepackage[dvips]{graphics,epsfig} 
\usepackage{amssymb,amsmath,amsthm,mathrsfs} 
\usepackage{graphicx}
\usepackage{lineno}

\usepackage[authoryear,sort&compress]{natbib}
\usepackage{url}
\usepackage{enumerate}
\usepackage{colordvi,color,pspicture}
\usepackage{psfrag}
\usepackage{rotating}
\usepackage{verbatim}

\graphicspath{{Figures/}}

\newcommand{\I}{\mathbf{I}}
\newcommand{\la}{\text{and}}
\newcommand{\lo}{\text{or}}

\newcommand{\E}{\mathbf{E}}

\newcommand{\U}{\mathcal{U}}

\newcommand{\g}{\gamma}
\newcommand{\G}{\Gamma}
\newcommand{\V}{\mathcal{V}}
\newcommand{\A}{\mathcal{A}}
\newcommand{\Y}{\mathcal{Y}}
\newcommand{\y}{\mathsf{y}}
\newcommand{\X}{\mathcal{X}}

\newcommand{\mE}{\mathcal{E}}
\newcommand{\mI}{\mathcal{I}}
\newcommand{\NY}{N_{\Y}}

\newcommand{\NPE}{N_{PE}}
\newcommand{\TY}{T(\Y_3)}
\newcommand{\N}{\mathcal{N}}
\newcommand{\R}{\mathbb{R}}

\newcommand{\msP}{\mathscr{P}}

\newcommand{\ve}{\varepsilon}
\newcommand{\Var}{\mathbf{Var}}
\newcommand{\Cov}{\mathbf{Cov}}

\DeclareMathOperator{\argsup}{argsup}

\DeclareMathOperator{\PAE}{PAE}
\DeclareMathOperator{\HLAE}{HLAE}

\theoremstyle{plain}
\newtheorem{theorem}{Theorem}[section]
\newtheorem{lemma}[theorem]{Lemma}

\newtheorem{corollary}[theorem]{Corollary}

\theoremstyle{definition}

\theoremstyle{remark}

\newtheorem{remark}[theorem]{Remark}

\begin{document}

\title{Technical Report \# KU-EC-09-5:\\
Relative Edge Density of the Underlying Graphs Based on
Proportional-Edge Proximity Catch Digraphs for Testing Bivariate Spatial Patterns}
\author{
Elvan Ceyhan\thanks{Address:
Department of Mathematics, Ko\c{c} University, 34450 Sar{\i}yer, Istanbul, Turkey.
e-mail: elceyhan@ku.edu.tr, tel:+90 (212) 338-1845, fax: +90 (212) 338-1559.
}
}
\date{\today}
\maketitle

\begin{abstract}
\noindent
The use of data-random graphs in
statistical testing of spatial patterns is introduced recently.
In this approach, a
random directed graph is constructed from the data using the
relative positions of the points from various classes.
Different random graphs result from different definitions of the
proximity region associated with each data point and different
graph statistics can be employed for pattern testing.
The approach used in this article is based on underlying graphs of
a family of data-random digraphs which is determined by a
family of parameterized proximity maps.
The relative edge density of the AND- and OR-underlying graphs is used
as the summary statistic, providing an alternative to the
relative arc density and domination number of the digraph employed previously.
Properly scaled,
relative edge density of the underlying graphs is a $U$-statistic,
facilitating analytic study of its asymptotic
distribution using standard $U$-statistic central limit theory.
The approach is illustrated with an application to the testing of
bivariate spatial clustering patterns of segregation and association.
Knowledge of the asymptotic distribution allows evaluation of the Pitman
asymptotic efficiency,
hence selection of the proximity map
parameter to optimize efficiency.
Asymptotic efficiency and Monte Carlo simulation analysis indicate that
the AND-underlying version is better (in terms of power and efficiency)
for the segregation alternative,
while the OR-underlying version is better for the association alternative.
The approach presented here is also
valid for data in higher dimensions.
\end{abstract}


\vspace{.25 in}



\newpage


\section{Introduction}
\label{sec:intro}
Classification and clustering have received considerable attention
in the statistical literature.
In this article, a graph-based approach for testing
bivariate spatial clustering patterns is introduced.
The analysis of spatial point patterns in natural populations has been
extensively studied and have important implications in epidemiology,
population biology, and ecology.
The patterns of points from one class with respect to points from other classes,
rather than the pattern of points from one-class with respect to the ground, are investigated.
The spatial relationships among two or more classes have important implications especially
for plant species.
See, for example, \cite{pielou:1961} and \cite{dixon:1994, dixon:EncycEnv2002}.

The goal of this article is to derive the asymptotic distribution of
the relative edge density of underlying graphs based on a
particular digraph family and use it to test the spatial pattern of complete
spatial randomness against spatial segregation or association.
Complete spatial randomness (CSR) is roughly defined as the lack of
spatial interaction between the points in a given study area.
Segregation is the pattern in which points of one class tend to
cluster together, i.e., form one-class clumps.
In association, the points of one class tend to occur more frequently around points from
the other class.
For convenience and generality, we call the different types of
points ``classes", but the class can be replaced by any
characteristic of an observation at a particular location. For
example, the pattern of spatial segregation has been investigated
for plant species (Diggle 1983), age classes of plants (\cite{hamill:1986})
and sexes of dioecious plants (\cite{nanami:1999}).

In recent years,
the use of mathematical graphs has also gained popularity in spatial analysis (\cite{roberts:2000}).
In spatial pattern analysis
graph theoretic tools provide a way to move beyond Euclidean metrics for spatial analysis.
For example, graph-based approaches have been proposed to determine paths among
habitats at various scales and dispersal movement distances, and
balance data requirements with information content (\cite{fall:2007}).
Although only recently introduced to landscape ecology, graph theory is well suited to
ecological applications concerned with connectivity or movement (\cite{minor:2007}).
However, conventional graphs do not explicitly maintain geographic reference,
reducing utility of other geo-spatial information.
\cite{fall:2007} introduce spatial graphs that
integrate a geometric reference system that ties patches and paths
to specific spatial locations and spatial dimensions
thereby preserving the relevant spatial information.
After a graph is constructed using spatial data,
usually the scale is lost (see for instance, \cite{su:2007}).
Many concepts in spatial ecology depend on the idea of spatial adjacency which
requires information on the close vicinity of an object.
Graph theory conveniently can be used to express
and communicate adjacency information allowing one to compute
meaningful quantities related to spatial point pattern.
Adding vertex and edge properties to graphs extends the problem domain to network modeling (\cite{keitt:2007}).
\cite{wu:2008} propose a new measure based on graph theory and spatial interaction,
which reflects intra-patch and inter-patch relationships by
quantifying contiguity within patches and potential contiguity among patches.
\cite{friedman:1983} also propose a graph-theoretic method
to measure multivariate association,
but their method is not designed to analyze spatial interaction
between two or more classes;
instead it is an extension of generalized correlation coefficient
(such as Spearman's $\rho$ or Kendall's $\tau$)
to measure multivariate (possibly nonlinear) correlation.

A new type of spatial clustering test
using directed graphs (i.e., digraphs) which is based on
the relative positions of the data points from various classes has also been developed recently.
Data-random digraphs are directed graphs in which each vertex corresponds to a data point,
and directed edges (i.e., arcs) are defined in terms of some bivariate function on the data.
For example, nearest neighbor digraphs are defined by placing
an arc between each vertex and its nearest neighbor.
\cite{priebe:2001} introduced the class cover catch digraphs (CCCDs) in $\R$ and gave the exact and
the asymptotic distribution of the domination number of the CCCDs.
\cite{devinney:2002a}, \cite{marchette:2003}, \cite{priebe:2003b},
\cite{priebe:2003a}, and \cite{devinney:2006} applied the concept in higher dimensions and
demonstrated relatively good performance of CCCDs in classification.
Their methods involve \emph{data reduction} (i.e., \emph{condensing}) by
using approximate minimum dominating sets as \emph{prototype sets}
(since finding the exact minimum dominating set is an NP-hard
problem in general --- e.g., for CCCD in multiple dimensions --- (see \cite{devinney:2006}).
Furthermore the exact and the asymptotic
distribution of the domination number of the CCCDs are not
analytically tractable in multiple dimensions.
For the domination number of CCCDs for one-dimensional data, a SLLN result is proved in \cite{devinney:2002b},
and this result is extended by \cite{wiermanSLLN:2008}; furthermore,
a generalized SLLN result is provided by \cite{wiermanSLLN:2008},
and a CLT is also proved by \cite{xiangCLT:2009}.
The asymptotic distribution of the domination number of CCCDs for non-uniform data
in $\mathbb{R}$ is also calculated in a rather general setting (\cite{ceyhan:dom-num-CCCD-NonUnif}).
\cite{ceyhan:Phd-thesis} generalized CCCDs to what is called \emph{proximity catch digraphs} (PCDs).
The first PCD family is introduced by \cite{ceyhan:CS-JSM-2003};
the parametrized version of this PCD is developed by \cite{ceyhan:arc-density-CS}
where the relative arc density of the PCD is calculated and used for spatial pattern analysis.
\cite{ceyhan:dom-num-NPE-SPL} introduced another digraph family called \emph{proportional edge PCDs}
and calculated the asymptotic distribution of its domination number and
used it for the same purpose.
The relative arc density of this PCD family is also computed
and used in spatial pattern analysis (\cite{ceyhan:arc-density-PE}).
\cite{ceyhan:dom-num-NPE-MASA} derived the asymptotic distribution
of the domination number of proportional-edge PCDs for two-dimensional uniform  data.

The underlying graphs based on digraphs are obtained by replacing arcs in the digraph
by edges based on bivariate relations.
If symmetric arcs are replaced by edges, then we obtain the AND-underlying graph;
and if all arcs are replaced by edges without allowing multi-edges,
then we obtain the OR-underlying graph.
The statistical tool utilized in this article is the asymptotic theory of $U$-statistics.
Properly scaled, we demonstrate that the relative
edge density of the underlying graphs of proportional-edge PCDs is a $U$-statistic,
which has asymptotic normality by the general central limit theory of $U$-statistics.
For the digraphs introduced by \cite{priebe:2001},
whose relative arc density is also of the $U$-statistic form,
the asymptotic mean and variance of the relative density is not analytically tractable,
due to geometric difficulties encountered.
However, for the PCDs introduced in \cite{ceyhan:CS-JSM-2003}, \cite{ceyhan:arc-density-PE}, and \cite{ceyhan:arc-density-CS},
the relative arc density has tractable asymptotic mean and variance.

We define the underlying graphs of proportional-edge PCDs and
their relative edge density in Section \ref{sec:underlying},
provide the asymptotic distribution of the relative edge density under the null hypothesis
in Section \ref{sec:null-dist-edge-density},
and
describe the alternatives of segregation and association in Section \ref{sec:alt-seg-assoc}.
We prove the consistency of the relative edge density in Section \ref{sec:consistency},
and
provide Pitman asymptotic efficiency in Section \ref{sec:PAE}.
We present the Monte Carlo simulation analysis for finite sample performance
in Section \ref{sec:monte-carlo},
in particular, provide the Monte Carlo power analysis under segregation
in Section \ref{sec:power-seg},
and
under association in Section \ref{sec:power-assoc}.
We treat the multiple triangle case in Section \ref{sec:multiple-triangle-case},
provide extension to higher dimensions in Section \ref{sec:extend-high-dim}.
We provide the discussion and conclusions in Section \ref{sec:discussion},
and the tedious calculations and long proofs are deferred to the Appendix.

\section{Relative Edge Density of Underlying Graphs}
\label{sec:underlying}

\subsection{Preliminaries}
\label{sec:preliminaries}
The main difference between a graph and a
digraph is that edges are directed in digraphs,
hence are called arcs.
So the arcs are denoted as ordered pairs
while edges are denoted as unordered pairs.
The \emph{underlying graph} of a
digraph is the graph obtained by replacing each arc $uv \in \A$ or
each symmetric arc, $\{uv,\,vu\} \subset \A$ by the edge $(u,v)$.
The former underlying graph will be referred as the \emph{OR-underlying graph},
while the latter as the \emph{AND-underlying graph}.
That is,
the AND-underlying graph for digraph $D=(\V, \A)$ is the graph
$G_{\la}(D)=(\V, \mE_{\la})$ where $\mE_{\la}$ is the set of edges
such that  $(u,v)\in \mE_{\la}$ iff $uv \in \A$ \emph{and} $vu \in \A$.
The OR-underlying graph for $D=(\V, \A)$ is the graph
$G_{\lo}(D)=(\V, \mE_{\lo})$ where $\mE_{\lo}$ is the set of edges such that
$(u,v)\in \mE_{\lo}$ iff $uv \in \A$ \emph{or} $vu \in \A$.

The relative edge density of a graph $G=(\V,\mE)$ of order $|\V| = n$,
denoted $\rho(G)$, is defined as
$$\rho(G) = \frac{2\,|\mE|}{n(n-1)}$$
where $|\cdot|$ denotes the set cardinality function (\cite{janson:2000}).
Thus $\rho(G)$ represents the ratio of the number of edges
in the graph $G$ to the number of edges in the complete graph of order $n$, which is $n(n-1)/2$.


Let $(\Omega,\mathcal{M})$ be a measurable space and
consider $N:\Omega \rightarrow \wp(\Omega)$,
where $\wp(\cdot)$ represents the power set functional.
Then given $\Y_m \subset \Omega$,
the {\em proximity map}
$\NY(\cdot)$ associates with each point $x \in \Omega$
a {\em proximity region} $\NY(x) \subseteq \Omega$.
The $\G_1$-region $\G_1(\cdot,N):\Omega \rightarrow \wp(\Omega)$
associates the region $\G_1(x,\NY):=\{z \in \Omega: x \in \NY(z)\}$ with each point $x \in \Omega$.
If $X_1,X_2,\ldots,X_n$ are $\Omega$-valued random variables,
then the $\NY(X_i)$ (and $\G_1(X_i,\NY)$), $i=1,2,\ldots,n$ are random sets.
If the $X_i$ are independent and identically distributed,
then so are the random sets $\NY(X_i)$ (and $\G_1(X_i,\NY)$).

Consider the data-random PCD $D$
with vertex set $\V=\{X_1,X_2,\ldots,X_n\}$
and arc set $\A$ defined by
$X_iX_j \in \A \iff X_j \in \NY(X_i)$.
The {\em AND-underlying graph}, $G_{\la}$, of $D$ with the vertex set $\V$ and
the edge set $\mE_{\la}$ is defined by $(X_i,X_j) \in \mE_{\la}$ iff
$X_iX_j \in \A \text{ and } X_jX_i \in \A$.
Likewise, the {\em OR-underlying graph}, $G_{\lo}$, of $D$ with the vertex set $\V$
and the edge set $\mE_{\lo}$ is defined by
$(X_i,X_j) \in \mE_{\lo} \iff X_iX_j \in \A \text{ or } X_jX_i \in \A$.
Then
$\left( X_i,X_j \right) \in \mE_{\la}$ iff $X_j \in \NY(X_i) \text{ and } X_i \in \NY(X_j)$
iff $X_j \in \NY(X_i) \text{ and } X_j \in \G_1(X_i,\NY)$ iff $X_j \in \NY(X_i) \cap  \G_1(X_i,\NY)$.
Similarly, $\left( X_i,X_j \right) \in \mE_{\lo}$ iff $X_j \in \NY(X_i) \cup  \G_1(X_i,\NY)$.
Since the random digraph $D$ depends on
the (joint) distribution of the $X_i$ and on
the map $\NY$, so do the underlying graphs.
The adjective {\em proximity} --- for the catch digraph $D$ and for the map $\NY$ ---
comes from thinking of the region $\NY(x)$ as representing those points in $\Omega$ ``close'' to $x$
(\cite{toussaint:1980} and \cite{jaromczyk:1992}).

\subsection{Relative Edge Density of the AND-Underlying Graphs}
\label{sec:AND-edge-density}
The relative edge density of $G_{\la}(D)$, the AND-underlying graph based on digraph $D$,
is denoted as $\rho_{\la}(D)$.
For $X_i \stackrel{iid}{\sim}F$, $\rho_{\la}(D)$ is a $U$-statistic,
$$\rho_{\la}(D)=\frac{2}{n\,(n-1)}\sum\hspace*{-0.1 in}\sum_{i < j \hspace*{0.25 in}}   \hspace*{-0.1 in} \,h^{\la}_{ij}$$
where
\begin{eqnarray*}
h^{\la}_{ij} & = & h_{\la}(X_i,X_j;N)=\I(\left( X_i,X_j \right)\in \mE_{\la})= \I(X_iX_j \in \A)\cdot \I(X_jX_i \in \A)\\
&=&\I(X_i \in N(X_j))\cdot \I(X_j \in N(X_i))=\I(X_j \in N(X_i)\cap \G_1\left( X_i,N \right)).
\end{eqnarray*}
is the number of symmetric arcs between $X_i$ and $X_j$ in $D$ or
number of edges between $X_i$ and $X_j$ in $G_{\la}(D)$.
Note that $h^{\la}_{ij}$ is a symmetric kernel with finite variance since $0 \le h_{\la}(X_i,X_j;N)\le 1$.
Moreover, $\rho_{\la}(D)$ is a random
variable that depends on $n$, $F$, and $N(\cdot)$ (i.e., $\mathcal Y$).
But $\E\left[ \rho_{\la}(D) \right]$ only depends on $F$ and $N(\cdot)$.
Then
\begin{equation}
\label{eqn:E[rho-and-D]}
0\le \E\left[ \rho_{\la}(D) \right]=\frac{2}{n\,(n-1)}\sum\hspace*{-0.1 in}\sum_{i < j \hspace*{0.25 in}}
\hspace*{-0.1 in} \,\E[h^{\la}_{ij}]=\E\left[ h^{\la}_{12} \right]=\mu_{\la}(N)
\end{equation}
where $\E\left[ h^{\la}_{12} \right]=P(X_1X_2 \in \A\, \, , \,\, X_2X_1 \in \A)=
P(X_2 \in N(X_1)\cap \G_1(X_1,N))=\mu_{\la}(N)$ is the \emph{symmetric arc probability}.
Note that $\mu_{\la}(N)=P(X_j \in N(X_i)\cap \G_1\left( X_i,N \right))$ for $i \not= j$.
Furthermore,
\begin{equation}
\label{eqn:Var[rho-and-D]}
0 \le \Var\left[ \rho_{\la}(D) \right]=
\frac{4}{n^2\,(n-1)^2}\Var\left[\sum\hspace*{-0.1 in}\sum_{i < j \hspace*{0.25 in}}   \hspace*{-0.1 in} h^{\la}_{ij}\right].
\end{equation}
Expanding this expression,
we have
$$\Var\left[ \rho_{\la}(D) \right]=\frac{2}{n\,(n-1)}\Var\left[ h^{\la}_{12} \right]+
\frac{4\,(n-2)}{n\,(n-1)} \, \Cov\left[ h^{\la}_{12},h^{\la}_{13} \right].$$
Let $A_{ij}$ be the event that $\{X_iX_j \in \A\} =\{X_j \in N(X_i)\}$,
then $h^{\la}_{ij}=\I(A_{ij})\cdot \I(A_{ji})=\I(A_{ij} \cap A_{ji})$.
In particular, $h^{\la}_{12}=\I(A_{12})\cdot \I(A_{21})=\I(A_{12} \cap A_{21})$.
Then
$$\Var\left[ h^{\la}_{12} \right]=
\E\left[ \left( h^{\la}_{12} \right)^2 \right]-\left(\E\left[ h^{\la}_{12} \right]\right)^2=
\E\left[ \left( h^{\la}_{12} \right)^2 \right]-\left(\mu_{\la}(N)\right)^2.$$
Furthermore,
$\E\left[ (h^{\la}_{12})^2 \right] = \E\left[ (\I(A_{12}\cap A_{21}))^2 \right]=\E[(\I(A_{12}\cap A_{21})]=\mu_{\la}(N)$.
So
$$\Var\left[ h^{\la}_{12} \right]=\mu_{\la}(N)-\left[\mu_{\la}(N)\right]^2.$$
Moreover,
$$\Cov\left[ h^{\la}_{12},h^{\la}_{13} \right]=
\E\left[ h^{\la}_{12}.h^{\la}_{13} \right]-\E\left[ h^{\la}_{12} \right]\E\left[ h^{\la}_{13} \right]$$

where $\E\left[ h^{\la}_{12} \right]=\E\left[ h^{\la}_{13} \right]=\mu_{\la}(N)$ and,
\begin{eqnarray*}
\E\left[ h^{\la}_{12}.h^{\la}_{13} \right] & = & \E[\I(A_{12} \cap A_{21})\,(\I(A_{13}\cap A_{31})] = \E[(\I(A_{12} \cap A_{21} \cap A_{13}\cap A_{31})] \\
                & = & P(X_2 \in N(X_1)\cap \G_1(X_1,N)\, , \, X_3 \in N(X_1)\cap \G_1(X_1,N))\\
& = & P(\{X_2,X_3\} \subset N(X_1)\cap \G_1(X_1,N)).
\end{eqnarray*}
Thus
$$\Cov\left[ h^{\la}_{12},h^{\la}_{13} \right]=P(\{X_2,X_3\} \subset N(X_1)\cap \G_1(X_1,N))-\left[\mu_{\la}(N)\right]^2.$$

\subsubsection{The Joint Distribution of $\left( h^{\la}_{12},h^{\la}_{13} \right)$}
By definition $\left( h^{\la}_{12},h^{\la}_{13} \right)$ is a discrete random
variable with four possible values:
$$\left( h^{\la}_{12},h^{\la}_{13} \right)\in \{(0,0),(1,0),(0,1),(1,1)\}.$$
Then finding the joint distribution of $\left( h^{\la}_{12},h^{\la}_{13} \right)$
is equivalent to finding the joint probability mass function of
$\left( h^{\la}_{12},h^{\la}_{13} \right)$.

First, note that
$$\left( h^{\la}_{12},h^{\la}_{13} \right)=(0,0)\text{ iff }
h^{\la}_{12}=h^{\la}_{13}=0 \text{ iff }\I(A_{12} \cap A_{21})=\I(A_{13} \cap A_{31})=0 \text{ iff }$$
$$\I(X_2 \not\in N(X_1)\cap \G_1(X_1,N))=\I(X_3 \not\in N(X_1)\cap \G_1(X_1,N))=1 \text{ iff }$$
$$\I(X_2 \in \TY\setminus N(X_1)\cap \G_1(X_1,N))=\I(X_3 \in \TY\setminus N(X_1)\cap \G_1(X_1,N))=1 \text{ iff }$$
$$\I(\{X_2,X_3 \} \subset \TY\setminus [N(X_1)\cap \G_1(X_1,N)])=1.$$

Hence $P(\left( h^{\la}_{12},h^{\la}_{13} \right)=(0,0))=P(\{X_2,X_3 \} \subset
\TY\setminus [N(X_1)\cap \G_1(X_1,N)])$.

Next, note that $\left( h^{\la}_{12},h^{\la}_{13} \right)=(1,1)$ iff
$h^{\la}_{12}=h^{\la}_{13}=1$. So
$P\left(\left( h^{\la}_{12},h^{\la}_{13} \right)=(1,1) \right)=\E\left[ h^{\la}_{12}.h^{\la}_{13} \right]$.

Furthermore, by symmetry
$P\left(\left( h^{\la}_{12},h^{\la}_{13} \right)=(0,1)\right)=P(\left( h^{\la}_{12},h^{\la}_{13} \right)=(1,0)).$

Hence
\begin{eqnarray*}
P\left(\left( h^{\la}_{12},h^{\la}_{13} \right)=(0,1)\right)
&=&P\left(\left( h^{\la}_{12},h^{\la}_{13} \right)=(1,0)\right)\\
&=&\frac{1}{2}\,\left[1-\left(P\left(\left( h^{\la}_{12},h^{\la}_{13} \right)=(0,0)\right)+
    P\left(\left( h^{\la}_{12},h^{\la}_{13} \right)=(1,1)\right)\right)\right].
\end{eqnarray*}

\subsection{Relative Edge Density of OR-Underlying Graphs}
\label{sec:OR-edge-density}
The relative edge density of $G_{\lo}(D)$, the OR-underlying graph of digraph $D$,
is denoted as $\rho_{\lo}(D)$.
For $X_i \stackrel{iid}{\sim}F$, $\rho_{\lo}(D)$ is a $U$-statistic,
$$\rho_{\lo}(D)=\frac{2}{n\,(n-1)}\sum\hspace*{-0.1 in}\sum_{i < j \hspace*{0.25 in}}   \hspace*{-0.1 in} \,h^{\lo}_{ij}$$
where
\begin{eqnarray*}
h^{\lo}_{ij} = h_{\lo}(X_i,X_j;N) &=& \I(\left( X_i,X_j \right)\in \mE_{\lo}) = \max(\I(X_iX_j \in \A), \I(X_jX_i \in \A)) \\
                & = & \I(X_i \in N(X_i)\cup \G_1\left( X_i,N \right)).
\end{eqnarray*}
is the number of edges between $X_i$ and $X_j$ in $G_{\lo}(D)$.
Note that $h^{\lo}_{ij}$ is a symmetric kernel with finite variance since
$0 \le h_{\lo}(X_i,X_j;N)\le 1$.  Moreover, $\rho_{\lo}(D)$ is a
random variable that depends on $n$, $F$, and $N(\cdot)$ (i.e.,  $\mathcal Y$).
But $\E[\rho_{\lo}(D)]$ does only depend on $F$ and $N(\cdot)$.
Then
$$0\le \E[\rho_{\lo}(D)]=
\frac{2}{n\,(n-1)}\sum\hspace*{-0.1 in}\sum_{i < j \hspace*{0.25 in}}\hspace*{-0.1 in} \,\E[h^{\lo}_{ij}]=
\E\left[ h^{\lo}_{12} \right]=P(X_2 \in N(X_1)\cup \G_1(X_1,N))$$
where $\E\left[ h^{\lo}_{12} \right]=P(X_1X_2 \in \A \vee X_2X_1 \in \A)$ which
we denote as $\mu_{\lo}(N)$ for brevity of notation.

Similar to the AND-underlying case,
$$0 \le \Var[\rho_{\lo}(D)]=\frac{2}{n\,(n-1)}\Var\left[ h^{\lo}_{12} \right]+
\frac{4\,(n-2)}{n\,(n-1)} \, \Cov[h^{\lo}_{12},h^{\lo}_{13}].$$
Notice that $h^{\lo}_{12}=\max(\I(A_{12}),\I(A_{21}))=\I(A_{12} \cup A_{21})$.
Then
$$\Var\left[ h^{\lo}_{12} \right]=
\E\left[ (h^{\lo}_{12})^2 \right]-\left( \E\left[ h^{\lo}_{12} \right] \right)^2=
\mu_{\lo}(N)-[\mu_{\lo}(N)]^2$$
since $\E\left[ h^{\lo}_{12} \right]=\mu_{\lo}(N)$ and
$$\E\left[ (h^{\lo}_{12})^2 \right] = \E\left[(\I(A_{12}\cup A_{21}))^2\right]=\E[\I(A_{12}\cup A_{21})] = \mu_{\lo}(N).$$

Furthermore,
$$\Cov[h^{\lo}_{12},h^{\lo}_{13}]=\E\left[ h^{\lo}_{12}.h^{\lo}_{13} \right]-\E\left[ h^{\lo}_{12} \right]\E[h^{\lo}_{13}]$$

where $\E\left[ h^{\lo}_{12} \right]=\E[h^{\lo}_{13}]=\mu_{\lo}(N)$, and
\begin{eqnarray*}
\E\left[ h^{\lo}_{12}.h^{\lo}_{13} \right] & = & \E[(\I(A_{12} \cup A_{21})\,(\I(A_{13}\cup A_{31})]= \E[(\I((A_{12} \cup A_{21})\cap (A_{13}\cup A_{31}))] \\
                & = & P(X_2 \in N(X_1)\cup \G_1(X_1,N)\, , \, X_3 \in N(X_1)\cup \G_1(X_1,N))\\
& = & P(\{X_2,X_3\} \subset N(X_1)\cup \G_1(X_1,N).
\end{eqnarray*}
So
$$\Cov[h^{\lo}_{12},h^{\lo}_{13}]=P(\{X_2,X_3\} \subset N(X_1)\cup \G_1(X_1,N))-[\mu_{\lo}(N)]^2.$$
\begin{remark}
Note that $h_{ij}=h^{\la}_{ij}+h^{\lo}_{ij}$, since
$\I(A_{ij})+\I(A_{ji})=\I(A_{ij} \cap A_{ji})+\I(A_{ij} \cup
A_{ji})$. $\square$
\end{remark}

\subsubsection{The Joint Distribution of $\left( h^{\lo}_{12},h^{\lo}_{13} \right)$}
Finding the joint distribution of $\left( h^{\lo}_{12},h^{\lo}_{13} \right)$ is
equivalent to finding the joint probability mass function of
$\left( h^{\lo}_{12},h^{\lo}_{13} \right)$, i.e.,  finding
$P(\left( h^{\lo}_{12},h^{\lo}_{13} \right)=(i,j)) \text{ for each } (i,j)\in
\{(0,0),(1,0),(0,1),(1,1)\}.$

First, note that
$$\left( h^{\lo}_{12},h^{\lo}_{13} \right)=(0,0)\text{ iff } h^{\lo}_{12}=h^{\lo}_{13}=0 \text{ iff } \I(A_{12} \cup A_{21})=\I(A_{13} \cup A_{31})=0 \text{ iff } $$
$$\I(X_2 \not\in N(X_1)\cup \G_1(X_1,N))=\I(X_3 \not\in N(X_1)\cup \G_1(X_1,N))=1 \text{ iff } $$
$$\I(\{X_2,X_3 \} \subset \TY\setminus [N(X_1)\cup \G_1(X_1,N)])=1. $$

Hence $P(\left( h^{\lo}_{12},h^{\lo}_{13} \right)=(0,0))=P(\{X_2,X_3 \} \subset
\TY\setminus [N(X_1)\cup \G_1(X_1,N)])$.

Next, note that $\left( h^{\lo}_{12},h^{\lo}_{13} \right)=(1,1)$ iff
$h^{\lo}_{12}=h^{\lo}_{13}=1$.
$P(\left( h^{\lo}_{12},h^{\lo}_{13} \right)=(1,1))=\E\left[ h^{\lo}_{12}.h^{\lo}_{13} \right]$.

By symmetry
$P(\left( h^{\lo}_{12},h^{\lo}_{13} \right)=(0,1))=P(\left( h^{\lo}_{12},h^{\lo}_{13} \right)=(1,0))$.
Hence
\begin{eqnarray*}
P(\left( h^{\lo}_{12},h^{\lo}_{13} \right)=(0,1)) & = & P(\left( h^{\lo}_{12},h^{\lo}_{13} \right)=(1,0)) \\
                & = & \frac{1}{2}\,\left(1-\left[P(\left( h^{\lo}_{12},h^{\lo}_{13} \right)=(0,0))+P(\left( h^{\lo}_{12},h^{\lo}_{13} \right)=(1,1))\right]\right).
\end{eqnarray*}

\subsection{Proportional-Edge Proximity Maps and the Associated Regions}
\label{sec:prop-edge}
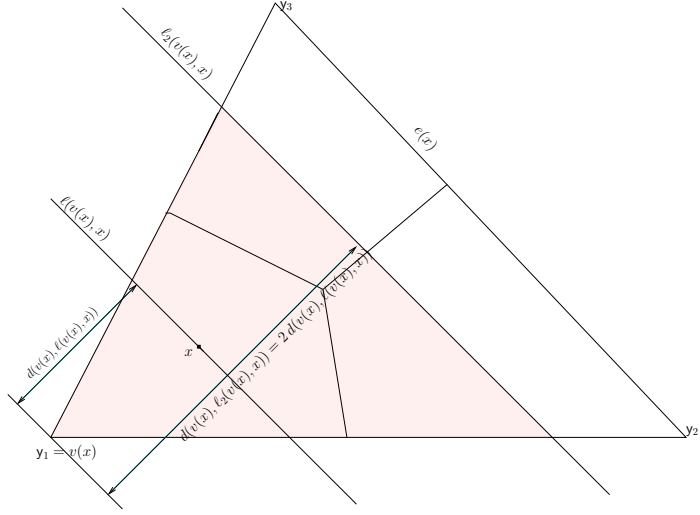
\begin{figure} [ht]
    \centering
   \scalebox{.4}{\input{Nofnu2.pstex_t}}
   \caption{Construction of proportional-edge proximity region, $\NY^{r=2}(x)$ (shaded region) for an $x \in R(\y)1)$. }
\label{fig:ProxMapDef1}
    \end{figure}

\begin{figure} [ht]
    \centering
    \scalebox{.4}{\input{Gammaofnu2.pstex_t}}
    \caption{Construction of the $\G_1$-region, $\G_1^{r=2}(x)$ (shaded region) for an $x \in R(\y)1)$. }
\label{fig:ProxMapDef2}
    \end{figure}
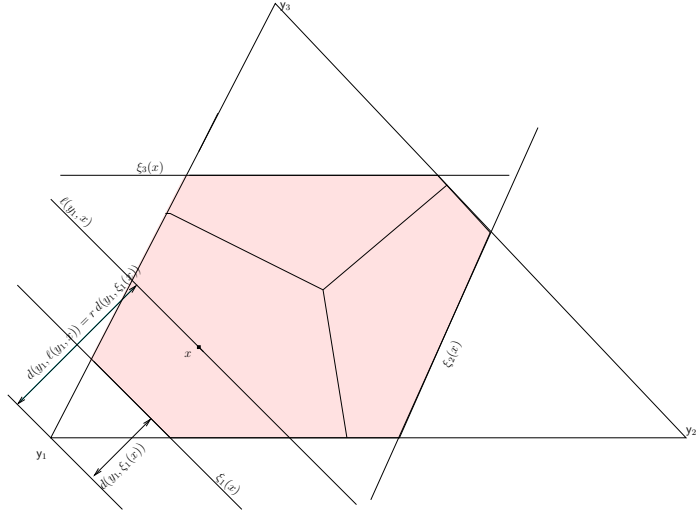

Let $\Omega = \mathbb{R}^2$ and $\Y_3 = \{\y_1,\y_2,\y_3\} \subset \mathbb{R}^2$ be three non-collinear points.
Denote by $\TY$ the triangle (including the interior) formed by these three points.
For $r \in [1,\infty]$
define $\NPE^r(x)$ to be the {\em proportional-edge}  proximity map with parameter $r$
and $\G_1^r(x):=\G_1\left( x,\NPE^r \right)$ to be the corresponding $\G_1$-region as follows;
see also Figures  \ref{fig:ProxMapDef1} and \ref{fig:ProxMapDef2}.
Let ``vertex regions'' $R(\y_1)$, $R(\y_2)$, $R(\y_3)$
partition $\TY$ using segments from the
center of mass of $\TY$ to the edge midpoints.
For $x \in \TY \setminus \Y_3$, let $v(x) \in \Y_3$ be the
vertex whose region contains $x$; $x \in R(v(x))$.
If $x$ falls on the boundary of two vertex regions,
or at the center of mass, we assign $v(x)$ arbitrarily.
Let $e(x)$ be the edge of $\TY$ opposite $v(x)$.
Let $\ell(v(x),x)$ be the line parallel to $e(x)$ through $x$.
Let $d(v(x),\ell(v(x),x))$ be the Euclidean (perpendicular) distance from $v(x)$ to $\ell(v(x),x)$.
For $r \in [1,\infty)$ let $\ell_r(v(x),x)$ be the line parallel to $e(x)$
such that $d(v(x),\ell_r(v(x),x)) = rd(v(x),\ell(v(x),x))$ and $d(\ell(v(x),x),\ell_r(v(x),x)) < d(v(x),\ell_r(v(x),x))$.
Let $T_r(x)$ be the triangle similar to
and with the same orientation as $\TY$
having $v(x)$ as a vertex and $\ell_r(v(x),x)$ as the opposite edge.
Then the proportional-edge proximity region
$\NPE^r(x)$ is defined to be $T_r(x) \cap \TY$.

Furthermore, let $\xi_i(x)$ be the line such that $\xi_i(x)\cap \TY \not=\emptyset$ and
$r\,d(\y_i,\xi_i(x))=d(\y_i,\ell(\y_i,x))$ for $i=1,2,3$.
Then $\G_1^r(x)\cap R(\y_i)=\{z \in R(\y_i): d(\y_i,\ell(\y_i,z)) \ge d(\y_i,\xi_i(x)\}$, for $i=1,2,3$.
Hence  $\G_1^r(x)=\bigcup_{i=1}^3 (\G_1^r(x)\cap R(\y_i))$.
Notice that $r \ge 1$ implies $x \in \NPE^r(x)$ and $x \in \G_1^r(x)$.
Furthermore,
$\lim_{r \rightarrow \infty} \NPE^r(x) = \TY$ for all $x \in \TY \setminus \Y_3$,
and so we define $\NY^{\infty}(x) = \TY$ for all such $x$.
For $x \in \Y_3$, we define $\NPE^r(x) = \{x\}$ for all $r \in [1,\infty]$.
Then,
for $x \in R(\y_i)$
$\lim_{r \rightarrow \infty} \G_1^r(x) = \TY \setminus \{\y_j,\y_k\}$
for distinct $i,j$, and $k$.

Notice that $X_i \stackrel{iid}{\sim} F$,
with the additional assumption
that the non-degenerate two-dimensional
probability density function $f$ exists
with support in $\TY$,
implies that the special cases in the construction of $\NPE^r$
---$X$ falls on the boundary of two vertex regions,
or at the center of mass, or $X \in \Y_3$ ---
occur with probability zero.
Note that for such an $F$, $\NY(x)$ is a triangle a.s.
and $\G_1(x)$ is a convex or nonconvex polygon.

\subsection{Relative Edge Density of the Underlying Graphs of Proportional-Edge PCDs}
\label{sec:relative-density-PEPCD}
Consider the underlying graphs of the data-random PCD
$D$ with vertex set $\V=\{X_1,X_2,\ldots,X_n\}$ and arc set $\A$
defined by $\left( X_i,X_j \right) \in \A \iff X_j \in \NY(X_i)$.
Recall that $\left( X_i,X_j \right) \in \mE_{\la}$ iff $X_j \in \NPE^r(X_i) \cap
\G_1\left( X_i,\NPE^r \right)$ and $\left( X_i,X_j \right) \in \mE_{\lo}$ iff $X_j \in \NPE^r(X_i)
\cup  \G_1\left( X_i,\NPE^r \right)$.

Let $h^{\la}_{ij}(r):=h_{\la}(X_i,X_j;\NPE^r)=\I(X_j \in \NPE^r(X_i)\cap \G^r_1\left( X_i \right))$
and $h^{\lo}_{ij}(r):=h_{\lo}(X_i,X_j;\NPE^r)=\I(X_j \in \NPE^r(X_i)\cup \G^r_1\left( X_i \right))$
for $i \not= j$.
The random variable $\rho^{\la}_n(r) := \rho_{\la}(\X_n;h,\NPE^r)$
depends on $n$ explicitly, and on $F$ and $\NPE^r$ implicitly.
The expectation $\E\left[\rho^{\la}_n(r)\right]$, however, is independent of $n$
and depends on only $F$ and $\NPE^r$.
Let $\mu_{\la}(r):=\E\left[ h^{\la}_{12}(r) \right]$
and
$\nu_{\la}(r):=\Cov\left[ h^{\la}_{12}(r),h^{\la}_{13}(r) \right]$.
Then
\begin{eqnarray}
0 \le \E\left[ \rho^{\la}_n(r) \right] = \E\left[ h^{\la}_{12}(r) \right] \le 1.
\end{eqnarray}
The variance $\Var\left[ \rho^{\la}_n(r) \right]$ simplifies to
\begin{eqnarray}
0 \leq
  \Var\left[ \rho^{\la}_n(r) \right] =
     \frac{2}{n(n-1)} \Var\left[ h^{\la}_{12}(r) \right] +
     \frac{4\,(n-2)}{n(n-1)} \Cov\left[ h^{\la}_{12}(r),h^{\la}_{13}(r) \right]
  \leq 1/4.
\end{eqnarray}
A central limit theorem for $U$-statistics (\cite{lehmann:1999}) yields
\begin{eqnarray}
\sqrt{n}\left(\rho^{\la}_n(r)-\mu_{\la}(r) \right)
\stackrel{\mathcal{L}}{\longrightarrow}
\N\left(0,4\,\nu_{\la}(r) \right)
\end{eqnarray}
provided $\nu_{\la}(r) > 0$.
The asymptotic variance of $\rho^{\la}_n(r)$, $4\,\nu_{\la}(r)$,
depends on only $F$ and $\NPE^r$.
Thus we need determine only $\mu_{\la}(r)$ and
$\nu_{\la}(r)$ in order to obtain the normal
approximation
\begin{eqnarray}
\label{eqn:asy-norm-and}
\rho^{\la}_n(r) \stackrel{\text{approx}}{\sim}
\N\left(\mu_{\la}(r),\frac{4\,\nu_{\la}(r)}{n}\right).
\end{eqnarray}

The above paragraph holds for
$\rho^{\lo}_n(r)=\rho_{\lo}(\X_n;h,\NPE^r)$ also with $\rho^{\la}_n(r)$ is
replaced by $\rho^{\lo}_n(r)$, $h^{\la}_{12}(r)$ and $h^{\la}_{13}(r)$ are
replaced by $h^{\lo}_{12}$ and $h^{\lo}_{13}$, respectively.


For $r=1$, $\NPE^{r=1}(x) \cap \G_1^{r=1}(x)=\ell(v(x),x)$
which has zero $\R^2$-Lebesgue measure.
Then we have
$\E\left[ \rho^{\la}_n(r=1) \right] = \E\left[ h^{\la}_{12}(r=1) \right]=
\mu_{\la}(r=1)=
P(X_2 \in \NPE^{r=1}(X_1) \cap \G_1^{r=1}(X_1))=0$.
Similarly,
$P(\{X_2,X_3\} \subset \NPE^{r=1}(X_1) \cap \G_1^{r=1}(X_1))=0$.
Thus, $\nu_{\la}(r=1)=0$.
Furthermore, for $r=\infty$,
$\NPE^{r=\infty}(x) \cap \G_1^{r=\infty}(x)=\TY$ for all $x \in \TY \setminus \Y_3$.
Then $\E\left[ \rho^{\la}_n(r=\infty) \right] =
\E\left[ h^{\la}_{12}(r=\infty) \right]=
\mu_{\la}(r=\infty)=
P(X_2 \in \NPE^{r=\infty}(X_1) \cap \G_1^{r=\infty}(X_1)=
P(X_2 \in \TY)=1$.
Similarly, $P(\{X_2,X_3\} \subset \NPE^{r=\infty}(X_1) \cap \G_1^{r=\infty}(X_1))=1$.
Hence $\nu_{\la}(r=\infty)=0$.
Therefore, the CLT result in Equation \eqref{eqn:asy-norm-and} holds only for $r \in (1,\infty)$.
Furthermore, $\rho^{\la}_n(r=1)=0$ a.s. and $\rho^{\la}_n(r=\infty)=1$ a.s.

For $r=1$, $\NPE^{r=1}(x) \cup \G_1^{r=1}(x)$
has positive $\R^2$-Lebesgue measure.
Then
$P(\{X_2,X_3\} \subset \NPE^{r=1}(X_1) \cup \G_1^{r=1}(X_1))>0$.
Thus, $\nu_{\lo}(r=1) \not= 0$.
On the other hand, for $r=\infty$,
$\NPE^{r=\infty}(X_1) \cup \G_1^{r=\infty}(X_1))=\TY$ for all $X_1 \in \TY$.
Then $\E\left[ \rho^{\lo}_n(r=\infty) \right] = \E\left[ h^{\lo}_{12}(r=\infty) \right]=
P(X_2 \in \NPE^{r=\infty}(X_1) \cup \G_1^{r=\infty}(X_1))=
\mu_{\lo}(r=\infty)= P(X_2 \in \TY)=1$.
Similarly, $P(\{X_2,X_3\} \subset \NPE^{r=\infty}(X_1) \cup \G_1^{r=\infty}(X_1))=1$.
Hence $\nu_{\lo}(r=\infty)=0$.
Therefore, the CLT result for the OR-underlying case holds only for $r \in [1,\infty)$.
Moreover $\rho^{\lo}_n(r=\infty)=1$ a.s.

\begin{remark}
\label{rem:arc-density}
\textbf{Relative Arc Density of PCDs:}
The relative arc density of the digraph $D$ is denoted as $\rho(D)$.
For $X_i \stackrel{iid}{\sim}F$, $\rho(D)$ is also shown to be a $U$-statistic (\cite{ceyhan:arc-density-PE}),
$$\rho(D)=\frac{1}{n\,(n-1)}\sum\hspace*{-0.1 in}\sum_{i < j \hspace*{0.25 in}}   \hspace*{-0.1 in} \,h_{ij}$$
where
$h_{ij} = h(X_i,X_j;N)=\I(X_iX_j \in \A)=\I(X_j \in N(X_i)))$
is the number of arcs between $X_i$ and $X_j$ in $D$.
Here
$$0\le \E\left[ \rho(D) \right]=\frac{1}{n\,(n-1)}\sum\hspace*{-0.1 in}\sum_{i < j \hspace*{0.25 in}}
\hspace*{-0.1 in} \,\E[h_{ij}]=\E\left[ h_{12} \right]/2.$$
Furthermore,
Moreover,
$$\Cov\left[ h_{12},h_{13} \right]=
P(\{X_2,X_3\} \subset N(X_1))-\left[\E\left[ h_{12} \right]\right]^2.$$

Let $h_{ij}(r):=h(X_i,X_j;\NPE^r)=\I(X_j \in \NPE^r(X_i))$ for $i \not= j$
and the random variable $\rho_n(r) := \rho(\X_n;h,\NPE^r)$.
Let $\mu(r):=\E\left[ \rho_n(r) \right]$
and
$\nu(r):=\Cov\left[ h_{12}(r),h_{13}(r) \right]$.
A central limit theorem for $U$-statistics (\cite{lehmann:1999}) yields
\begin{eqnarray}
\sqrt{n}\left(\rho_n(r)-\mu(r)\right)
\stackrel{\mathcal{L}}{\longrightarrow}
\N\left(0,\nu(r)\right)
\end{eqnarray}
provided $\nu(r) > 0$.
The explicit forms of asymptotic mean $\mu(r)$
and variance $\nu(r)$ are provided in \cite{ceyhan:arc-density-PE}.
$\square$
\end{remark}

\section{Relative Edge Density under Null and Alternative Patterns}
\label{sec:rel-edge-null-alternative}

\subsection{Null Distribution of Relative Edge Density}
\label{sec:null-dist-edge-density}
The null hypothesis is generally some form of
{\em complete spatial randomness};
thus we consider
$$H_o: X_i \stackrel{iid}{\sim} \mathcal{U}(\TY).$$
If it is desired to have the sample size be a random variable,
we may consider a spatial Poisson point process on $\TY$ as our null hypothesis.

We first present a ``geometry invariance" result which will simplify our subsequent analysis
by allowing us to consider the special case of the equilateral triangle.
Let $\rho^{\la}_n(r) := \rho_{\la}(n;\mathcal{U}(\TY),\NPE^r)$ and
$\rho^{\lo}_n(r) := \rho_{\lo}(n;\mathcal{U}(\TY),\NPE^r)$.
\begin{theorem}
\label{thm:geo-inv-NYr-under}
\textbf{Geometry Invariance:}
Let $\Y_3 = \{\y_1,\y_2,\y_3\} \subset \mathbb{R}^2$
be three non-collinear points.
For $i=1,2,\ldots,n$
let $X_i \stackrel{iid}{\sim} F = \mathcal{U}(\TY)$,
the uniform distribution on the triangle $\TY$.
Then for any $r \in [1,\infty]$
the distribution of $\rho^{\la}_n(r)$ and $\rho^{\lo}_n(r)$
is independent of $\Y_3$,
and hence the geometry of $\TY$.
\end{theorem}

\noindent \textbf{Proof:}
A composition of translation, rotation, reflections, and scaling
will take any given triangle $T_o = T(\y_1,\y_2,\y_3)$
to the ``basic'' triangle $T_b = T((0,0),(1,0),(c_1,c_2))$
with $0 < c_1 \le \frac{1}{2}$, $c_2 > 0$ and $(1-c_1)^2+c_2^2 \le 1$,
preserving uniformity.
The transformation $\phi: \mathbb{R}^2 \rightarrow \mathbb{R}^2$
given by $\phi(u,v) = \left( u+\frac{1-2\,c_1}{\sqrt{3}}\,v,\frac{\sqrt{3}}{2\,c_2}\,v \right)$
takes $T_b$ to
the equilateral triangle
$T_e = T\left((0,0),(1,0),\left(1/2,\sqrt{3}/2\right)\right)$.
Investigation of the Jacobian shows that $\phi$
also preserves uniformity.
Furthermore, the composition of $\phi$ with the rigid motion transformations and scaling
maps
     the boundary of the original triangle $T_o$
  to the boundary of the equilateral triangle $T_e$,
     the median lines of $T_o$
  to the median lines of $T_e$,
and  lines parallel to the edges of $T_o$
  to lines parallel to the edges of $T_e$.
Since the joint distribution of any collection of the $h^{\la}_{ij}(r)$ and $h^{\lo}_{ij}(r)$
involves only probability content of unions and intersections
of regions bounded by precisely such lines,
and the probability content of such regions is preserved since uniformity is preserved,
the desired result follows.
$\blacksquare$

Based on Theorem \ref{thm:geo-inv-NYr-under},
for our proportional-edge proximity map and the uniform null hypothesis,
we may assume that
$\TY$ is a standard equilateral triangle
with $\Y_3 = \{(0,0),(1,0),(1/2,\sqrt{3}/2)\}$
henceforth.

In the case of this (proportional-edge proximity map, uniform null hypothesis) pair,
the asymptotic null distribution of
$\rho^{\la}_n(r)$ and
$\rho^{\lo}_n(r)$
as a function of $r$ can be derived.
Recall that
$\mu_{\la}(r)=\E\left[ h^{\la}_{12}(r) \right]=P(X_2 \in \NPE^r(X_1)\cap \G_1^r(X_1))=\mu_{\la}(r)$ and
$\mu_{\lo}(r)=\E\left[ h^{\lo}_{12} \right]=P(X_2 \in \NPE^r(X_1)\cup \G_1^r(X_1))=\mu_{\lo}(r)$
are the probability of an edge occurring between any two vertices in the AND- and OR-underlying graphs, respectively.

\begin{theorem}
\label{thm:asy-norm-under}
\textbf{Asymptotic Normality:}
For $r \in (1,\infty)$,
$$\sqrt{n}\left(\rho^{\la}_n(r)-\mu_{\la}(r)\right)\Big/\sqrt{4\,\nu_{\la}(r)}
\stackrel{\mathcal{L}}{\longrightarrow}  \mathcal{N}(0,1)$$
and
for $r \in [1,\infty)$,
$$\sqrt{n}\left(\rho^{\lo}_n(r)-\mu_{\lo}(r)\right)\Big/\sqrt{4\,\nu_{\lo}(r)}
\stackrel{\mathcal{L}}{\longrightarrow}  \mathcal{N}(0,1).$$
\end{theorem}
where
\begin{eqnarray}
\label{eqn:Asymean_and}
\mu_{\la}(r) =
 \begin{cases}
-{\frac{1}{54}}\,{\frac{(-1+r)(5\,r^5-148\,r^4+245\,r^3-178\,r^2-232\,r+128)}{r^2(r+2)(r+1)}} &\text{for} \quad r \in [1,4/3), \\
  -{\frac{1}{216}}\,{\frac{101\,r^5-801\,r^4+1302\,r^3-732\,r^2-536\,r+672}{r(r+2)(r+1)}}                                 &\text{for} \quad r \in [4/3,3/2), \\
  \frac{1}{8}\,{\frac{r^8-13\,r^7+30\,r^6+148\,r^5-448\,r^4+264\,r^3+288\,r^2-368\,r+96}{r^4(r+2)(r+1)}}  &\text{for} \quad r \in [3/2,2), \\
  \frac{(r^3+3\,r^2-2+2\,r)(-1+r)^2}{r^4(r+1)} &\text{for} \quad r \in [2,\infty),
 \end{cases}
\end{eqnarray}
\begin{eqnarray}
\label{eqn:Asymean_or}
\mu_{\lo}(r) =
\begin{cases}
           \frac{47\,r^6-195\,r^5+860\,r^4-846\,r^3-108\,r^2+720\,r-256}{108\,r^2(r+2)(r+1)}&\text{for} \quad r \in [1,4/3),\\
           \frac{175\,r^5-579\,r^4+1450\,r^3-732\,r^2-536\,r+672}{216\,r\,(r+2)(r+1)} &\text{for} \quad r \in [4/3,3/2),\\
           -\frac{3\,r^8-7\,r^7-30\,r^6+84\,r^5-264\,r^4+304\,r^3+144\,r^2-368\,r+96}{8\,r^4(r+1)(r+2)} &\text{for} \quad r \in [3/2,2),\\
           \frac{r^5+r^4-6\,r+2}{r^4(r+1)} &\text{for} \quad r \in [2,\infty),\\
         \end{cases}
\end{eqnarray}
\begin{equation}
\label{eqn:Asyvar_and}
\nu_{\la}(r)=\sum_{i=1}^{11}\vartheta^{\la}_i(r)\,\I(\mI_i),
\end{equation}
\begin{equation}
\label{eqn:Asyvar_or}
\nu_{\lo}(r)=\sum_{i=1}^{11}\vartheta^{\lo}_i(r)\,\I(\mI_i)
\end{equation}
where $\vartheta^{\la}_i(r)$ and $\vartheta^{\lo}_i(r)$ are provided in Appendix Sections 1 and 2,
and the derivations of $\mu_{\la}(r)$ and $\nu_{\la}(r)$ are provided in Appendix 3,
while those of $\mu_{\lo}(r)$ and $\nu_{\lo}(r)$ are provided in Appendix 4.

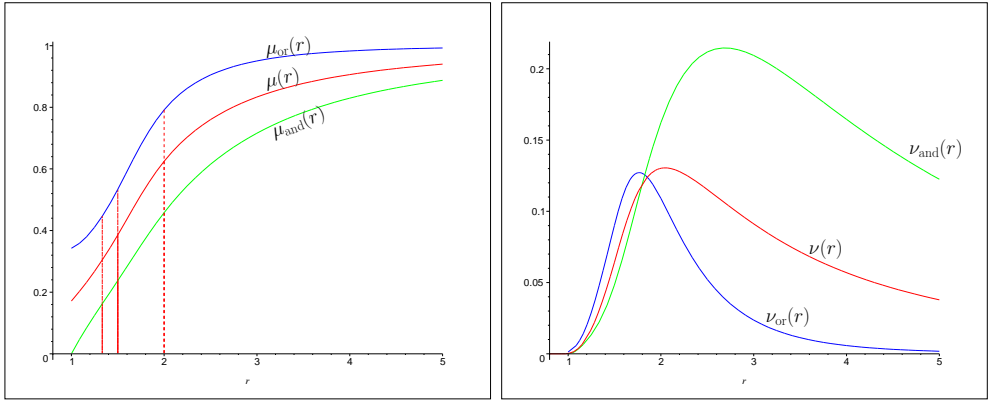
\begin{figure} [ht]
\centering
\scalebox{.3}{\input{means3.pstex_t}}
\scalebox{.3}{\input{var3.pstex_t}}
\caption{Result of Theorem \ref{thm:asy-norm-under}:
asymptotic null means $\mu(r)$, $\mu_{\la}(r)$, and $\mu_{\lo}(r)$ (left) and
variances $\nu(r)$, $4\,\nu_{\la}(r)$, and $4\,\nu_{\lo}(r)$ (right),
from Equations \eqref{eqn:Asymean_and}, \eqref{eqn:Asymean_or}, and
\eqref{eqn:Asyvar_and}, \eqref{eqn:Asyvar_or}, respectively.
Some values of note:
$\mu(1) = 37/216$, $\mu_{\la}(1) =0$, and $\mu_{\lo}(1) = 37/108$,
$\lim_{r \rightarrow \infty} \mu(r)=
\lim_{r \rightarrow \infty} \mu_{\la}(r)=\lim_{r \rightarrow \infty} \mu_{\lo}(r) = 1$,
$\nu_{\la}(r=1)=0$ and $\lim_{r \rightarrow \infty}\nu_{\la}(r)=0$,
$\nu_{\lo}(r=1)=1/3240$ and $\lim_{r \rightarrow \infty}\nu_{\lo}(r)=0$,
and $\argsup_{r \in [1,\infty]} \nu(r) \approx 2.045$ with
$\sup_{r \in [1,\infty]} \nu(r) \approx .1305$,
$\argsup_{r \in [1,\infty]}4\, \nu_{\la}(r) \approx 2.69$ with
$\sup_{r \in [1,\infty]} \nu_{\la}(r) \approx .0537$,
$\argsup_{r \in [1,\infty]}\nu_{\lo}(r) \approx 1.765$ with
$\sup_{r \in [1,\infty]}\nu_{\lo}(r) \approx .0318$.}
\label{fig:asymptotics}
\end{figure}

Notice that $\mu_{\la}(r=1)=0$ and $\lim_{r\rightarrow \infty}\mu_{\la}(r)=1$ (at rate $O(r^{-1})$);
and $\mu_{\lo}(r=1)=37/108$ and $\lim_{r\rightarrow \infty}\mu_{\lo}(r)=1$ (at rate $O(r^{-1})$).


To illustrate the limiting distribution, for example, $r=2$ yields
$$\frac{\sqrt{n}(\rho^{\la}_n(2)-\mu_{\la}(2))}{\sqrt{4\,\nu_{\la}(2)}}=
\sqrt{\frac{362880\,n}{58901}} \left(\rho^{\la}_n(2) -
\frac{11}{24}\right)\stackrel{\mathcal{L}}{\longrightarrow}
\N(0,1)$$ and
$$\frac{\sqrt{n}(\rho^{\lo}_n(2)-\mu_{\lo}(2))}{\sqrt{4\,\nu_{\lo}(2)}}=
\sqrt{\frac{120960\,n}{13189}} \left(\rho^{\lo}_n(2) -
\frac{19}{24}\right)\stackrel{\mathcal{L}}{\longrightarrow}
\N(0,1)$$ or equivalently,
$$\rho^{\la}_n(2) \stackrel{\text{approx}}{\sim} \N\left(\frac{11}{24},\frac{58901}{362880\,n}\right)
\text{  and  } \rho^{\lo}_n(2) \stackrel{\text{approx}}{\sim}
\N\left(\frac{19}{24},\frac{13189}{120960\,n}\right).$$

By construction of the underlying graphs,
there is a natural ordering of the means of relative arc and edge densities.
\begin{lemma}
The means of the relative edge densities and arc density have the following ordering:
$\mu_{\la}(r) < \mu(r) < \mu_{\lo}(r)$ for all $r \in [1,\infty)$.
Furthermore, for $r=\infty$ we have $\mu_{\la}(r) = \mu(r) = \mu_{\lo}(r)=1$.
\end{lemma}

\noindent \textbf{Proof:}
Recall that $\mu_{\la}(r) = \E[\rho_n^{\la}(r)]=P(X_2 \in \NPE^r(X_1) \cap \G_1^r(X_1))$,
$\mu(r) = \E[\rho_n(r)]=P(X_2 \in \NPE^r(X_1))$, and
$\mu_{\lo}(r) = \E[\rho_n^{\lo}(r)]=P(X_2 \in \NPE^r(X_1) \cup \G_1^r(X_1))$.
And $\NPE^r(X_1) \cap \G_1^r(X_1) \subseteq \NPE^r(X_1) \subseteq \NPE^r(X_1) \cup \G_1^r(X_1)$
with probability 1 for all $r \ge 1$ with equality holding for $r=\infty$ only.
Then the desired result follows.
See also Figure \ref{fig:asymptotics}.
$\blacksquare$

Note that the above lemma holds for all $X_i$ that has a continuous distribution on $\TY$.
There is also a stochastic ordering for the relative edge and
arc densities as follows.
\begin{theorem}
For sufficiently small $r$, $\rho_n^{\la}(r) <^{ST} \rho_n(r) <^{ST} \rho_n^{\lo}(r)$ as $n \rightarrow \infty$.
\end{theorem}

\noindent \textbf{Proof:}
Above we have proved that $\mu_{\la}(r) < \mu(r) < \mu_{\lo}(r)$ for all $r \in [1,\infty)$.
For small $r$ ($r \le \widehat r \approx 1.8$) the asymptotic variances have the same ordering,
$4\,\nu_{\la}(r) < \nu(r) <4\,\nu_{\lo}(r)$.
Since $\rho_n^{\la}(r), \, \rho_n(r),\, \rho_n^{\lo}(r)$ are asymptotically normal,
then the desired result follows.
See also Figure \ref{fig:asymptotics}.
$\blacksquare$

\begin{figure}[ht]
\centering
\psfrag{Density}{\Huge{density}}
\rotatebox{-90}{ \resizebox{2.1 in}{!}{ \includegraphics{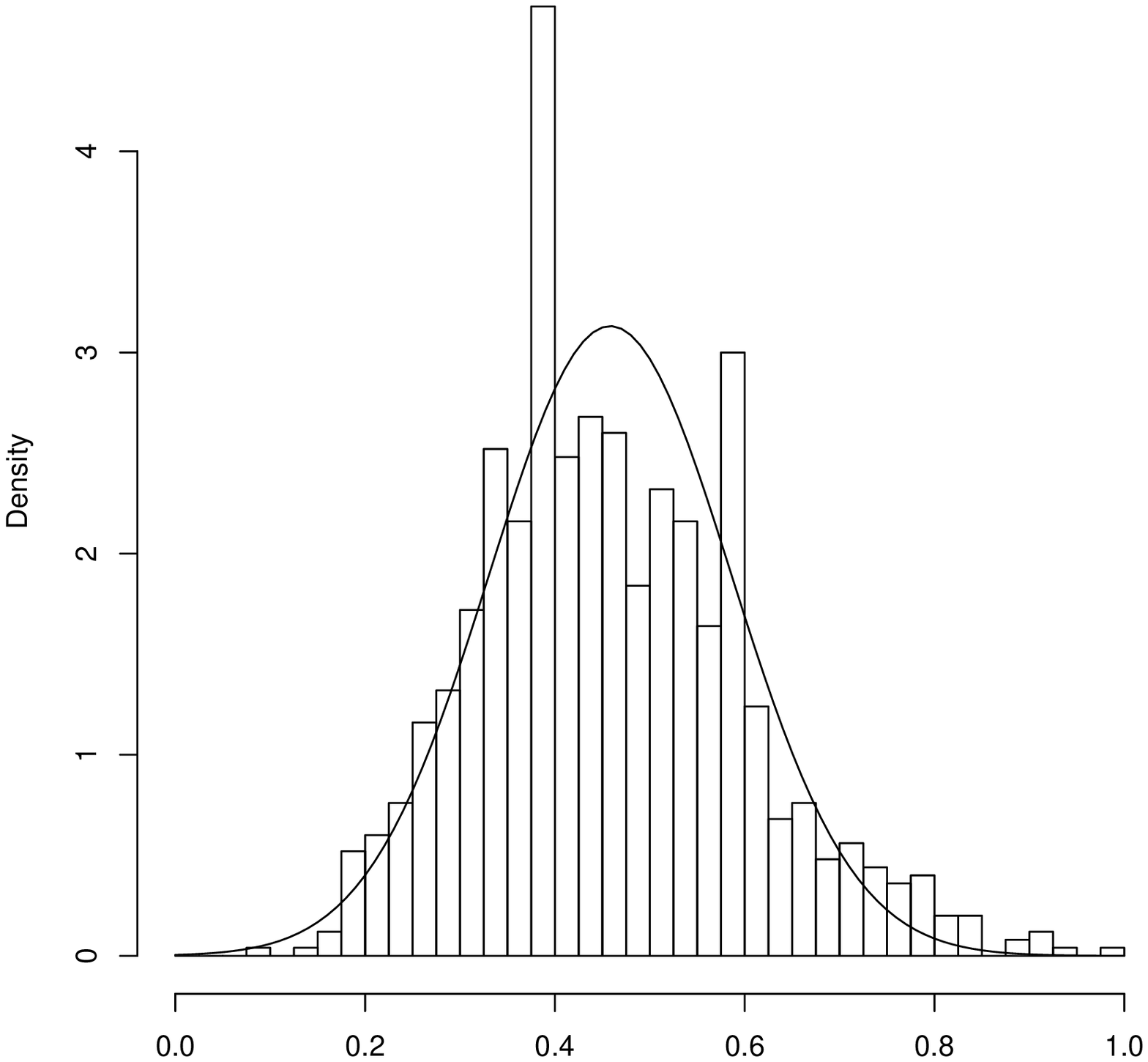} } }
\rotatebox{-90}{ \resizebox{2.1 in}{!}{ \includegraphics{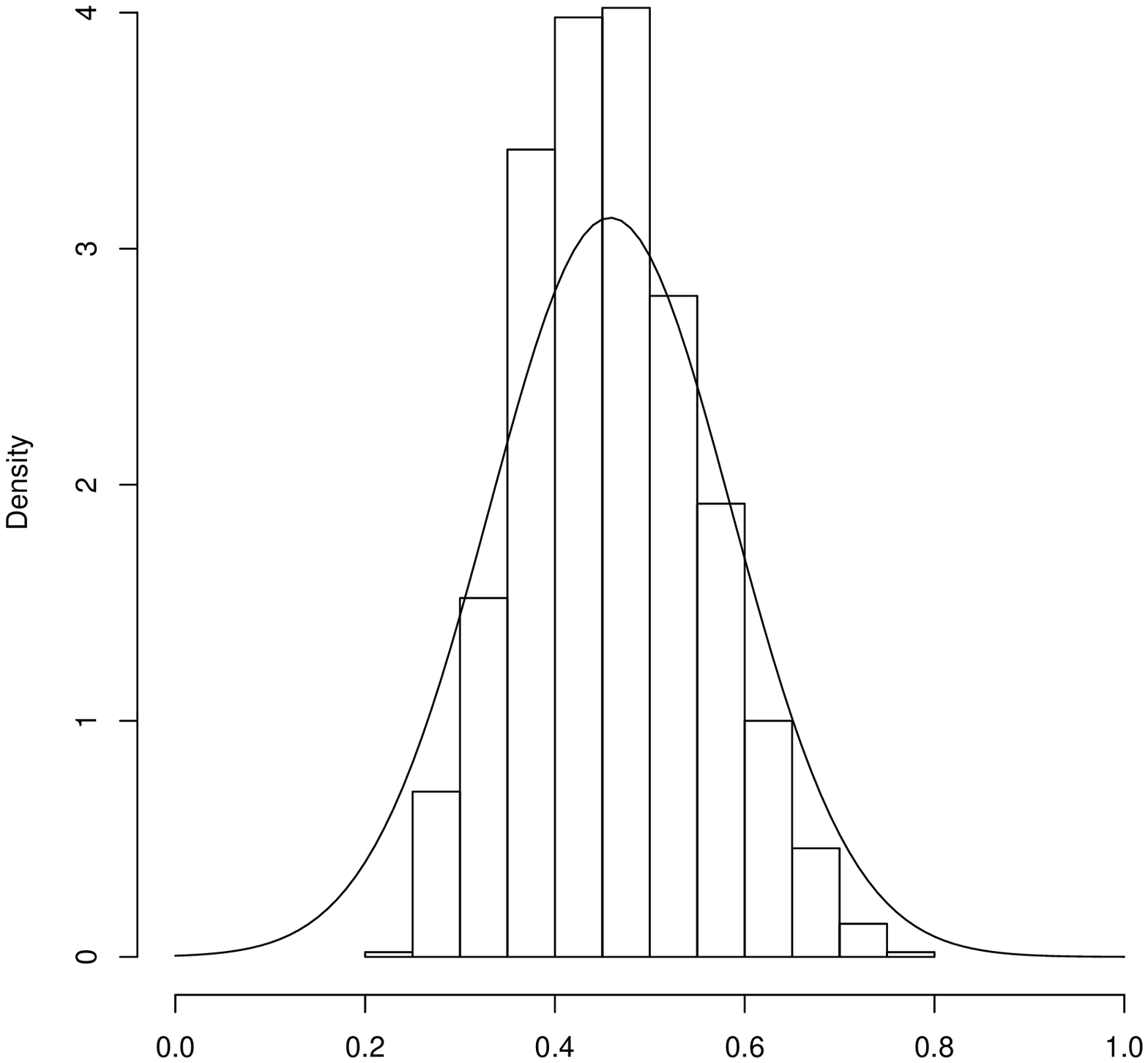} } }
\rotatebox{-90}{ \resizebox{2.1 in}{!}{ \includegraphics{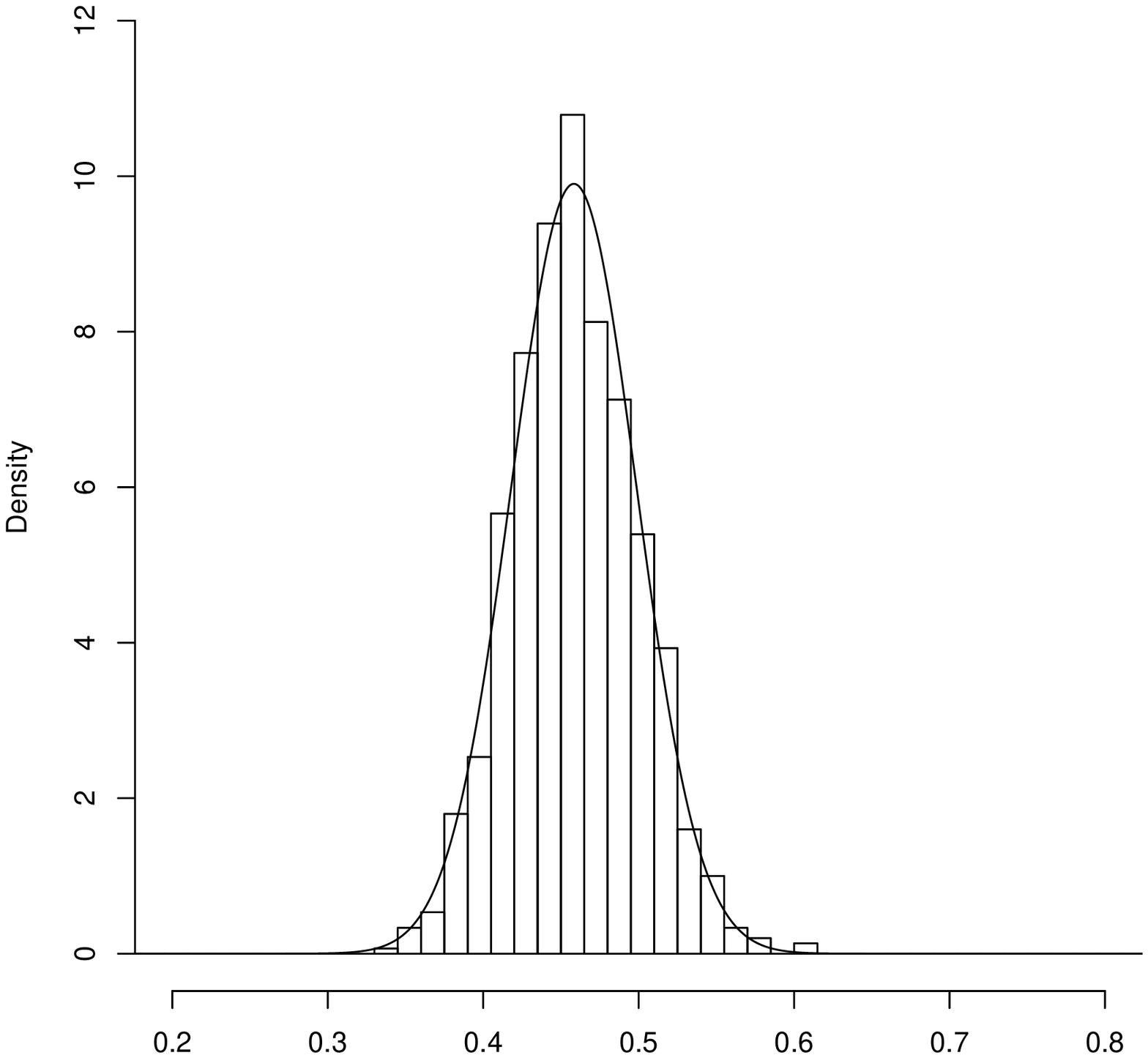} } }
\caption{
\label{fig:NormApprox_1}
Depicted are $\rho^{\la}_n(2) \stackrel{\text{approx}}{\sim} \mathcal{N}\left(\frac{11}{24},\frac{58901}{362880\,n}\right)$
for $n=10,\,20,\,100$ (left to right).
Histograms are based on 1000 Monte Carlo replicates.
Solid lines are the corresponding normal densities.
Notice that the vertical axes are differently scaled.
}
\end{figure}

\begin{figure}[ht]
\centering
\psfrag{Density}{\Huge{density}}
\rotatebox{-90}{ \resizebox{2.1 in}{!}{ \includegraphics{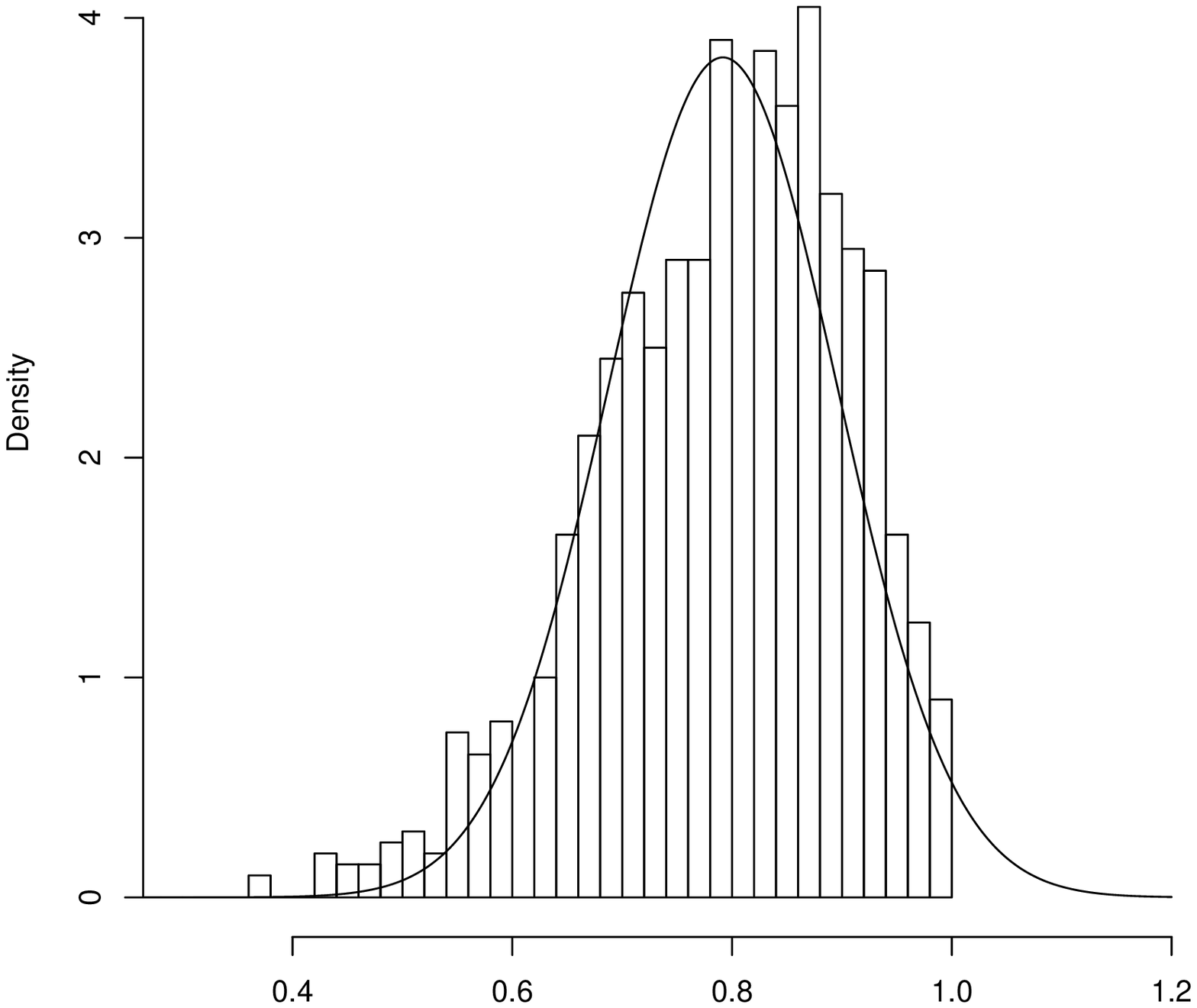} } }
\rotatebox{-90}{ \resizebox{2.1 in}{!}{ \includegraphics{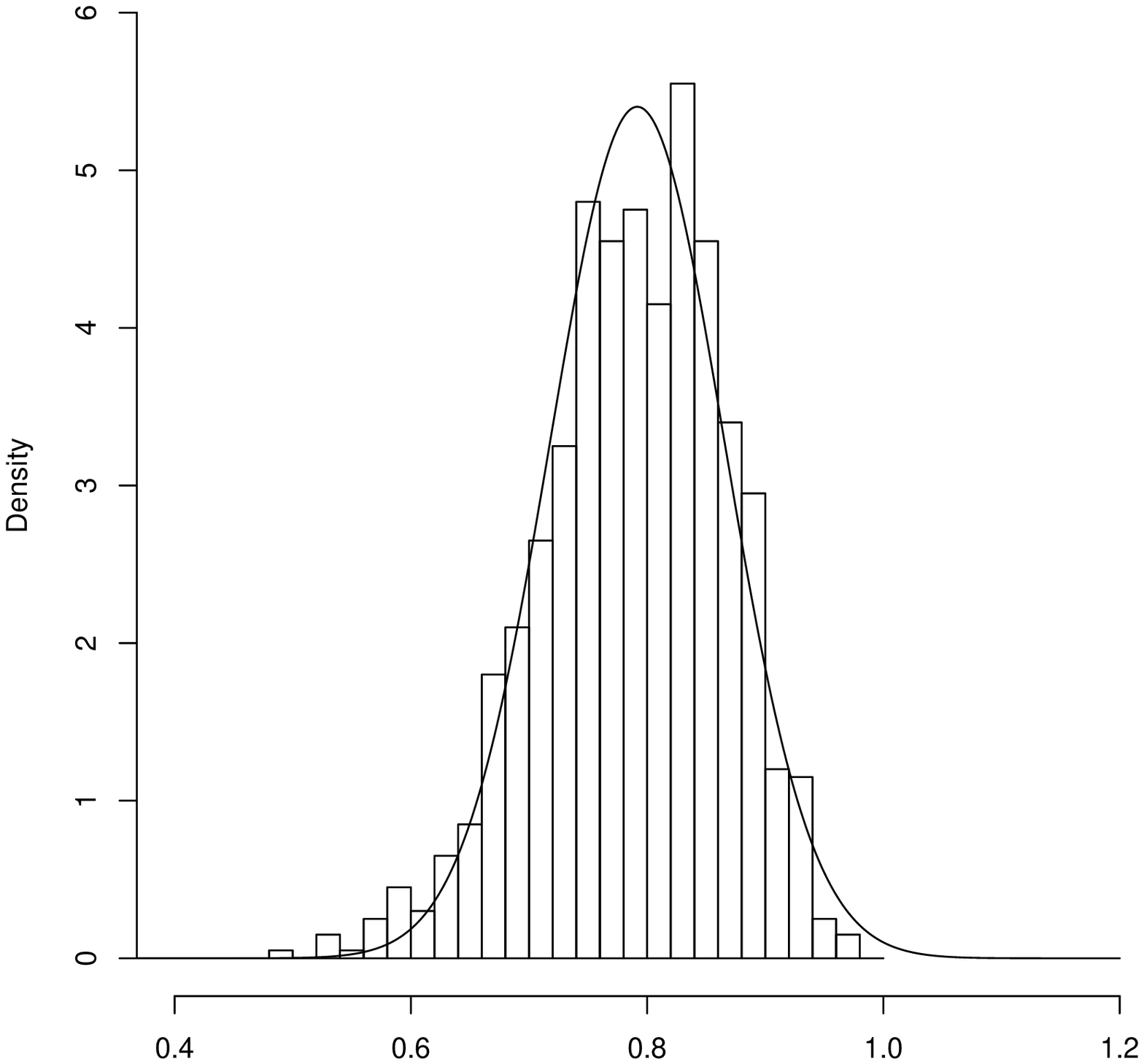} } }
\rotatebox{-90}{ \resizebox{2.1 in}{!}{ \includegraphics{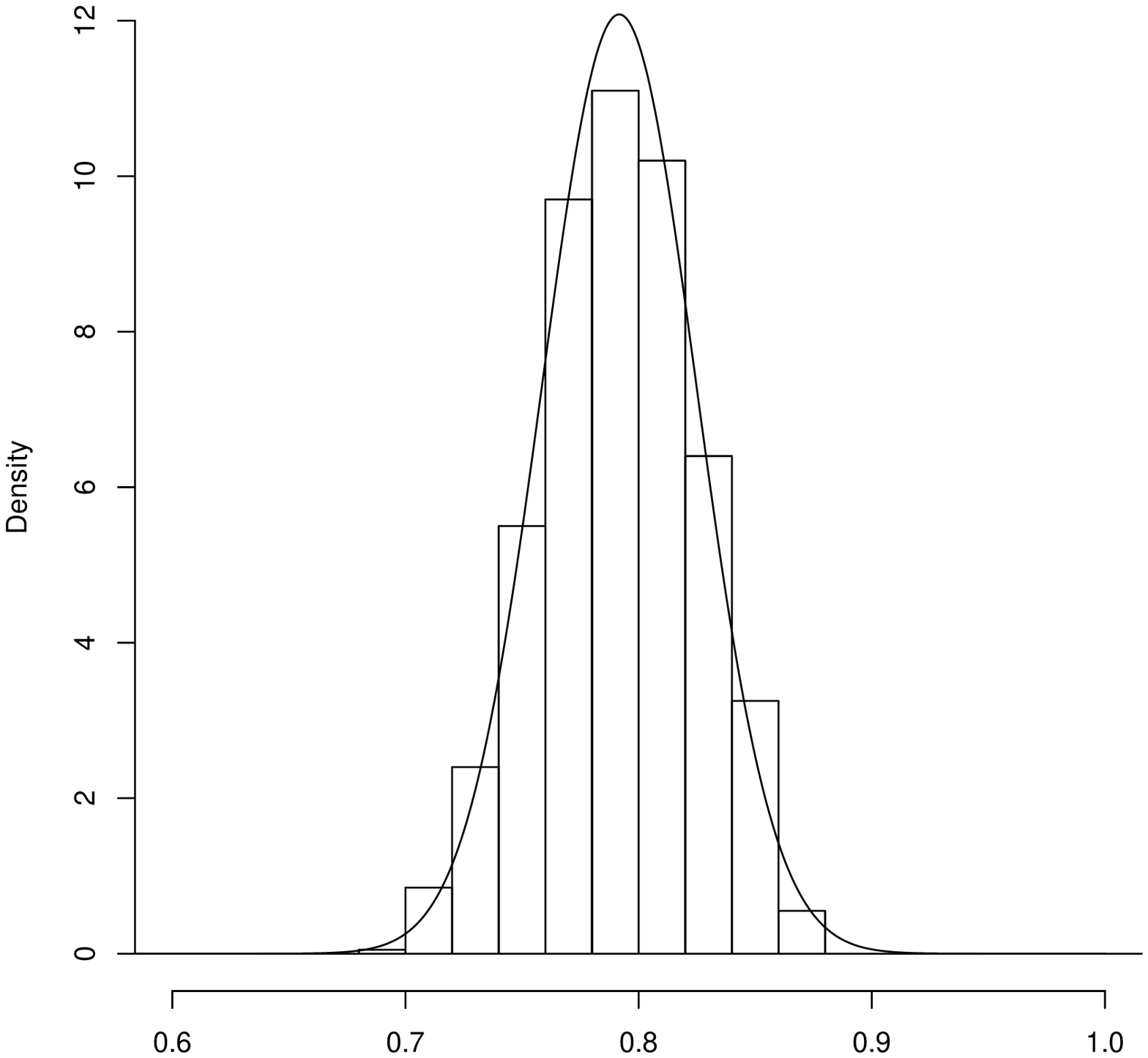} } }
\caption{
\label{fig:NormApprox2}
Depicted are $\rho^{\lo}_n(2) \stackrel{\text{approx}}{\sim} \mathcal{N}\left(\frac{19}{24},\frac{13189}{120960\,n}\right)$
for $n=10,20,100$ (left to right).
Histograms are based on 1000 Monte Carlo replicates.
Solid lines are the corresponding normal densities.
Notice that the vertical axes are differently scaled.
}
\end{figure}

\begin{figure}[ht]
\centering
\psfrag{Density}{\Huge{density}}
\rotatebox{-90}{ \resizebox{2.1 in}{!}{ \includegraphics{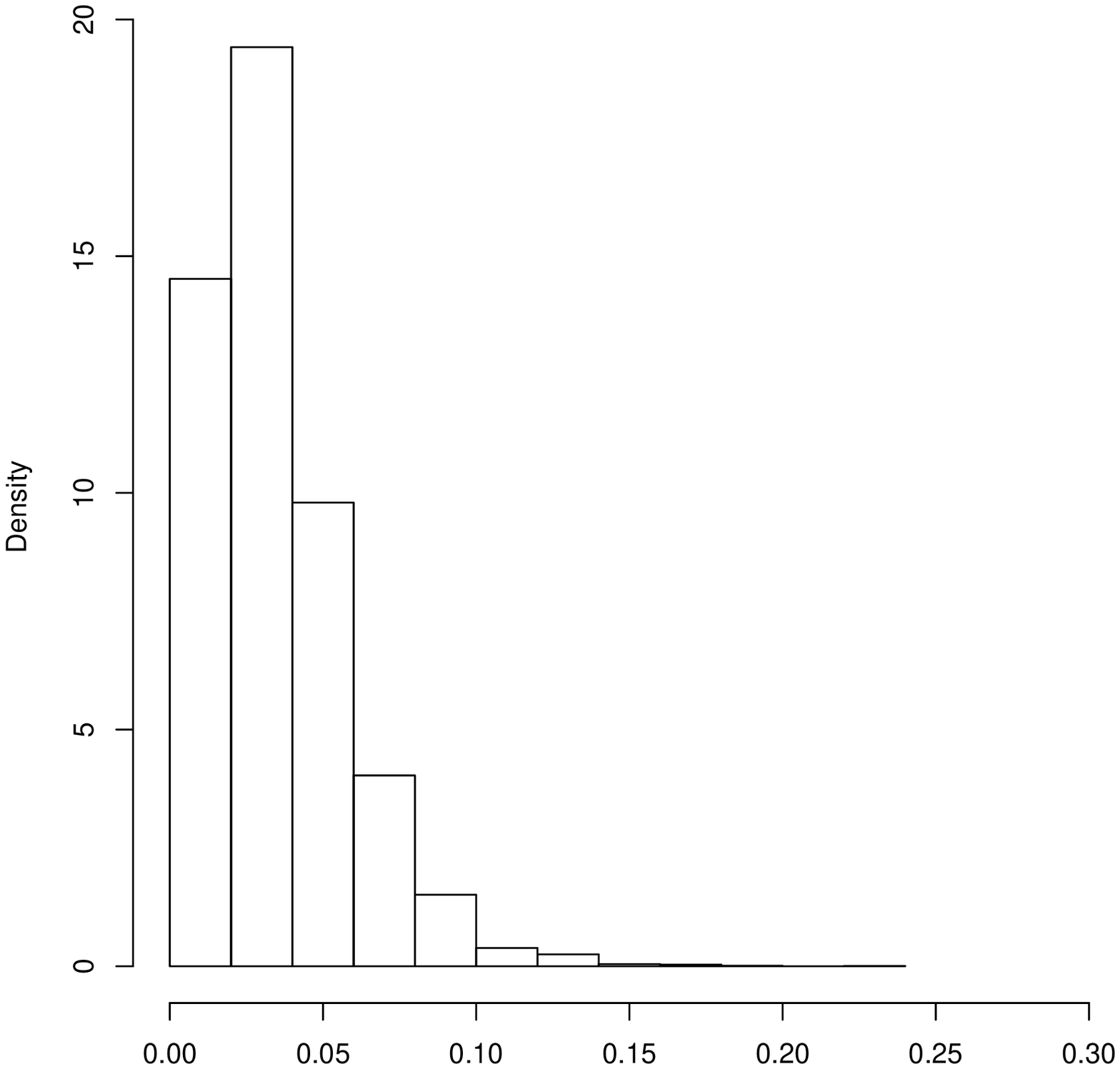} } }
\rotatebox{-90}{ \resizebox{2.1 in}{!}{ \includegraphics{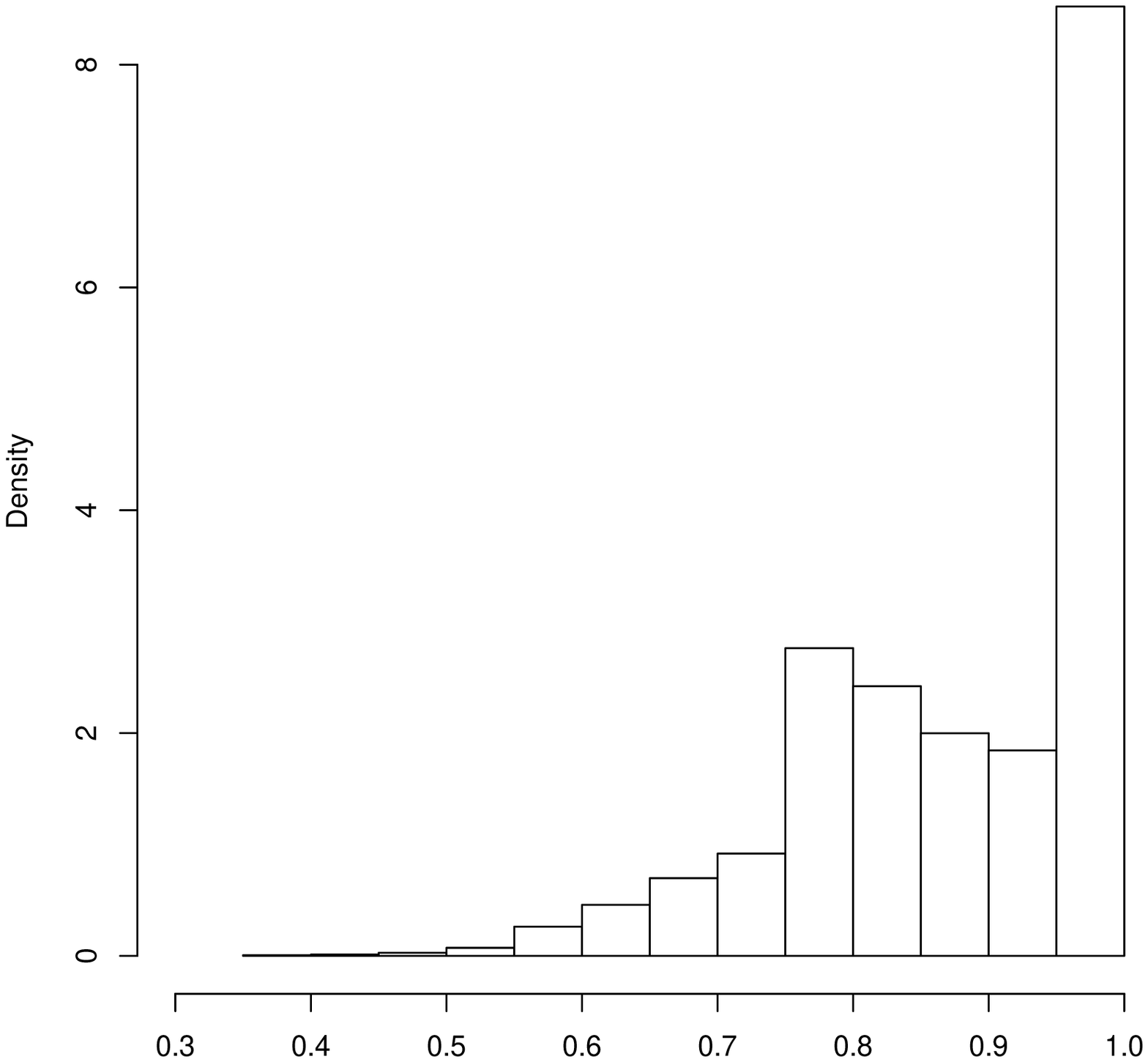} } }
\caption{
\label{fig:ANDskew}
Depicted are the histograms for 10000 Monte Carlo replicates of $\rho^{\la}_{10}(1.05)$ (left)
and $\rho^{\la}_{10}(5)$ (right) indicating severe small sample skewness for extreme values of $r$.
Notice that the vertical axes are differently scaled.
}
\end{figure}

\begin{figure}[ht]
\centering
\psfrag{Density}{\Huge{density}}
\rotatebox{-90}{ \resizebox{2.1 in}{!}{ \includegraphics{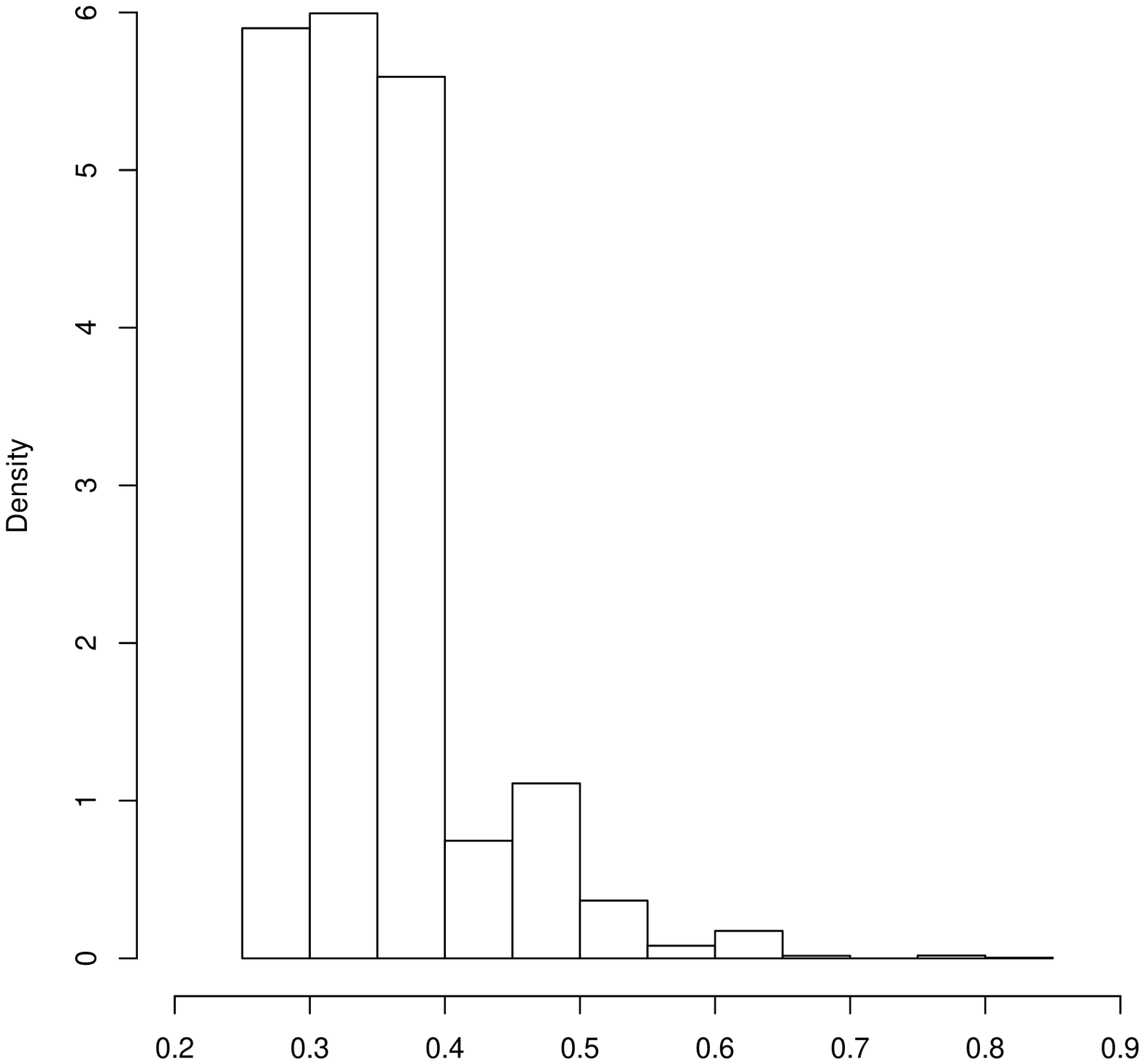} } }
\rotatebox{-90}{ \resizebox{2.1 in}{!}{ \includegraphics{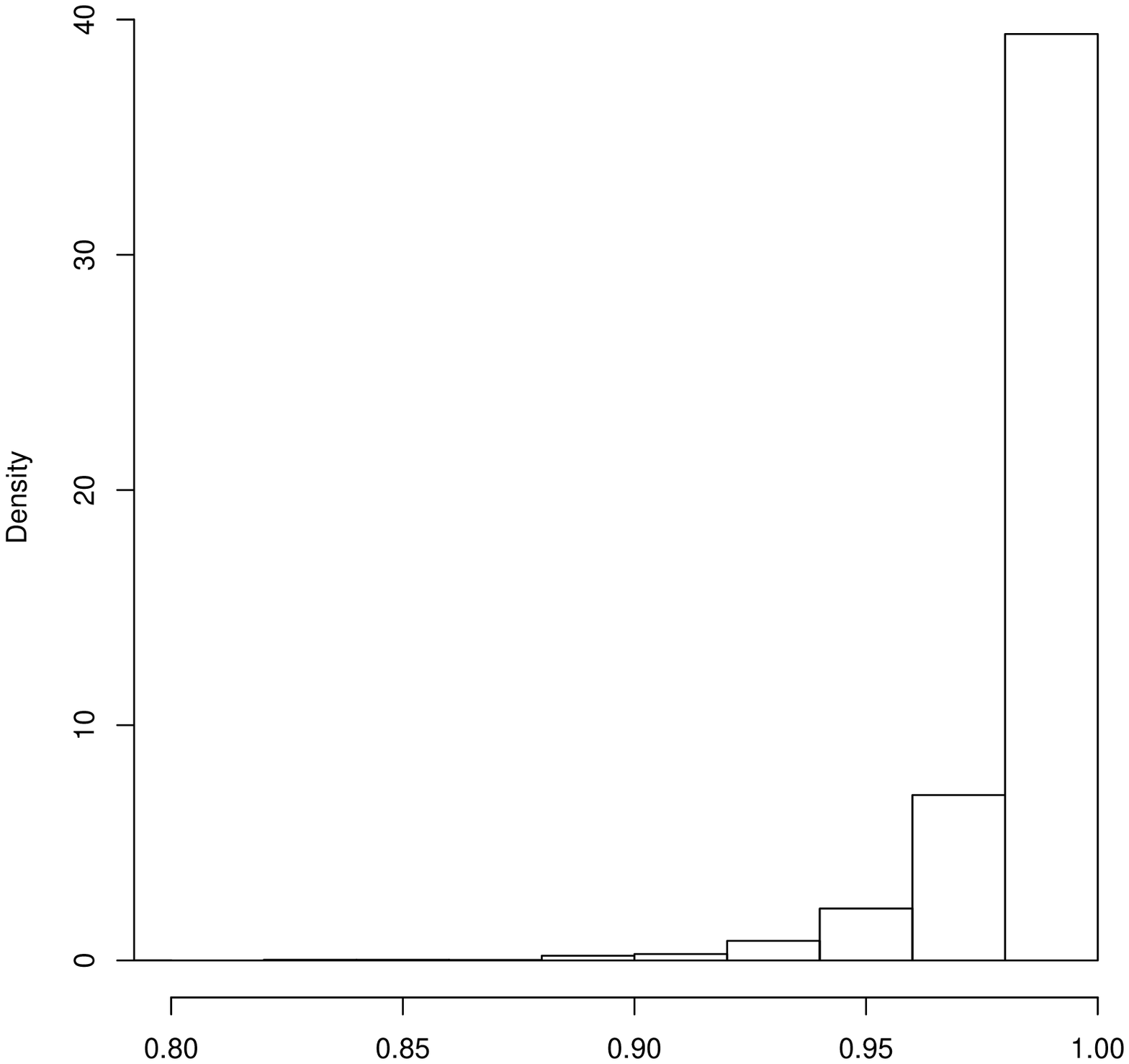} } }
\caption{
\label{fig:ORskew}
Depicted are the histograms for 10000 Monte Carlo replicates of $\rho^{\lo}_{10}(1)$ (left)
and $\rho^{\lo}_{10}(5)$ (right) indicating severe small sample skewness for extreme values of $r$.
Notice that the vertical axes are differently scaled.
}
\end{figure}

Figures \ref{fig:NormApprox_1} and \ref{fig:NormApprox2}
indicate that, for $r=2$,
the normal approximation is accurate even for small $n$
although kurtosis may be indicated for $n=10$ in the AND-underlying case,
and skewness may be indicated for $n=10$ in the OR-underlying case.
Figures \ref{fig:ANDskew} and \ref{fig:ORskew} demonstrate,
however, that severe skewness obtains for some values of $n$, $r$.
The finite sample variance and skewness may be derived analytically
in much the same way as was $4\,\nu_{\la}(r)$
(and $4\,\nu_{\lo}(r)]$)
for the asymptotic variance.
In fact,
the exact distribution of $\rho^{\la}_n(r)$ (and $\rho^{\lo}_n(r)$)
is, in principle, available
by successively conditioning on the values of the $X_i$.
Alas,
while the joint distribution of $h^{\la}_{12}(r),h^{\la}_{13}(r)$
(and $h^{\lo}_{12}(r),h^{\lo}_{13}(r)$) is available,
the joint distribution of $\{h^{\la}_{ij}(r)\}_{1 \leq i < j \leq n}$ (and $\{h^{\lo}_{ij}(r)\}_{1 \leq i < j \leq n}$),
and hence the calculation for the exact distribution of $\rho^{\la}_n(r)$ (and $\rho^{\lo}_n(r)$),
is extraordinarily tedious and lengthy for even small values of $n$.

\vspace*{0.05 in}

Let $\g_n(r)$ be the domination number of the proportional-edge PCD
based on $\X_n$ which is a random sample from $\U(\TY)$.
Additionally, let $\g^\la_n(r)$ and $\g^\lo_n(r)$ be
the domination number of the AND- and OR-underlying graphs
based on the proportional-edge PCD, respectively.
Then we have the following stochastic ordering for the domination numbers.
\begin{theorem}
For all $r \in [1,\infty)$ and $n>1$,
$\g^\lo_n(r) <^{ST} \g_n(r) <^{ST} \g^\la_n(r)$.
\end{theorem}

\noindent \textbf{Proof:}
For all $x \in \TY $,
we have $\NPE^r(x)\cap \G^r_1\left( x \right) \subseteq
\NPE^r(x) \subseteq \NPE^r(x) \cup \G^r_1\left( x \right)$.
For $X \sim \U(\TY)$, we have $\NPE^r(X)\cap \G^r_1\left( X \right) \subsetneq
\NPE^r(X) \subsetneq \NPE^r(X) \cup \G^r_1\left( X \right)$ a.s.
Moreover,
$\g_n(r)=1$ iff $\X_n \subset \NPE^r(X_i)$ for some $i$;
$\g^\la_n(r)=1$ iff $\X_n \subset \NPE^r(X_i)\cap \G^r_1\left( X_i \right)$ for some $i$;
and
$\g^\lo_n(r)=1$ iff $\X_n \subset \NPE^r(X_i)\cup \G^r_1\left( X_i \right)$ for some $i$.
So it follows that
$P(\g^\la_n(r)=1) < P(\g_n(r)=1) < P(\g^\lo_n(r)=1)$.
In a similar fashion,
we have $P(\g^\la_n(r)\le 2) < P(\g_n(r)\le 2) < P(\g^\lo_n(r) \le 2)$.
Since $P(\g_n(r)\le 3)=1$ (\cite{ceyhan:dom-num-NPE-SPL}),
it follows that $P(\g^\lo_n(r) \le 3)=1$ also holds
as $P(\g_n(r)\le 3) < P(\g^\lo_n(r) \le 3)$.
Hence the desired stochastic ordering follows.
$\blacksquare$

Note the stochastic ordering in the above theorem
holds for any continuous distribution $F$ with support
being in $\TY$.

\subsection{Alternatives: Segregation and Association}
\label{sec:alt-seg-assoc}
The phenomenon known as {\em segregation}
involves observations from different classes
having a tendency to repel each other
--- in our case, this means the $X_i$ tend to fall away from all elements of $\Y_3$.
{\em Association} involves observations from different classes
having a tendency to attract one another,
so that the $X_i$ tend to fall near an element of $\Y_3$.
See, for instance, \cite{dixon:1994} and \cite{coomes:1999}.

We define two simple classes of alternatives,
$H^S_{\ve}$ and $H^A_{\ve}$
with $\ve \in \left( 0,\sqrt{3}/3 \right)$,
for segregation and association, respectively.
For $\y \in \Y_3$,
let $e(\y)$ denote the edge of $\TY$ opposite vertex $\y$,
and for $x \in \TY$
let $\ell_\y(x)$ denote the line parallel to $e(\y)$ through $x$.
Then define
$T(\y,\ve) = \{x \in \TY: d(\y,\ell_\y(x)) \leq \ve\}$.
Let $H^S_{\ve}$ be the model under which
$X_i \stackrel{iid}{\sim} \mathcal{U}(\TY \setminus \cup_{\y \in \Y_3} T(\y,\ve))$
and $H^A_{\ve}$ be the model under which
$X_i \stackrel{iid}{\sim} \mathcal{U}(\cup_{\y \in \Y_3} T(\y,\sqrt{3}/3 - \ve))$.
Thus the segregation model excludes the possibility of
any $X_i$ occurring near a $\y_j$,
and the association model requires
that all $X_i$ occur near a $\y_j$.
The $\sqrt{3}/3 - \ve$ in the definition of the
association alternative is so that $\ve=0$
yields $H_o$ under both classes of alternatives.

\begin{remark}
These definitions of the alternatives
are given for the standard equilateral triangle.
The geometry invariance result of Theorem \ref{thm:geo-inv-NYr-under}
still holds under the alternatives $H^S_{\ve}$ and $H^A_{\ve}$.
In particular, the segregation alternative with $\ve \in \left( 0,\sqrt{3}/4 \right)$ in the standard equilateral triangle
corresponds to the case that in an arbitrary triangle, $\delta \times 100\%$ of the area is carved
away as forbidden from the vertices using line segments parallel to the opposite edge
where $\delta = 4\ve^2$ (which implies $\delta \in (0,3/4)$).
But the segregation alternative with $\ve \in \left( \sqrt{3}/4,\sqrt{3}/3 \right)$ in the standard equilateral triangle
corresponds to the case that in an arbitrary triangle, $\delta \times 100\%$ of the area is carved
away as forbidden around the vertices using line segments parallel to the opposite edge
where $\delta = 1-4 \left(1-\sqrt{3}\ve \right)^2$ (which implies $\delta \in (3/4,1)$).
This argument is for the segregation alternative;
a similar construction is available for the association alternative.
$\square$
\end{remark}

The asymptotic normality of the relative edge density under the alternatives follows as in the null case.
\begin{theorem}
\label{thm:asy-norm-under-alt}
\textbf{Asymptotic Normality under the Alternatives:}
Let $\mu_{\la}(r,\ve)$ be the mean and $\nu_{\la}(r,\ve)$ be the variance of $\rho^{\la}_n(r)$
under the alternatives for $r \in [1,\infty)$ and $\ve \in \left( 0,\sqrt{3}/3 \right)$.
Then under $H^S_{\ve}$ and $H^A_{\ve}$,
$\sqrt{n}(\rho^{\la}_n(r)-\mu_{\la}(r,\ve)) \stackrel {\mathcal L}{\longrightarrow} \N(0,4\,\nu_{\la}(r,\ve))$
for the values of the pair $(r,\ve)$ for which $\nu_{\la}(r,\ve)>0$.
A similar result holds for $\rho_n^{\lo}(r)$.
\end{theorem}

\noindent \textbf{Proof:}
Under the alternatives, i.e.,  $\ve>0$ ,
$\rho^{\la}_n(r)$ is a $U$-statistic
with the same symmetric kernel $h^{\la}_{ij}(r)$ as in the null case.
Let $\E_{\ve}[\cdot]$ be the expectation with respect to the uniform distribution under the
alternatives with $\ve \in \left( 0,\sqrt{3}/3 \right)$.
The mean $\mu_{\la}(r,\ve)=\E_{\ve}\left[\rho^{\la}_n(r)\right] = \E_{\ve}\left[ h^{\la}_{12}(r) \right]$,
now a function of both $r$ and $\ve$, is again in $[0,1]$. The asymptotic variance,
$4\,\nu_{\la}(r,\ve)=4\,\Cov\left[ h^{\la}_{12}(r),h^{\la}_{13}(r) \right]$, also a function of both $r$ and $\ve$,
is bounded above by $1/4$, as before.
Thus asymptotic normality obtains provided $\nu_{\la}(r,\ve) > 0$;
otherwise $\rho^{\la}_n(r)$ is degenerate.
Then under $H^S_{\ve}$, $\nu_{\la}(r,\ve)>0$ for $(r,\ve)$ in
  $(1,\sqrt{3}/(2 \ve)) \times (0,\sqrt{3}/4]$ or
  $(1,\sqrt{3}/\ve-2) \times (\sqrt{3}/4,\sqrt{3}/3)$,
and under $H^A_{\ve}$, $\nu_{\la}(r,\ve)>0$ for $(r,\ve)$ in
 $(1,\infty) \times \left( 0,\sqrt{3}/3 \right)$.
Also under $H^S_{\ve}$, $\nu_{\lo}(r,\ve)>0$ for $(r,\ve)$ in
  $[1,\sqrt{3}/(2 \ve)) \times (0,\sqrt{3}/4]$ or
  $[1,\sqrt{3}/\ve-2) \times (\sqrt{3}/4,\sqrt{3}/3)$,
and under $H^A_{\ve}$, $\nu_{\lo}(r,\ve)>0$ for $(r,\ve)$ in
 $(1,\infty) \times \left( 0,\sqrt{3}/3 \right)$ or $\{1\} \times (0,\sqrt{3}/12)$.
$\blacksquare$

Notice that for the association class of alternatives
any $r \in (1,\infty)$ yields asymptotic normality
for all $\ve \in \left( 0,\sqrt{3}/3 \right)$ in both AND- and OR-underlying cases,
while for the segregation class of alternatives
only $r=1$ yields this universal asymptotic normality in the OR-underlying case,
and such an $\ve$ does not exist for the AND-underlying case.

The relative edge density of the underlying graphs based on the PCD
is a test statistic for the segregation/association alternative;
rejecting for extreme values of $\rho^{\la}_n(r)$ is appropriate
since under segregation we expect $\rho^{\la}_n(r)$ to be large,
while under association we expect $\rho^{\la}_n(r)$ to be small. The same holds for $\rho^{\lo}_n(r)$.
Using the test statistics
\begin{equation}
R^{\la}_n(r) = \sqrt{n} \left(\rho^{\la}_n(r) - \mu_{\la}(r)\right)\Big/\sqrt{4\,\nu_{\la}(r)}, \text{ and }
R^{\lo}_n(r) = \sqrt{n} \left(\rho^{\lo}_n(r) - \mu_{\lo}(r)\right)\Big/\sqrt{4\,\nu_{\lo}(r)}
\end{equation}
for AND- and OR-underlying cases, respectively,
the asymptotic critical value
for the one-sided level $\alpha$ test against segregation
is given by
\begin{equation}
z_{\alpha} = \Phi^{-1}(1-\alpha)
\end{equation}
where $\Phi(\cdot)$ is the standard normal distribution function.
The test rejects for $R^{\la}_n(r)>z_{\alpha}$ against segregation.
Against association,
the test rejects for $R^{\la}_n(r)<z_{1-\alpha}$.
The same holds for the test statistic $R^{\lo}_n(r)$.

\section{Asymptotic Performance of Relative Edge Density}
\label{sec:asymptotic}

\subsection{Consistency}
\label{sec:consistency}

\begin{theorem}
The test against $H^S_{\ve}$ which rejects for
$R^{\la}_n(r)>z_{\alpha}$ and the test against $H^A_{\ve}$ which
rejects for $R^{\la}_n(r)<z_{1-\alpha}$ are consistent for $r \in
(1,\infty)$ and $\ve \in \left( 0,\sqrt{3}/3 \right)$.
The same holds for
$R^{\lo}_n(r)$ with $r\in [1,\infty)$.
\end{theorem}

\noindent \textbf{Proof:}
Since the variance of the asymptotically normal test statistic,
under both the null and the alternatives,
converges to 0 as $n \rightarrow \infty$
(or might be zero for $n < \infty$),
it remains to show that the mean under the null, $\mu_{\la}(r)=\E\left[\rho^{\la}_n(r)\right]$,
is less than (greater than) the mean under the alternative,
$\mu_{\la}(r,\ve)=\E\left[\rho^{\la}_n(r)\right]$ against segregation (association) for $\ve > 0$.
Whence it will follow that power converges to 1 as $n \rightarrow \infty$.
Let $P_{\ve}(\cdot)$ be the probability with respect to the uniform distribution under the
alternatives with $\ve \in \left( 0,\sqrt{3}/3 \right)$.
Then against segregation, we have
{\small
\begin{eqnarray*}
\mu_{\la}(r) &=& P_0(X_2 \in \NPE^r(X_1) \cap \G_1^r(X_1))\\
&=& P_0(X_2 \in \NPE^r(X_1) \cap \G_1^r(X_1), X_1 \in \cup_{\y \in \Y_3} T(\y,\ve))+
P_0(X_2 \in \NPE^r(X_1) \cap \G_1^r(X_1), X_1 \in \TY
\setminus \cup_{\y \in \Y_3} T(\y,\ve))\\
&=& P_0(X_2 \in \NPE^r(X_1) \cap \G_1^r(X_1)| X_1 \in \cup_{\y \in \Y_3} T(\y,\ve))\,
P_0( X_1 \in \cup_{\y \in \Y_3} T(\y,\ve))\\
& &+P_0(X_2 \in \NPE^r(X_1)| X_1 \in \TY \setminus
\cup_{\y \in \Y_3} T(\y,\ve))\,P_0(X_1 \in \TY \setminus \cup_{\y \in \Y_3} T(\y,\ve))\\
&<&P_0(X_2 \in \NPE^r(X_1) \cap \G_1^r(X_1)| X_1 \in \cup_{\y \in \Y_3} T(\y,\ve))\,p_1 +
P_{\ve}(X_2 \in \NPE^r(X_1) \cap \G_1^r(X_1)| X_1 \in \TY \setminus \cup_{\y \in \Y_3} T(\y,\ve))\,p_2\\
&=& \E_0(I(X_2 \in \NPE^r(X_1) \cap \G_1^r(X_1))| X_1 \in \cup_{\y \in \Y_3} T(\y,\ve))\,p_1\\
& & +\E_{\ve}(I(X_2 \in \NPE^r(X_1) \cap \G_1^r(X_1))| X_1 \in
\TY \setminus \cup_{\y \in \Y_3} T(\y,\ve))\,p_2
\end{eqnarray*}
}
where $p_1=P_0(X_1 \in \cup_{\y \in \Y_3} T(\y,\ve))$ and
$p_2=P_0(X_1 \in \TY \setminus \cup_{\y \in \Y_3} T(\y,\ve))=1-p_1$.
Then $$\mu_{\la}(r,\ve)>\mu_{\la}(r)\,\frac{(1-p_1)}{p_2}=\mu_{\la}(r).$$
Likewise, we have $\mu_{\la}(r,\ve)=\E_{\ve}\left[\rho^{\la}_n(r)\right]<\E[\rho_n(r)]=\mu_{\la}(r)$, for association.

The consistency follows for the OR-underlying case in a similar fashion.
$\blacksquare$

\subsection{Pitman Asymptotic Efficiency}
\label{sec:PAE}
Pitman asymptotic efficiency (PAE)
provides an investigation of ``local asymptotic power''
--- local about $H_o$.
This involves the limit as $\ve \rightarrow 0$ as well as
the limit as $n \rightarrow \infty$.
A detailed discussion of PAE can be found in \cite{kendall:1979} and \cite{eeden:1963}.
For segregation or association alternatives with the AND-underlying graphs,
the PAE is given by $\displaystyle \frac{\left( (\mu_{\la})^{(k)}(r,\ve=0) \right)^2}{\nu_{\la}(r)}$
where $(\mu_{\la})^{(k)}(r,\ve=0)$ is the $k^{th}$ derivative with respect to $\ve$
so that $(\mu_{\la})^{(k)}(r,\ve=0) \not=0$ but $(\mu_{\la})^{(k-1)}(r,\ve=0)=0$ for $k=1,2,\ldots$.
Likewise the same holds for the OR-underlying case.
Then under segregation alternative $H^S_{\ve}$, the PAE is given by
$$
\PAE^S_\la (r) =  \frac{\left( (\mu_{\la})^{\prime\prime}(r,\ve=0) \right)^2}{\nu_{\la}(r)}
\text{ and }
\PAE^S_\lo (r) =  \frac{\left( (\mu_{\lo})^{\prime\prime}(r,\ve=0) \right)^2}{\nu_{\lo}(r)}
$$
since $(\mu_{\la})^{\prime}(r,\ve=0) = 0$ and $(\mu_{\lo})^{\prime}(r,\ve=0) = 0$.
Under association alternative
$H^A_{\ve}$ is
$$
\PAE^A_\la (r) = \frac{\left( (\mu_{\la})^{\prime\prime}(r,\ve=0) \right)^2}{\nu_{\la}(r)}
\text{ and }
\PAE^A_\lo (r) = \frac{\left( (\mu_{\lo})^{\prime\prime}(r,\ve=0) \right)^2}{\nu_{\lo}(r)}
$$
since $(\mu_{\la})^{\prime}(r,\ve=0) =(\mu_{\lo})^{\prime}(r,\ve=0) = 0$.
Equations (\ref{eqn:Asyvar_and}) and (\ref{eqn:Asyvar_or})  provide the denominators;
the numerators require a bit of additional work,
but $\mu_{\la}(r,\ve)$ and $\mu_{\lo}(r,\ve)$ are available for small enough $\ve$,
which is all we need here.
See Appendix 5 for explicit forms of $\mu_{\la}(r,\ve)$
and $\mu_{\lo}(r,\ve)$ for segregation and association,
and the derivations of $\mu_{\la}(r,\ve)$ and $\mu_{\lo}(r,\ve)$
are provided in Appendix 6.

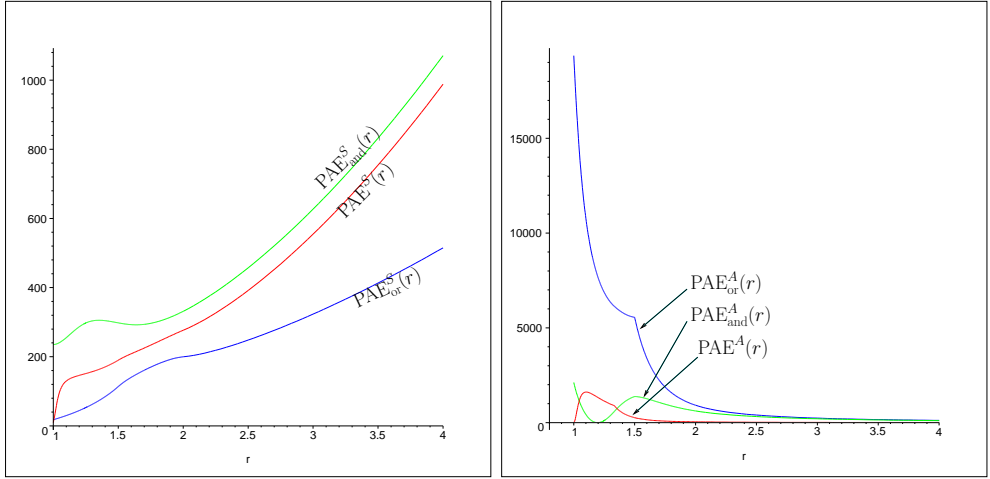
\begin{figure}[ht]
\centering
\scalebox{.3}{\input{pae_seg_und.pstex_t}}
\scalebox{.3}{\input{pae_agg_und.pstex_t}}
\caption{ \label{fig:PAECurves}
Pitman asymptotic efficiency against segregation (left) and association (right)
as a function of $r$.
Some values of note:
$\PAE^S(r=1) = 160/7$,
$\PAE^S_{\la}(r=1) = 4000/17$,
$\PAE^S_{\lo}(r=1) = 160/9$,
$\lim_{r \rightarrow \infty} \PAE^S(r)=\lim_{r \rightarrow \infty}
\PAE^S_\la (r)=\lim_{r \rightarrow \infty} \PAE^S_\lo (r) = \infty$,
and $\PAE^S_\la (r)$ has a local supremum at $\approx 1.35$.
Also
$\PAE^A(r=1) = 0$, $\PAE^A_{\la}(r=1) =\PAE^A_{\lo}(r=1)=\infty$,
$\lim_{r \rightarrow \infty} \PAE^A(r)=
\lim_{r \rightarrow \infty} \PAE^A_\la(r)=
\lim_{r \rightarrow \infty} \PAE^A_\la(r) = 0$,
$\argsup_{r \in [1,\infty]} \PAE^A(r) \approx 1.1$, and
$\PAE^A_\la(r)$ has a local supremum at $r =1.5$ and
a local infimum at $r \approx 1.2$
}
\end{figure}

Let $\PAE^S(r)$ and $\PAE^A(r)$ denote the PAE score against the segregation
and association alternatives, respectively, for the relative arc density
of the PCD based on $\NPE^r$ (see \cite{ceyhan:arc-density-PE} more detail).
Figure \ref{fig:PAECurves}
presents the PAE as a function of $r$ for both segregation and association in the digraph,
AND, and OR-underlying graph cases.
For large $n$ and small $\ve$,
PAE analysis suggests
choosing $r$ large for testing against segregation in all three cases and
choosing $r$ small for testing against association, arbitrarily close to 1
for the AND- and OR-underlying cases, but around 1.1 for the digraph case.
Furthermore, in segregation, $\PAE^S_\lo (r) < \PAE^S(r) <\PAE^S_\la (r)$,
suggesting the use of AND-underlying version.
Under association, $\max\left(\PAE^S_{\la}(r),\PAE^S(r) < \PAE^S_\lo (r)\right)$
implying the use of OR-underlying version.

\begin{remark}
\label{rem:HLAE}
{\bf Hodges-Lehmann Asymptotic Efficiency:}
Hodges-Lehmann asymptotic efficiency (HLAE) (\cite{hodges:1956}) is given by
$$
\HLAE(\rho^{\la}_n(r),\ve):=\frac{(\mu_{\la}(r,\ve)-\mu_{\la}(r))^2}{\nu_{\la}(r,\ve)}.
$$
Unlike PAE, HLAE does not involve the limit as $\ve \rightarrow 0$.
Since this requires the mean and, especially,
the asymptotic variance of $\rho^{\la}_n(r)$ {\it under the alternative},
we avoid the explicit investigation of HLAE.
HLAE for OR-underlying graphs can be defined similarly.
The ordering of HLAE seems to be the same as that of PAE. $\square$
\end{remark}

\begin{remark}
\label{rem:asymptotic-power-function}
{\bf Asymptotic Power Function Analysis:}
The asymptotic power function (\cite{kendall:1979})
allows investigation of power as a function of $r$, $n$, and $\ve$
using the asymptotic critical value and an appeal to normality.
Under a specific segregation alternative $H^S_{\ve}$,
the asymptotic power function for AND-underlying graphs is given by
\begin{eqnarray*}
\Pi_\la^S(r,n,\ve) =
1-\Phi\left(\frac{z_{\alpha}\,\sqrt{\nu_{\la}(r)}}{\sqrt{\nu_{\la}(r,\ve)}}+
\frac{\sqrt{n}\,(\mu_{\la}(r)-\mu_{\la}(r,\ve))}{\sqrt{\nu_{\la}(r,\ve)}}\right).
\end{eqnarray*}
  Under $H^A_{\ve}$,
we have
\begin{eqnarray*}
\Pi_\la^A(r,n,\ve) =
\Phi\left(\frac{z_{1-\alpha}\,\sqrt{\nu_{\la}(r)}}{\sqrt{\nu_{\la}(r,\ve)}}+
\frac{\sqrt{n}\,(\mu_{\la}(r)-\mu_{\la}(r,\ve))}{\sqrt{\nu_{\la}(r,\ve)}}\right).
\end{eqnarray*}
For OR-underlying graphs, the asymptotic power functions,
$\Pi_\lo^S(r,n,\ve)$ and $\Pi_\lo^A(r,n,\ve)$, are defined similarly.
However it is not investigated in this article. $\square$
\end{remark}

\section{Monte Carlo Simulation Analysis for Finite Sample Performance}
\label{sec:monte-carlo}
We implement the Monte Carlo simulations under the above described
null and alternatives for $r \in \{1,\,11/10,\,6/5,\,4/3,\,\sqrt{2},\\
\,3/2,\,2,\,3,\,5\}$.

\subsection{Monte Carlo Power Analysis under Segregation}
\label{sec:power-seg}
\begin{figure}[ht]
\centering
\psfrag{kernel density estimate}{\Huge{kernel density estimate}}
\psfrag{relative density}{\Huge{relative edge density}}
\rotatebox{-90}{ \resizebox{2.5 in}{!}{ \includegraphics{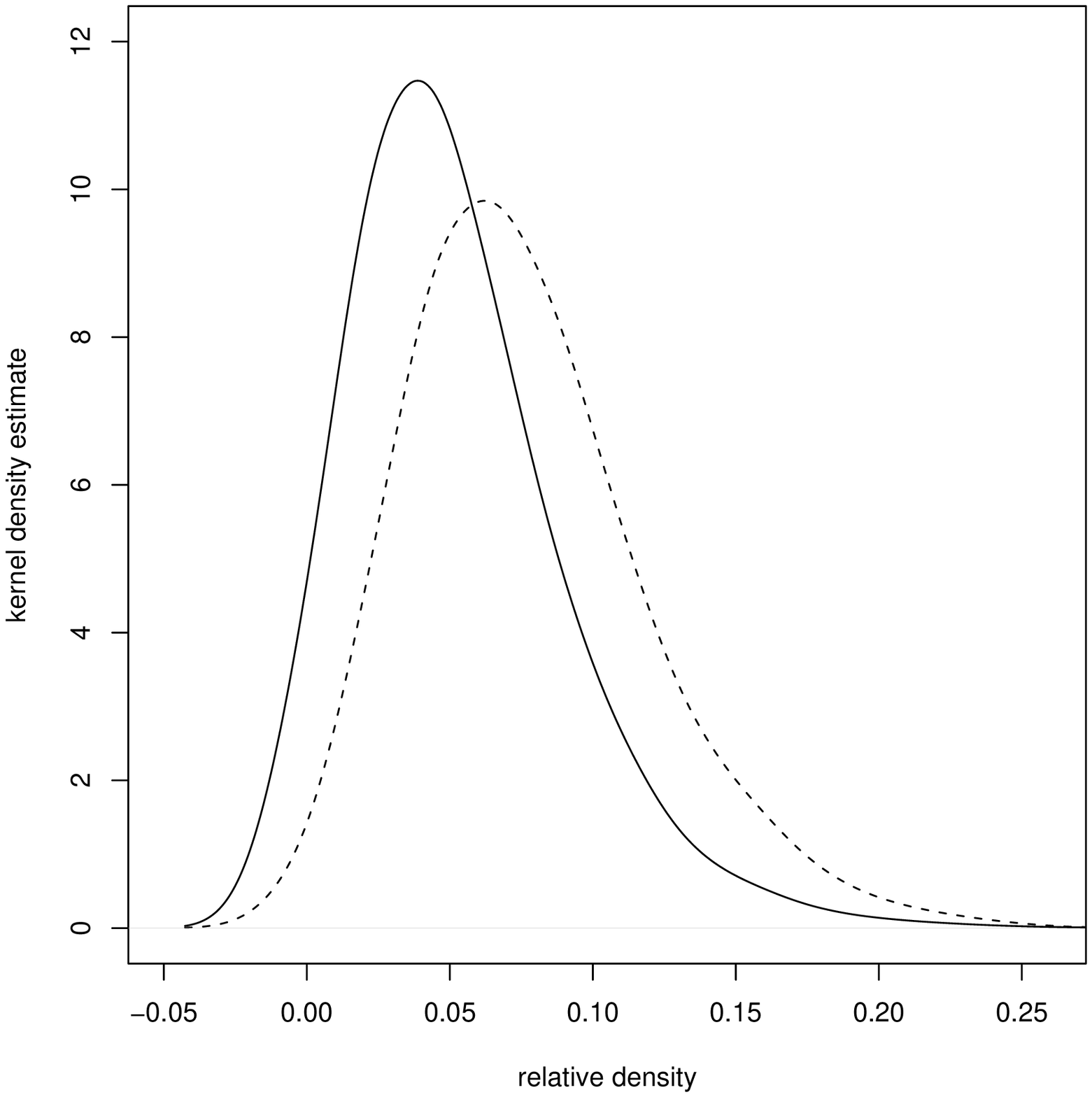} } }
\rotatebox{-90}{ \resizebox{2.5 in}{!}{ \includegraphics{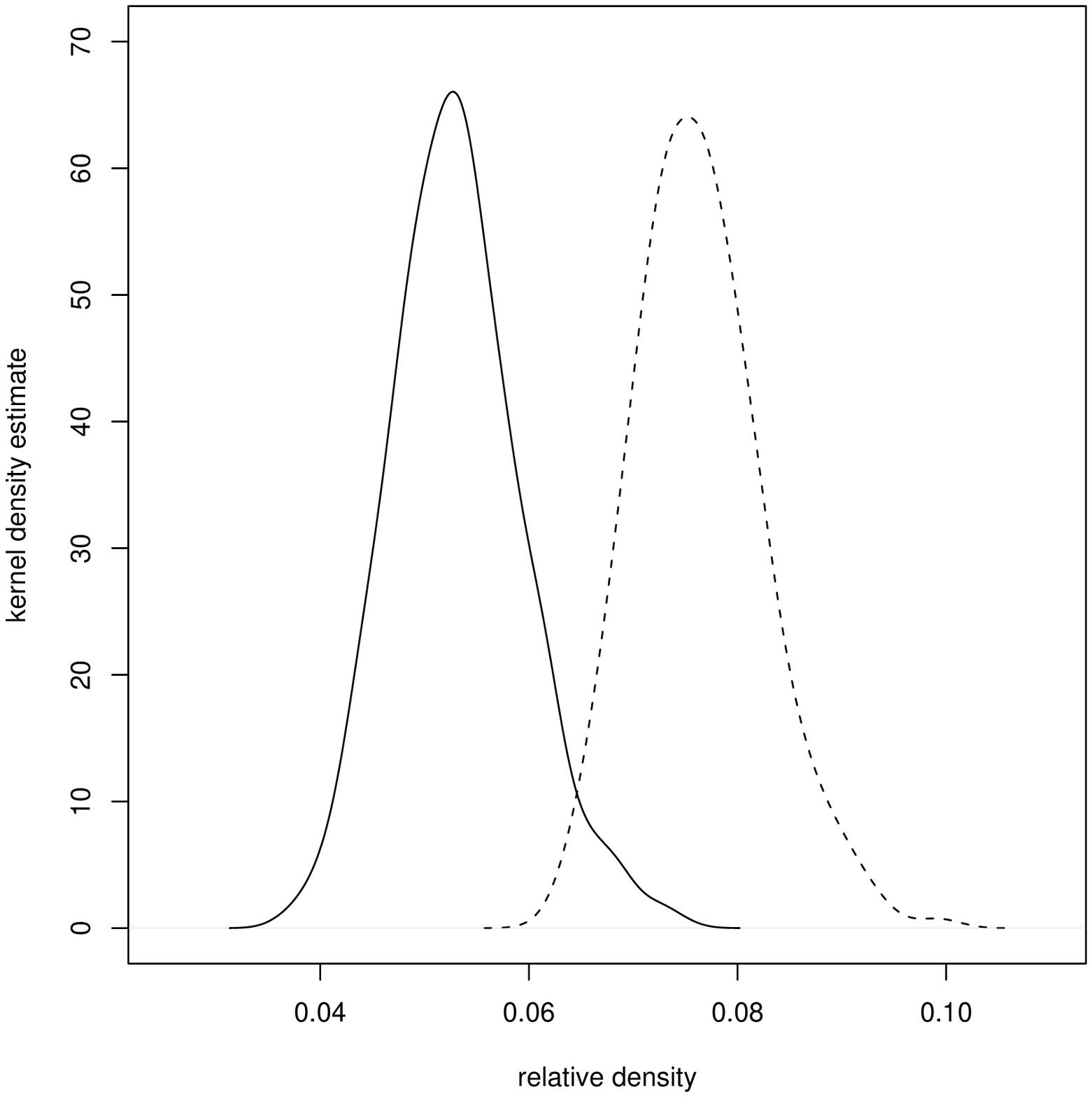} } }
\rotatebox{-90}{ \resizebox{2.5 in}{!}{ \includegraphics{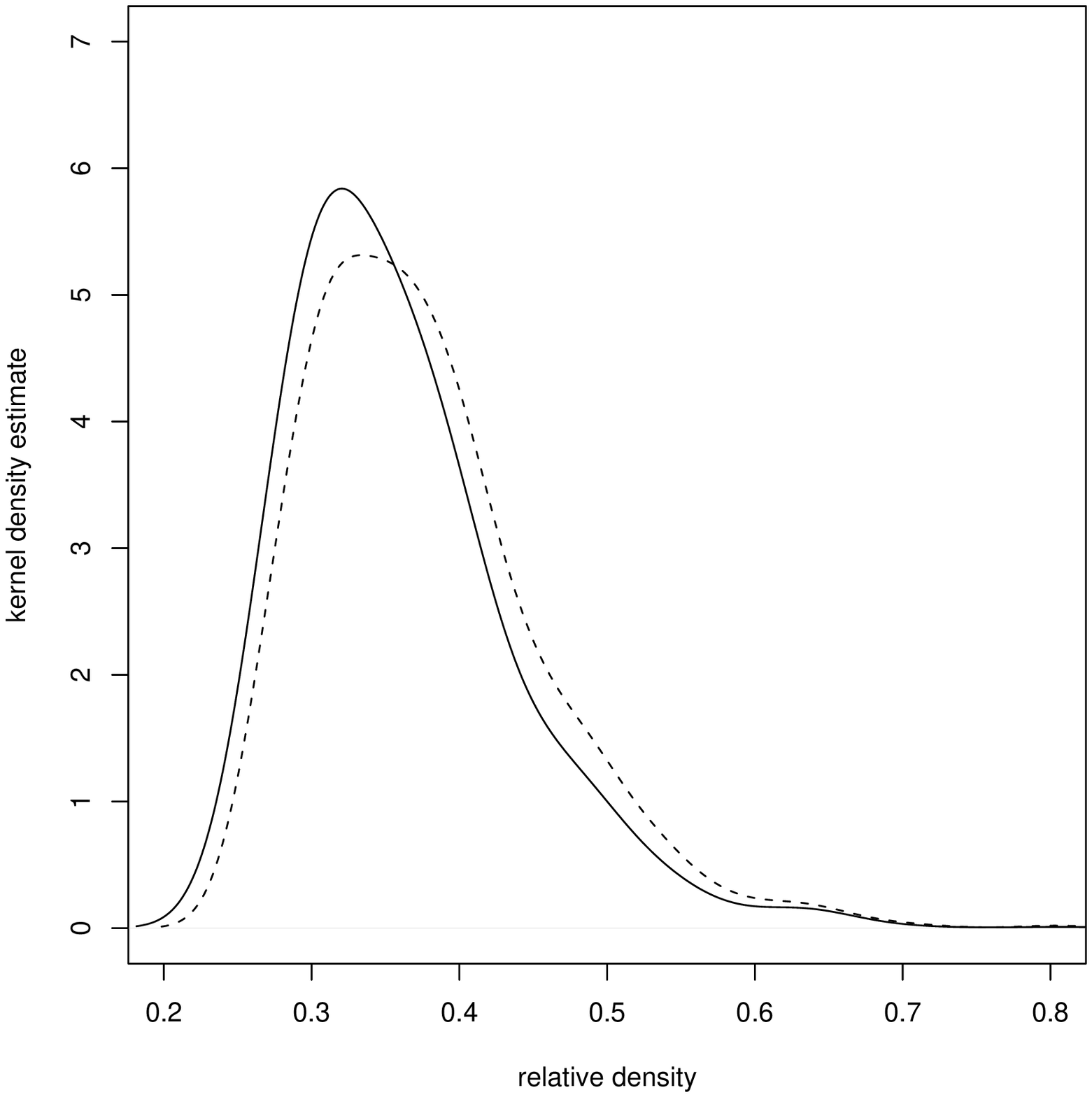} } }
\rotatebox{-90}{ \resizebox{2.5 in}{!}{ \includegraphics{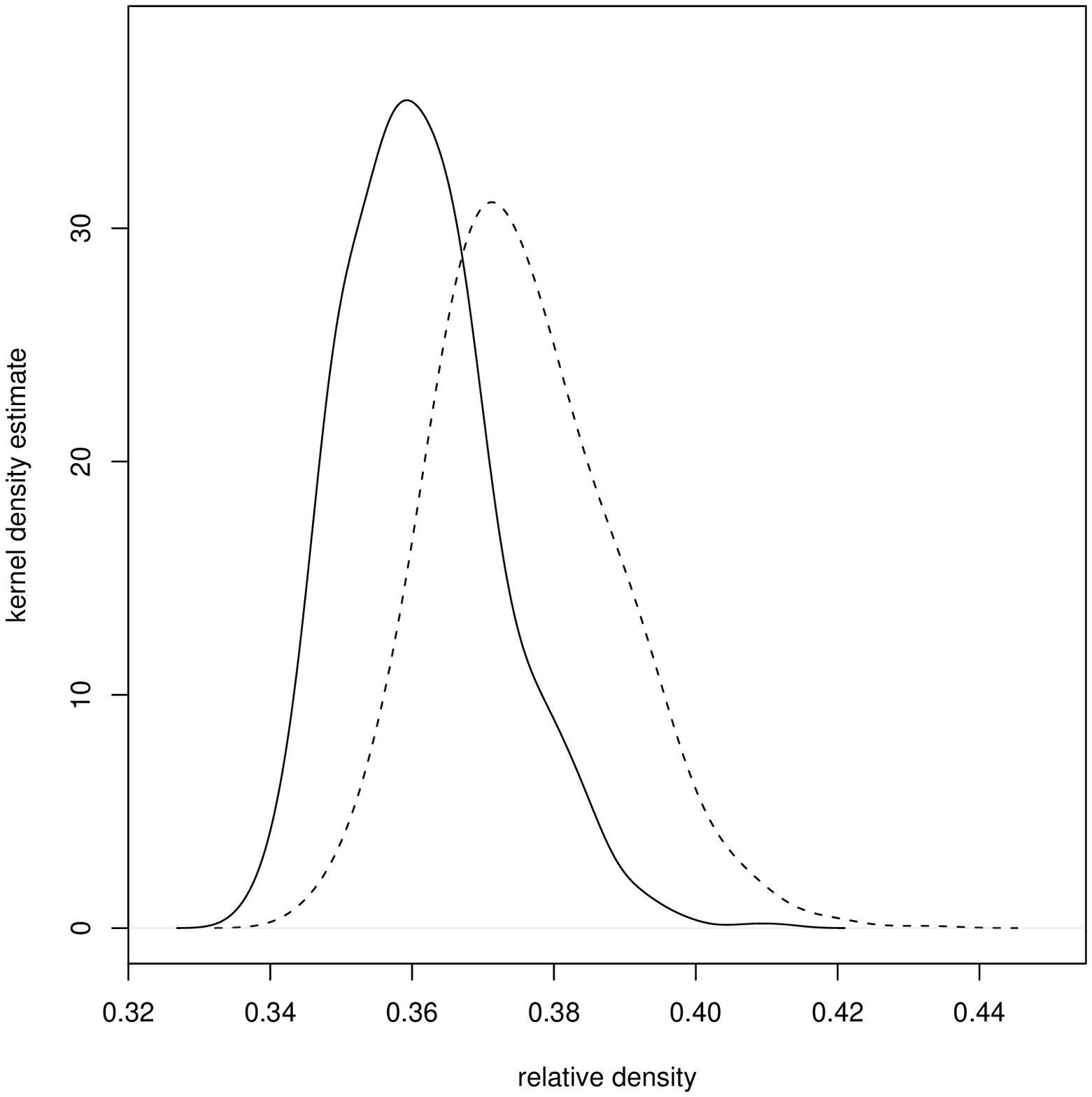} } }
\caption{ \label{fig:SegSimKernel}
Two Monte Carlo experiments against the segregation alternatives $H^S_{\sqrt{3}/8}$.
Depicted are kernel density estimates of $\rho^{\la}_n(1.1)$ for
$n=10$ (top left) and $n=100$ (top right) and $\rho^{\lo}_n(1.1)$ for
$n=10$ (bottom left) and $n=100$ (bottom right) under the null (solid) and alternative (dashed) cases.
}
\end{figure}

In Figure \ref{fig:SegSimKernel}, we present a Monte Carlo
investigation against the segregation alternative $H^S_{\sqrt{3}/8}$
for $r=1.1$ and $n=10$ (left) and $n=100$ (right).
The empirical power estimates are calculated based on the Monte Carlo critical values.
Let $\widehat{\beta}^S_{mc}\left(\rho^{\la}_n(r)\right)$ and
$\widehat{\beta}^S_{mc}\left(\rho^{\lo}_n(r)\right)$ stand
for the corresponding empirical power estimates
for the AND- and OR-underlying cases.
With $n=10$, the
null and alternative probability density functions for
$\rho^{\la}_{10}(1.1)$ and $\rho^{\lo}_{10}(1.1)$ are very similar,
implying small power (10,000 Monte Carlo replicates yield empirical power values
$\widehat{\beta}^S_{mc}\left(\rho^{\la}_{10}\right) = 0.1318$ and
$\widehat{\beta}^S_{mc}\left(\rho^{\lo}_{10}\right) = 0.0539$).
Among the 10000 Monte Carlo replicates under $H_o$,
we find the $95^{th}$ percentile value and use it as the Monte Carlo critical value
at $.05$ level for the segregation alternative,
and use $5^{th}$ percentile value for the association alternative.
With $n=100$,
there is more separation between null and alternative probability
density functions in the underlying cases where separation is much
less emphasized in the OR-underlying case;
1000 Monte Carlo replicates yield $\widehat{\beta}^S_{mc}\left(\rho^{\la}_{100}\right) = 0.994$
and $\widehat{\beta}^S_{mc}\left(\rho^{\lo}_{100}\right) = 0.298$
where the empirical power estimates are based on Monte Carlo critical values.
Notice also that the probability density functions are skewed right
for $n=10$ in both underlying cases, while approximate normality
holds for $n=100$.

\begin{figure}[ht]
\centering
\psfrag{power}{\Huge{power}}
\psfrag{r}{\Huge{$r$}}
\rotatebox{-90}{ \resizebox{2.5 in}{!}{ \includegraphics{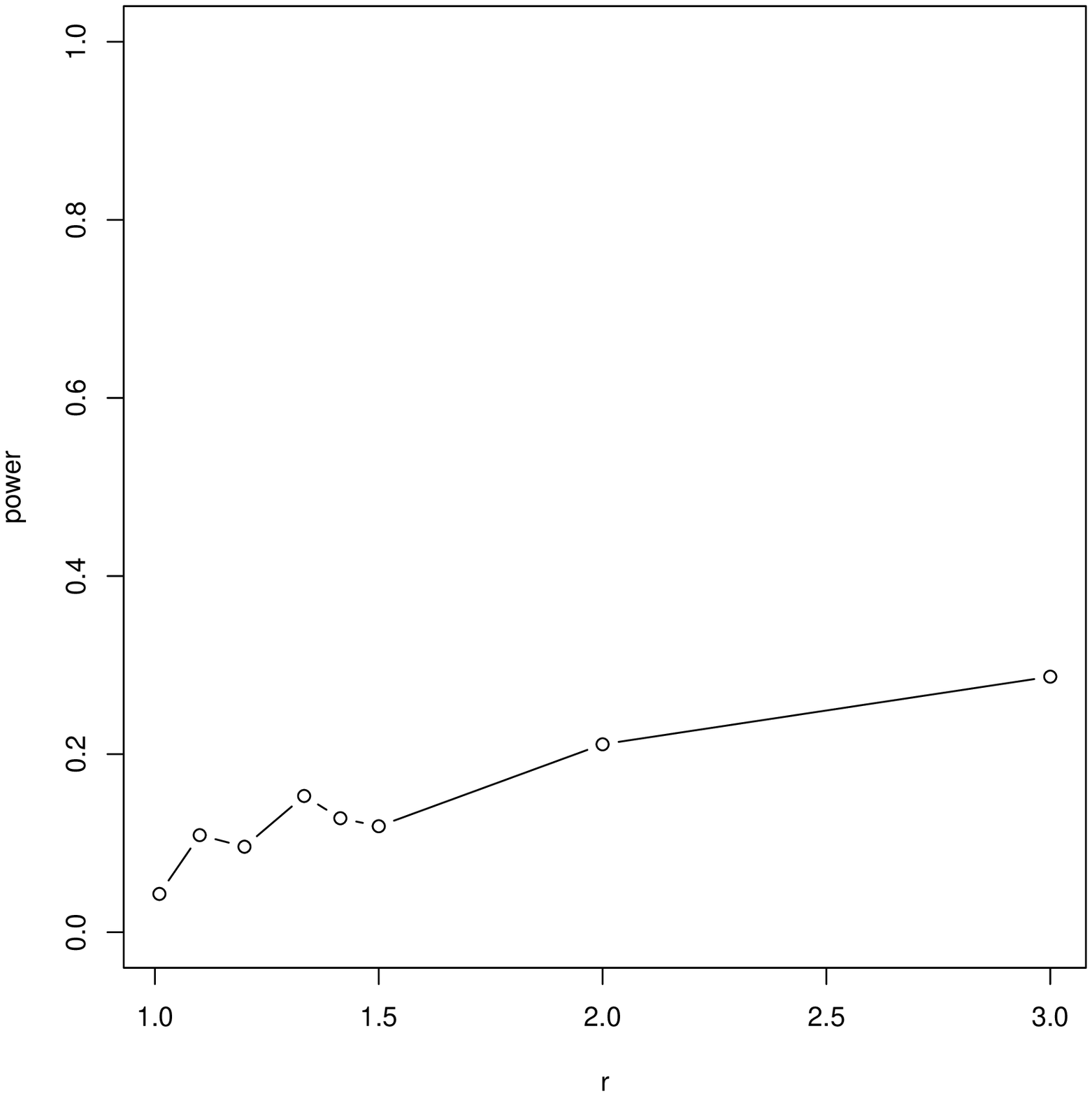} } }
\rotatebox{-90}{ \resizebox{2.5 in}{!}{ \includegraphics{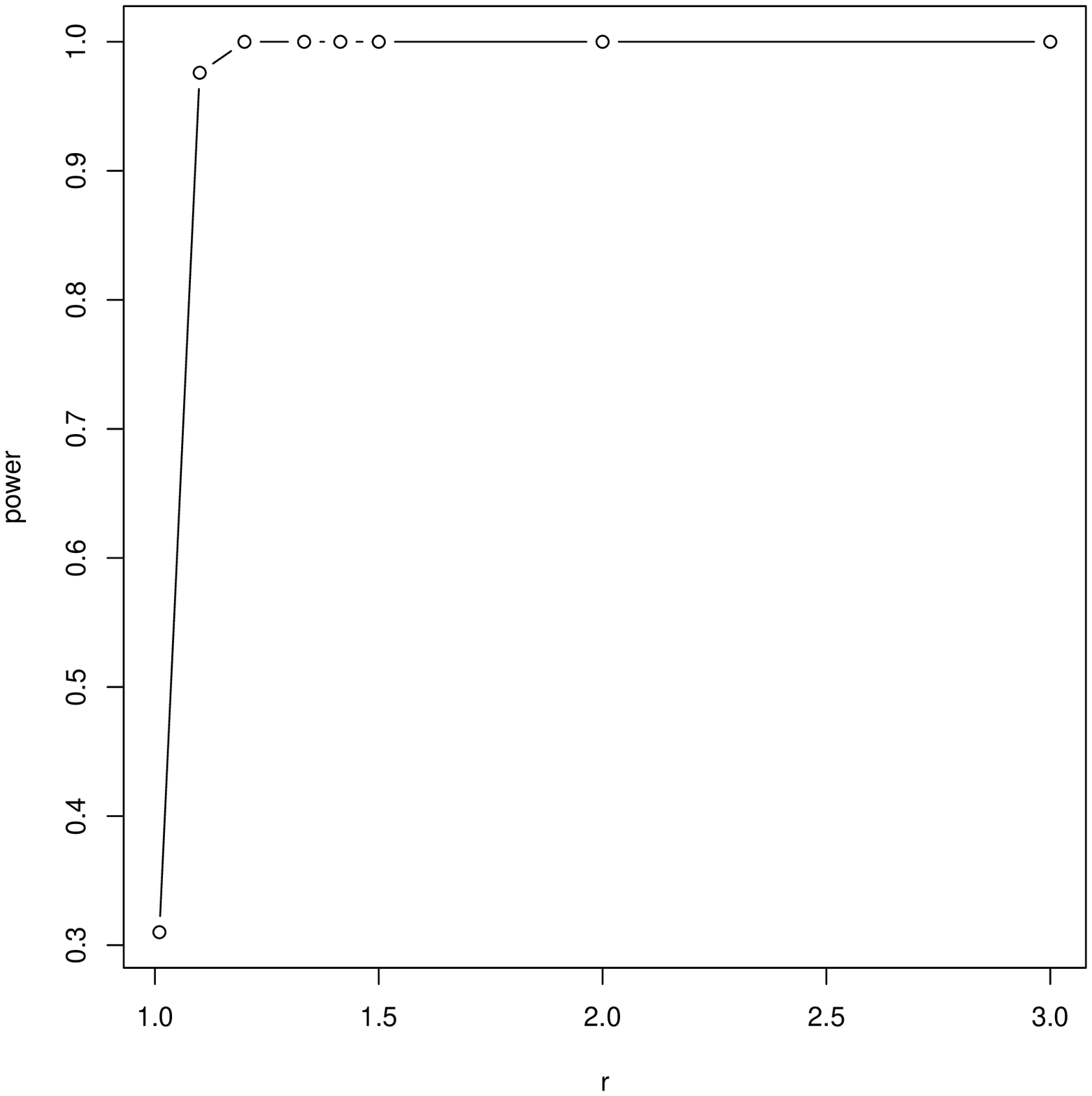} } }
\rotatebox{-90}{ \resizebox{2.5 in}{!}{ \includegraphics{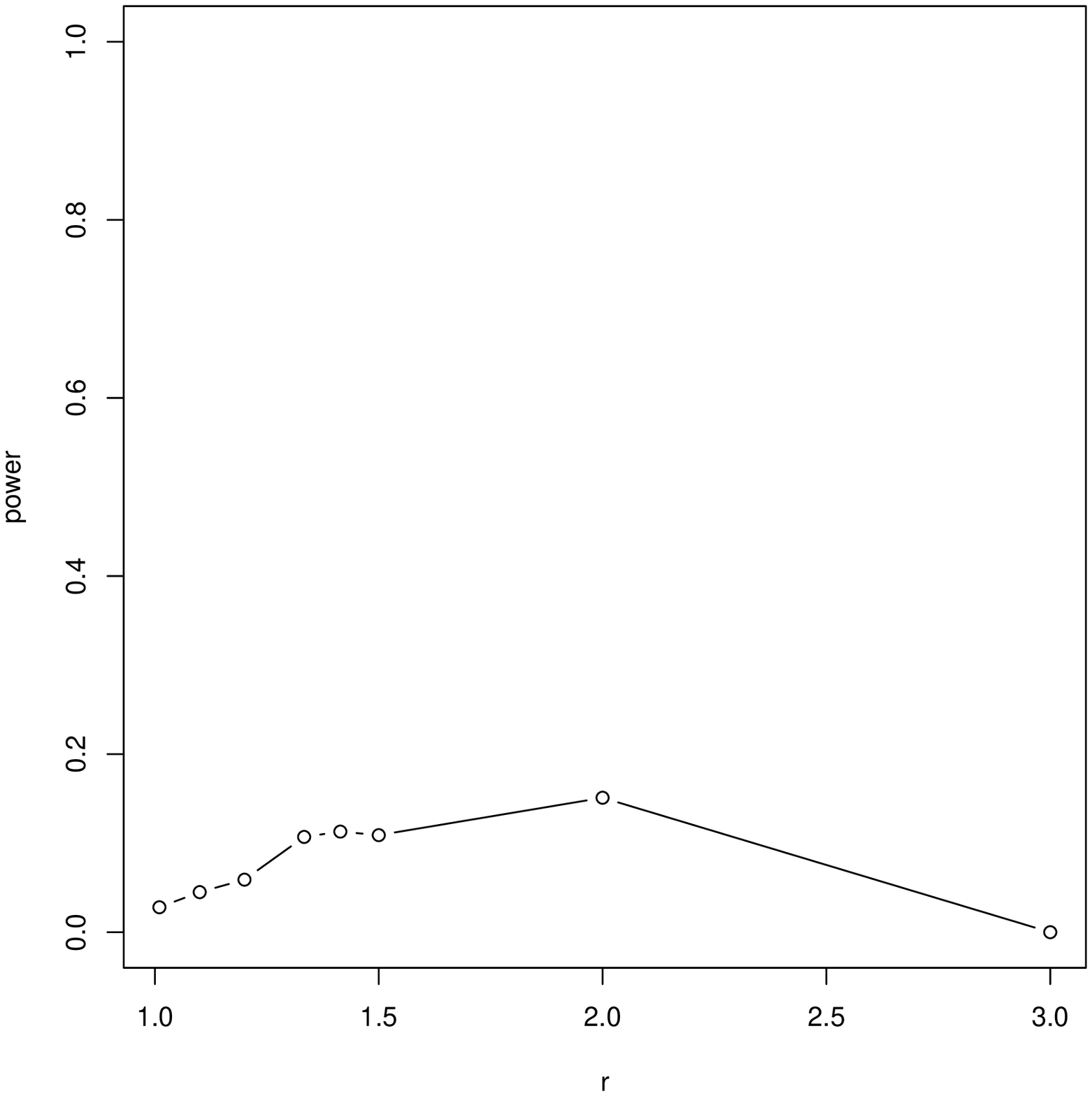} } }
\rotatebox{-90}{ \resizebox{2.5 in}{!}{ \includegraphics{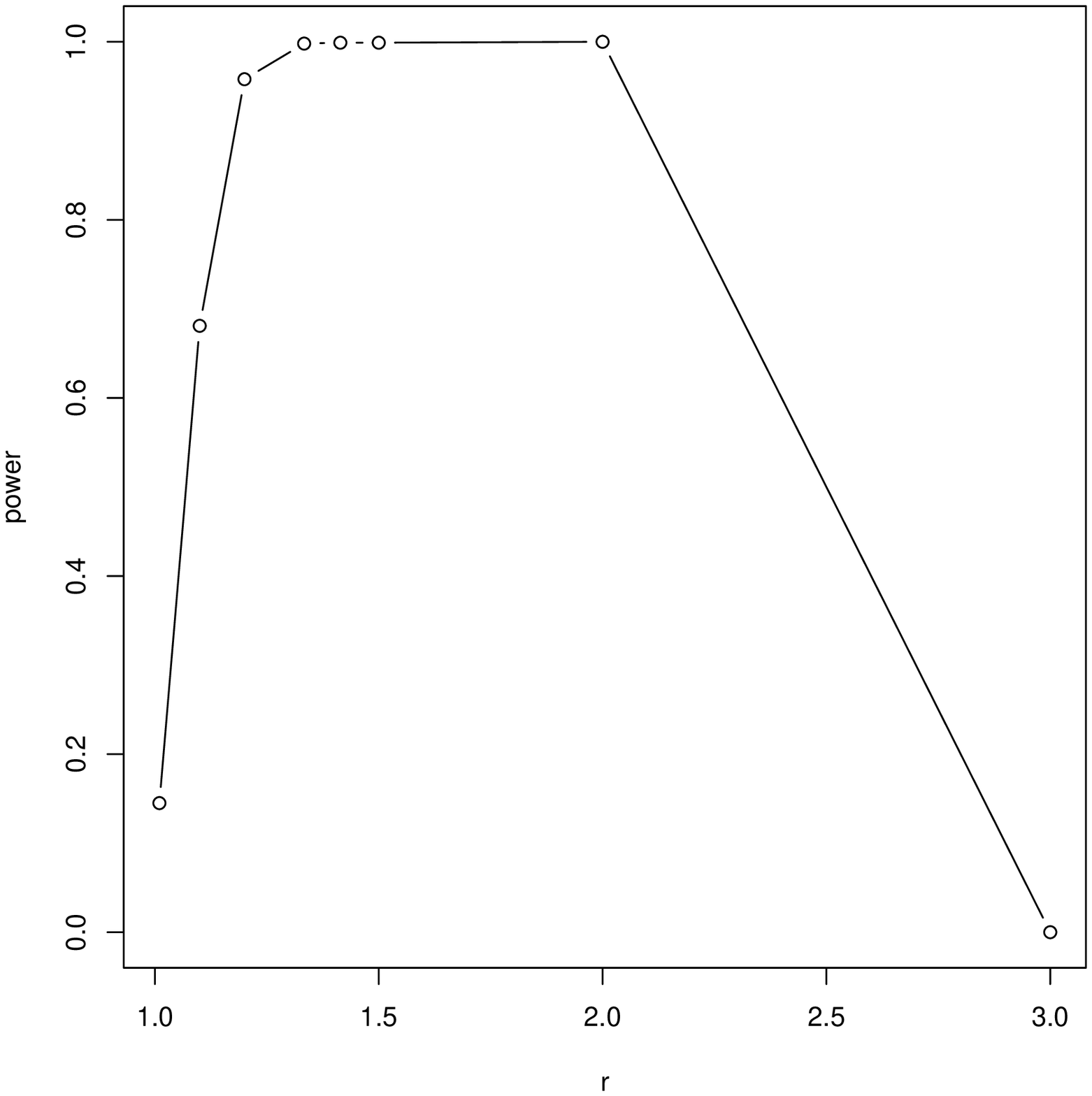} } }
\caption{
\label{fig:SegSimPower}
Empirical power estimates based on Monte Carlo critical values as a function of $r$
against segregation alternatives with the AND-underlying case (top two)
and OR-underlying case (bottom two); in both cases, we have
$H^S_{\sqrt{3}/8}$ (left)
and
$H^S_{\sqrt{3}/4}$ (right)
for $n=10$ and $N_{mc}=1000$ Monte Carlo replicates.
}
\end{figure}

For a given alternative and sample size
we may consider optimizing the empirical power of the test
as a function of the proximity factor $r$.
Figure \ref{fig:SegSimPower} presents
a Monte Carlo investigation of empirical power based on
Monte Carlo critical values against
$H^S_{\sqrt{3}/8}$
and
$H^S_{\sqrt{3}/4}$
as a function of $r$ for $n=10$ with 1000 replicates.
The corresponding empirical power estimates
are given in Table \ref{tab:emp-val-S-Under}.
Our Monte Carlo estimates of $r^*_{\ve}$,
the value of $r$ which maximizes the power against $H^S_{\ve}$,
are
$r^*_{\sqrt{3}/8} = 3$
and
$r^*_{\sqrt{3}/4} \in [4/3,3]$ in the AND-underlying case,
and
$r^*_{\sqrt{3}/8} = 2$
and
$r^*_{\sqrt{3}/4} \in [4/3,2]$ in the OR-underlying case.
That is, more severe segregation (larger $\ve$)
suggests a smaller choice of $r$ in both cases.
For both $\ve$ values, smaller $r$ values are suggested
in the OR-underlying case compared to the AND-underlying case.

\begin{table}[ht]
\centering
\begin{tabular}{|c|c|c|c|c|c|c|c|c|}
\hline
\multicolumn{9}{|c|}{$n=10$ and $N_{mc}=1000$ AND-underlying case} \\
\hline
$r$  & 1 & 11/10 & 6/5 & 4/3 & $\sqrt{2}$ & 3/2 & 2 & 3  \\
\hline
$\widehat{C}^S_{mc}$ & $0.0\bar{2}$ & $0.\bar{1}$ & .2 & $0.2\bar 8$ &  $0.3\bar 5$ & $0.4\bar 2$ & $0.7\bar{3}$ & $0.9\bar{7}$\\
\hline
$\widehat{\alpha}^S_{mc}(n)$      & 0.023 & 0.048 & 0.035 & 0.044 & 0.040 & 0.036 & 0.031 & 0.039\\
\hline
$\widehat{\beta}^S_{mc}(\sqrt{3}/8)$ & 0.043 & 0.109 & 0.096 & 0.153 & 0.128 & 0.119 & 0.211 & 0.287\\
\hline
$\widehat{\beta}^S_{mc}(\sqrt{3}/4)$ & 0.000 & 0.98 & 1 & 1 & 1 & 1 & 1 & 1\\
\hline
\multicolumn{9}{|c|}{$n=10$ and $N_{mc}=1000$ OR-underlying case} \\
\hline
$r$  & 1 & 11/10 & 6/5 & 4/3 & $\sqrt{2}$ & 3/2 & 2 & 3  \\
\hline
$\widehat{C}^S_{mc}$ & $0.4\bar 8$ & $0.4\bar 8$ & $0.5\bar 3$ & $0.6\bar 2$ & $0.6\bar 8$ & $0.7\bar 3$ & $0.9\bar 5$ & 1.00\\
\hline
$\widehat{\alpha}^S_{mc}(n)$      & 0.030 & 0.045 & 0.049 & 0.043 & 0.037 & 0.043 & 0.034 & 0.000\\
\hline
$\widehat{\beta}^S_{mc}(\sqrt{3}/8)$ & 0.028 & 0.045 & 0.059 & 0.107 & 0.113 & 0.109 & 0.151 & 0.000\\
\hline
$\widehat{\beta}^S_{mc}(\sqrt{3}/4)$ & 0.145 & 0.681 & 0.958 & 0.998 & 0.999 & 0.999 & 1.000 & 0.000\\
\hline
\end{tabular}
\caption{ \label{tab:emp-val-S-Under}
The Monte Carlo critical values, $\widehat{C}^S_{mc}$,
empirical significance levels, $\widehat{\alpha}^S_{mc}(n)$,
and empirical power estimates, $\widehat{\beta}^S_{mc}$, based on the Monte Carlo critical values under
$H^S_{\sqrt{3}/8}$ and $H^S_{\sqrt{3}/4}$, $N_{mc}=1000$, and $n=10$ at
$\alpha=.05$.}
\end{table}

\begin{figure}[ht]
\centering
\psfrag{power}{\Huge{power}}
\psfrag{r}{\Huge{$r$}}
\rotatebox{-90}{ \resizebox{2.5 in}{!}{ \includegraphics{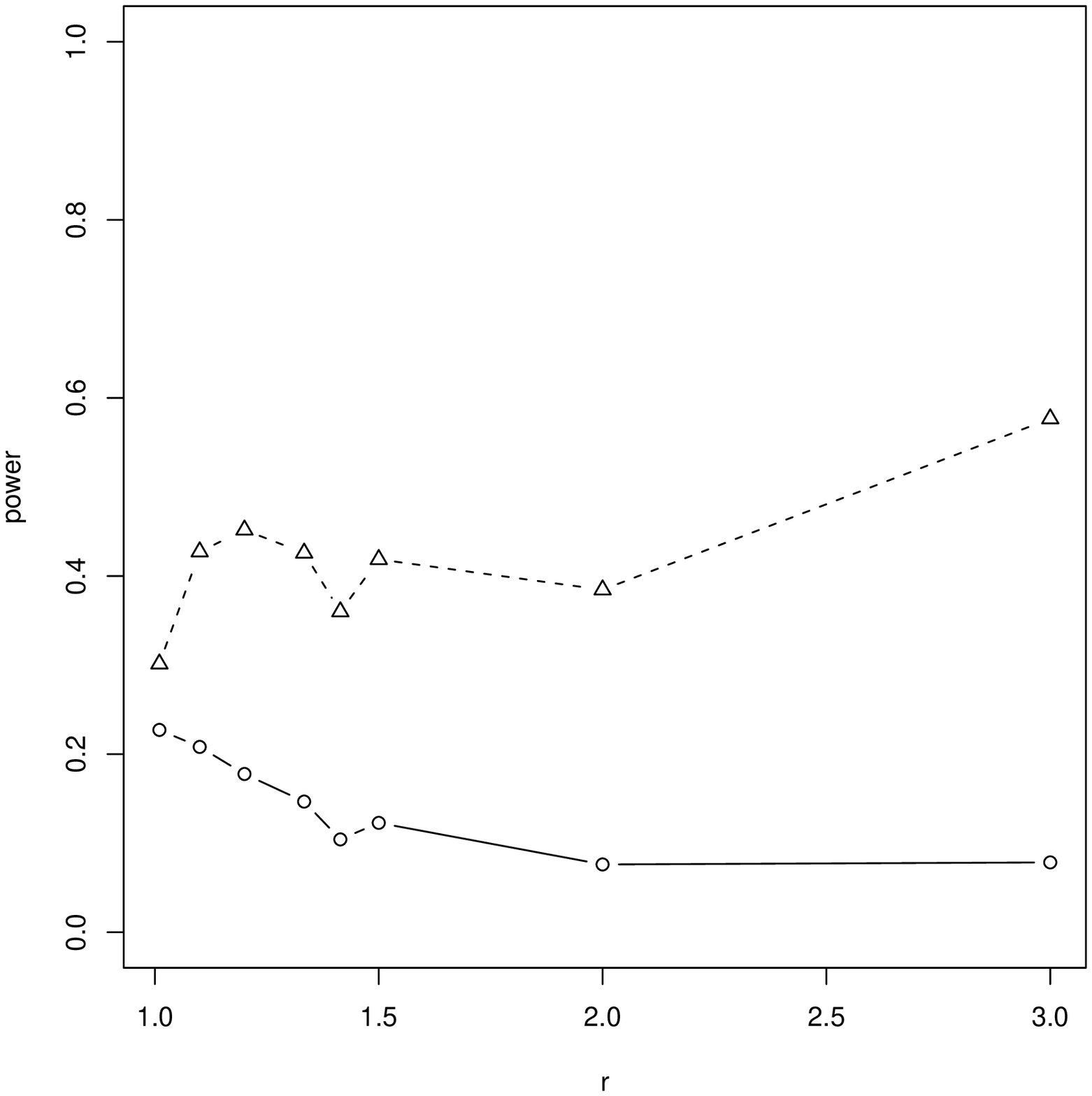} } }
\rotatebox{-90}{ \resizebox{2.5 in}{!}{ \includegraphics{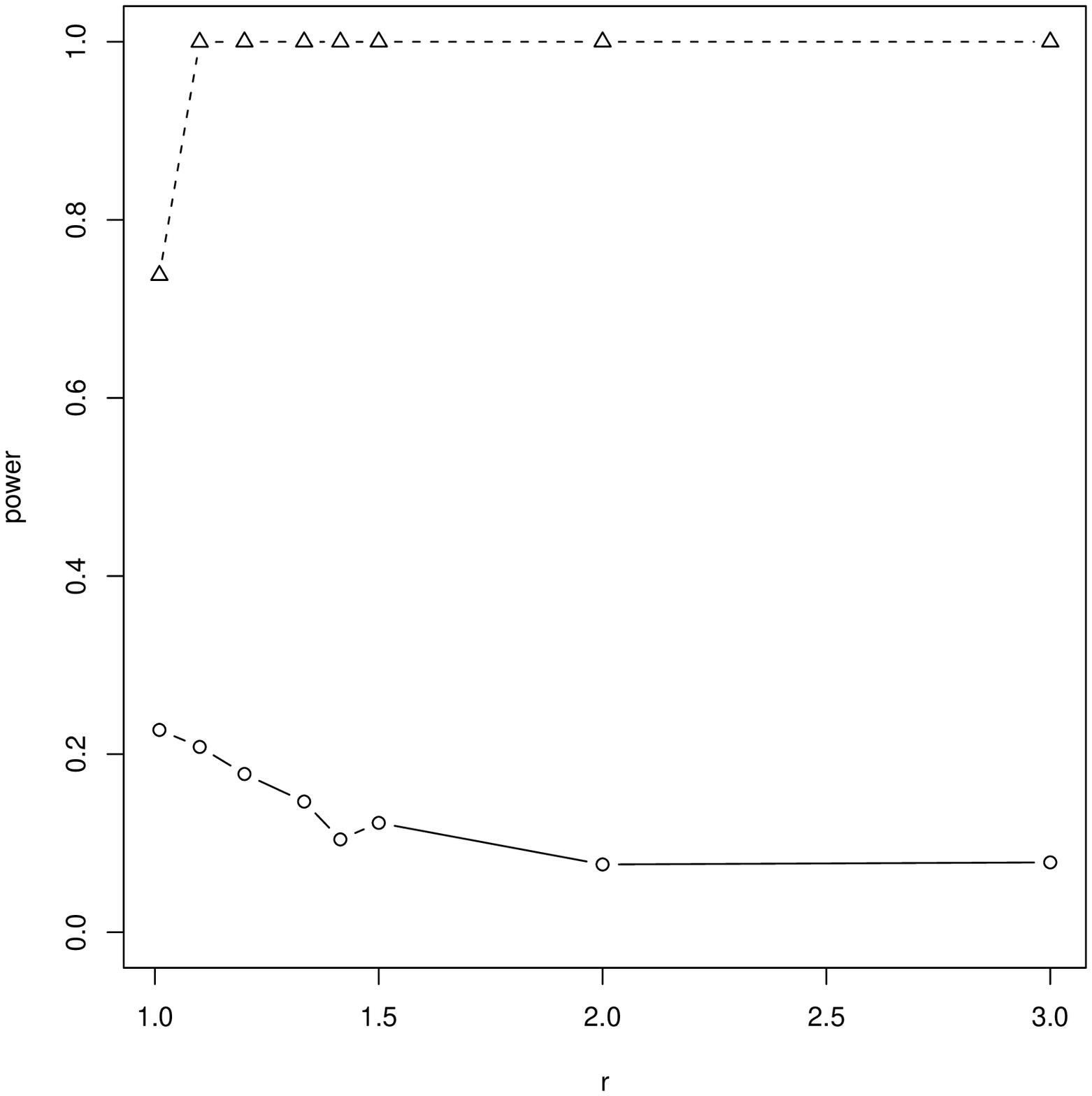} } }
\rotatebox{-90}{ \resizebox{2.5 in}{!}{ \includegraphics{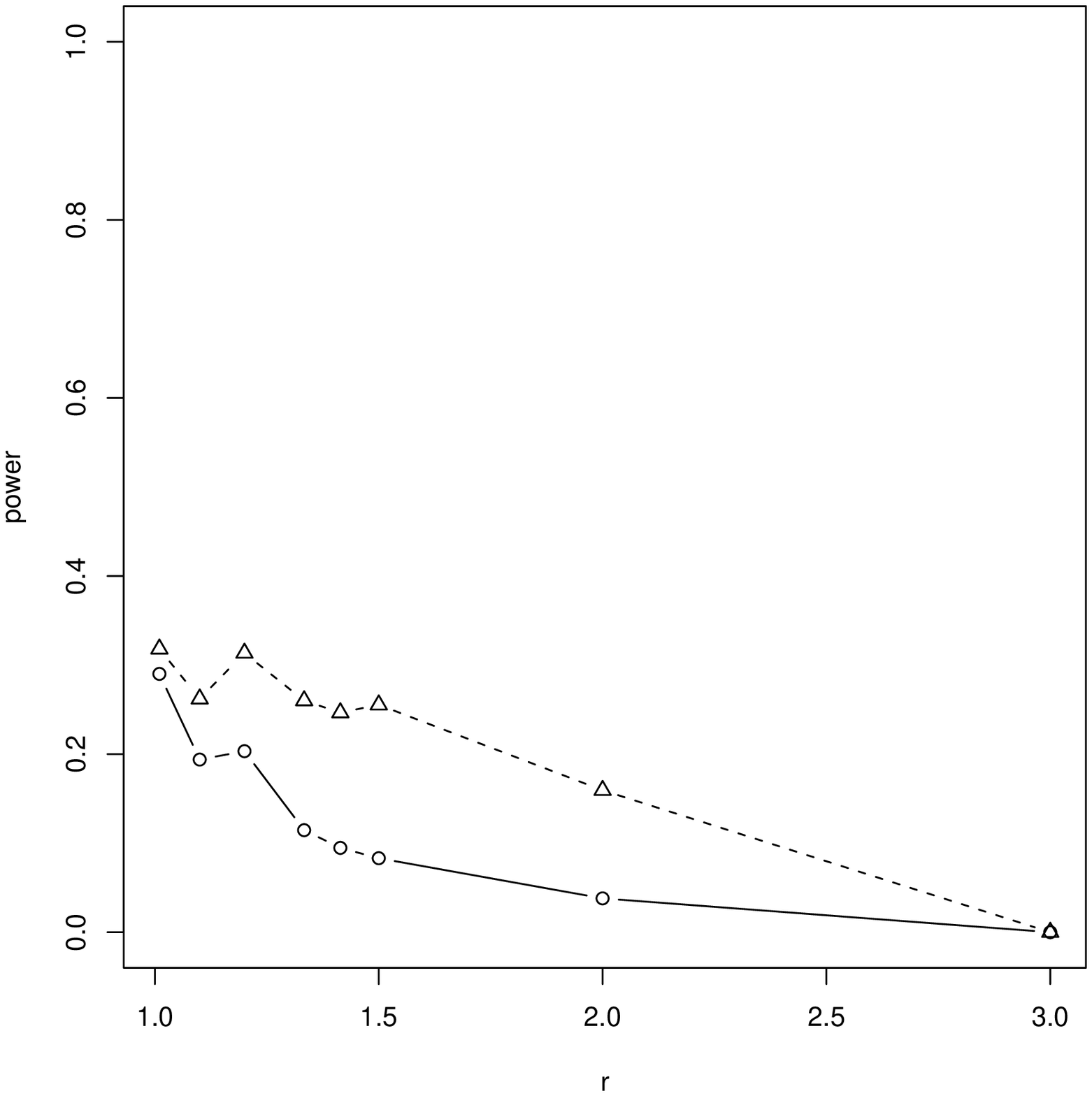} } }
\rotatebox{-90}{ \resizebox{2.5 in}{!}{ \includegraphics{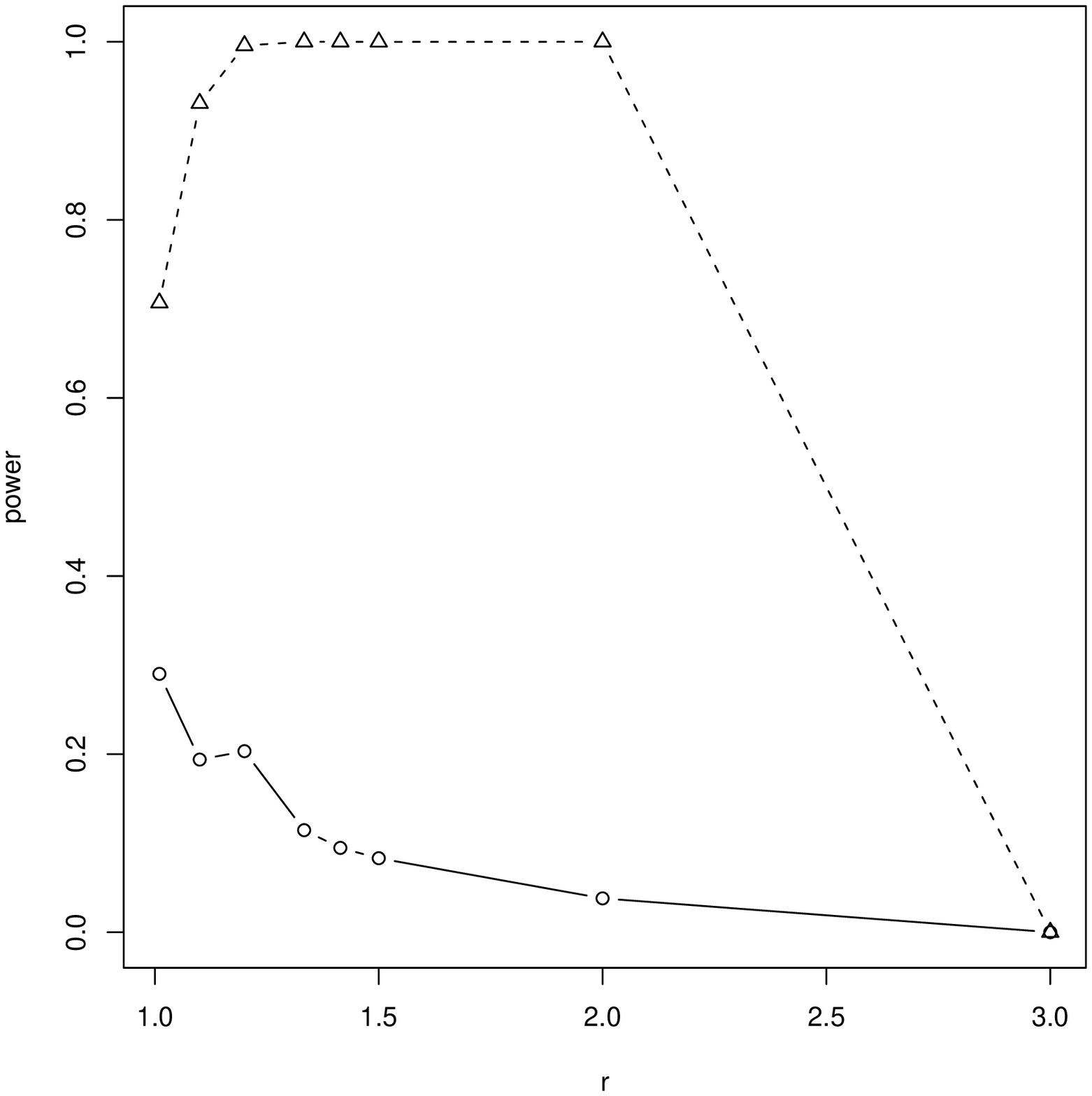} } }
\caption{ \label{fig:SegSimPowerCurve}
The empirical size (circles joined with solid lines) and power estimates (triangles with dotted lines)
based on the asymptotic critical value against segregation
alternatives in the AND-underlying case (top two) and the OR-underlying case (bottom two);
in both cases,
$H^S_{\sqrt{3}/8}$ (left)
and
$H^S_{\sqrt{3}/4}$ (right)
as a function of $r$, for $n=10$ and $N_{mc}=10000$.
}
\end{figure}

For a given alternative and sample size we may consider analyzing
the power of the test --- using the asymptotic critical value--- as
a function of the proximity factor $r$.
Let $\widehat{\alpha}_n(r)$ denote the empirical significance levels
and $\widehat{\beta}_n(r)$ empirical power estimates
based on the asymptotic critical value.
Figure \ref{fig:SegSimPowerCurve} presents a Monte Carlo investigation of
empirical power based on asymptotic critical value
against $H^S_{\sqrt{3}/8}$ and $H^S_{\sqrt{3}/4}$ as a
function of $r$ for $n=10$.
The corresponding empirical power estimates
are given in Table \ref{tab:asy-emp-val-S-Under}.
In the AND-underlying case, the empirical significance level, $\widehat{\alpha}_{n=10}(r)$, is
closest to $.05$ for $r=2$ and $3$ which have the empirical power
$\widehat{\beta}_{10}(2) = .3846$ and $\widehat{\beta}_{10}(3) = .5767$
for $\ve=\sqrt{3}/8$, and $\widehat{\beta}_{10}(2)
=\widehat{\beta}_{10}(3) =1$ for $\ve=\sqrt{3}/4$.
In the OR-underlying case, the empirical significance level,
$\widehat{\alpha}_{n=10}(r)$, is closest to $.05$ for $r=2$ ---larger for
all $r$ values --- which have the empirical power
$\widehat{\beta}_{10}(2) = .1594$ for $\ve=\sqrt{3}/8$, and
$\widehat{\beta}_{10}(2)=1$ for $\ve=\sqrt{3}/4$.
So, for small sample sizes, moderate values of $r$ is more appropriate for
normal approximation, as they yield the desired significance level,
and the more severe the segregation, higher the power estimate.
Furthermore, the AND-underlying version seems to perform better than
the OR-underlying version for segregation alternatives.

\begin{table}[ht]
\centering
\begin{tabular}{|c|c|c|c|c|c|c|c|c|}
\hline
\multicolumn{9}{|c|}{$n=10$ and $N_{mc}=10000$ AND-underlying case} \\
\hline
$r$  & 1 & 11/10 &6/5 & 4/3 & $\sqrt{2}$ & 3/2 & 2 & 3 \\
\hline
$\widehat{\alpha}_S(n)$ & 0.2272 & 0.2081 & 0.1777 & 0.1467 & 0.1042 & 0.1228 & 0.0761 & 0.0784 \\
\hline
$\widehat{\beta}^S_{n}(r,\sqrt{3}/8)$ & 0.3014 & 0.4273 & 0.4518 & 0.4259 & 0.3600 & 0.4187 & 0.3846 & 0.5767\\
\hline
$\widehat{\beta}^S_{n}(r,\sqrt{3}/4)$ & 0.6519 & 0.9985 & 1.0000 & 1.0000 & 1.0000 & 1.0000 & 1.0000 & 1.0000  \\
\hline
\multicolumn{9}{|c|}{$n=10$ and $N_{mc}=10000$ OR-underlying case} \\
\hline
$r$  & 1 & 11/10 &6/5 & 4/3 & $\sqrt{2}$ &3/2 & 2 & 3 \\
\hline
$\widehat{\alpha}_S(n)$ & 0.2901 & 0.1939 & 0.2033 & 0.1146 & 0.0947 & 0.0831 & 0.0380 & 0.0000 \\
\hline
$\widehat{\beta}^S_{n}(r,\sqrt{3}/8)$ & 0.3182 & 0.2621 & 0.3135 & 0.2601 & 0.2466 & 0.2554 & 0.1594 & 0.0000\\
\hline
$\widehat{\beta}^S_{n}(r,\sqrt{3}/4)$ & 0.7069 & 0.9310 & 0.9958 & 1.0000 & 1.0000 & 0.9999 & 1.0000 & 0.0000 \\
\hline
\end{tabular}
\caption{ \label{tab:asy-emp-val-S-Under}
The empirical significance levels, $\widehat{\alpha}_S(n)$, and
empirical power values, $\widehat{\beta}^S_{n}(r,\ve)$,
based on asymptotic critical values under $H^S_{\ve}$ for
$\ve=\sqrt{3}/8,\,\sqrt{3}/4$, $N_{mc}=10000$, and $n=10$ at
$\alpha=.05$.}
\end{table}

\subsection{Monte Carlo Power Analysis under Association}
\label{sec:power-assoc}

\begin{figure}[ht]
\centering
\psfrag{kernel density estimate}{\Huge{kernel density estimate}}
\psfrag{relative density}{\Huge{relative edge density}}
\rotatebox{-90}{ \resizebox{2.5 in}{!}{ \includegraphics{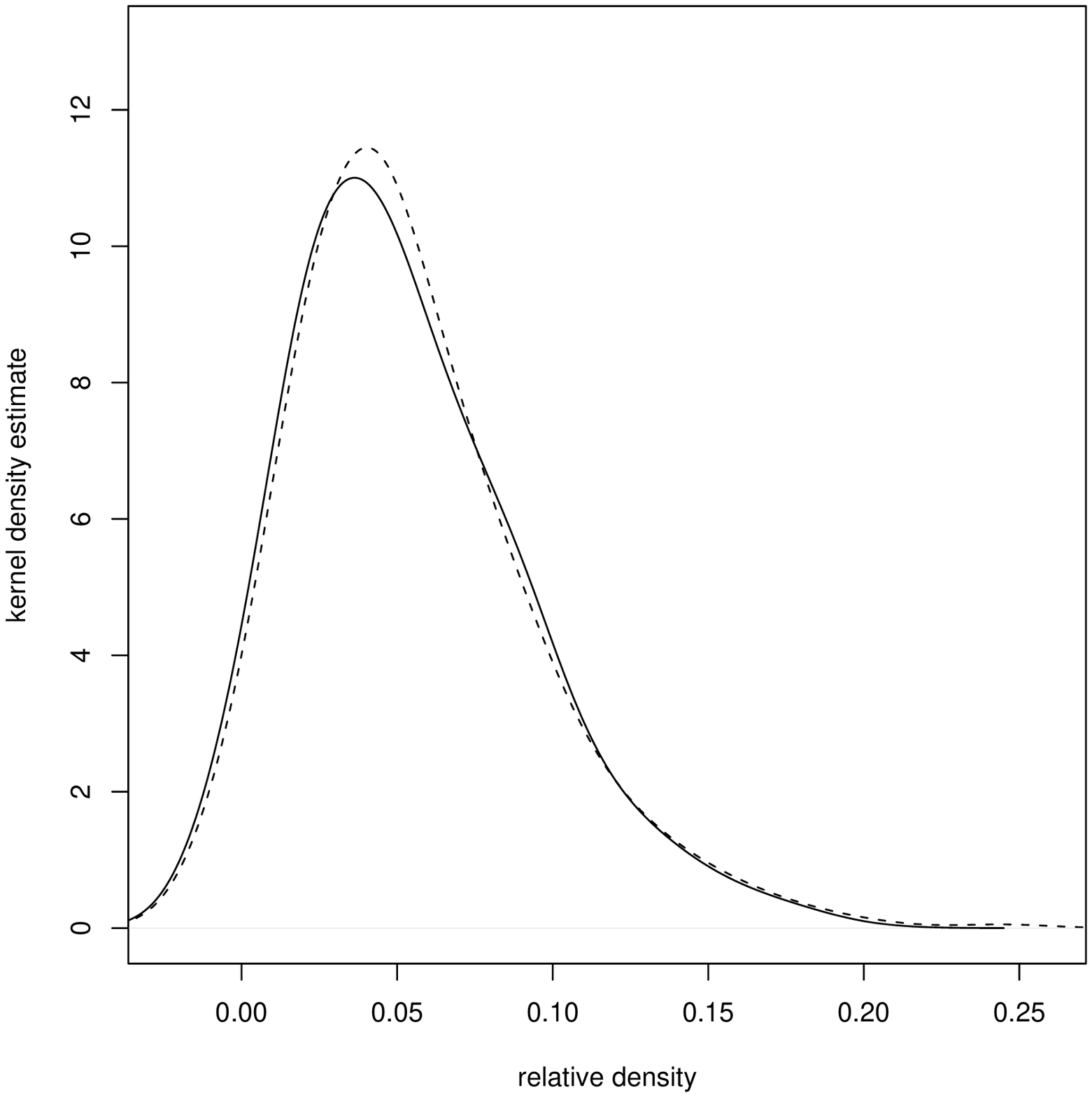} } }
\rotatebox{-90}{ \resizebox{2.5 in}{!}{ \includegraphics{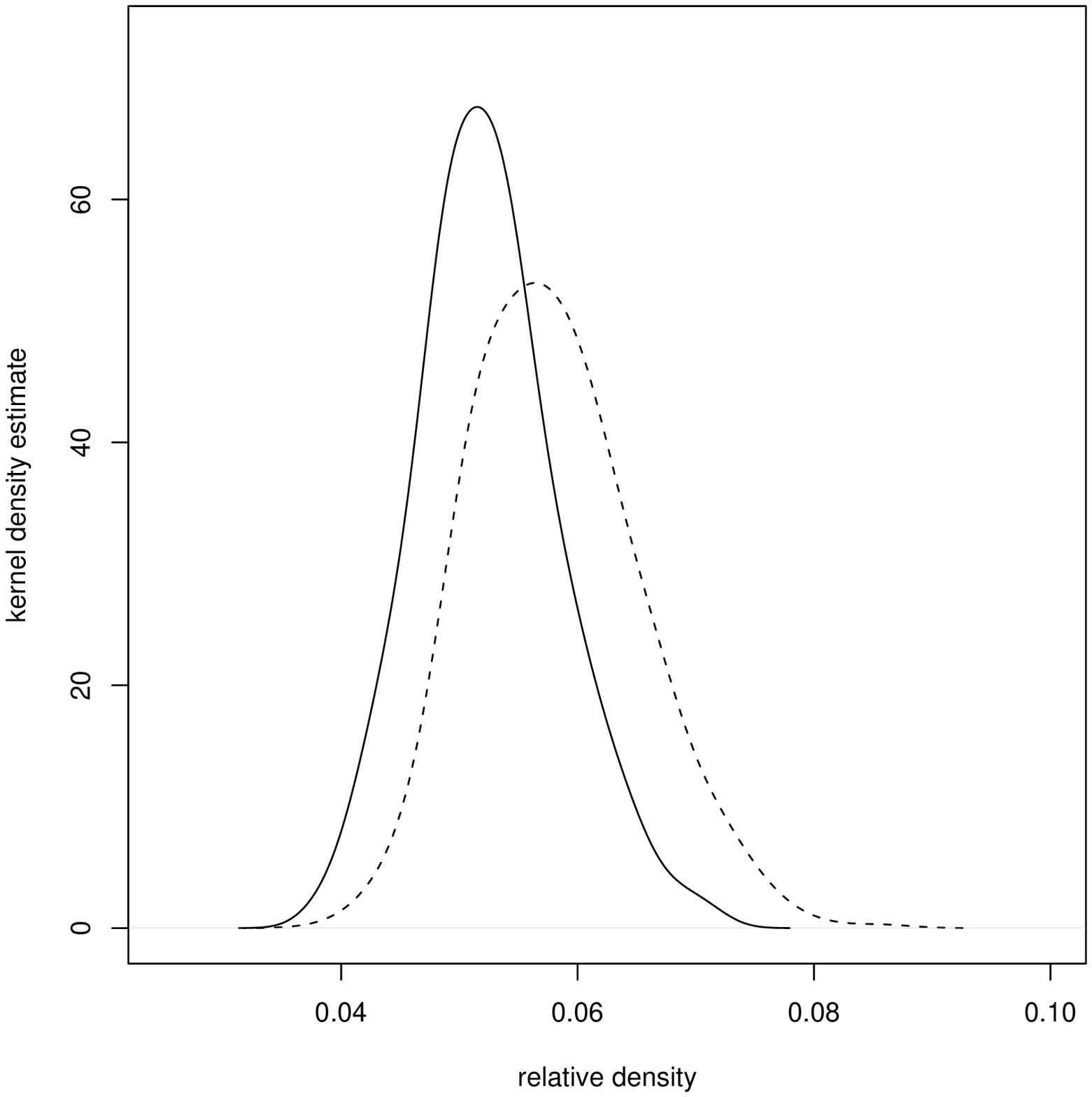} } }
\rotatebox{-90}{ \resizebox{2.5 in}{!}{ \includegraphics{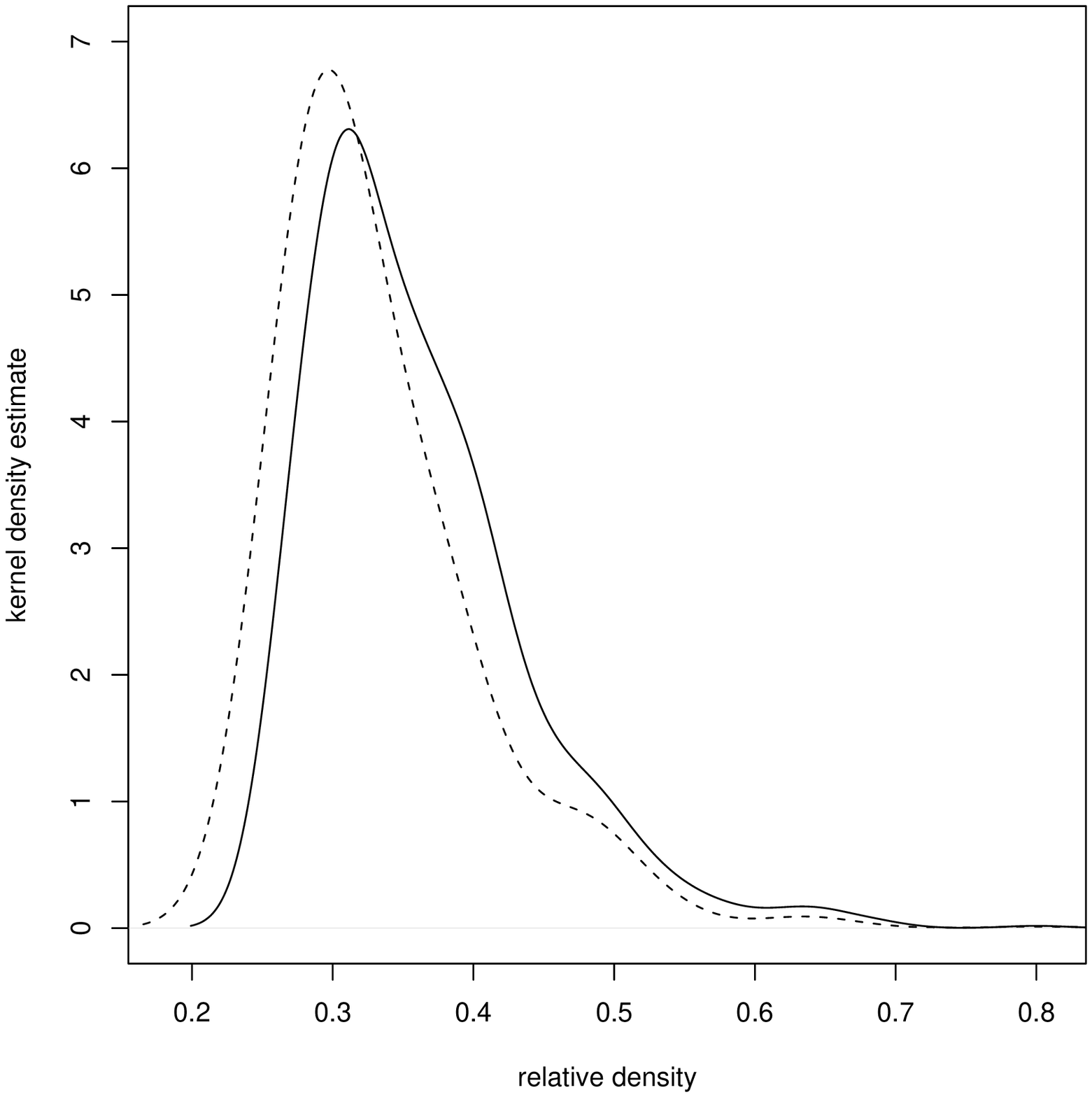} } }
\rotatebox{-90}{ \resizebox{2.5 in}{!}{ \includegraphics{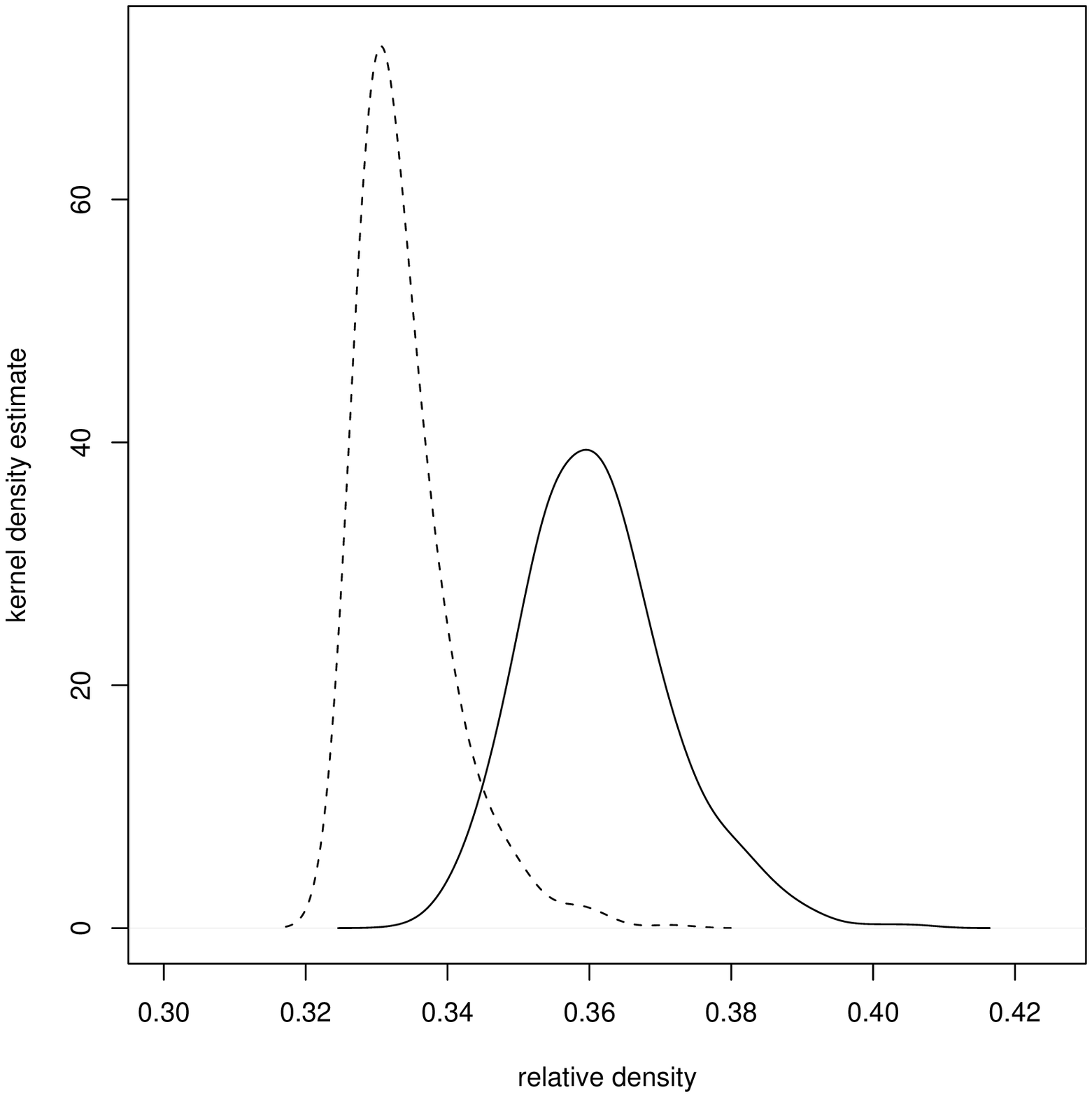} } }
\caption{ \label{fig:AggSimKernel}
Two Monte Carlo experiments against the association alternative $H^A_{\sqrt{3}/12}$.
Depicted are kernel density estimates of $\rho^{\la}_n(1.1)$ for
$n=10$ (top left) and $n=100$ (top right) and $\rho^{\lo}_n(1.1)$ for
$n=10$ (bottom left) and $n=100$ (bottom right) under the null (solid) and alternative (dashed).
}
\end{figure}

In Figure \ref{fig:AggSimKernel}, we present a Monte Carlo
investigation against the association alternative
$H^A_{\sqrt{3}/12}$ for $r=1.1$ and $n=10$ (left) and $n=100$ (right).
The empirical power estimates are calculated based on the Monte Carlo critical values
Let $\widehat{\beta}^A_{mc}\left(\rho^{\la}_n(r)\right)$ and
$\widehat{\beta}^A_{mc}\left(\rho^{\lo}_n(r)\right)$ stand
for the corresponding empirical power estimates
for the AND- and OR-underlying cases.
As above, with $n=10$, the null and alternative probability
density functions for $\rho^{\la}_{10}(1.1)$ and
$\rho^{\lo}_{10}(1.1)$ are very similar, implying small power---
in fact, virtually no power---
(10,000 Monte Carlo replicates yield the following empirical power estimates
based on Monte Carlo critical values:
$\widehat{\beta}_{mc}^A\left(\rho^{\la}_{10}\right) = 0.0$ and
$\widehat{\beta}_{mc}^A\left(\rho^{\lo}_{10}\right) = 0.0$).
With $n=100$,
there is more separation between null and alternative probability
density functions in the underlying cases where separation is much
less emphasized in the AND-underlying case; for this case,
1000 Monte Carlo replicates yield the following empirical power estimates
based on Monte Carlo critical values:
$\widehat{\beta}_{mc}^A\left(\rho^{\la}_{100}\right) = 0.0\overline{09}$ and
$\widehat{\beta}_{mc}^A\left(\rho^{\lo}_{100}\right) = 0.939$.
Notice also that the probability density functions are skewed right for $n=10$ in
both underlying cases, with more skewness in OR-underlying  case,
while approximate normality holds for $n=100$ for both cases.

\begin{figure}[ht]
\centering
\psfrag{power}{\Huge{power}}
\psfrag{r}{\Huge{$r$}}
\rotatebox{-90}{ \resizebox{2.5 in}{!}{ \includegraphics{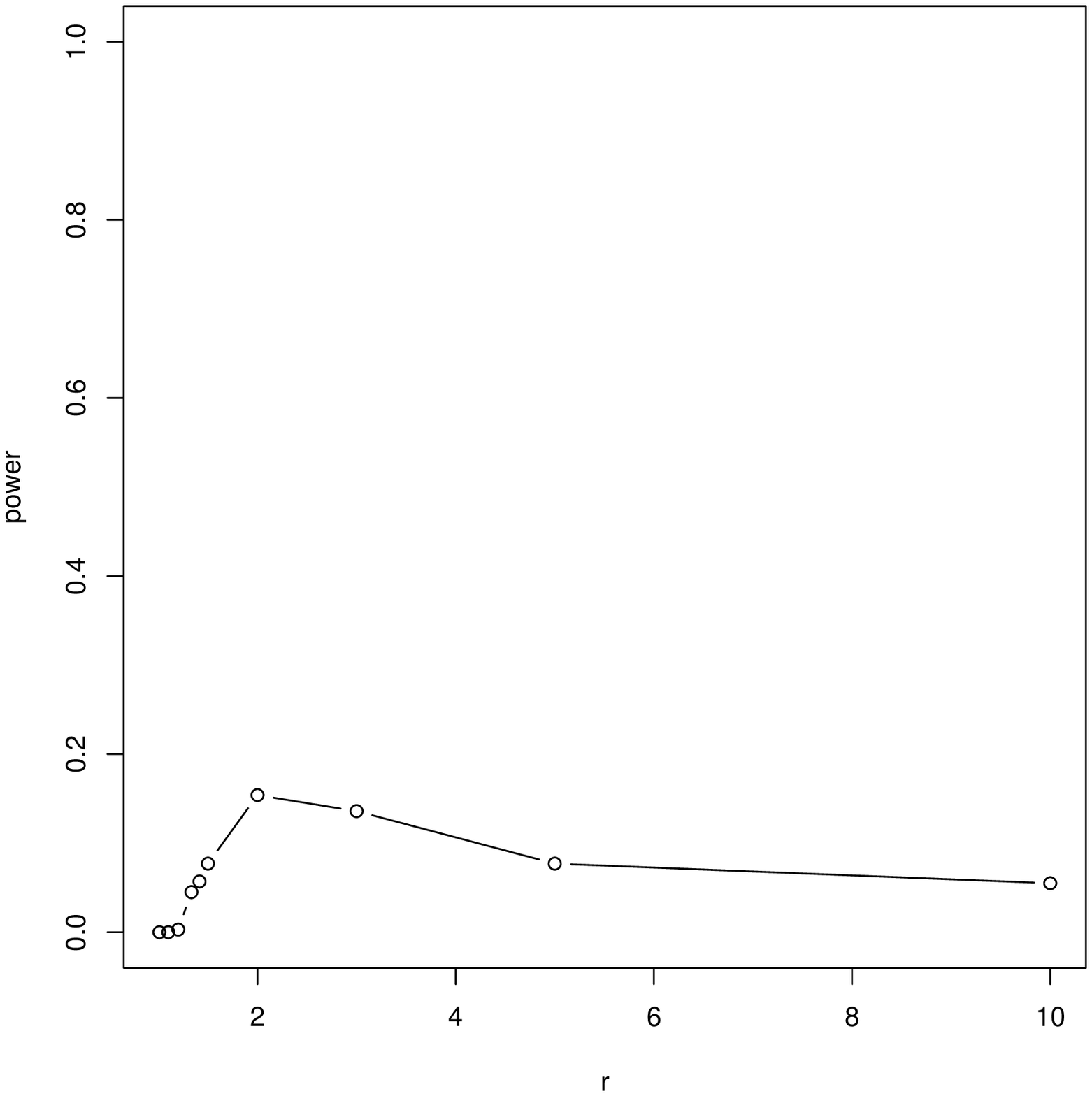} } }
\rotatebox{-90}{ \resizebox{2.5 in}{!}{ \includegraphics{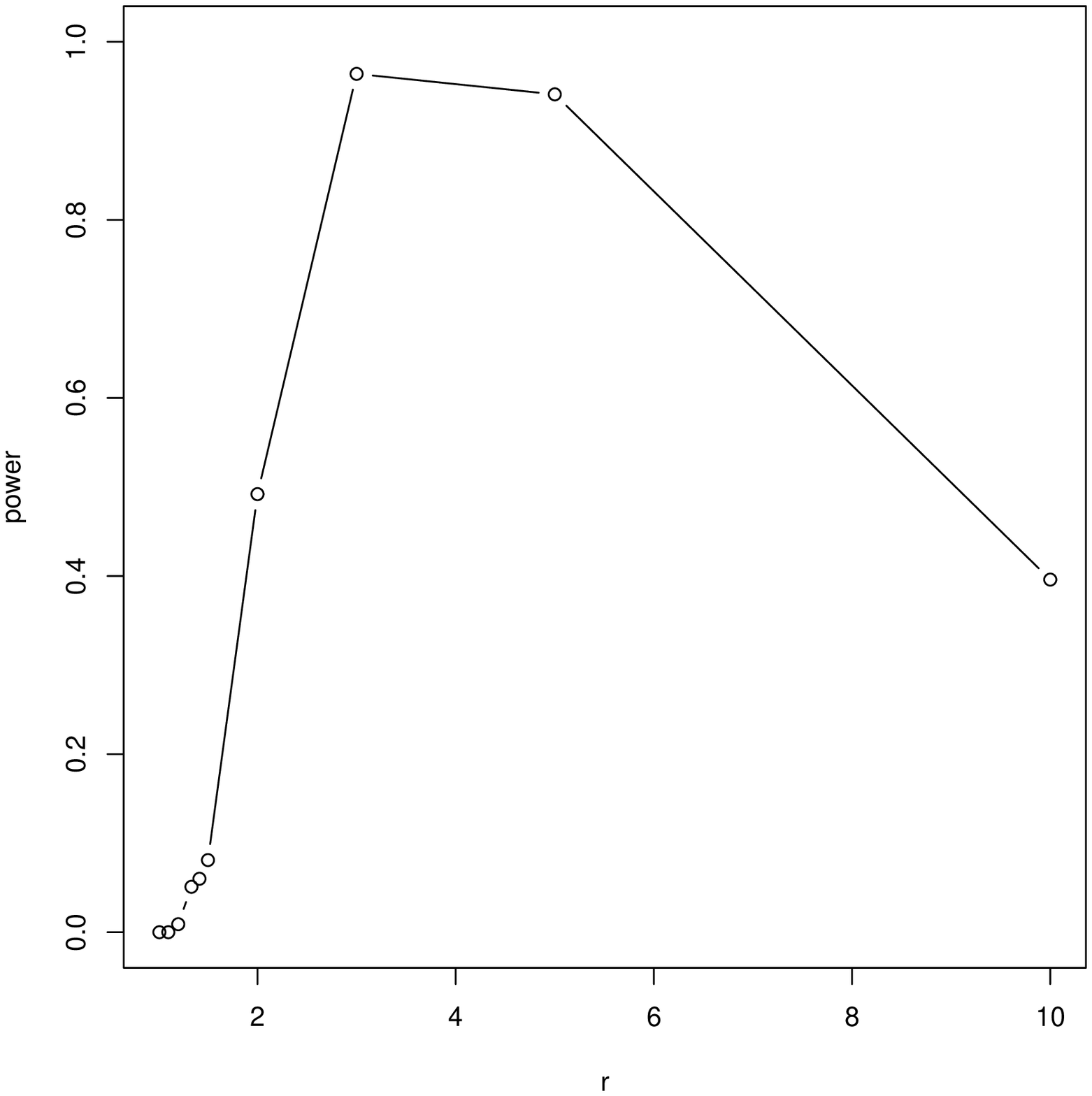} } }
\rotatebox{-90}{ \resizebox{2.5 in}{!}{ \includegraphics{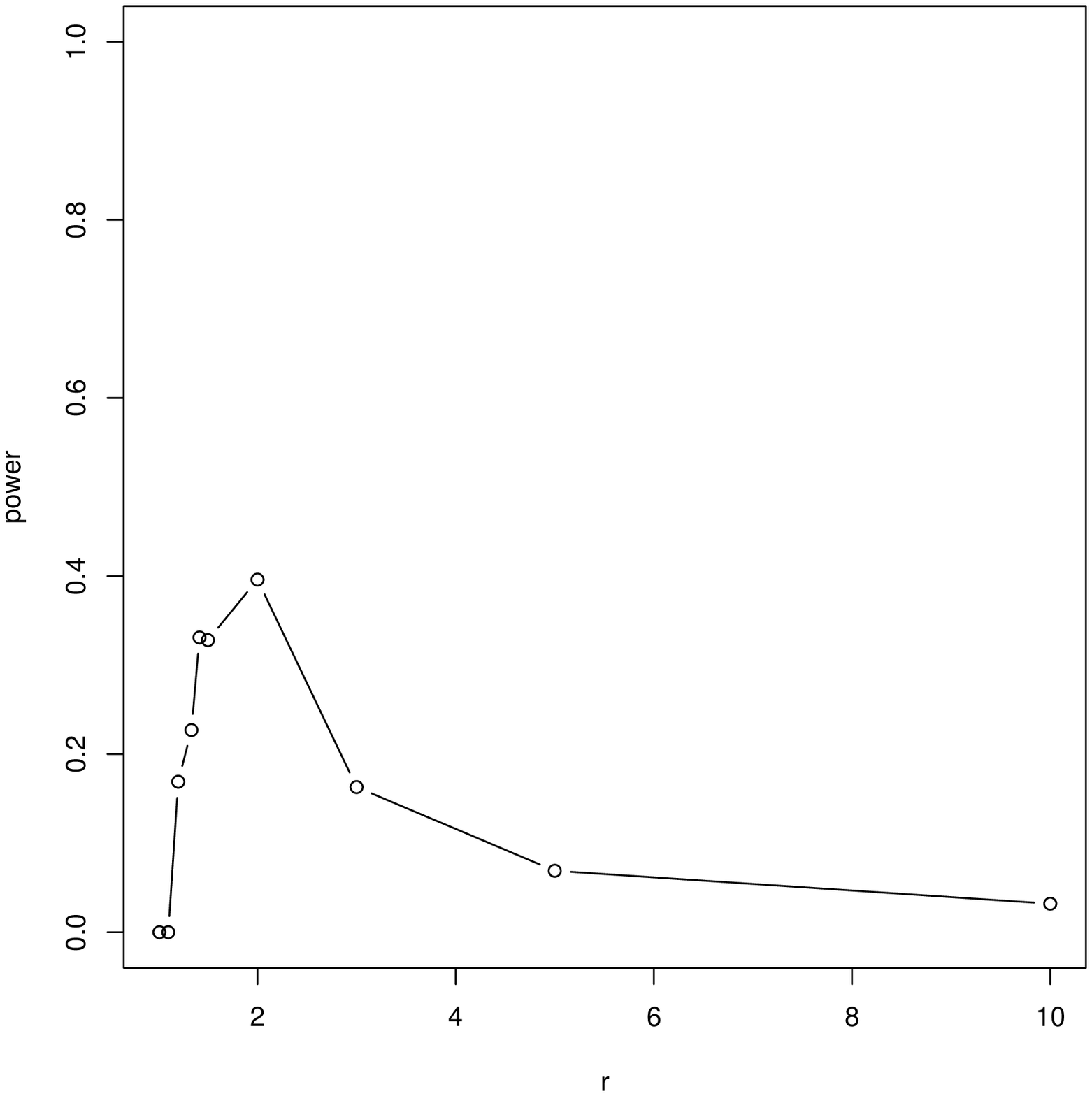} } }
\rotatebox{-90}{ \resizebox{2.5 in}{!}{ \includegraphics{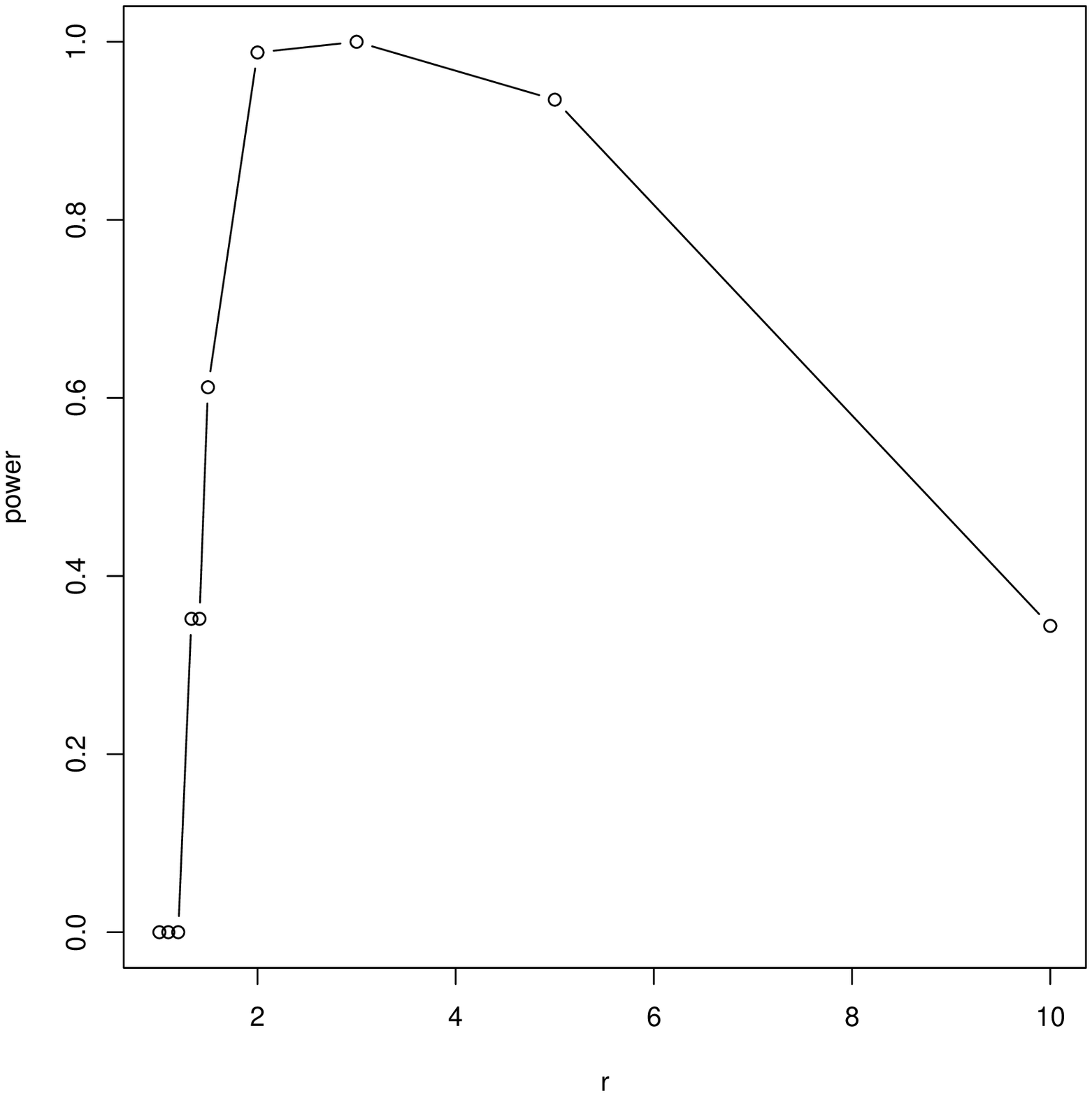} } }
\caption{ \label{fig:AggSimMCPower}
Empirical power estimates based on Monte Carlo critical values against the association alternatives
with the AND-underlying case (top two) and OR-underlying case (bottom two), in both cases,
$H^A_{\sqrt{3}/12 }$ (left)
and
$H^A_{5\,\sqrt{3}/24}$ (right)
as a function of $r$, for $n=10$ and $N_{mc}=1000$.
}
\end{figure}

In Figure \ref{fig:AggSimMCPower}, we also present a Monte Carlo
investigation of empirical power based on Monte Carlo critical values
against $H^A_{\sqrt{3}/12}$ and $H^A_{5\,\sqrt{3}/24}$
as a function of $r$ for $n=10$ with 1000 replicates.
The corresponding empirical power estimates are presented in Table \ref{tab:emp-val-A-Under}.
Our Monte Carlo estimates of $r^*_{\ve}$ are
$r^*_{\sqrt{3}/12} = 2$ and $r^*_{5\,\sqrt{3}/24} =3$ in both underlying cases.
That is, more severe association (larger
$\ve$) suggests a larger choice of $r$ in both cases.

\begin{table}[ht]
\centering
\begin{tabular}{|c|c|c|c|c|c|c|c|c|c|c|}
\hline
\multicolumn{11}{|c|}{$n=10$ and $N_{mc}=1000$ AND-underlying case} \\
\hline
$r$  & 1 & 11/10 & 6/5 & 4/3 & $\sqrt{2}$ & 3/2 & 2 & 3 & 5 & 10 \\
\hline
$\widehat{C}^A_{mc}$ & 0.0 & 0.0  & $0.0\bar 2$ &  $0.0\bar 6$ & $0.0\bar 8$ & $0.\bar{1}$ & $0.2\bar{4}$ & $0.4\bar{6}$ & $0.6\bar{8}$ & $0.8\bar{2}$\\
\hline
$\widehat{\alpha}^A_{mc}(n)$      & 0.000 & 0.000 & 0.005 & 0.030 & 0.027 & 0.037 & 0.038 & 0.043 & 0.048 & 0.041\\
\hline
$\widehat{\beta}^A_{mc}(\sqrt{3}/12)$ & 0.000 & 0.000 & 0.003 & 0.045 & 0.057 & 0.077 & 0.154 & 0.136 & 0.077  & 0.055\\
\hline
$\widehat{\beta}^A_{mc}(5\,\sqrt{3}/24)$ & 0.000 & 0.000  &0.009 & 0.051 & 0.060 & 0.081 & 0.492 & 0.964 & 0.941 & 0.396\\
\hline
\multicolumn{11}{|c|}{$n=10$ and $N_{mc}=1000$ OR-underlying case} \\
\hline
$r$  & 1 & 11/10 & 6/5 & 4/3 & $\sqrt{2}$ & 3/2 & 2 & 3 & 5 & 10 \\
\hline
$\widehat{C}^A_{mc}$ & $0.2\bar 6$ & $0.2\bar 6$ & $0.2\bar 8$ & $0.3\bar 1$ & $0.\bar 3$ & $0.3\bar 5$ & $0.6$ & $0.8\bar 4$ & $0.9\bar 5$ & 1.00\\
\hline
$\widehat{\alpha}^A_{mc}(n)$      & 0.000 & 0.000 & 0.040 & 0.045 & 0.049 & 0.042 & 0.049  & 0.044 & 0.022 & 0.019\\
\hline
$\widehat{\beta}^A_{mc}(\sqrt{3}/12)$ & 0.000 & 0.000 & 0.169 & 0.227 & 0.331 & 0.328 & 0.396 & 0.163 & 0.069 & 0.032\\
\hline
$\widehat{\beta}^A_{mc}(5\,\sqrt{3}/24)$ & 0.000 & 0.000 & 0.000 & 0.352 & 0.352 & 0.612 & 0.988 & 1.000 & 0.935 & 0.344 \\
\hline
\end{tabular}
\caption{ \label{tab:emp-val-A-Under}
Monte Carlo critical values, $\widehat{C}^A_{mc}$,
empirical significance levels, $\widehat{\alpha}^A_{mc}(n)$,
and empirical power estimates, $\widehat{\beta}^A_{mc}$, based on Monte Carlo critical values under
$H^A_{\sqrt{3}/12 }$
and
$H^A_{5\,\sqrt{3}/24}$,
$N_{mc}=1000$, and $n=10$ at
$\alpha=.05$.}
\end{table}

\begin{figure}[ht]
\centering
\psfrag{power}{\Huge{power}}
\psfrag{r}{\Huge{$r$}}
\rotatebox{-90}{ \resizebox{2.5 in}{!}{ \includegraphics{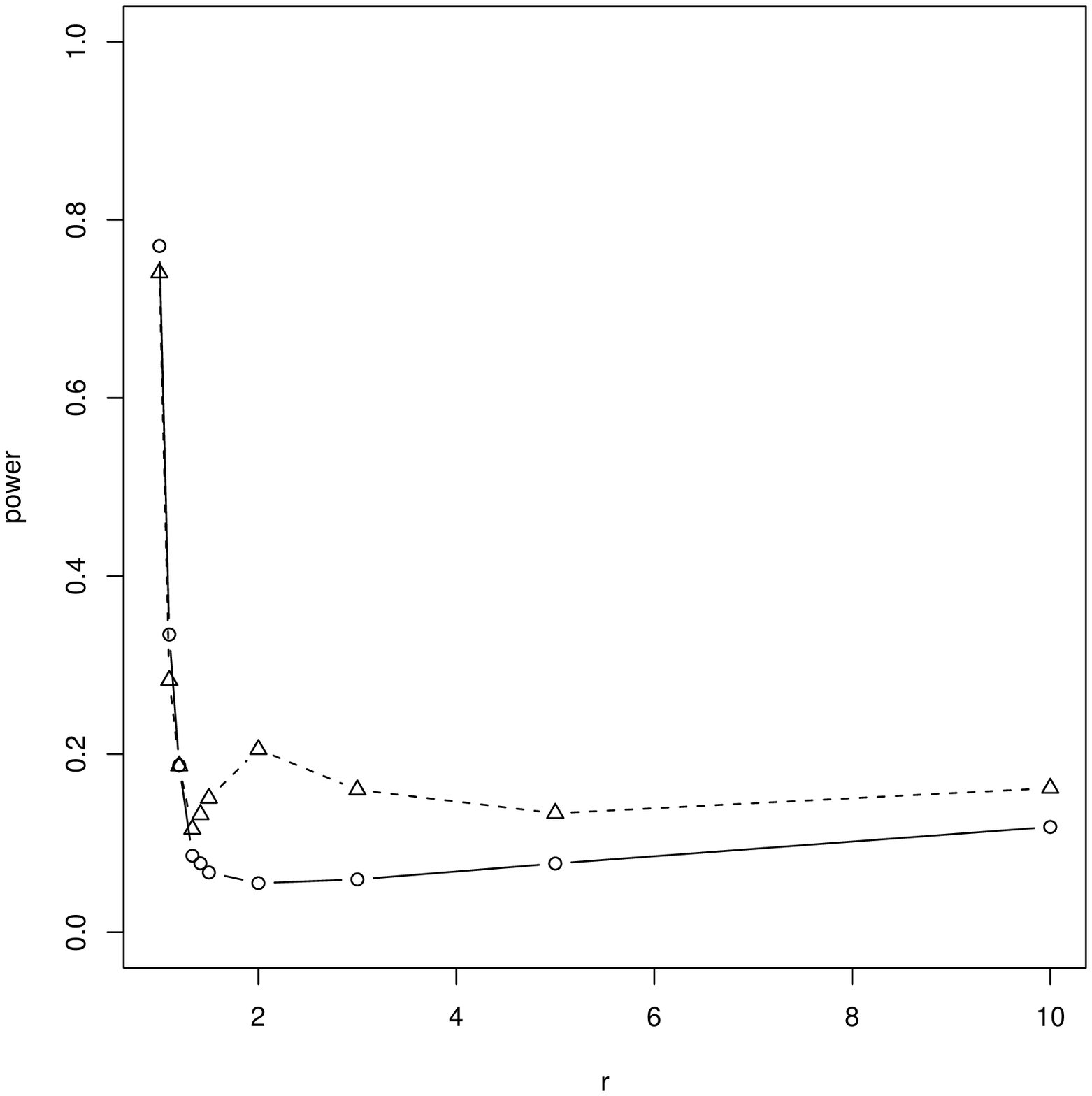} } }
\rotatebox{-90}{ \resizebox{2.5 in}{!}{ \includegraphics{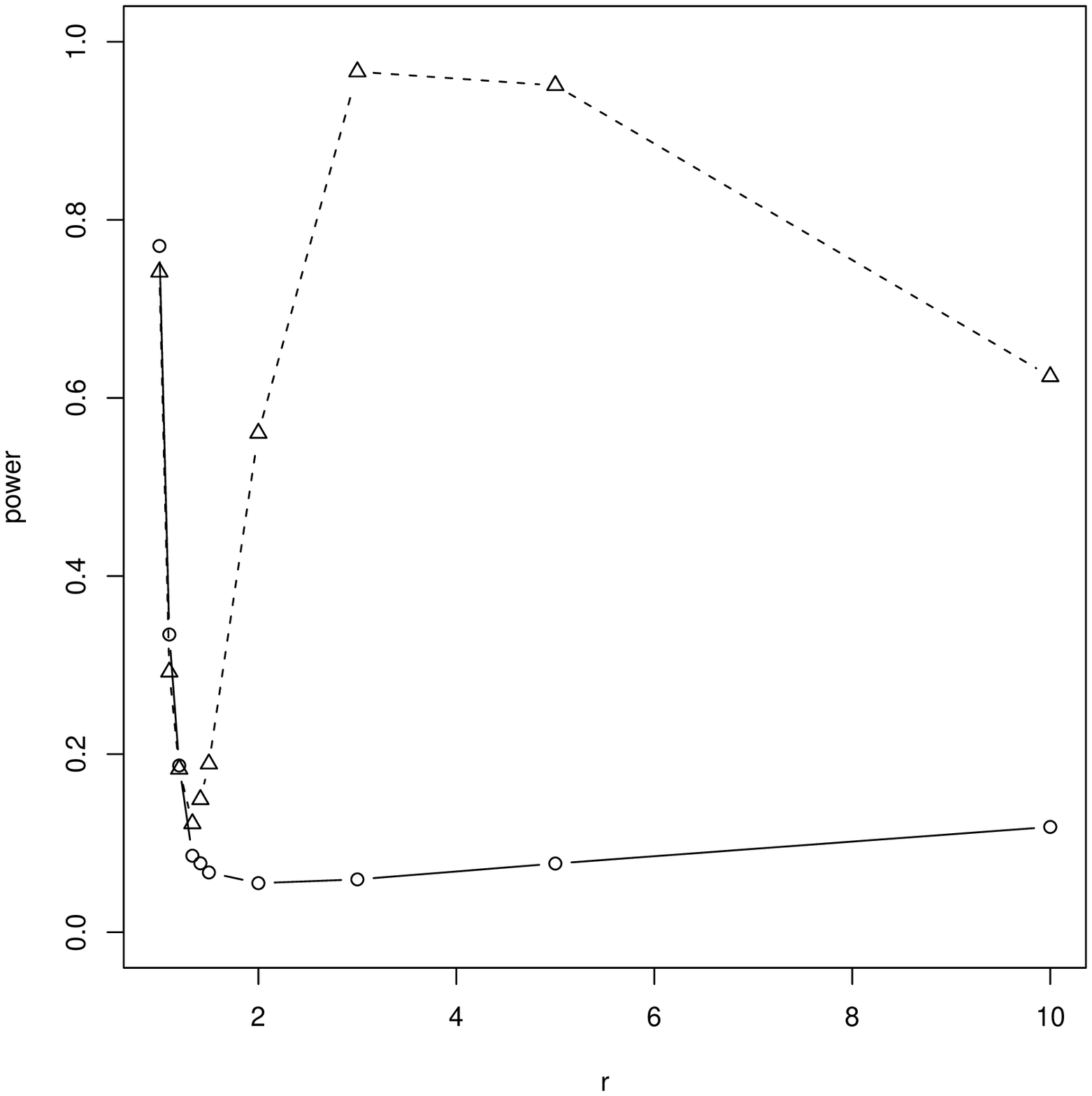} } }
\rotatebox{-90}{ \resizebox{2.5 in}{!}{ \includegraphics{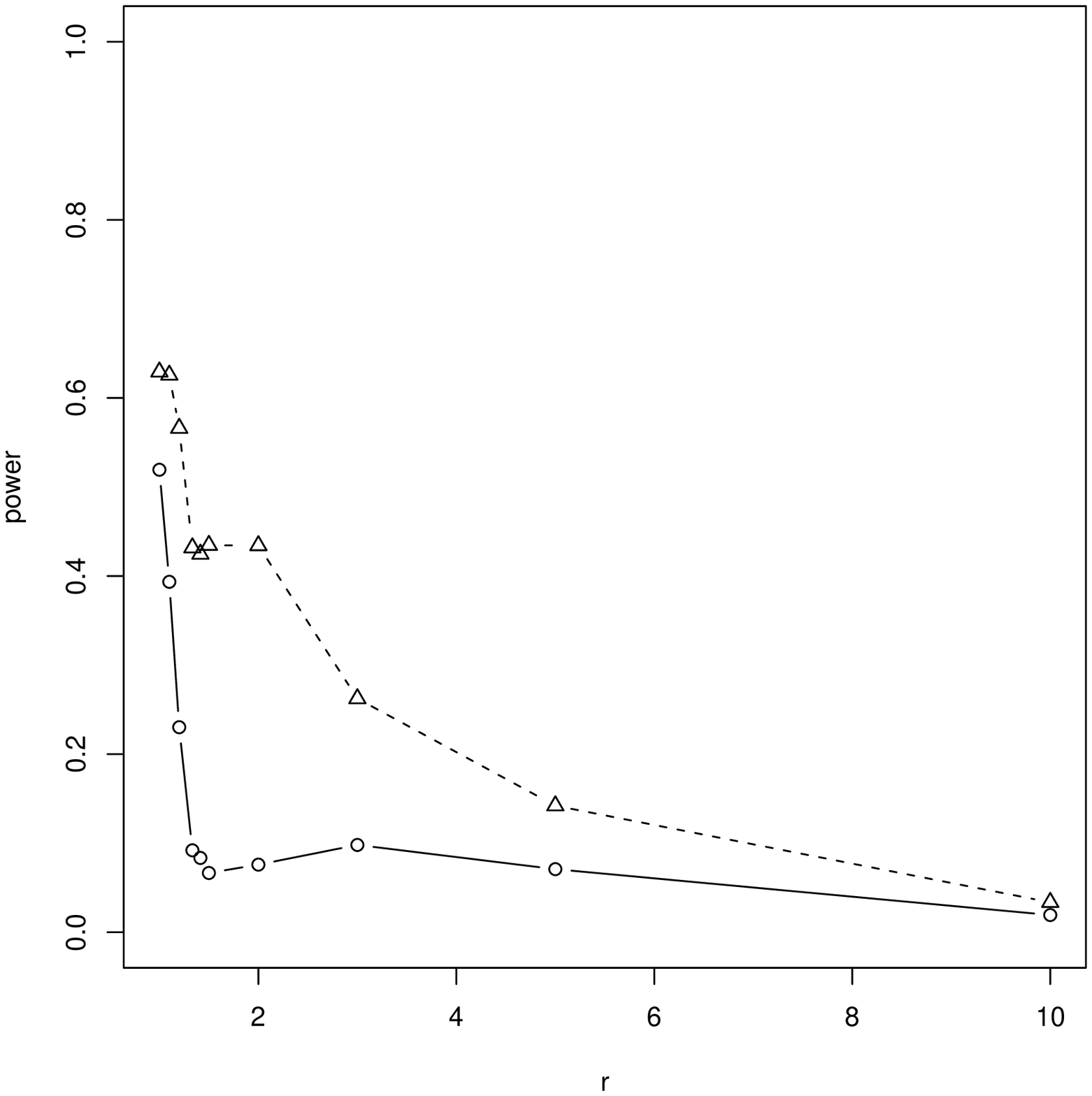} } }
\rotatebox{-90}{ \resizebox{2.5 in}{!}{ \includegraphics{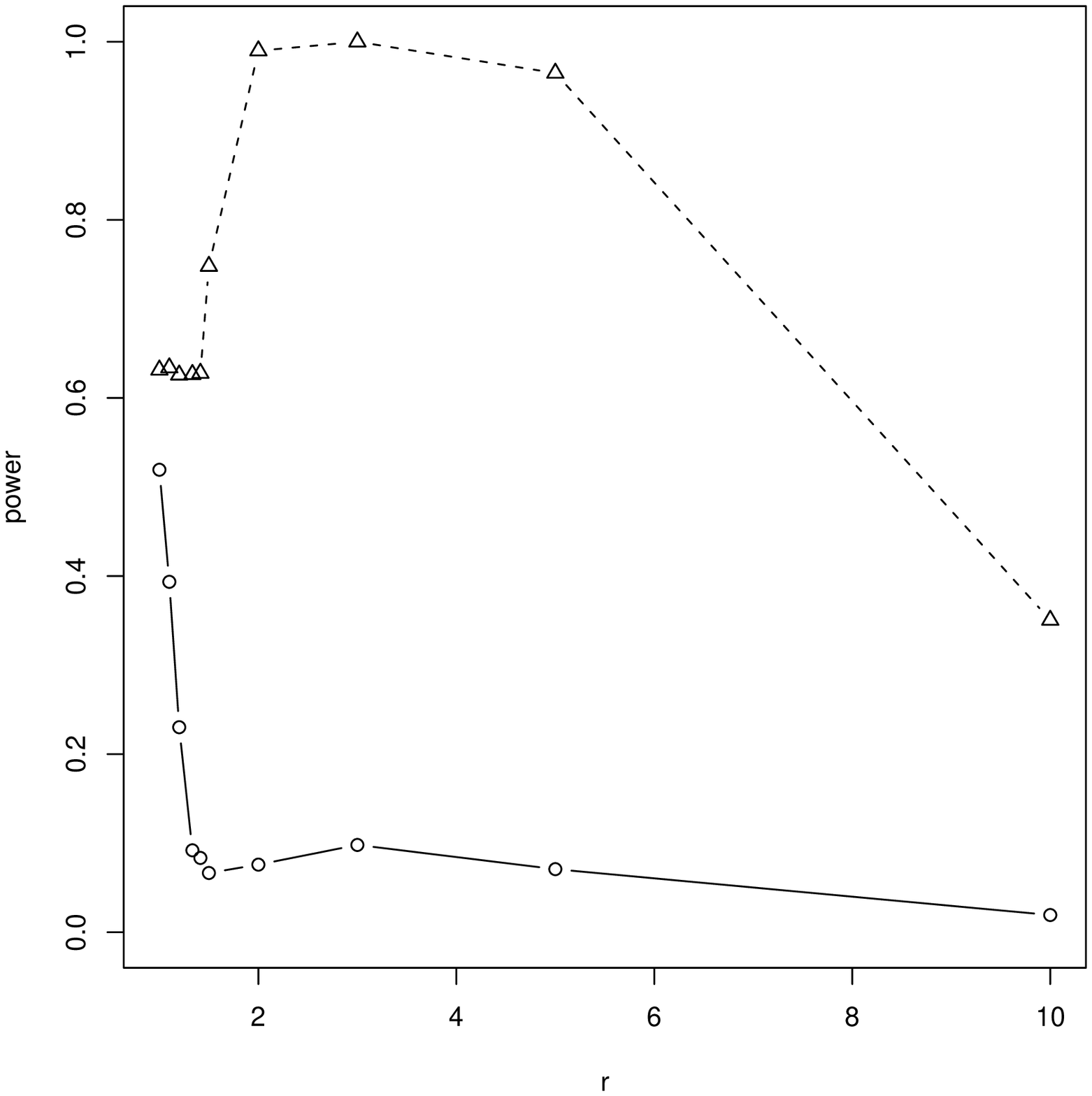} } }
\caption{ \label{fig:AggSimAsyPower}
The empirical size (circles joined with solid lines) and power estimates (triangles with dotted lines)
based on the asymptotic critical value against association
alternatives in the AND-underlying case (top two) and the OR-underlying case (bottom two),
in both cases,
$H^A_{\sqrt{3}/12}$ (left)
and
$H^A_{5\sqrt{3}/24}$ (right)
as a function of $r$, for $n=10$ and $N_{mc}=10000$.
}
\end{figure}

In Figure \ref{fig:AggSimAsyPower}, we present a Monte Carlo
investigation of power based on asymptotic critical values
against $H^A_{\sqrt{3}/12}$ and $H^A_{5\,\sqrt{3}/24}$ as a function of $r$ for $n=10$.
In the AND-underlying case,  the empirical
significance level, $\widehat{\alpha}_{n=10}(r)$, is about $.05$ for $r=2$
and $3$ which have the empirical power $\widehat{\beta}_{10}(2)
\approx .2$ with maximum power at $r=2$ for
$\ve=\sqrt{3}/12$, and $\widehat{\beta}_{10}(3) =1$ for
$\ve=5\,\sqrt{3}/24$.
In the OR-underlying case,  the
empirical significance level, $\widehat{\alpha}_{n=10}(r)$, is closest to
$.05$ for $r=1.5$ which have the empirical power
$\widehat{\beta}_{10}(1.5) \approx .45$ for
$\ve=\sqrt{3}/12$, and $\widehat{\beta}_{10}(1.5) =1$ for
$\ve=5\,\sqrt{3}/24$.
So, for small sample sizes, moderate
values of $r$ is more appropriate for normal approximation, as they
yield the desired significance level, and the more severe the
association, higher the power estimate.
Furthermore, the OR-underlying
version seems to perform better than the AND-underlying version for
association alternatives.
The empirical significance levels, and empirical power
$\widehat{\beta}^S_n(r,\ve)$ values based on asymptotic critical values under
$H^A_{\ve}$ for $\ve=\sqrt{3}/12,\,5\,\sqrt{3}/24$
are given in Table \ref{tab:asy-emp-val-A-Under}.

\begin{table}[ht]
\centering
\begin{tabular}{|c|c|c|c|c|c|c|c|c|c|c|}
\hline
\multicolumn{11}{|c|}{$n=10$ and $N_{mc}=1000$ AND-underlying case} \\
\hline
$r$  & 1 & 11/10 &6/5 & 4/3 & $\sqrt{2}$ & 3/2 & 2 & 3 & 5 & 10 \\
\hline
$\widehat{\alpha}_A(n)$ & 0.7707 & 0.3343 & 0.1872 & 0.0859 & 0.0774 & 0.0671 & 0.0551 & 0.0593 & 0.0771 & 0.1182 \\
\hline
$\widehat{\beta}^A_{n}(r,\sqrt{3}/12)$ & 0.7406 & 0.2829 & 0.1869 & 0.1156 & 0.1323 & 0.1506 & 0.2053 & 0.1599 & 0.1336 & 0.1618\\
\hline
$\widehat{\beta}^A_{n}(r,5\,\sqrt{3}/24)$ & 0.7415 & 0.2923 & 0.1833 & 0.1220 & 0.1491 & 0.1891 & 0.5605 & 0.9664 & 0.9510 & 0.6241  \\
\hline
\multicolumn{11}{|c|}{$n=10$ and $N_{mc}=1000$ OR-underlying case} \\
\hline
$r$  & 1 & 11/10 &6/5 & 4/3 & $\sqrt{2}$ &3/2 & 2 & 3 &5 & 10 \\
\hline
$\widehat{\alpha}_A(n)$ & 0.5194 & 0.3935 & 0.2302 & 0.0920 & 0.0834 & 0.0665 & 0.0759 & 0.0980 & 0.0708 & 0.0193 \\
\hline
$\widehat{\beta}^A_{n}(r,\sqrt{3}/12)$ & 0.6293 & 0.6258 & 0.5661 & 0.4318 & 0.4247 & 0.4346 & 0.4343 & 0.2624 & 0.1421 & 0.0336\\
\hline
$\widehat{\beta}^A_{n}(r,5\,\sqrt{3}/24)$ & 0.6315 & 0.6340 & 0.6259 &  0.6265 & 0.6279 & 0.7480 & 0.9900 & 1.0000 & 0.9649 & 0.3505 \\
\hline
\end{tabular}
\caption{ \label{tab:asy-emp-val-A-Under}
The empirical significance level and empirical power estimates based on asymptotic critical values under $H^A_{\ve}$ for
$\ve=\sqrt{3}/12,\,5\,\sqrt{3}/24$, $N_{mc}=10000$, and $n=10$ at
$\alpha=.05$.}
\end{table}

\section{Multiple Triangle Case}
\label{sec:multiple-triangle-case}
Suppose $\Y_m$ is a finite collection of $m>3$ points in $\mathbb{R}^2$.
Consider the Delaunay triangulation (assumed to exist) of $\Y_m$.
Let $T_i$ denote the $i^{th}$ Delaunay triangle,
$J_m$ denote the number of triangles, and
$C_H(\Y_m)$ denote the convex hull of $\Y_m$.
We wish to investigate
$H_o: X_i \stackrel{iid}{\sim} \mathcal{U}(C_H(\Y_m))$
against segregation and association alternatives using
the relative edge densities of the associated underlying graphs.
The underlying graphs are constructed using the PCD $D$,
which is constructed using $\NPE^{r}(\cdot)$ as described in Section \ref{sec:prop-edge},
where for $X_i \in T_j$,
the three points in $\Y_m$ defining the
Delaunay triangle $T_j$ are used as $\Y_{[j]}$.
We consider various versions of the relative edge density
as a test statistic in the multiple triangle case.

\subsection{First Version of Relative Edge Density in the Multiple Triangle Case}
\label{sec:version-I-mult-tri}
For $J_m>1$,
as in Section \ref{sec:relative-density-PEPCD},
let $\rho^{\la}_{I,n}(r)=2\,\left|\mE_{\la}\right|/(n\,(n-1))$
and $\rho^{\lo}_n(r)=2\,\left|\mE_{\lo}\right|/(n\,(n-1))$.
Let $\mE^\la_i$ be the number of edges and
$\rho^{\la}_{{}_{[i]}}(r)$ be the relative edge density for triangle $i$ in the AND-underlying case,
and $\mE^\lo_i$ and $\rho^{\lo}_{{}_{[i]}}(r)$ be similarly defined for OR-underlying case.
Let $n_i$ be the number of $X$ points in $T_i$ for $i=1,2,\ldots,J_m$.
Letting $w_i = A(T_i) / A(C_H(\Y_m))$ with $A(\cdot)$ being the area functional,
we obtain the following as a corollary to Theorem \ref{thm:asy-norm-under}.

\begin{corollary}
\label{cor:MT-asy-norm-NYr}
The asymptotic null distribution for $\rho^{\la}_{I,n}(r)$ conditional on $\Y_m$
for $r \in (1,\infty)$
is given by
\begin{equation}
\sqrt{n}\left(\rho^{\la}_{I,n}(r)-\widetilde \mu_{\la}(r)\right)
\stackrel{\mathcal L}{\longrightarrow}\\
\mathcal{N}
 \left(
   0,
   4\,\widetilde \nu_{\la}(r)
 \right),
\end{equation}
where $\widetilde \mu_{\la}(r)=\mu_{\la}(r) \left(\sum_{i=1}^{J_m}w_i^2\right)$
and
$\widetilde \nu_{\la}(r)=
\left[  \nu_{\la}(r) \left(\sum_{i=1}^{J_m}w_i^3 \right)+
\left( \mu_{\la}(r) \right)^2\left(\sum_{i=1}^{J_m}w_i^3-\left(\sum_{j=1}^{J_m}w_i^2 \right)^2\right) \right]$
with $\mu_{\la}(r)$ and $\nu_{\la}(r)$ being as in Equations (\ref{eqn:Asymean_and}) and (\ref{eqn:Asyvar_and}),
respectively.
The asymptotic null distribution of $\rho^{\lo}_{I,n}(r)$ with $r \in [1,\infty)$ is similar.
\end{corollary}

The Proof is provided in Appendix 7.
By an appropriate application of the Jensen's Inequality, we see
that $\sum_{i=1}^{J_m}w_i^3 \ge \left(\sum_{i=1}^{J_m}w_i^2 \right)^2.$
So the covariance above is zero iff $\nu_{\la}(r)=0$ and
$\sum_{i=1}^{J_m}w_i^3=\left(\sum_{i=1}^{J_m}w_i^2 \right)^2$, so
asymptotic normality may hold even though $\nu_{\la}(r)=0$.
The same holds for the OR-underlying case.

Under the segregation (association) alternatives with
$\delta \times 100 \%$ where $\delta =4\,\ve^2/3$ around the
vertices of each triangle is forbidden (allowed),
we obtain the above asymptotic distribution of $\rho^{\la}_{I,n}(r)$ with
$\mu_{\la}(r)$ being replaced by $\mu_{\la}(r,\ve)$ and
$\nu_{\la}(r)$ by $\nu_{\la}(r,\ve)$.
The OR-underlying case is similar.

\subsection{Other Versions of Relative Edge Density in the Multiple Triangle Case}
\label{sec:version-II-mult-tri}

Let $\displaystyle \Xi^{\la}_n(r):=\sum_{i=1}^{J_m}\frac{n_i\,(n_i-1)}{n\,(n-1)} \rho^{\la}_{{}_{[i]}}(r)$.
Then $\Xi^{\la}_n(r) = \rho^{\la}_{I,n}(r)$,
since $\displaystyle \Xi^{\la}_n(r)=\sum_{i=1}^{J_m}\frac{n_i\,(n_i-1)}{n\,(n-1)} \rho^{\la}_{{}_{[i]}}(r)=
\frac{\sum_{i=1}^{J_m} 2\,|\mE^\la_i|}{n\,(n-1)}=\frac{2\,\left|\mE_{\la}\right|}{n\,(n-1)}=\rho^{\la}_{I,n}(r)$.
Similarly, $\Xi^{\lo}_n(r) = \rho^{\lo}_n(r)$.

Furthermore, let $\widehat{\Xi}^{\la}_n:=
\sum_{i=1}^{J_m}w_i^2\,\rho^{\la}_{{}_{[i]}}(r)$ where $w_i$ is as above.
So $\widehat{\Xi}^{\la}_n$ a mixture of $\rho^{\la}_{{}_{[i]}}(r)$'s.
Then since $\rho^{\la}_{{}_{[i]}}(r)$'s are asymptotically independent,
$\Xi^{\la}_n(r), \, \rho^{\la}_{I,n}(r)$ are asymptotically normal;
i.e., for large $n$ their distribution is approximately
$\mathcal N\left( \widetilde \mu_{\la}(r),4\,\widetilde \nu_{\la}(r)/n\right)$.
A similar result holds for the OR-underlying case.

In Section \ref{sec:version-I-mult-tri},
the denominator of $\rho^{\la}_{I,n}(r)$ has
$n(n-1)/2$ as the maximum number of edges possible.
However, by definition,
given the $n_i$'s
we can at most have a graph with $J_m$ complete components,
each with order $n_i$ for $i=1,2,\ldots,J_m$.
Then the maximum number of edges possible is $n_t:=\sum_{i=1}^{J_m}n_i\,(n_i-1)/2$
which suggests another version of relative edge density:
$\displaystyle \rho^{\la}_{II,n}(r):=\frac{\left|\mE_{\la}\right|}{n_t}$.
Then $\displaystyle \rho^{\la}_{II,n}(r)=\frac{\sum_{i=1}^{J_m} |\mE^\la_i|}{n_t}=
\sum_{i=1}^{J_m}\frac{n_i\,(n_i-1)}{2\,n_t}\,\rho^{\la}_{{}_{[i]}}(r)$.
Since $\frac{n_i\,(n_i-1)}{2\,n_t} \ge 0$ for each $i$,
and $\displaystyle \sum_{i=1}^{J_m}\frac{n_i\,(n_i-1)}{2\,n_t}=1$,
$\rho^{\la}_{II,n}(r)$ is a mixture of $\rho^{\la}_{{}_{[i]}}(r)$'s.

\begin{theorem}
\label{thm:MT-asy-norm-II}
The asymptotic null distribution for $\rho^{\la}_{II,n}(r)$ conditional on $\Y_m$
for $r \in (1,\infty)$
is given by
\begin{equation}
\sqrt{n}\left(\rho^{\la}_{II,n}(r)-\breve \mu_{\la}(r)\right)
\stackrel{\mathcal L}{\longrightarrow}\\
\mathcal{N}
 \left(
   0,
   4\,\breve \nu_{\la}(r)
 \right),
\end{equation}
where $\breve \mu_{\la}(r)=\mu_{\la}(r)$
and
$\breve \nu_{\la}(r)=
\left[  \nu_{\la}(r) \left(\sum_{i=1}^{J_m}w_i^3 \right)\Big/\left(\sum_{i=1}^{J_m}w_i^2 \right)^2 \right]$
with $\mu_{\la}(r)$ and $\nu_{\la}(r)$ being as in Equations (\ref{eqn:Asymean_and}) and (\ref{eqn:Asyvar_and}),
respectively.
The asymptotic null distribution of $\rho^{\lo}_{II,n}(r)$ with $r \in [1,\infty)$ is similar.
\end{theorem}

Proof is provided in Appendix 8.
Notice that the covariance $\breve \nu_{\la}(r)$ is zero iff $\nu_{\la}(r)=0$,
Under the segregation (association) alternatives,
we obtain the above asymptotic distribution of $\rho^{\la}_{II,n}(r)$ with
$\mu_{\la}(r)$ being replaced by $\mu_{\la}(r,\ve)$ and
$\nu_{\la}(r)$ by $\nu_{\la}(r,\ve)$.
The OR-underlying case is similar.

\begin{remark}
\label{rem:comp-versions-mult-tri}
\textbf{Comparison of Versions of Relative Edge Density in the Multiple Triangle Case:}
Among the versions of the relative edge density we considered,
$\Xi^{\la}_n(r) = \rho^{\la}_{I,n}(r)$ for all $n>1$,
and
$\widehat{\Xi}^{\la}_n$ and $\rho^{\la}_{I,n}(r)$
are asymptotically equivalent (i.e., they have the
same asymptotic distribution in the limit).
However, $\rho^{\la}_{I,n}(r)$ and $\rho^{\la}_{II,n}(r)$
do not have the same distribution for finite or infinite $n$.
But we have $\rho^{\la}_{I,n}(r)=\frac{2\,n_t}{n(n-1)}\rho^{\la}_{II,n}(r)$
and $\widetilde \mu_{\la}(r) < \breve \mu_{\la}(r)=\mu_{\la}(r)$,
since $\sum_{i=1}^{J_m}w_i^2<1$.
Furthermore,
since $\frac{2\,n_t}{n(n-1)}=\sum_{i=1}^{J_m} \frac{n_i(n_i)}{n(n-1)}
\longrightarrow \sum_{i=1}^{J_m} w_i^2$,
we have
$\lim_{n_i \rightarrow \infty} \Var[\sqrt{n}\rho^{\la}_{I,n}(r)]=
\left(\sum_{i=1}^{J_m} w_i^2\right)^2 \lim_{n_i \rightarrow \infty} \Var[\sqrt{n}\rho^{\la}_{I,n}(r)]$
Hence $\breve \nu_{\la}(r) \ge \widetilde \nu_{\la}(r)$.
Therefore,
we choose $\rho^{\la}_{I,n}(r)$
for further analysis in the multiple triangle case.
Moreover, asymptotic normality might hold for $\rho^{\la}_{I,n}(r)$
even if $\nu_{\la}(r)=0$.
$\square$
\end{remark}

\subsection{Power Analysis for the Multiple Triangle Case}
\label{sec:power-mult-tri}
Let $S^{\la}_n(r):=\rho^{\la}_{I,n}(r)$
and
$S^{\lo}_n(r):=\rho^{\lo}_{I,n}(r)$.
Thus in the case of $J_m>1$ (i.e., $m > 3$),
we have a (conditional) test of $H_o: X_i \stackrel{iid}{\sim} \U(C_H(\Y_m))$
which once again rejects against segregation for large values of $S^{\la}_n(r)$ and rejects
against association for small values of $S^{\la}_n(r)$.
The same holds for $S^{\lo}_n(r)$.

\begin{figure}[ht]
\centering
\rotatebox{-00}{ \resizebox{2.in}{!}{ \includegraphics{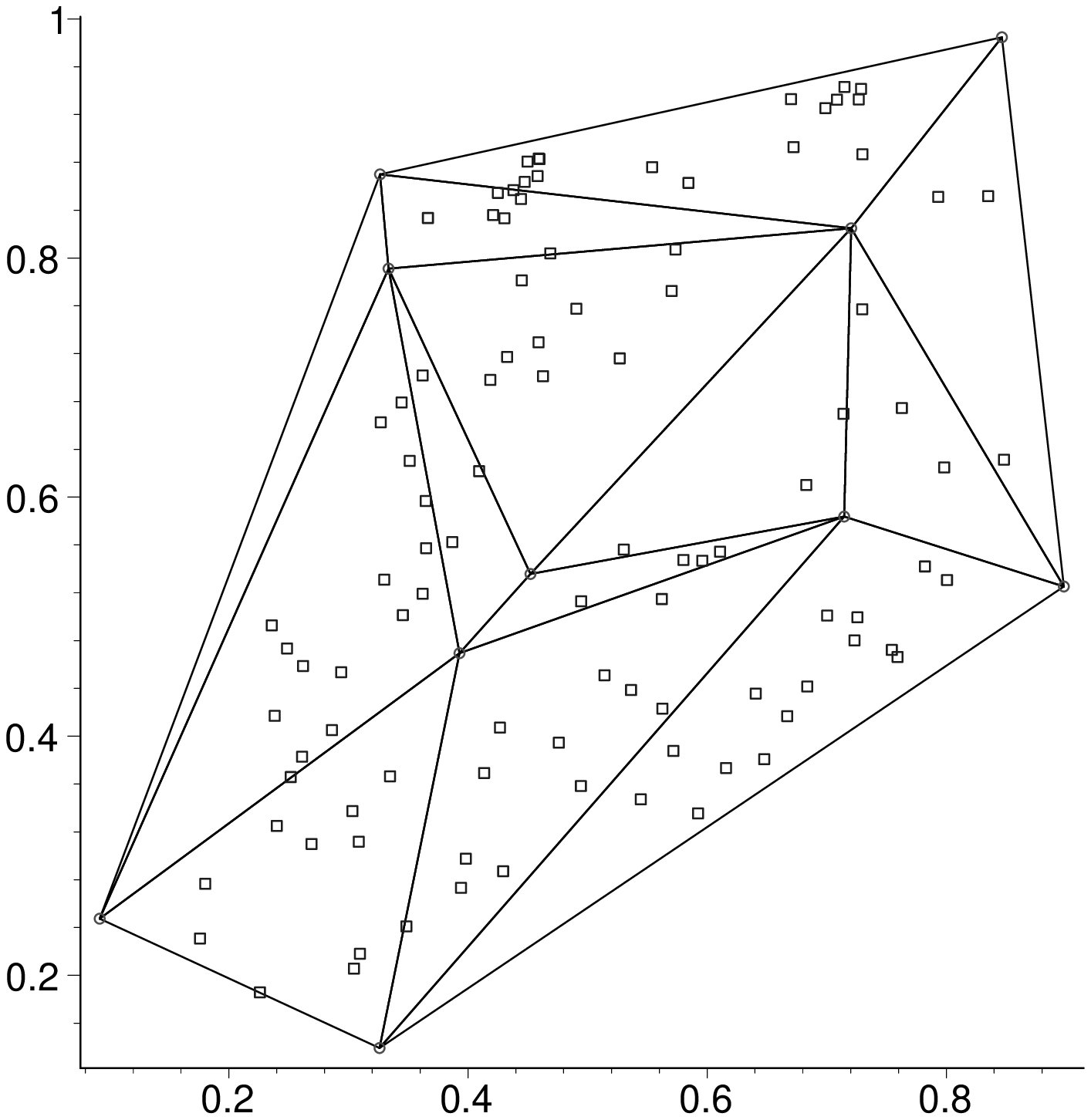} } }
\rotatebox{-00}{ \resizebox{2.in}{!}{ \includegraphics{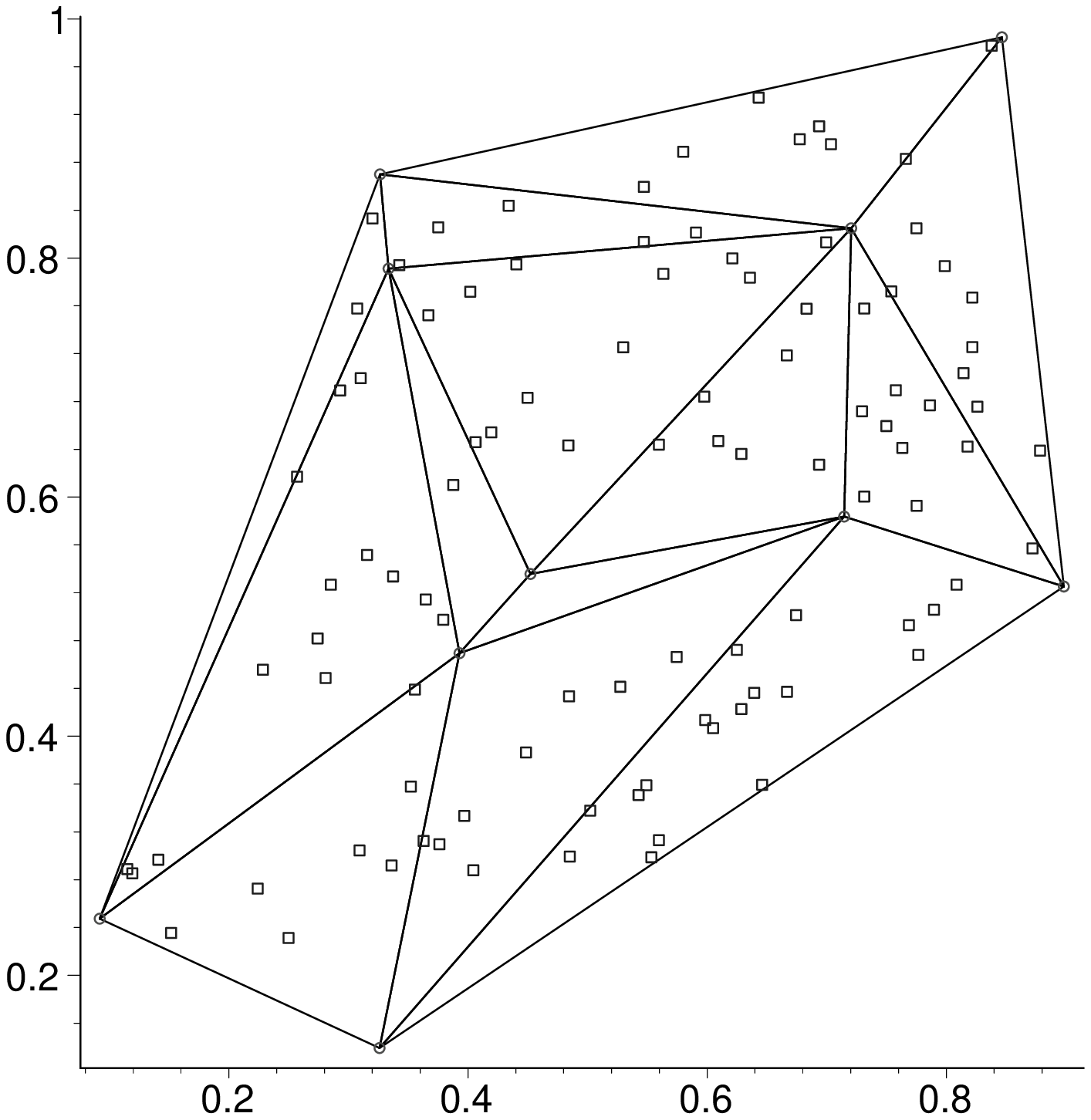} } }
\rotatebox{-00}{ \resizebox{2.in}{!}{ \includegraphics{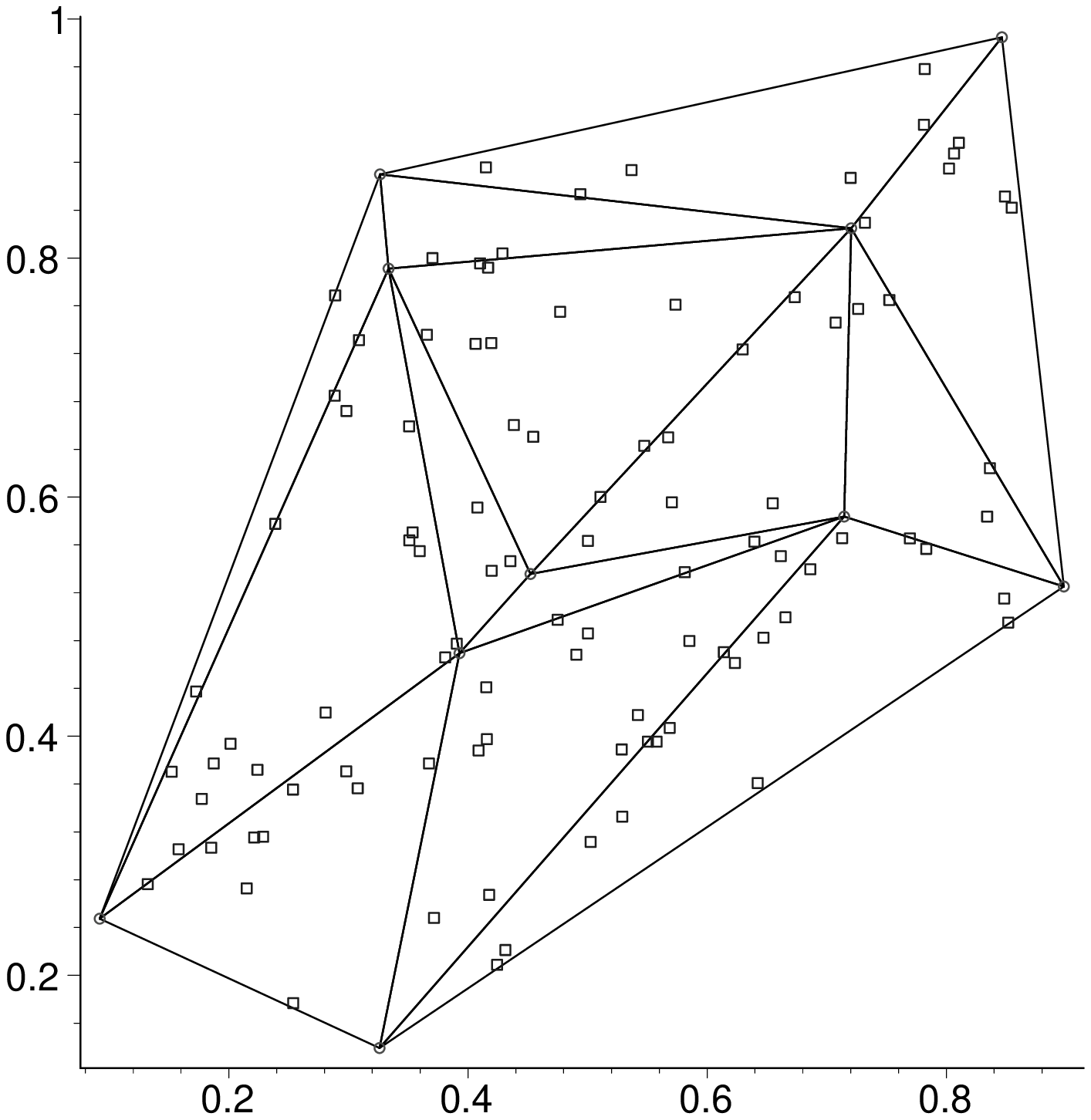} } }
\caption{\label{fig:deldata1}
Realization of segregation (left), $H_o$ (middle), and association (right) for $|\Y_m|=10$ and $n=100$.
}
\end{figure}

\begin{figure}[ht]
\centering
\rotatebox{-00}{ \resizebox{2.in}{!}{ \includegraphics{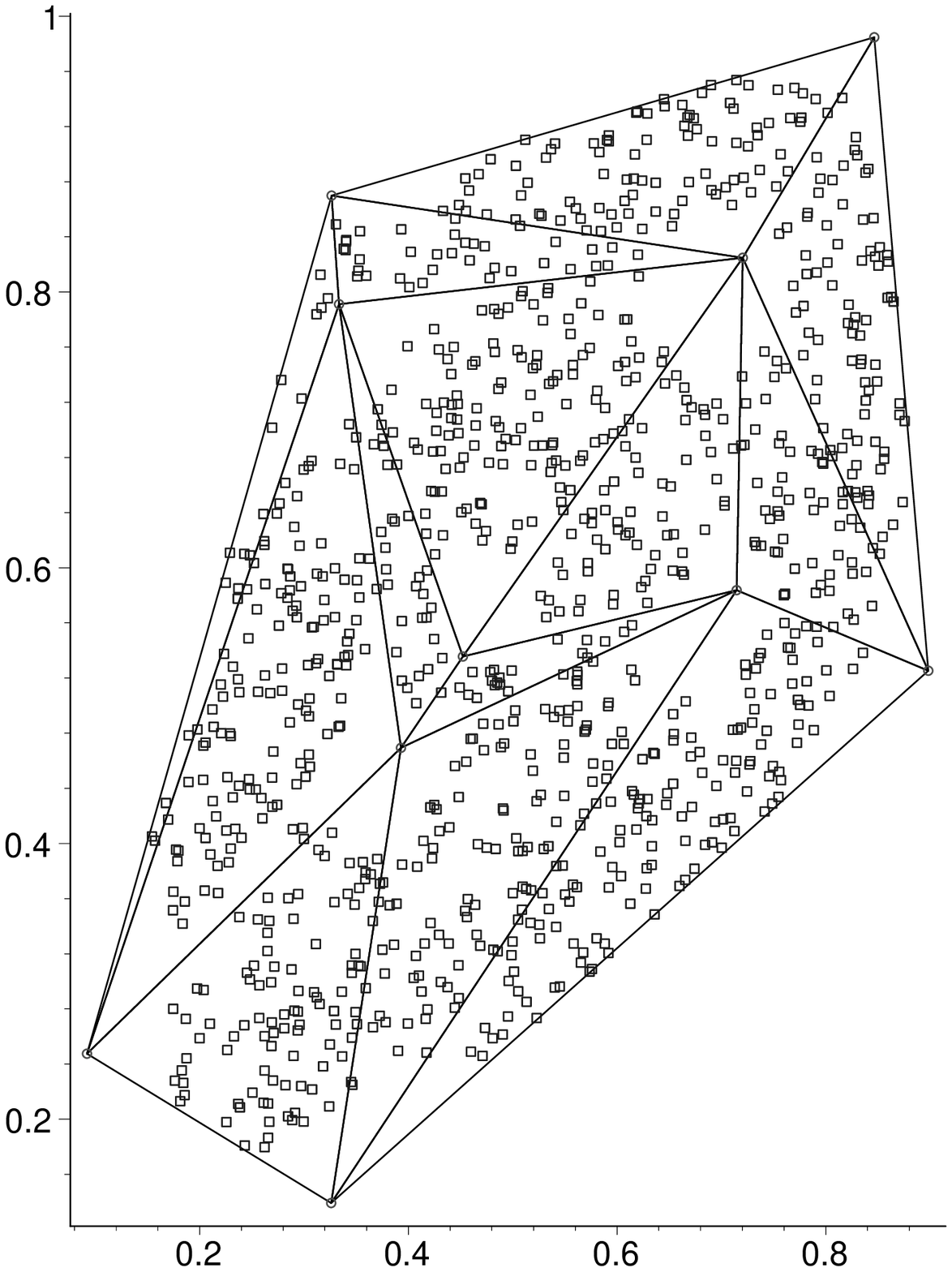} } }
\rotatebox{-00}{ \resizebox{2.in}{!}{ \includegraphics{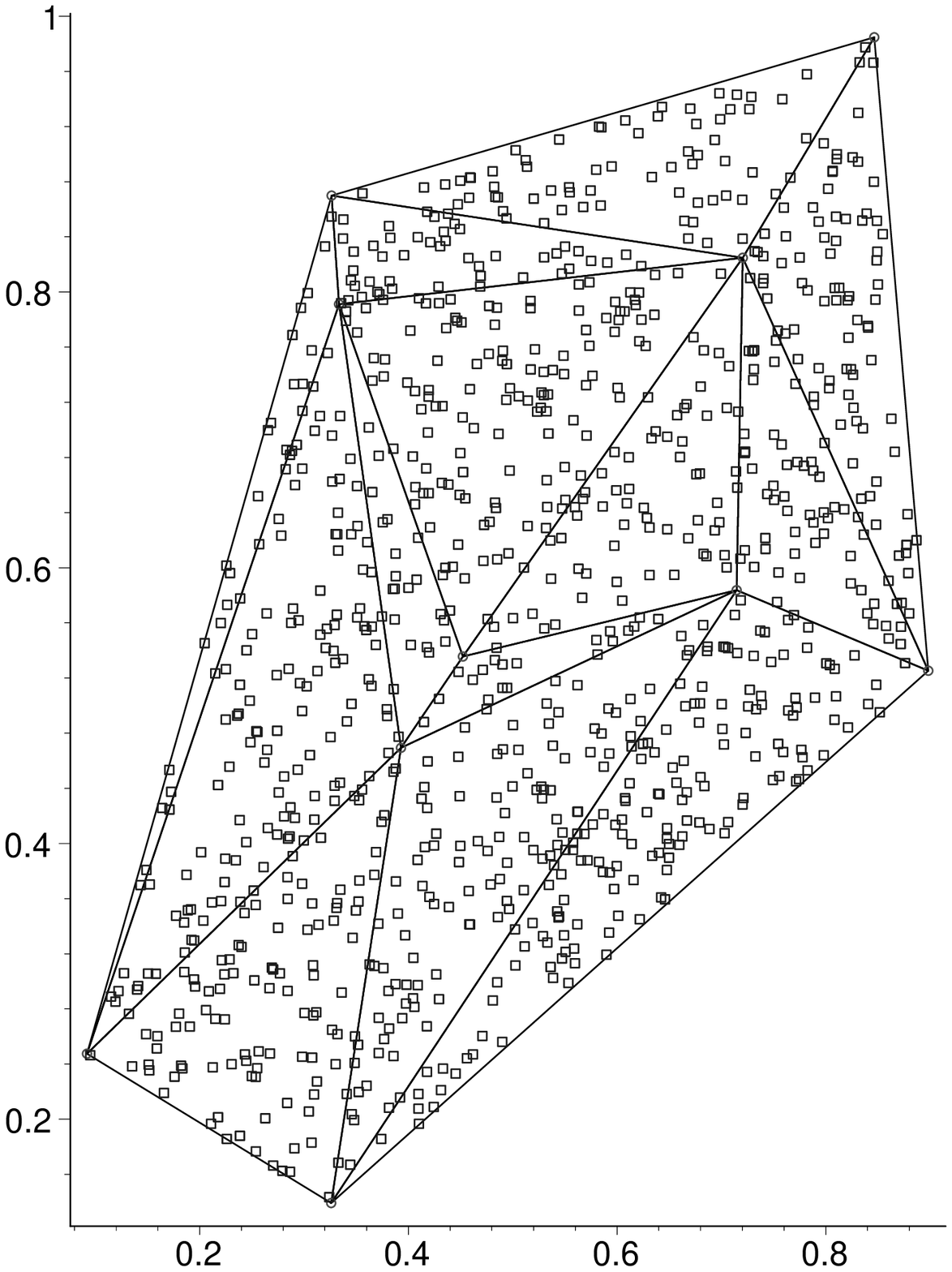} } }
\rotatebox{-00}{ \resizebox{2.in}{!}{ \includegraphics{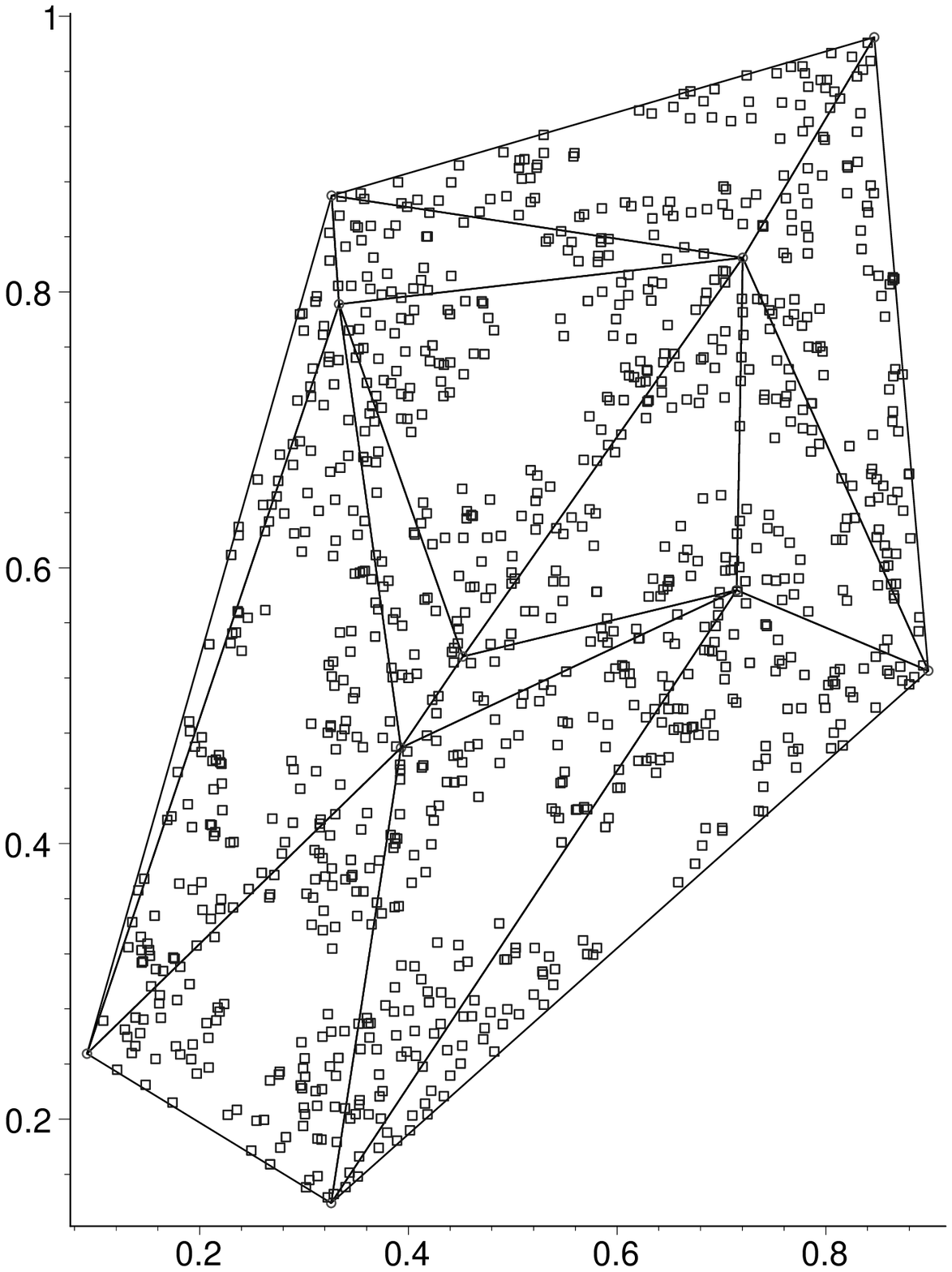} } }
\caption{\label{fig:deldata2}
Realization of segregation (left), $H_o$ (middle), and association (right) for $|\Y_m|=10$ and $n=1000$.
}
\end{figure}

Depicted in Figures \ref{fig:deldata1} and \ref{fig:deldata2} are
the realizations of 100 and 1000 observations, respectively,
independent identically distributed according to the segregation
with $\delta=1/16$, null, and association with $\delta=1/4$ (from
left to right) for $|\Y_m|=10$ and $J_{10}=13$.

With $n=100$ , for the null realization, the $p$-value is greater
than 0.1 for all $r$ except $r=1,\,4/3,\,\sqrt{2}$ for both
alternatives in the AND-underlying case, and for all $r$  values and
both alternatives in the OR-underlying case.
For the segregation
realization with $\delta=1/16$, we obtain $p < 0.018$ for all $r$
values except $r=1$ in the AND-underlying case and $p < 0.02$ for
all $r$ values in the OR-underlying case.
For the association
realization with $\delta=1/4$, we obtain $p < 0.043$ for $r=2,\,3$
in the AND-underlying case and $p < 0.05$ for
$r=4/3,\,\sqrt{2},\,1.5,\,2$ in the OR-underlying case.

With $n=1000$, in the AND-underlying case under the null distribution,
$p >.05$ for all $r$ values relative to segregation and association.
Under segregation with $\delta=1/16$,
$p<.01$ for all $r$ values considered.
Under association with $\delta=1/4$,
$p<.01$ for $r \in \{4/3,\sqrt{2},1.5,2,3,5 \}$ and
$p>.05$ for the other $r$ values considered.
In the OR-underlying case under the null distribution,
$p >.05$ for all $r$ values relative to segregation and association.
Under segregation with $\delta=1/16$,
$p<.01$ for $r \in \{1.1,1.2,4/3,\sqrt{2},1.5,2,3,5 \}$ and
$p>.05$ for the other $r$ values considered.
Under association with $\delta=1/4$,
$p<.01$ for $r \in \{1.1,1.2,4/3,\sqrt{2},1.5,2,3 \}$ and
$p>.05$ for the other $r$ values considered.

\begin{figure}[ht]
\centering
\psfrag{power}{\Huge{power}}
\psfrag{r}{\Huge{$r$}}
\rotatebox{-90}{ \resizebox{2.5 in}{!}{\includegraphics{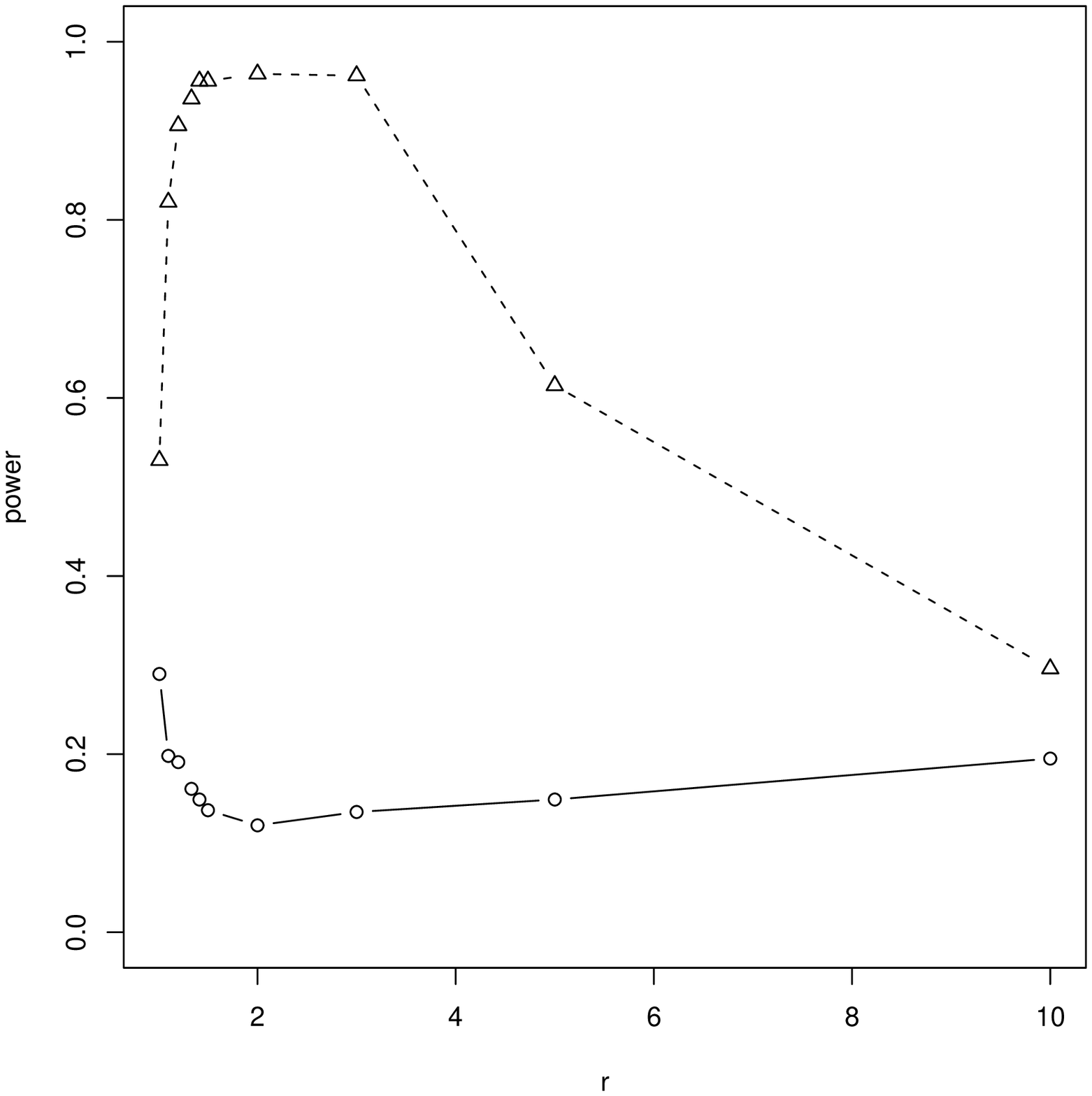} } }
\rotatebox{-90}{\resizebox{2.5 in}{!}{ \includegraphics{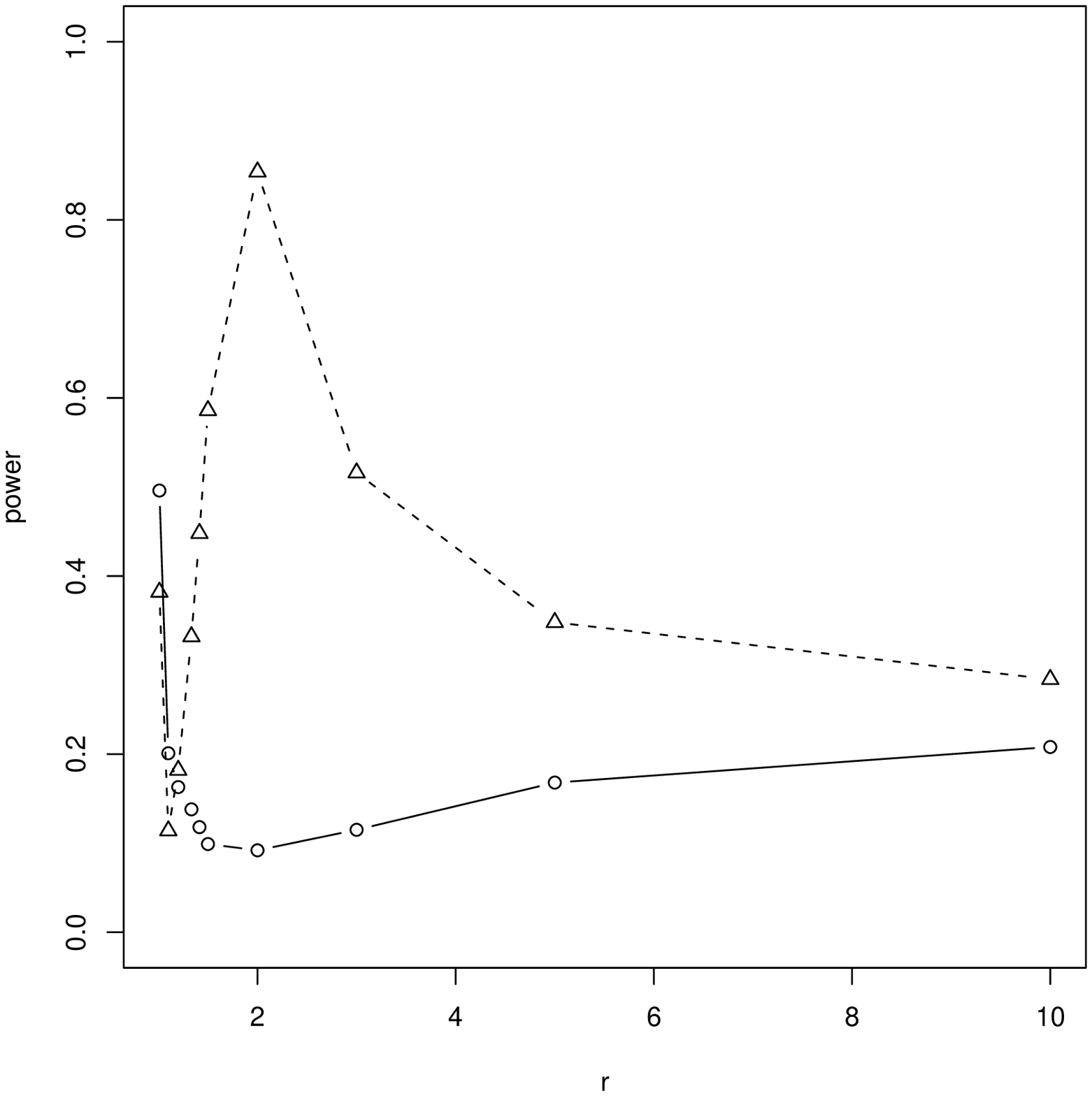} } }
\rotatebox{-90}{ \resizebox{2.5 in}{!}{\includegraphics{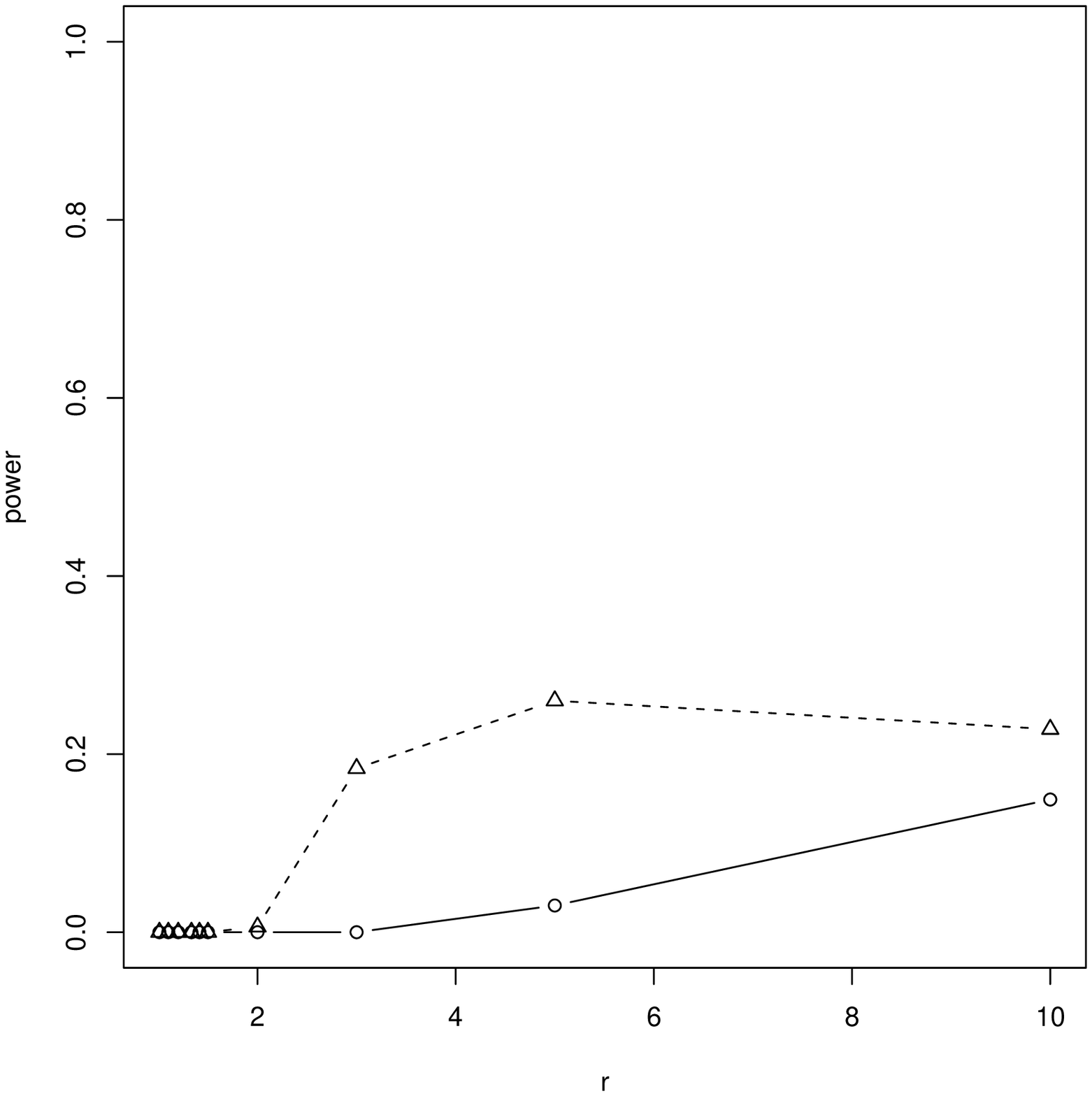} } }
\rotatebox{-90}{\resizebox{2.5 in}{!}{ \includegraphics{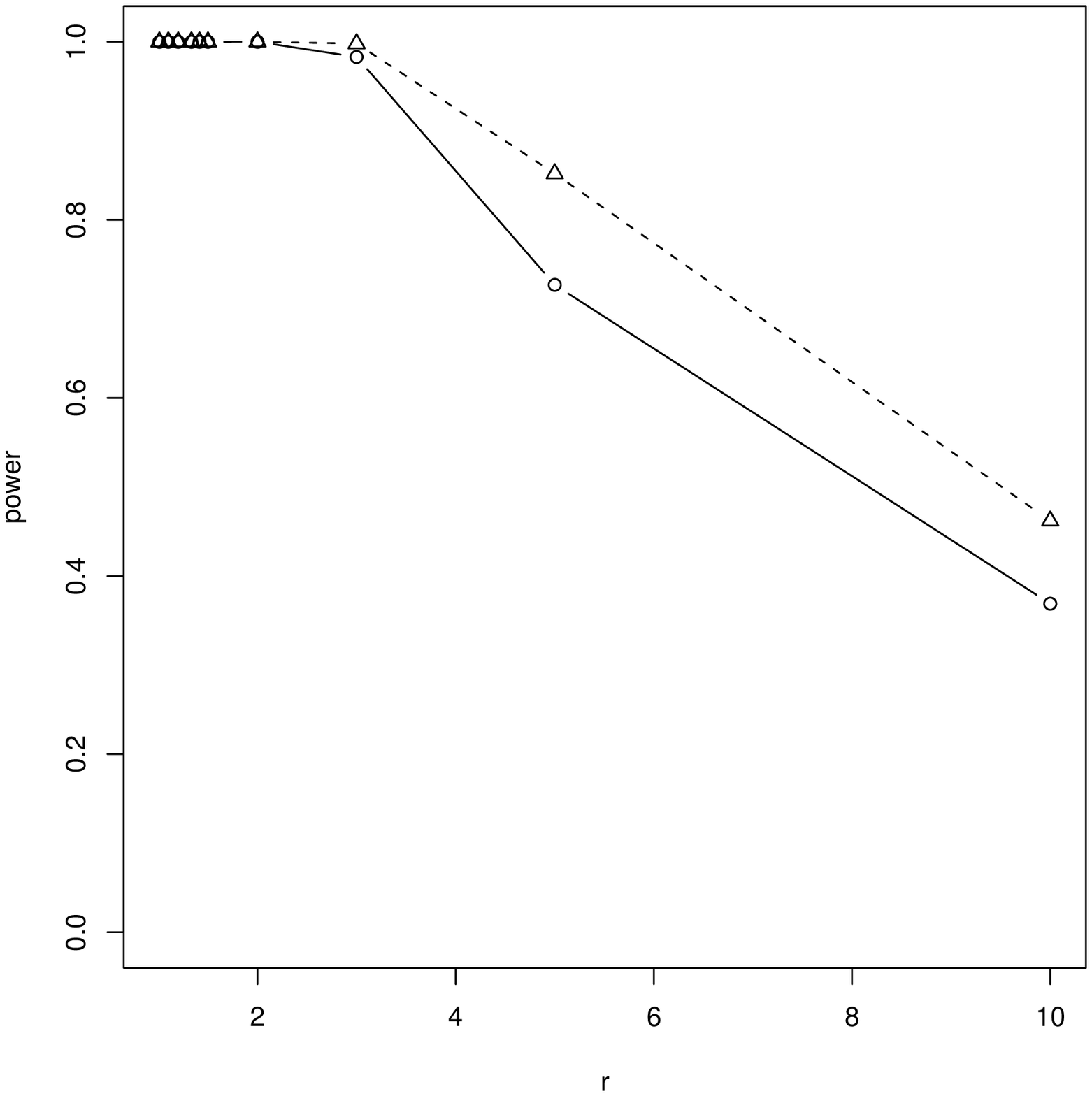} } }
\caption{
\label{fig:MT-asy-pow1}
The empirical size (circles joined with solid lines) and power estimates (triangles with dotted lines)
based on the asymptotic critical value for the AND-underlying case (top) and the
OR-underlying case (bottom) in the multiple triangle case, in both cases, $H^S_{\sqrt{3}/8}$
(left) and $H^A_{\sqrt{3}/12}$ (right) as a function of $r$, for $n=100$.
}
\end{figure}

\begin{figure}[ht]
\centering
\psfrag{power}{\Huge{power}}
\psfrag{r}{\Huge{$r$}}
\rotatebox{-90}{ \resizebox{2.5 in}{!}{ \includegraphics{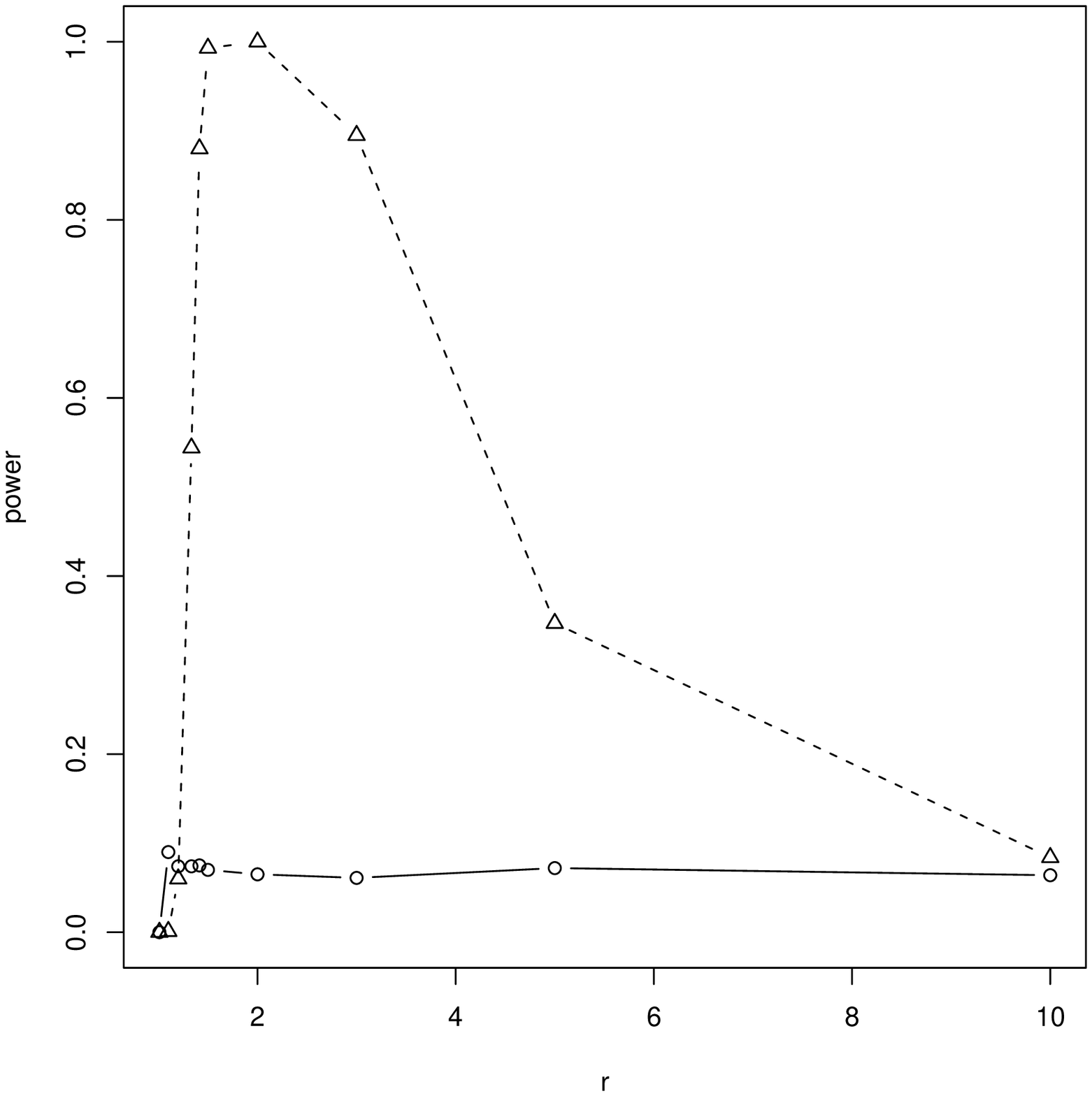} } }
\rotatebox{-90}{\resizebox{2.5 in}{!}{ \includegraphics{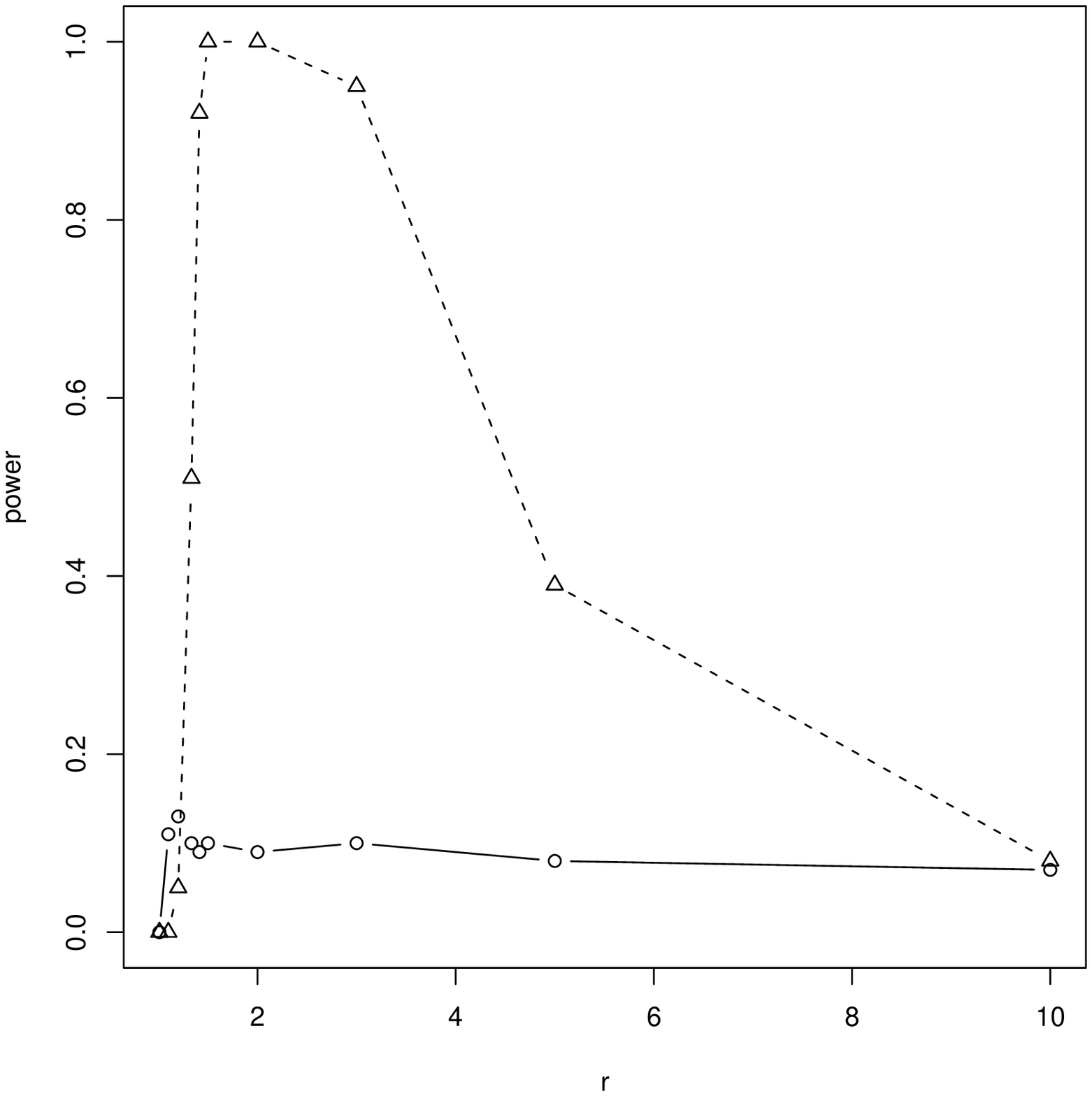} } }
\rotatebox{-90}{ \resizebox{2.5 in}{!}{\includegraphics{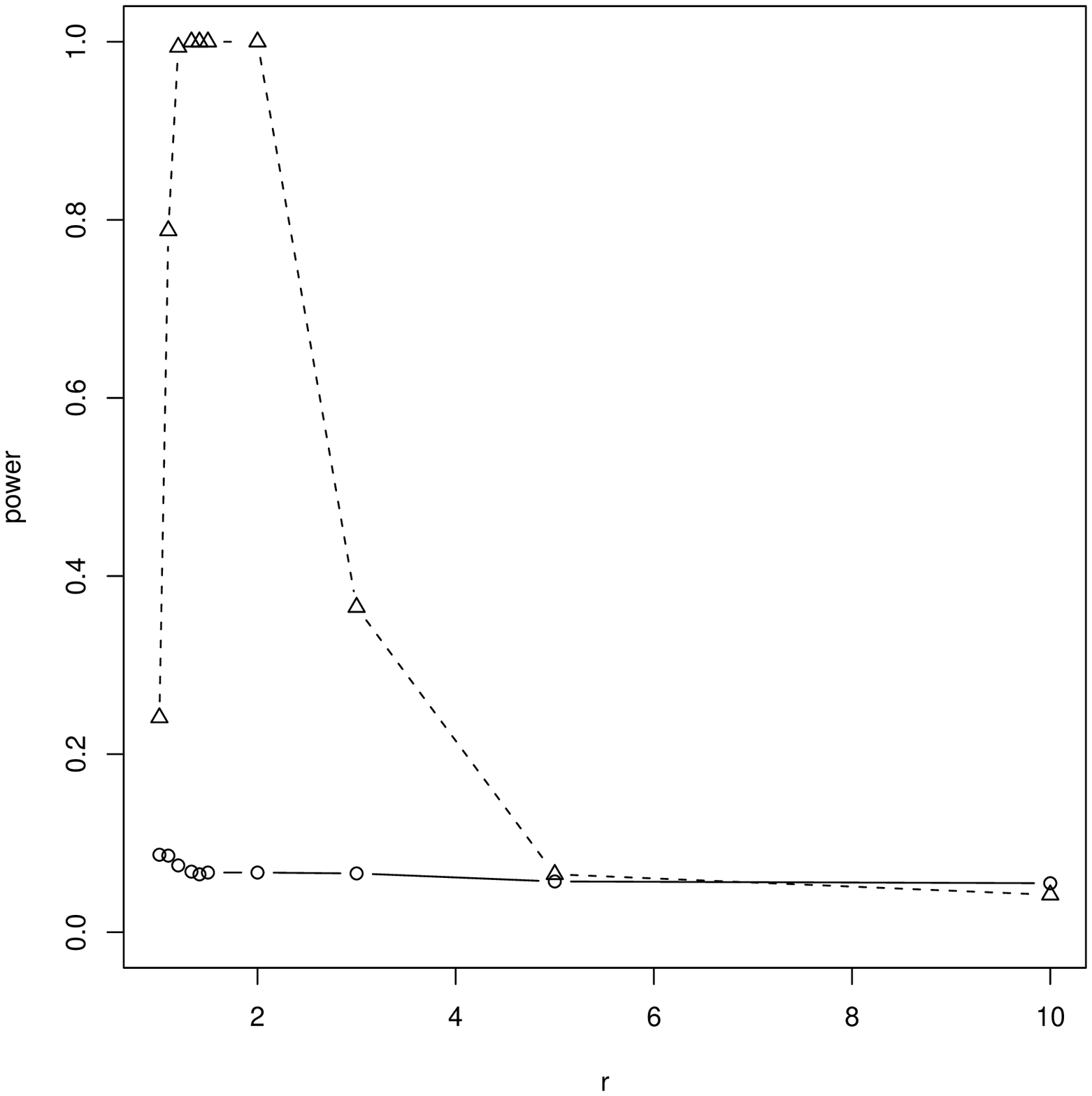} } }
\rotatebox{-90}{\resizebox{2.5 in}{!}{ \includegraphics{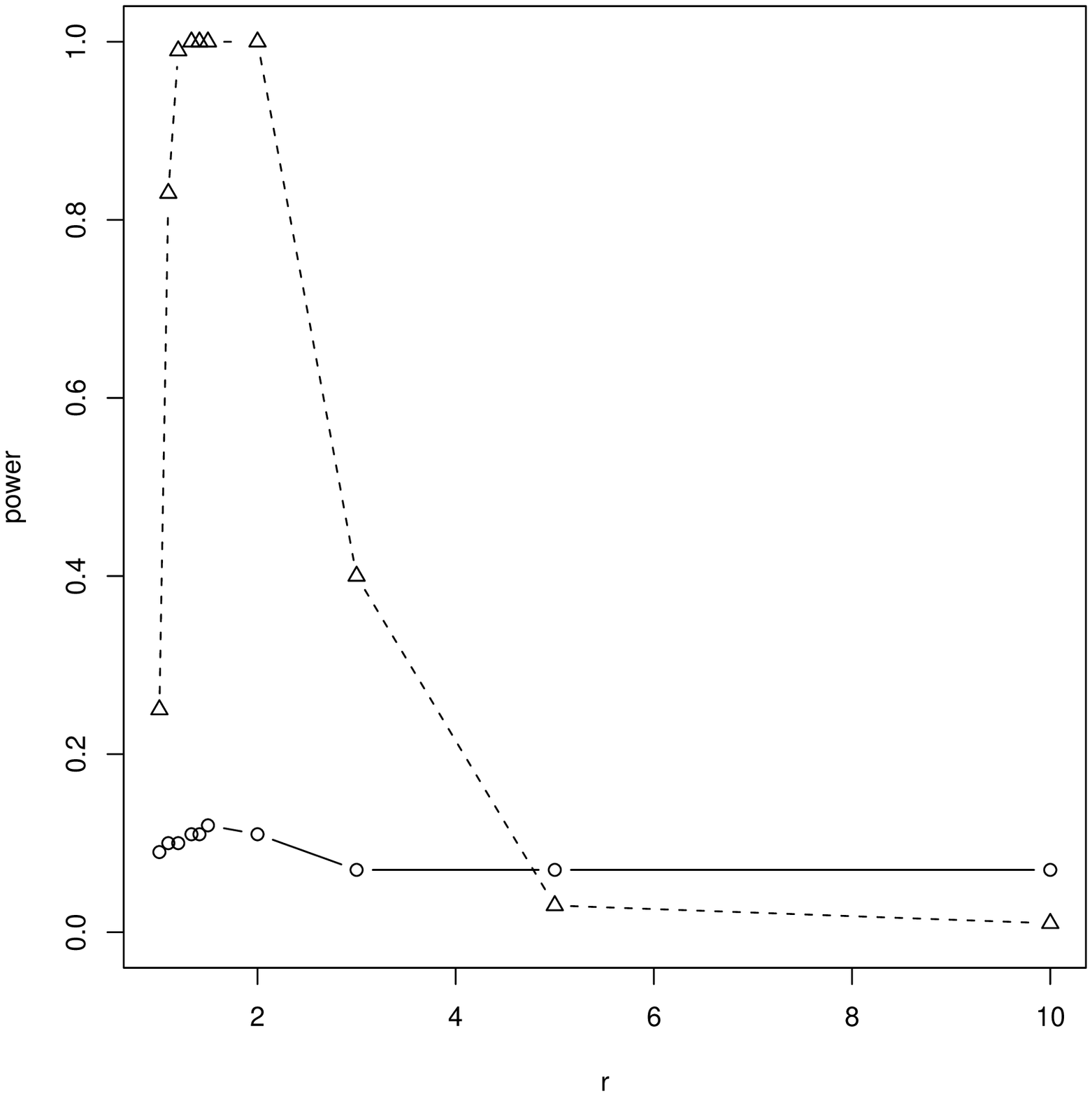} } }
\caption{
\label{fig:MT-asy-pow2}
The empirical size (circles joined with solid lines) and power estimates (triangles with dotted lines)
based on the asymptotic critical value for the AND-underlying case (top) and the
OR-underlying case (bottom) in the multiple triangle case, in both cases, $H^S_{\sqrt{3}/8}$
(left) and $H^A_{\sqrt{3}/12}$ (right) as a function of $r$, for $n=500$.
}
\end{figure}

We repeat the null realization $1000$ times for $n=100$ and find the estimated
significance level above $0.05$ for the AND-underlying case relative
to both alternatives with smallest being $0.12$ at $r=2$ relative to
segregation and $0.099$ at $r=2$ relative to association.
The associated empirical size and power estimates are presented in Figures \ref{fig:MT-asy-pow1} and \ref{fig:MT-asy-pow2}.
These results indicate that $n=100$ (i.e., the average number of points per triangle being about 8)
is not enough for the normal approximation in the AND-underlying case.
For the OR-underlying case the estimated significance level relative to segregation is
closest to $0.05$ is $0.03$ at $r=5$ and all much different at other $r$ values.
The estimated significance level relative to association
are larger than $0.25$ for all $r$ values.
Again the number of points per triangle is not large enough for normal approximation.
With $n=500$ (i.e., the average number of points per triangle being about 40),
the estimated significance levels get closer to $0.05$,
however they still are all above $0.05$,
hence for moderate sample sizes, the tests using the
relative edge density of the underlying graphs are liberal in rejecting $H_o$.
The empirical power analysis suggests the choice of $r=2$
---a moderate $r$ value---for both alternatives in both underlying cases.
Note also that AND-underlying case seems to perform better for segregation.

The PAE is given for $J_m=1$ in Section \ref{sec:PAE}.
For $J_m>1$, the analysis will depend both the number of triangles as well as the sizes of the triangles.
So the optimal $r$ values suggested for the $J_m=1$ case
does not necessarily hold for $J_m>1$,
so it needs to be updated, given the $\Y_m$ points.
The conditional test presented here is appropriate when the $\Y_m$ are fixed.
An unconditional version requires the joint distribution of the number and size of Delaunay triangles
when $\Y_m$ is, for instance, a Poisson point pattern.
Alas, this joint distribution is not available (\cite{okabe:2000}).

\subsection{Extension to Higher Dimensions}
\label{sec:extend-high-dim}
The extension to $\mathbb{R}^d$ for $d > 2$ is straightforward.
Let $\Y_{d+1} = \{\y_1,\y_2,\ldots,\y_{d+1}\}$ be $d+1$ non-coplanar points.
Denote the simplex formed by these $d+1$ points as $\mathfrak S (\Y_{d+1})$.
A simplex is the simplest polytope in $\mathbb{R}^d$
having $d+1$ vertices, $d\,(d+1)/2$ edges and $d+1$ faces of dimension $(d-1)$.
For $r \in [1,\infty]$, define the proportional-edge proximity map as follows.
Given a point $x$ in $\mathfrak S (\Y_{d+1})$,
let $y := \arg\min_{y \in \Y_{d+1}} \mbox{volume}(Q_y(x))$
where $Q_y(x)$ is the polytope with vertices being the $d\,(d+1)/2$ midpoints of the edges,
the vertex $\y$ and $x$.
That is, the vertex region for vertex $v$ is the polytope with vertices
given by $v$ and the midpoints of the edges.
Let $v(x)$ be the vertex in whose region $x$ falls.
(If $x$ falls on the boundary of two vertex regions or at the center of mass, we assign $v(x)$ arbitrarily.)
Let $\varphi(x)$ be the face opposite to vertex $v(x)$,
and $\eta(v(x),x)$ be the hyperplane parallel to $\varphi(x)$ which contains $x$.
Let $d(v(x),\eta(v(x),x))$ be the (perpendicular) Euclidean distance from $v(x)$ to $\eta(v(x),x)$.
For $r \in [1,\infty)$, let $\eta_r(v(x),x)$ be the hyperplane parallel to $\varphi(x)$
such that $d(v(x),\eta_r(v(x),x))=r\,d(v(x),\eta(v(x),x))$ and $d(\eta(v(x),x),\eta_r(v(x),x))< d(v(x),\eta_r(v(x),x))$.
Let $\mathfrak S_r(x)$ be the polytope similar to and with the same orientation as $\mathfrak S$ having $v(x)$ as a vertex and $\eta_r(v(x),x)$ as the opposite face.
Then the proportional-edge proximity region $\NPE^r(x):=\mathfrak S_r(x) \cap \mathfrak S(\Y_{d+1})$.
Furthermore, let $\zeta_i(x)$ be the hyperplane such that
$\zeta_i(x) \cap \mathfrak S(\Y_{d+1}) \not=\emptyset$ and $r\,d(\y_i,\zeta_i(x))=d(\y_i,\eta(\y_i,x))$
for $i=1,2,\ldots,d+1$.
Then $\G_1^r(x)\cap R(\y_i)=\{z \in R(\y_i): d(\y_i,\eta(\y_i,z)) \ge d(\y_i,\zeta_i(x)\}$, for $i=1,2,3$.
Hence  $\G_1^r(x)=\cup_{j=1}^{d+1} (\G_1^r(x)\cap R(\y_i))$.  Notice that $r \ge 1$ implies $x \in \NPE^r(x)$ and $x \in \G_1^r(x)$.

Theorem 1 generalizes,
so that any simplex $\mathfrak S$ in $\mathbb{R}^d$
can be transformed into a regular polytope
(with edges being equal in length and faces being equal in volume)
preserving uniformity.
Delaunay triangulation becomes Delaunay tessellation in $\mathbb{R}^d$,
provided no more than 4 points being cospherical
(lying on the boundary of the same sphere).
In particular, with $d=3$, the general simplex is a tetrahedron
(4 vertices, 4 triangular faces and 6 edges),
which can be mapped into a regular tetrahedron
(4 faces are equilateral triangles) with vertices
$(0,0,0)\,(1,0,0)\,(1/2,\sqrt{3}/2,0),\,(1/2,\sqrt{3}/4,\sqrt{3}/2)$.

Asymptotic normality of the $U$-statistic and consistency of the tests hold for $d>2$ in both underlying cases.

\section{Discussion and Conclusions}
\label{sec:discussion}
In this article, we consider the asymptotic distribution of the
relative edge density of the underlying graphs based on
(parametrized) proportional-edge proximity catch digraphs (PCDs),
for testing bivariate spatial point patterns of segregation and association.
To our knowledge the PCD-based methods are the only graph theoretic
methods for testing spatial patterns in literature
(\cite{ceyhan:dom-num-NPE-SPL}, \cite{ceyhan:arc-density-PE}, and \cite{ceyhan:arc-density-CS}).
The proportional-edge PCDs lend themselves for such a purpose,
because of the geometry invariance property for uniform data on Delaunay triangles.
Let the two samples of sizes $n$ and $m$ be from classes $\X$ and $\Y$, respectively,
with $\X$ points being used as the vertices of the PCDs and $\Y$ points
being used in the construction of Delaunay triangulation.
For the relative density approach to be appropriate,
$n$ should be much larger compared to $m$.
This implies that $n$ tends to infinity while $m$ is assumed to be fixed.
That is, the difference in the relative abundance of the two classes
should be large for this method.
Such an imbalance usually confounds the
results of other spatial interaction tests.
Furthermore,
we can perform Monte Carlo randomization to remove the conditioning on $\Y_m$.

Previously, \cite{ceyhan:arc-density-PE} employed the relative (arc) density
of the proportional-edge PCDs for testing bivariate spatial patterns.
In this work, we consider the AND- and OR-underlying graphs based on this PCD;
in particular, we demonstrate that relative edge density of these underlying PCDs
is a $U$-statistic, and employing asymptotic normality of $U$-statistics,
we derive the asymptotic distribution of the relative edge density.
We then use relative edge density as a test statistic for testing segregation and association.

The null hypothesis is assumed to be CSR of $\X$ points,
i.e., the uniformness of $\X$ points in the convex hull of $\Y$ points.
Although we have two classes here, the null pattern is not the CSR independence,
since for finite $m$,
we condition on $m$ and the areas of the Delaunay triangles based on $\Y$ points
as long as they are not co-circular.

There are many types of parametrizations for the alternatives.
The particular parametrization of the alternatives in this article
is chosen so that the distribution of the relative edge density under the alternatives
would be geometry invariant
(i.e., independent of the geometry of the support triangles).
The more natural alternatives (i.e.,
the alternatives that are more likely to be found in practice)
can be similar to or might be approximated by our parametrization.
Because in any segregation alternative,
the $\X$ points will tend to be further away from $\Y$ points
and in any association alternative $\X$ points will tend to cluster around the $\Y$ points.
And such patterns can be detected by the test
statistics based on the relative edge density,
since under segregation (whether it is parametrized as in Section \ref{sec:alt-seg-assoc}
or not) we expect them to be larger,
and under association (regardless of the parametrization)
they tend to be smaller.

Our Monte Carlo simulation analysis and asymptotic efficiency analysis based on Pitman
asymptotic efficiency reveals that AND-underlying graph has
better power performance against segregation compared to the digraph and OR-underlying version.
On the other hand, OR-underlying graph has better
power performance against association compared to the digraph
and AND-underlying version.
When the number of $\X$ points per triangle is less than 30,
we recommend the use Monte Carlo randomization, otherwise
we recommend the use of normal approximation as $n \rightarrow \infty$.
Furthermore,
when testing against segregation we recommend the parameter $r \approx 2$,
while for testing against association we recommend the parameters $r \in (2,3)$
as they exhibit the better performance in terms of size and power.

\section*{Acknowledgments}
This work was partially sponsored by the Defense Advanced Research Projects Agency as administered
by the Air Force Office of Scientific Research under contract DOD F49620-99-1-0213 and
by Office of Naval Research Grant N00014-95-1-0777 and by TUBITAK Kariyer Project Grant 107T647.


\clearpage
\section*{APPENDIX}

\subsection*{Appendix 1: The Variance of Relative Edge Density for the AND-Underlying Graph Version:}
\begin{figure}
\centering
\psfrag{r}{\normalsize{$r$}}
\epsfig{figure=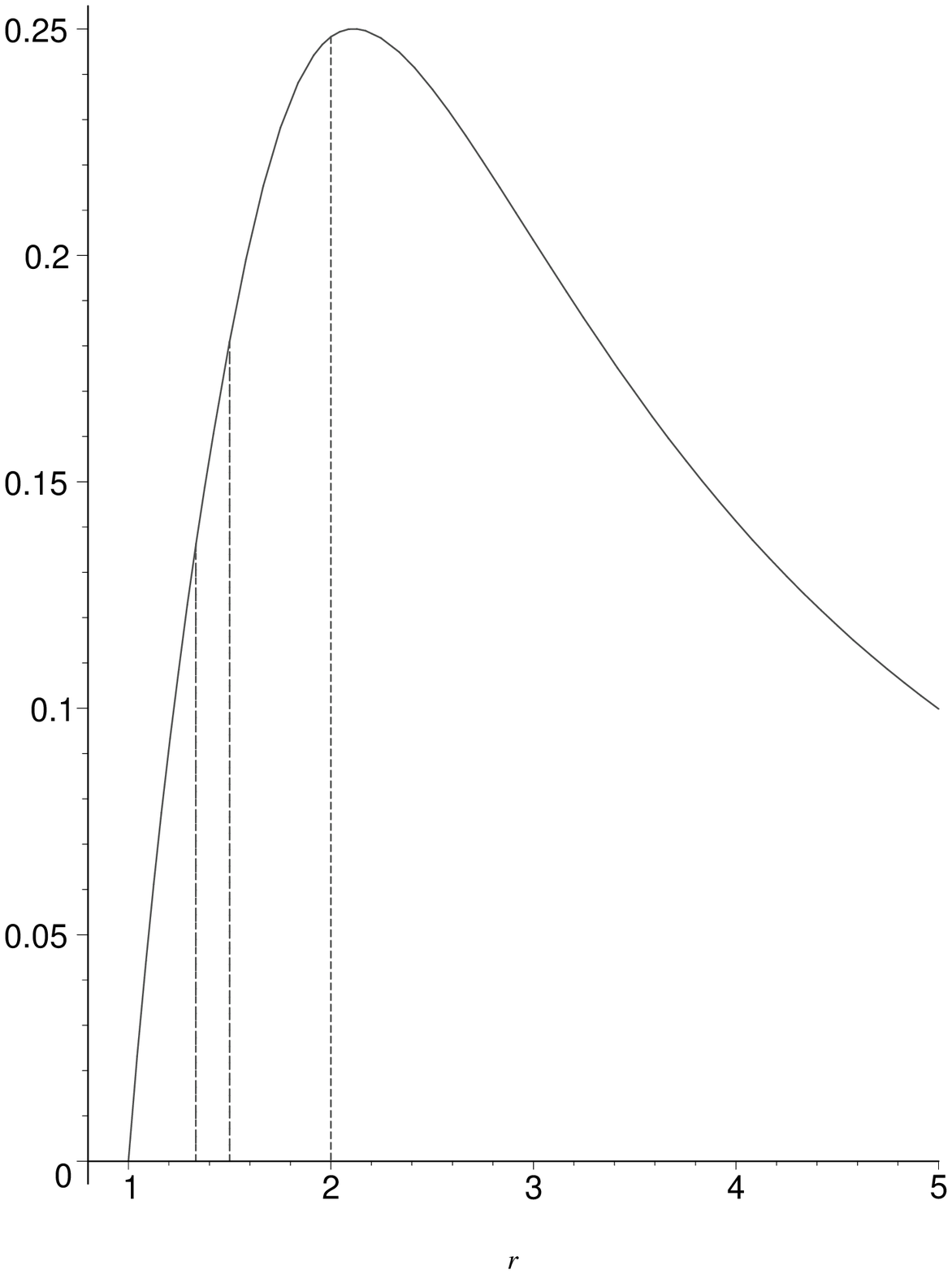, height=150pt, width=200pt}
\epsfig{figure=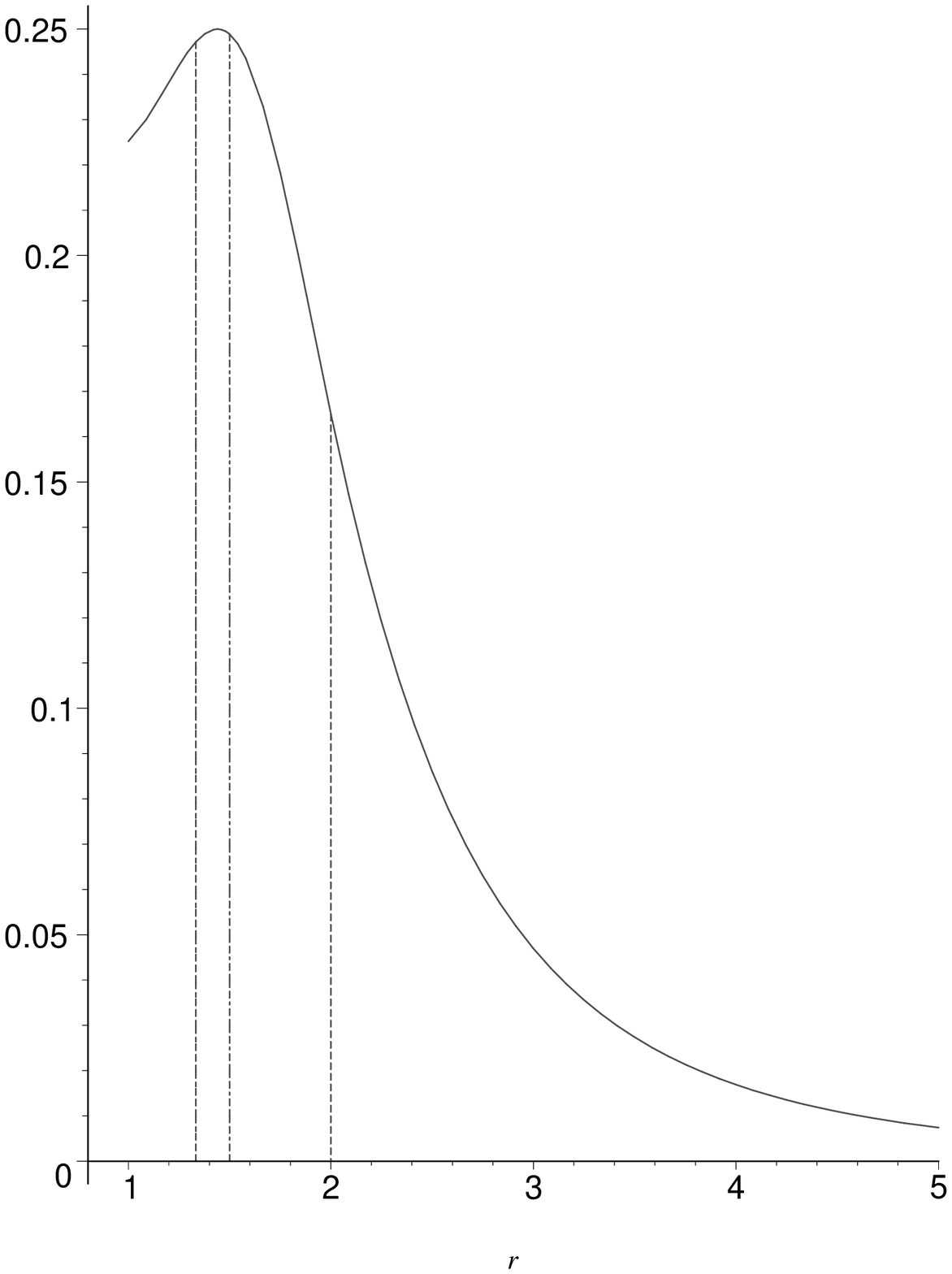, height=150pt, width=200pt}\\
\caption{
\label{var of nu under-graph}
$\Var\left[ h^{\la}_{12}(r) \right]$ (left)
and
$\Var\left[ h^{\lo}_{12}(r) \right]$ (right)
as a function of $r$ for $r \in [1,5] $.}
\end{figure}

The variance term is
$$
\Var\left[ h^{\la}_{12}(r) \right]=
\varphi^{\la}_{1,1}(r)\I(r \in [1,4/3))+
\varphi^{\la}_{1,2}(r)\I(r \in [4/3,3/2))+
\varphi^{\la}_{1,3}(r)\I(r \in [3/2,2))+
\varphi^{\la}_{1,4}(r)\I(r \in [2,\infty))
$$
where
$\varphi^{\la}_{1,1}(r)=-{\frac{(5\,r^6-153\,r^5+393\,r^4-423\,r^3-54\,r^2+360\,r-128)(447\,r^4-261\,r^3+54\,r^2+5\,r^6-153\,r^5+360\,r-128)}{2916\,r^4(r+2)^2(r+1)^2}}$,\\
$\varphi^{\la}_{1,2}(r)=-{\frac{(101\,r^5-801\,r^4+1302\,r^3-732\,r^2-536\,r+672)(1518\,r^3-84\,r^2-104\,r+101\,r^5-801\,r^4+672)}{46656\,r^2(r+2)^2(r+1)^2}}$,\\
$\varphi^{\la}_{1,3}(r)=-{\frac{(r^8-13\,r^7+30\,r^6+148\,r^5-448\,r^4+264\,r^3+288\,r^2-368\,r+96)(22\,r^6+124\,r^5-464\,r^4+r^8-13\,r^7+264\,r^3+288\,r^2-368\,r+96)}{64\,r^8(r+2)^2(r+1)^2}}$,\\
$\varphi^{\la}_{1,4}(r)={\frac{(r^5+r^4-3\,r^3-3\,r^2+6\,r-2)(3\,r^3+3\,r^2-6\,r+2)}{r^8(r+1)^2}}$.
See Figure \ref{var of nu under-graph}.
Note that $\Var_{\la}(r=1)=0$ and $\lim_{r \rightarrow \infty}\Var_{\la}(r)=0$ (at rate $O(r^{-2})$),
and $\argsup_{r \in [1,\infty)} \Var_{\la}(r) \approx 2.1126$ with $\sup \Var_{\la}(r)=.25$.

\begin{figure}
\centering
\psfrag{r}{\normalsize{$r$}}
\epsfig{figure=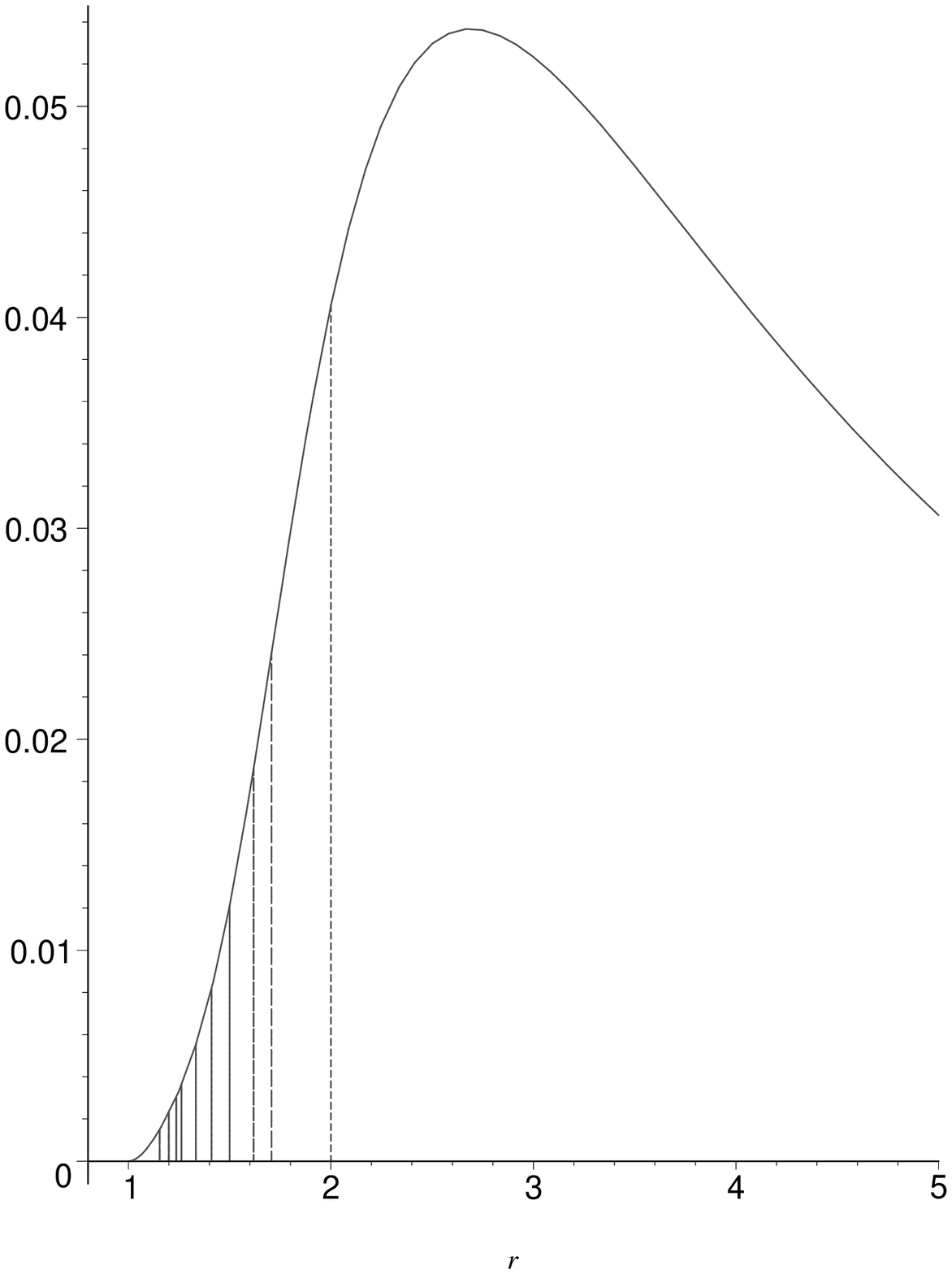, height=150pt, width=200pt}
\epsfig{figure=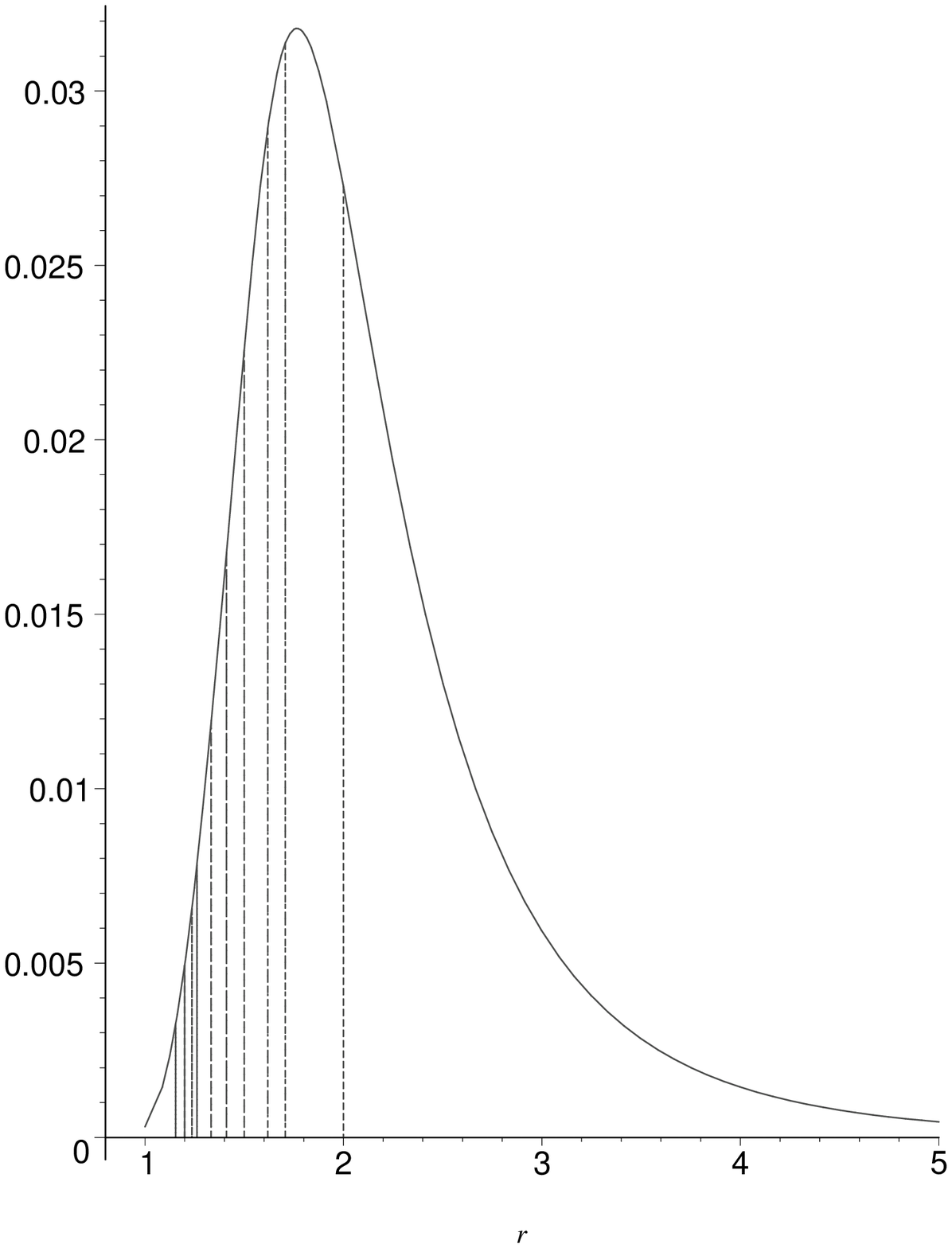, height=150pt, width=200pt}\\
\caption{
\label{cov of nu under-graph}
$\nu_{\la}(r)=\Cov\left[ h^{\la}_{12}(r),h^{\la}_{13}(r) \right]$ (left)
and
$\nu_{\lo}(r)=\Cov\left[ h^{\lo}_{12}(r),h^{\lo}_{13}(r) \right]$ (right)
as a function of $r$ for $r \in [1,5] $.}
\end{figure}

Moreover,
$$\nu_{\la}(r):=\Cov\left[ h^{\la}_{12}(r),h^{\la}_{13}(r) \right]=\sum_{i=1}^{11}\vartheta^{\la}_i(r)\,\I(\mI_i)$$
where
{\small
\begin{multline*}
\vartheta^{\la}_1(r)=-\frac{1}{58320\,(2\,r^2+1)(r+2)^2(r+1)^3r^6}\, ((r-1)^2(972\,r^{19}+8748\,r^{18}+44456\,r^{17}+140328\,r^{16}+121371\,r^{15}\\
-412117\,r^{14}-27145\,r^{13}-4503501\,r^{12}+1336147\,r^{11}+10640999\,r^{10}-982009\,r^9-6677105\,r^8-2274458\,r^7\\
-1150162\,r^6+249126\,r^5+1232530\,r^4+1234372\,r^3+226776\,r^2-184944\,r-81920))
\end{multline*}
\begin{multline*}
\vartheta^{\la}_2(r)=-\frac{1}{116640\,(2\,r^2+1)(r+2)^2(r+1)^3r^6}\,(486\,r^{21}+3402\,r^{20}-269\,r^{19}-45155\,r^{18}-118850\,r^{17}+443518\,r^{16}\\
+3251855\,r^{15}-13836295\,r^{14}+13434672\,r^{13}+11140788\,r^{12}-27667544\,r^{11}+13293088\,r^{10}+7159710\,r^9-\\
13013598\,r^8 +4185440\,r^7+3262952\,r^6+586636\,r^5-1616444\,r^4-680120\,r^3-55952\,r^2+219936\,r+49152)
\end{multline*}
\begin{multline*}
\vartheta^{\la}_3(r)=-\frac{1}{116640\,(2\,r^2+1)(r+2)^2(r+1)^3r^6}\,(486\,r^{21}+3402\,r^{20}-269\,r^{19}-45155\,r^{18}-118850\,r^{17}+443518\,r^{16}\\
+2751855\,r^{15}-13736295\,r^{14}+18084672\,r^{13}+8770788\,r^{12}-43009544\,r^{11}+24604048\,r^{10}+27137438\,r^9-30889822\,r^8\\
-2832544\,r^7+11101160\,r^6-4168820\,r^5+2364868\,r^4+2305864\,r^3-3041936\,r^2+219936\,r+49152)
\end{multline*}
\begin{multline*}
\vartheta^{\la}_4(r)=-\frac{1}{58320\,(r+2)^3(r^2-2)(2\,r^2+1)(r+1)^3r^6}\,(3632\,r^{22}+25632\,r^{21}-60328\,r^{20}-441888\,r^{19}+1353430\,r^{18}\\
-297666\,r^{17}-4791125\,r^{16}+12849927\,r^{15}-10894618\,r^{14}-26295324\,r^{13}+62283823\,r^{12}-2280753\,r^{11}-81700012\,r^{10}\\
+32551926\,r^9+39974410\,r^8-11284026\,r^7-5806580\,r^6-9167580\,r^5-2004944\,r^4+4646688\,r^3+1931776\,r^2-489024\,r-98304)
\end{multline*}
\begin{multline*}
\vartheta^{\la}_5(r)=\vartheta^{\la}_6(r)=-\frac{1}{58320\,(r+2)^3(2\,r^2+1)(r^2+1)(r+1)^3r^6}\,(3632\,r^{22}+25632\,r^{21}-49432\,r^{20}-364992\,r^{19}+958940\,r^{18}\\
-1167012\,r^{17}+1200518\,r^{16}+5424126\,r^{15}-23566328\,r^{14}+23837088\,r^{13}+11797395\,r^{12}-41623065\,r^{11}+39261953\,r^{10}\\
-8239197\,r^9-30178496\,r^8+27901506\,r^7-4936170\,r^6+61038\,r^5+4719720\,r^4-5513952\,r^3+340736\,r^2+23328\,r+65536)
\end{multline*}
\begin{multline*}
\vartheta^{\la}_7(r)=\frac{1}{466560\,(r+2)^3(2\,r^2+1)(r^2+1)(r+1)^3r^5}\,(1562\,r^{21}-11142\,r^{20}-103099\,r^{19}+2105697\,r^{18}-9774118\,r^{17}+\\
10220280\,r^{16}+27825711\,r^{15}-69243129\,r^{14}+81624200\,r^{13}-76052574\,r^{12}-65530400\,r^{11}+262451196\,r^{10}-178092280\,r^9\\
-69106464\,r^8+158439568\,r^7-97568688\,r^6+12246288\,r^5+17591952\,r^4-21111616\,r^3+15628032\,r^2-2545664\,r+993024)
\end{multline*}
\begin{multline*}
\vartheta^{\la}_8(r)=-\frac{1}{1920\,(r+2)^3(r^2+1)(2\,r^2+1)(r+1)^3r^{10}}\,(2\,r^{26}-30\,r^{25}-2395\,r^{23}+281\,r^{24}+8770\,r^{22}
+29528\,r^{21}-268053\,r^{20}+\\
245667\,r^{19}+2066216\,r^{18}-5313494\,r^{17}-1589216\,r^{16}+18512684\,r^{15}-18946136\,r^{14}-2665248\,r^{13}+22789584\,r^{12}-\\
32987760\,r^{11}+20482512\,r^{10}+13109584\,r^9-28084416\,r^8+17326976\,r^7-3864576\,r^6-4579328\,r^5+6666240\,r^4-3576320\,r^3\\
+635904\,r^2-116736\,r+61440)
\end{multline*}
\begin{multline*}
\vartheta^{\la}_9(r)=-\frac{1}{1920\,(r+2)^3(r^2+1)(2\,r^2+1)(r+1)^3r^{10}}\,(2\,r^{26}-30\,r^{25}-2395\,r^{23}281\,r^{24}+8258\,r^{22}
+31064\,r^{21}-262677\,r^{20}+\\
225443\,r^{19}+2052136\,r^{18}-5219030\,r^{17}-1608928\,r^{16}+18337836\,r^{15}-18837080\,r^{14}-2598688\,r^{13}+22736336\,r^{12}-\\
32858736\,r^{11}+20384720\,r^{10}+12930896\,r^9-27988416\,r^8+17416832\,r^7-3862784\,r^6-4575488\,r^5+6638848\,r^4-3603200\,r^3\\
+640512\,r^2-107520\,r+63488)
\end{multline*}
\begin{multline*}
\vartheta^{\la}_{10}(r)=-\frac{1}{1920\,(r+2)^3(r-1)(r+1)^3(2\,r^2-1)r^{10}}\,(2\,r^{25}+307\,r^{23}-32\,r^{24}-2612\,r^{22}
+11572\,r^{21}+21934\,r^{20}-328867\,r^{19}+\\
524994\,r^{18}+2446870\,r^{17}-8676180\,r^{16}-437020\,r^{15}+36944680\,r^{14}-40677696\,r^{13}-44860384\,r^{12}+106256352\,r^{11}-\\
15515040\,r^{10}-98636848\,r^9+66358080\,r^8+27142272\,r^7-42614272\,r^6+7781120\,r^5+7327232\,r^4-3388672\,r^3+\\
430592\,r^2-171008\,r+63488)
\end{multline*}
\begin{multline*}
\vartheta^{\la}_{11}(r)=\frac{1}{15\,(2\,r^2-1)(r+1)^3r^{10}}\,(30\,r^{13}+90\,r^{12}-127\,r^{11}-621\,r^{10}+320\,r^9+1568\,r^8-858\,r^7
-1370\,r^6+909\,r^5+\\
295\,r^4-292\,r^3+44\,r^2+6\,r-2)
\end{multline*}
}
\noindent
and
$\mI_1=[1,2/\sqrt{3}),\;
\mI_2=[2/\sqrt{3},6/5),\;
\mI_3=[6/5,\sqrt{5}-1),\;
\mI_4=[\sqrt{5}-1,(6+2\,\sqrt{2})/7),\;
\mI_5=[(6+2\,\sqrt{2})/7,4/3),\;
\mI_6=[4/3,(6+\sqrt{15})/7),\;
\mI_7=[(6+\sqrt{15})/7,3/2),\;
\mI_8=[3/2,(1+\sqrt{5})/2),\;
\mI_9=[(1+\sqrt{5})/2,1+1/\sqrt{2}),\;
\mI_{10}=[1+1/\sqrt{2},2),\;
\mI_{11}=[2,\infty)$.
See Figure \ref{cov of nu under-graph}.
Note that $\Cov_{\la}(r=1)=0$ and $\lim_{r \rightarrow \infty}\nu_{\la}(r)=0$ (at rate $O(r^{-2})$),
and $\argsup_{r \in [1,\infty)} \nu_{\la}(r)\approx 2.69$ with $\sup \nu_{\la}(r) \approx .0537$.

\subsection*{Appendix 2: The Variance of Relative Edge Density for the OR-Underlying Graph Version:}
The variance term is
$$
\Var\left[ h^{\lo}_{12}(r) \right]=
\varphi^{\lo}_{1,1}(r)\I(r \in [1,4/3))+
\varphi^{\lo}_{1,2}(r)\I(r \in [4/3,3/2))+
\varphi^{\lo}_{1,3}(r)\I(r \in [3/2,2))+
\varphi^{\lo}_{1,4}(r)\I(r \in [2,\infty))
$$
where
$\varphi^{\lo}_{1,1}(r)=-\frac{(47\,r^6-195\,r^5+860\,r^4-846\,r^3-108\,r^2+720\,r-256)(752\,r^4-1170\,r^3-324\,r^2+47\,r^6-195\,r^5+720\,r-256)}{11664\,r^4(r+2)^2(r+1)^2}$,\\
$\varphi^{\lo}_{1,2}(r)=-\frac{(175\,r^5-579\,r^4+1450\,r^3-732\,r^2-536\,r+672)(1234\,r^3-1380\,r^2-968\,r+175\,r^5-579\,r^4+672)}{46656\,r^2(r+2)^2(r+1)^2}$,\\
$\varphi^{\lo}_{1,3}(r)=-\frac{(3\,r^8-7\,r^7-30\,r^6+84\,r^5-264\,r^4+304\,r^3+144\,r^2-368\,r+96)(-22\,r^6+108\,r^5-248\,r^4+3\,r^8-7\,r^7+304\,r^3+144\,r^2-368\,r+96)}{64\,r^8(r+2)^2(r+1)^2}$,\\
\\
$\varphi^{\lo}_{1,4}(r)=2\,\frac{(r^5+r^4-6\,r+2)(3\,r-1)}{r^8(r+1)^2}$.
See Figure \ref{var of nu under-graph}.

Note that $\Var_{\lo}(r=1)=2627/11664$ and $\lim_{r \rightarrow \infty}\Var_{\lo}(r)=0$ (at rate $O(r^{-4})$),
and $\argsup_{r \in [1,\infty)} \Var_{\lo}(r) \approx 1.44$ with $\sup \Var_{\lo}(r) \approx .25$.

Moreover,
$$\nu_{\lo}(r):=\Cov[h^{\lo}_{12}(r),h^{\lo}_{13}(r)]=\sum_{i=1}^{11}\vartheta^{\lo}_i(r)\,\I(\mI_i)$$
where
{\small
\begin{multline*}
\vartheta^{\lo}_1(r)=-\frac{1}{58320\,(r^2+1)(2\,r^2+1)(r+1)^3(r+2)^3r^6}\,(1458\,r^{22}+13122\,r^{21}+50731\,r^{20}-84225\,r^{19}
-19193\,r^{18}-1823223\,r^{17}+\\
5576151\,r^{16}+2978697\,r^{15}-33432692\,r^{14}+37427862\,r^{13}+15883834\,r^{12}-60944766\,r^{11}+49876417\,r^{10}-1754523\,r^9-\\
36606859\,r^8+32338215\,r^7-10290256\,r^6-2234754\,r^5+7085471\,r^4-5608569\,r^3+1645826\,r^2-132876\,r+30824)
\end{multline*}
\begin{multline*}
\vartheta^{\lo}_2(r)=\vartheta^{\lo}_3(r)=-\frac{1}{116640\,(r^2+1)(2\,r^2+1)(r+1)^3(r+2)^3r^6}\,(1458\,r^{22}+13122\,r^{21}+62825\,r^{20}
-175011\,r^{19}+156014\,r^{18}-\\
3300900\,r^{17}+11053023\,r^{16}+5055135\,r^{15}-67685050\,r^{14}+75243552\,r^{13}+33155180\,r^{12}-120628524\,r^{11}+99831906\,r^{10}-\\
4883958\,r^9-74801558\,r^8+64360782\,r^7-19812000\,r^6-3667716\,r^5+14541630\,r^4-11254002\,r^3+3070468\,r^2-413208\,r+28880)
\end{multline*}
\begin{multline*}
\vartheta^{\lo}_4(r)=-\frac{1}{58320\,(r^2+1)(2\,r^2+1)(r^2-2)(r+2)^3(r+1)^3r^6}\,(972\,r^{24}+8748\,r^{23}+29590\,r^{22}-149106\,r^{21}-36820\,r^{20}-\\
986280\,r^{19}+5942884\,r^{18}+2883672\,r^{17}-47189711\,r^{16}+43450125\,r^{15}+85975304\,r^{14}-156173934\,r^{13}+27378901\,r^{12}+123606417\,r^{11}\\
-152209261\,r^{10}+64653597\,r^9+56621894\,r^8-88962768\,r^7+43754559\,r^6-5940597\,r^5-13006396\,r^4+17019366\,r^3-7037340\,r^2+\\
413208\,r-28880)
\end{multline*}
\begin{multline*}
\vartheta^{\lo}_5(r)=-\frac{1}{58320\,(r^2+1)(2\,r^2+1)(r+1)^3(r+2)^3r^6}\,(972\,r^{22}+8748\,r^{21}+31534\,r^{20}-131610\,r^{19}
+261546\,r^{18}-1552026\,r^{17}+\\
3745643\,r^{16}+4573731\,r^{15}-29416804\,r^{14}+26163354\,r^{13}+19600850\,r^{12}-43126062\,r^{11}+31497249\,r^{10}-7381467\,r^9-\\
22237963\,r^8+26778663\,r^7-9107024\,r^6-115074\,r^5+3136927\,r^4-5055609\,r^3+2292994\,r^2+14580\,r-1944)
\end{multline*}
\begin{multline*}
\vartheta^{\lo}_6(r)=\frac{1}{233280\,(r^2+1)(2\,r^2+1)(r+1)^3(r+2)^3r^6}\,(486\,r^{22}-7290\,r^{21}-181459\,r^{20}+1024401\,r^{19}
-2691213\,r^{18}+3921057\,r^{17}+\\
1844321\,r^{16}-33347697\,r^{15}+80028903\,r^{14}-29292735\,r^{13}-98093906\,r^{12}+125034492\,r^{11}-46658244\,r^{10}-57216612\,r^9+\\
88057996\,r^8-26383068\,r^7-12851392\,r^6+14179848\,r^5-8656508\,r^4+1593828\,r^3+134136\,r^2-58320\,r+7776)
\end{multline*}
\begin{multline*}
\vartheta^{\lo}_7(r)=\frac{1}{233280\,(r+2)^3(r^2+1)(2\,r^2+1)(r+1)^3(r-1)r^6}\,(486\,r^{23}-7776\,r^{22}-174169\,r^{21}
+1205860\,r^{20}-4656806\,r^{19}+\\
8763566\,r^{18}+7460036\,r^{17}-63559490\,r^{16}+91134324\,r^{15}+18516450\,r^{14}-122708655\,r^{13}+18577230\,r^{12}+80410332\,r^{11}-19357704\,r^{10}-\\
39129236\,r^9+75311048\,r^8-77449360\,r^7+4053376\,r^6+48283912\,r^5-40690240\,r^4+17736336\,r^3-4315680\,r^2+544320\,r-31104)
\end{multline*}
\begin{multline*}
\vartheta^{\lo}_8(r)=\frac{1}{960\,(r+2)^3(r^2+1)(2\,r^2+1)(r+1)^3r^8}\,(2\,r^{24}-30\,r^{23}-161\,r^{22}+107\,r^{21}+4137\,r^{20}
-10685\,r^{19}+8367\,r^{18}+\\
78713\,r^{17}-450859\,r^{16}+697707\,r^{15}+517846\,r^{14}-3723120\,r^{13}
+6565124\,r^{12}-1468692\,r^{11}-8695792\,r^{10}+9535720\,r^9-\\
6773160\,r^8+526744\,r^7+10691376\,r^6-7797264\,r^5+1137696\,r^4+523712\,r^3-2687872\,r^2+1701888\,r-245760)
\end{multline*}
\begin{multline*}
\vartheta^{\lo}_9(r)=\frac{1}{960\,(2\,r^2+1)(r+1)^2(r+2)^3(r^2+1)r^{10}}\,(2\,r^{25}-32\,r^{24}-129\,r^{23}+236\,r^{22}+4157\,r^{21}
-15610\,r^{20}+21289\,r^{19}+\\
67536\,r^{18}-511355\,r^{17}+1161830\,r^{16}-634128\,r^{15}-3001568\,r^{14}+9512164\,r^{13}-11014136\,r^{12}+2344968\,r^{11}+7126240\,r^{10}-\\
13850504\,r^9+14466592\,r^8-3823216\,r^7-4018976\,r^6+5155776\,r^5-4633984\,r^4+1959808\,r^3-244480\,r^2-3584\,r-1024)
\end{multline*}
\begin{multline*}
\vartheta^{\lo}_{10}(r)=\frac{1}{960\,(2\,r^2-1)(r+2)^3(r-1)(r+1)^2r^{10}}\,(2\,r^{24}-34\,r^{23}-101\,r^{22}+433\,r^{21}+5400\,r^{20}
-26982\,r^{19}+23049\,r^{18}+\\
166787\,r^{17}-717366\,r^{16}+1196092\,r^{15}+89468\,r^{14}-5130844\,r^{13}+12748688\,r^{12}-11274744\,r^{11}-12243496\,r^{10}+\\
33980568\,r^9-14886656\,r^8-19910592\,r^7+20667776\,r^6-1262208\,r^5-5402752\,r^4+2217088\,r^3-235776\,r^2-2560\,r-1024)
\end{multline*}
\begin{equation*}
\vartheta^{\lo}_{11}(r)=\frac{2}{15}\,{\frac
{180\,r^8-48\,r^7-648\,r^6+396\,r^5+214\,r^4-190\,r^3+39\,r^2-4\,r+1}{(2\,r^2-1)(r+1)^2r^{10}}}
\end{equation*}
}
\noindent
and
$\mI_1=[1,2/\sqrt{3}),\;
\mI_2=[2/\sqrt{3},6/5),\;
\mI_3=[6/5,\sqrt{5}-1),\;
\mI_4=[\sqrt{5}-1,(6+2\,\sqrt{2})/7),\;
\mI_5=[(6+2\,\sqrt{2})/7,4/3),\;
\mI_6=[4/3,(6+\sqrt{15})/7),\;
\mI_7=[(6+\sqrt{15})/7,3/2),\;
\mI_8=[3/2,(1+\sqrt{5})/2),\;
\mI_9=[(1+\sqrt{5})/2,1+1/\sqrt{2}),\;
\mI_{10}=[1+1/\sqrt{2},2),\;
\mI_{11}=[2,\infty)$.
See Figure \ref{cov of nu under-graph}.
Note that $\Cov_{\lo}(r=1)=1/3240$ and $\lim_{r \rightarrow \infty}\nu_{\lo}(r)=0$ (at rate $O(r^{-6})$),
and $\argsup_{r \in [1,\infty)} \nu_{\lo}(r) \approx 1.765$ with $\sup \nu_{\lo}(r) \approx .0318$.

\section*{Appendix 3: Derivation of $\mu_{\la}(r)$ and $\nu_{\la}(r)$ under the Null Case}
In the standard equilateral triangle,
let $\y_1=(0,0)$, $\y_2=(1,0)$, $\y_3=\bigl( 1/2,\sqrt{3}/2 \bigr)$,
$M_C$ be the center of mass,
$M_i$ be the midpoints of the edges $e_i$ for $i=1,2,3$.
Then $M_C=\bigl(1/2,\sqrt{3}/6\bigr)$,
$M_1=\bigl(3/4,\sqrt{3}/4 \bigr)$, $M_2=\bigl(1/4,\sqrt{3}/4\bigr)$, $M_3=(1/2,0)$.
%
%
Let $\X_n$ be a random sample of size $n$ from $\U(\TY)$.
For $x_1=(u,v)$, $\ell_r(x_1)=r\,v+r\,\sqrt{3}\,u-\sqrt{3}\,x.$
Next, let $N_1:=\ell_r(x_1)\cap e_3$ and $N_2:=\ell_r(x_1)\cap e_2$.

\subsection*{Derivation of $\mu_\la(r)$ in Theorem \ref{thm:asy-norm-under}}
First we find $\mu_\la(r)$ for $r \in (1,\infty)$.
Observe that, by symmetry,
$$\mu_\la(r)=P\bigl( X_2 \in \NPE^r(X_1) \cap \G_1^r(X_1) \bigr)=
6\,P\bigl( X_2 \in \NY^r(X_1) \cap \G_1^r(X_1), X_1 \in T_s \bigr)$$
where $T_s$ is the triangle with vertices $\y_1$, $M_3$, and $M_C$.
Let $\ell_s(r,x)$ be the line such that $r\,d(\y_1,\ell_s(r,x))=d(\y_1,e_1)$,
so $\ell_s(r,x)=\sqrt{3}\,(1/r-x)$.
Then if $x_1 \in T_s$ is above $\ell_s(r,x)$ then $\NPE^r(x_1)=\TY$,
otherwise, $\NPE^r(x_1)\subsetneq \TY$.

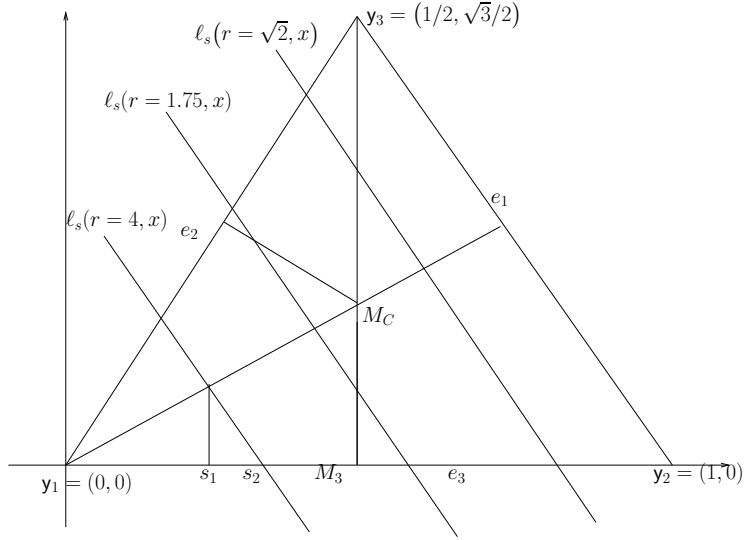
\begin{figure} [ht]
\centering
\scalebox{.4}{\input{ls_lam_cases.pstex_t}}
\caption{
\label{fig:ls-lam-cases}
The cases for relative position of $\ell_s(r,x)$ with various $r$ values.
These are the prototypes for various types of $\NPE^r(x_1)$.
}
\end{figure}

\begin{figure} [ht]
\centering
\scalebox{.27}{\input{G1ofxCase1.pstex_t}} 
\scalebox{.27}{\input{G1ofxCase2.pstex_t}}
\scalebox{.27}{\input{G1ofxCase3.pstex_t}}
\scalebox{.27}{\input{G1ofxCase4.pstex_t}}
\scalebox{.27}{\input{G1ofxCase5.pstex_t}}
\scalebox{.27}{\input{G1ofxCase6.pstex_t}}
\caption{
\label{fig:G_1-NYr-Cases-1}
The prototypes of the six cases of $\G^r_1\left(x\right)$ for $x \in T_s$ for $r \in [1,4/3)$.}
\end{figure}

\begin{figure}[ht]
\centering
\rotatebox{-90}{ \resizebox{8.5 in}{!}{\includegraphics{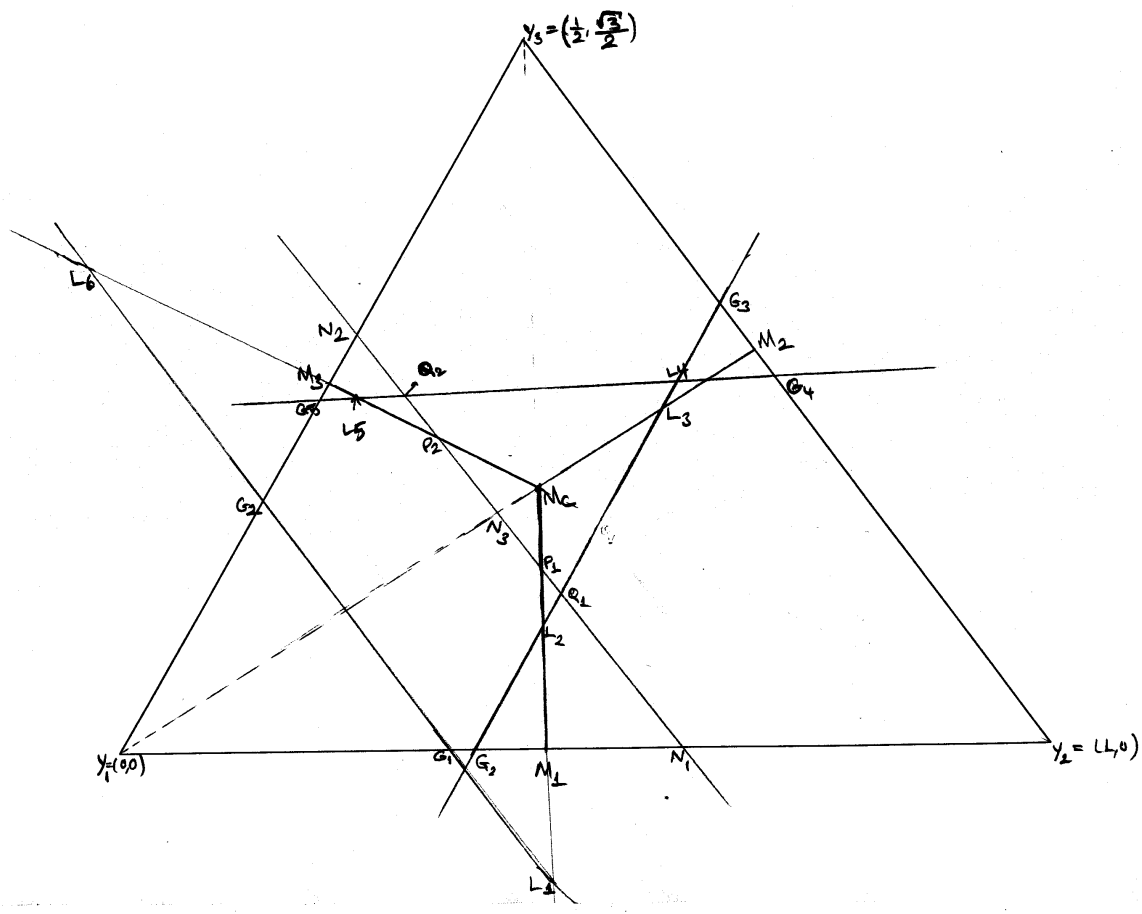}}}
\caption{
\label{fig:vertices-AND-OR}
An illustration of the vertices for possible types of $\NPE^r(x_1) \cap \G_1^r(x_1)$
for $x_1 \in T_s$.
}
\end{figure}

To compute $\mu_\la(r)$,
we need to consider various cases for $\NPE^r(X_1)$ and $\G_1^r(X_1)$ given $X_1=(x,y) \in T_s$.
See Figures \ref{fig:ls-lam-cases} and \ref{fig:G_1-NYr-Cases-1}.
For any $x=(u,v) \in T(\Y)$, $\G^r_1(x)$ is a convex or nonconvex polygon.
Let $\xi_i(r,x)$ be the line  between $x$ and the vertex $\y_i$
parallel to the edge $e_i$ such that $r\,d(\y_i,\xi_i(r,x))=d(\y_i,\ell_r(x)) \text{ for } i=1,2,3.$
Then $\G^r_1(x)\cap R(\y_i)$ is bounded by $\xi_i(r,x)$ and the median lines.
For $x=(u,v)$,
$\xi_1(r,x)=-\sqrt{3}\,x+(v+\sqrt{3}\,u)/r,\; \xi_2(r,x)=
(v+\sqrt{3}r\,(x-1)+\sqrt{3}(1-u))/r \text{ and }
\xi_3(r,x)=(\sqrt{3}(r-1)+2\,v)/(2\,r).$
For $r\in \Bigl[6/5,\sqrt{5}-1)$,
there are six cases regarding $\G_1^r(x)$ and one case for $\NPE^r(x)$.
See Figure \ref{fig:G_1-NYr-Cases-1} for the prototypes of these six cases of $\G_1\left( x,\NY^r \right)$.
For the AND-underlying version,
we determine the possible types of $\NPE^r(x_1) \cap \G_1^r(x_1)$
for $x_1 \in T_s$.
Depending on the location of $x_1$ and the value of the parameter $r$,
$\NPE^r(x_1) \cap \G_1^r(x_1)$ regions are polygons with various vertices.
See Figure \ref{fig:vertices-AND-OR} for the illustration of these vertices and
below for their explicit forms.

$G_1=\left( {\frac{\sqrt{3} y+3\, x}{3r}},0 \right)$,
$G_2=\left( -{\frac{\sqrt{3} y-3\,r+3-3\, x}{3r}},0\right) $,
$G_3=\left( -{\frac{\sqrt{3} y-6\,r+3-3\, x}{6r}},
-{\frac{\sqrt{3}\left( -\sqrt{3} y-3+3\, x \right) }{6r}}\right) $,
$G_4=\Bigl( {\frac{ \left( \sqrt{3}r+\sqrt{3}-2\, y \right) \sqrt{3}}{6 r}},\\
{\frac{\sqrt{3} \left( 3\,r -3+2\,\sqrt{3} y\right) }{6r}}\Bigr) $,
$G_5=\left( {\frac{ \left( \sqrt{3}r -\sqrt{3}+2\, y \right) \sqrt{3}}{6r}},
{\frac{\sqrt{3} \left( 3\,r -3+2\,\sqrt{3} y \right) }{6r}}\right) $,
$G_6=\left( {\frac{\sqrt{3} y+3\, x}{6r}},
{\frac{\sqrt{3} \left( \sqrt{3} y+3\, x \right) }{6r}}\right) $;

$P_1=\left( 1/2,\sqrt{3}/6\, \left( 2\,\sqrt{3}r\, y+6\,r\, x-3 \right) \right) $, and
$P_2=\left(-1/2+(\sqrt{3}r\, y+3\,r\, x)/2,
-\sqrt{3}/6\, \left( -3+\sqrt{3}r\, y+3\,r\, x \right)\right) $;

$L_1=\left( 1/2,{\frac{\sqrt{3} \left( 2\,\sqrt{3} y+6\, x-3\,r \right) }{6r}}\right) $,
$L_2=\left( 1/2,-{\frac{ \left( -2\,\sqrt{3} y-6+6\, x+3\,r \right) \sqrt{3}}{6r}}\right) $,
$L_3=\left( -{\frac{\sqrt{3} y-3\,r+3-3\, x}{2r}},
{\frac{\sqrt{3} \left(3\,r -\sqrt{3} y-3+3\, x \right) }{6r}}\right) $,
$L_4=\left( {\frac{3\,r -3+2\,\sqrt{3} y}{2r}},
{\frac{\sqrt{3} \left( 3\,r -3+2\,\sqrt{3} y \right) }{6r}}\right) $,
$L_5=\left( -{\frac{r -3+2\,\sqrt{3} y}{2r}},
{\frac{\sqrt{3} \left( 3\,r -3+2\,\sqrt{3} y \right) }{6r}}\right) $, and
$L_6=\left({\frac{-r+\sqrt{3} y+3\, x}{2r}},
-{\frac{\sqrt{3} \left( \sqrt{3} y+3\, x-3\,r \right) }{6r}}\right) $;
\noindent
$N_1=\left( \sqrt{3}r\, y/3+r\, x,0\right) $,
$N_2=\left( \sqrt{3}r\, y/6+r\, x/2,
\sqrt{3} \left( \sqrt{3} y/6+3\, x \right) r\right) $, and\\
$N_3=\left(\sqrt{3}r\, y/4+3\,r\, x/4,
\sqrt{3} \left( \sqrt{3} y/12+3\, x \right) r\right) $;
and
$Q_1=\left( {\frac{\sqrt{3}r^2 y+3\,r^2 x-\sqrt{
3} y+3\,r -3+3\, x}{6r}},
{\frac{ \left( \sqrt{3}r^2 y+3\,r^2 x+
\sqrt{3} y-3\,r+3-3\, x \right) \sqrt{3}}{6r}}\right) $,\\
and
$Q_2=\left( {\frac{2\,\sqrt{3}r^2 y+6\,r^2 x-3\,r+3-2\,\sqrt{3} y}{6r}},
{\frac{\sqrt{3} \left( 3\,r -3+2\,\sqrt{3} y \right) }{6r}}\right) $.

Let $\msP(a_1,a_2,\ldots, a_n)$ denote the polygon with vertices $a_1,a_2,\ldots, a_n $.
For $r \in \bigl[1,4/3\bigr)$,
there are 14 cases to consider for calculation of $\mu_{\la}(r)$ in the AND-underlying version.
Each of these cases correspond to the regions in Figure \ref{fig:cases-AND-OR},
where Case 1 corresponds to $R_i$ for $i=1,2,3,4$,
and Case $j$ for $j >1$ corresponds to $R_{j+3}$ for $j=1,2,\ldots,14$.
These regions are bounded by various combinations of the lines defined below.

Let $\ell_{am}(x)$ be the line joining $\y_1$ to $M_C$,
then $\ell_{am}(x)=\sqrt{3}x/3$.
Let also
$r_1(x)=\sqrt{3} \left( 2\,r+3\, x-3 \right)/3$,
$r_2(x)=\sqrt{3}/2-\sqrt{3}r/3$,
$r_3(x)=\left( 2\, x-2+r \right) \sqrt{3}/2$,
$r_4(x)=\sqrt{3}/2-\sqrt{3}r/4$,
$r_5(x)=-{\frac{\sqrt{3} \left( 2\,r\, x-1 \right) }{2r}}$,
$r_6(x)=-{\frac{\sqrt{3} \left( -2+3\,r\, x \right) }{3r}}$,
$r_7(x)=-{\frac{ \left( 1+r^2 x-r - x \right) \sqrt{3}}{r^2+1}}$,
$r_8(x)=-{\frac{ \left( r^2 x-1+ x \right) \sqrt{3}}{r^2-1}}$,
$r_9(x)=-{\frac{ \left( r^2 x-1 \right) \sqrt{3}}{r^2+2}}$,
$r_{10}(x)=-{\frac{ \left( -2\,r+2+r^2 x \right) \sqrt{3}}{-4+r^2}}$,
$r_{11}(x)=-{\frac{ \left( -2\,r+2-2\, x+r^2 x \right) \sqrt{3}}{r^2+2}}$,
$r_{12}(x)=-\left( 2\, x-r \right) \sqrt{3}/2$, and
$r_{13}(x)=-\left( -1+ x \right) \sqrt{3}/3$.
Furthermore, to determine the integration limits,
we specify the $x$-coordinate of the boundaries of these regions using $s_k$  for $k=0,1,\ldots,14$.
See also Figure \ref{fig:cases-AND-OR} for an illustration of these points
whose explicit forms are provided below.

$s_0=1-2\,r/3$,
$s_1=3/2-r$,
$s_2=3/(8\,r)$,
$s_3={\frac{-3\,r+2\,r^2+3}{6r}}$,
$s_4=1-r/2$,
$s_5={\frac{2\,r-r^2+1}{4r}}$,
$s_6=1/(2\,r)$,
$s_7=\frac{3}{2\, \left( 2\,r^2+1 \right)}$,
$s_8={\frac{9-3\,r^2+2\,r^3-2\,r}{6(r^2+1)}}$,
$s_9=1/\left( r+1 \right)$,
$s_{10}={\frac{-3\,r+2\,r^2+4}{6r}}$,
$s_{11}=3\,r/8$,
$s_{12}={\frac{6\,r-3\,r^2+4}{12r}}$,
$s_{13}=3/2-5\,r/6$, and
$s_{14}=r-1/2-r^3/8$.

%

Below, we compute $P(X_2 \in \NPE^r(X_1) \cap \G_1^r(X_1), X_1 \in T_s)$ for each of the 14 cases:
{\small
\noindent \textbf{Case 1:}
\begin{multline*}
P(X_2 \in \NPE^r(X_1) \cap \G_1^r(X_1), X_1 \in T_s)=
\left(\int_{0}^{s_2}\int_{0}^{\ell_{am}(x)} + \int_{s_2}^{s_6}\int_{0}^{r_5(x)}\right)
\frac{A(\msP(G_1,N_1,N_2,G_6))}{A(\TY)^2}dydx=\\
{\frac{ \left( r -1 \right)  \left( r+1 \right)
 \left( r^2+1 \right) }{64\,r^6}}
\end{multline*}
where
$A(\msP(G_1,N_1,N_2,G_6))=
\sqrt{3}/36\, \left( \sqrt{3} y+3\, x \right)^2r^2-
{\frac{\sqrt{3} \left( \sqrt{3} y+3\, x
\right)^2}{36\,r^2}}
$.

\noindent \textbf{Case 2:}
\begin{multline*}
P(X_2 \in \NPE^r(X_1) \cap \G_1^r(X_1), X_1 \in T_s)=
\left(\int_{s_5}^{s_6}\int_{r_5(x)}^{r_7(x)} + \int_{s_6}^{s_9}\int_{0}^{r_7(x)}\right)
\frac{A(\msP(G_1,N_1,P_2,M_3,G_6))}{A(\TY)^2}dydx=\\
{\frac{ \left( 9\,r^5+23\,r^4+24\,r^3+24\,r^2+13\,r+3 \right)
\left( r -1 \right)^4}{96\, r^6 \left( r+1 \right)^3}}
\end{multline*}
where
$A(\msP(G_1,N_1,P_2,M_3,G_6))=
-{\frac{\sqrt{3} \left( -4\,r^3\sqrt{3} y-12\,r^3 x+
2\,r^4\,y^2+4\,r^4\sqrt{3}y\, x+6\,r^4 x^2+3\,r^2+2\,y^2+
4\,\sqrt{3} y\, x+6\, x^2 \right) }{24\,r^2}}$.

\noindent \textbf{Case 3:}
\begin{multline*}
P(X_2 \in \NPE^r(X_1) \cap \G_1^r(X_1), X_1 \in T_s)=\\
\left( \int_{s_5}^{s_9}\int_{r_7(x)}^{r_3(x)} +
\int_{s_9}^{s_{12}}\int_{0}^{r_3(x)} +
\int_{s_{12}}^{1/2}\int_{0}^{r_6(x)} \right)
\frac{A(\msP(G_1,G_2,Q_1,P_2,M_3,G_6))}{A(\TY)^2}dydx=\\
{\frac{324\,r^{11}-1620\,r^{10}-618\,r^9+4626\,r^8+990\,r^7-2454\,r^6+2703\,r^5-
5571\,r^4-3827\,r^3+1455\,r^2+3072\,r+1024}
{7776\, \left( r+1 \right)^3r^6}}
\end{multline*}
where
$A(\msP(G_1,G_2,Q_1,P_2,M_3,G_6))=
-\Bigl[\sqrt{3} \bigl( -4\,\sqrt{3}r\, y-12\,x+4\,y^2+4\,r^2\,y^2-12\,r+9\,r^2+12
\,r\, x+4\,r^4\,y^2-12\, x^2r^2-24\,r^3 x+12\,r^4 x^2+8\,r^4\sqrt{
3} y\, x+12\, x^2+12\,r^2 x+6-8\,r^3\sqrt{3} y+4\,\sqrt{3} y+4\,\sqrt{3}r^2 y \bigr)\Bigr]\Big/\Bigl[24\,r^2\Bigr]
$.

\noindent \textbf{Case 4:}
\begin{multline*}
P(X_2 \in \NPE^r(X_1) \cap \G_1^r(X_1), X_1 \in T_s)=\\
\left(\int_{s_8}^{s_5}\int_{r_8(x)}^{r_2(x)}+
\int_{s_5}^{s_{10}}\int_{r_3(x)}^{r_2(x)}+
\int_{s_{10}}^{s_{12}}\int_{r_3(x)}^{r_6(x)}\right)
\frac{A(\msP(G_1,M_1,L_2,Q_1,P_2,M_3,G_6))}{A(\TY)^2}dydx=\\
\Bigl[512+138240\,r^7+3654\,r^{12}-255
\,r^8+43008\,r^3-12369\,r^2-86387\,r^4-193581
\,r^6+148224\,r^5-100608\,r^9+94802\,r^{10}-\\
35328\,r^{11}\Bigr]\Big/
\Bigl[7776\, \left( r^2+1 \right)^3r^6\Bigr]
\end{multline*}
where
$A(\msP(G_1,M_1,L_2,Q_1,P_2,M_3,G_6))=
-\Bigl[\sqrt{3} \bigl( 6\, x+3\,r^2-2\,\sqrt{3} y+2\,\sqrt{3}r^2 y+2\,r^4\,y^2-4\,
r^3\sqrt{3} y+4\,\sqrt{3} y\, x+2\,r^2 y^2+4\,r^4\sqrt{3} y\, x-6\, x^2 r^2-
12\,r^3 x+6\,r^4 x^2+6\,r^2 x-3 \bigr) \Bigr]\Big/\Bigl[12\,r^2\Bigr]
$.

\noindent \textbf{Case 5:}
\begin{multline*}
P(X_2 \in \NPE^r(X_1) \cap \G_1^r(X_1), X_1 \in T_s)=
\left(\int_{s_3}^{s_8}\int_{r_5(x)}^{r_2(x)} +
\int_{s_8}^{s_5}\int_{r_5(x)}^{r_8(x)} \right)
\frac{A(\msP(G_1,M_1,P_1,P_2,M_3,G_6))}{A(\TY)^2}dydx=\\
-{\frac{ \left( 177\,r^8-648\,r^7+570\,r^6-360\,r^5+28\,r^4-24\,r^3+174\,r^2+
72\,r+27 \right)  \left( -12\,r+7\,r^2+3 \right)^2}{7776\,\left( r^2+1 \right)^3r^6}}
\end{multline*}
where
$A(\msP(G_1,M_1,L_2,Q_1,P_2,M_3,G_6))=
-{\frac{\sqrt{3} \left( -4\,r^3\sqrt{3} y-12\,r^3 x+
3\,r^2+6\,r^4\sqrt{3} y\, x+9\,r^4 x^2+3\,r^4\,y^2+\,y^2+2
\,\sqrt{3} y\, x+3\, x^2 \right) }{12\,r^2}}
$.

\noindent \textbf{Case 6:}
\begin{multline*}
P(X_2 \in \NPE^r(X_1) \cap \G_1^r(X_1), X_1 \in T_s)=\\
\left( \int_{s_2}^{s_3}\int_{r_5(x)}^{\ell_{am}(x)} +
\int_{s_3}^{s_7}\int_{r_2(x)}^{\ell_{am}(x)} +
\int_{s_7}^{s_8}\int_{r_2(x)}^{r_8(x)} \right)
\frac{A(\msP(G_1,M_1,P_1,P_2,M_3,G_6))}{A(\TY)^2}dydx=\\
\Bigl[137472\,r^{18}-952704\,r^{17}+2792712\,r^{16}-5116608\,r^{15}+7057828\,r^{14}-
7725792\,r^{13}+7022682\,r^{12}-5484816\,r^{11}+\\
3631995\,r^{10}-2213712\,r^9+1213271\,r^8-578976\,r^7+292518\,r^6-
101952\,r^5+36612\,r^4-11664\,r^3+3051\,r^2-1296\,r+243\Bigr]\Big/\\
\Bigl[ \left(15552\, r^2+1 \right)^3 \left( 2\,r^2+1 \right)^3r^6\Bigr]
\end{multline*}
where
$A(\msP(G_1,M_1,P_1,P_2,M_3,G_6))=
-{\frac{\sqrt{3} \left( -4\,r^3\sqrt{3} y-12\,r^3 x+
3\,r^2+6\,r^4\sqrt{3} y\, x+9\,r^4 x^2+3\,r^4\,y^2+\,y^2+2
\,\sqrt{3} y\, x+3\, x^2 \right) }{12\,r^2}}
$.

\noindent \textbf{Case 7:}
\begin{multline*}
P(X_2 \in \NPE^r(X_1) \cap \G_1^r(X_1), X_1 \in T_s)=
\left(\int_{s_7}^{s_8}\int_{r_8(x)}^{r_9(x)}+
\int_{s_8}^{s_{10}}\int_{r_2(x)}^{r_9(x)}\right)
\frac{A(\msP(G_1,M_1,L_2,Q_1,P_2,M_3,G_6))}{A(\TY)^2}dydx=\\
-{\frac{ 4\left( 100\,r^{11}-408\,r^{10}+454\,r^9-564\,r^8+283\,r^7-108\,r^6-34\,r^5+
204\,r^4-r^3+132\,r^2+26\,r+24 \right)  \left( 2\,r -1 \right)^2 \left( r -1 \right)^2}
{243 \left( r^2+1 \right)^3r^3 \left( 2\,r^2+1 \right)^3}}
\end{multline*}
where
$A(\msP(G_1,M_1,L_2,Q_1,P_2,M_3,G_6))=
-\Bigl[\sqrt{3} \bigl( 6\, x+3\,r^2-2\,\sqrt{3} y+2\,\sqrt{3}r^2 y+2\,r^4\,y^2-
4\,r^3\sqrt{3} y+4\,\sqrt{3} y\, x+2\,r^2 y^2+4\,r^4\sqrt{3} y\, x-
6\,x^2 r^2-12\,r^3 x+6\,r^4 x^2+6\,r^2 x-3 \bigr) \Bigr] \Big / \Bigl[12\,r^2\Bigr]
$.

\noindent \textbf{Case 8:}
\begin{multline*}
P(X_2 \in \NPE^r(X_1) \cap \G_1^r(X_1), X_1 \in T_s)=
\left( \int_{s_{12}}^{s_{13}}\int_{r_6(x)}^{r_3(x)} +
\int_{s_{13}}^{1/2}\int_{r_6(x)}^{r_2(x)} \right)
\frac{A(\msP(G_1,G_2,Q_1,N_3,M_C,M_3,G_6))}{A(\TY)^2}dydx=\\
\Bigl[ \left( -2+r \right)  \bigl( 2369\,r^{11}-11342\,r^{10}+29934\,r^9-50340\,r^8+
54056\,r^7-51824\,r^6+48320\,r^5-20864\,r^4-640\,r^3\\
-1280\,r^2+512\,r+1024 \bigr) \Bigr] \Big/ \Bigl[15552\,r^6\Bigr]
\end{multline*}
where
$A(\msP(G_1,G_2,Q_1,N_3,M_C,M_3,G_6))=
-\Bigl[\sqrt{3} \bigl( 4\,\sqrt{3}r^2 y-12\, x-12\,r+5\,r^2+12\,r\, x+4\,y^2-
12\,x^2r^2+4\,r^2\,y^2+r^4\,y^2+2\,r^4\sqrt{3} y\, x-4\,r^3\sqrt{3} y+
6-12\,r^3 x+3\,r^4 x^2+12\, x^2+12\,r^2 x-4\,\sqrt{3}r\, y+4\,\sqrt{3} y \bigr) \Bigr] \Big/ \Bigl[24\,r^2\Bigr]
$.

\noindent \textbf{Case 9:}
\begin{multline*}
P(X_2 \in \NPE^r(X_1) \cap \G_1^r(X_1), X_1 \in T_s)=
\left( \int_{s_{10}}^{s_{12}}\int_{r_6(x)}^{r_2(x)} +
\int_{s_{12}}^{s_{13}}\int_{r_3(x)}^{r_2(x)} \right)
\frac{A(\msP(G_1,M_1,L_2,Q_1,N_3,M_C,M_3,G_6))}{A(\TY)^2}dydx=\\
-{\frac{ \left( 49\,r^8-168\,r^7+354\,r^6-528\,r^5+236\,r^4-96\,r^3-224\,r^2
+384\,r+64 \right)  \left( -12\,r+7\,r^2+4 \right)^2}{15552\,r^6}}
\end{multline*}
where
$A(\msP(G_1,M_1,L_2,Q_1,N_3,M_C,M_3,G_6))=
-\Bigl[\sqrt{3} \bigl( 8\,\sqrt{3} y\, x+4\,\sqrt{3}r^2 y+12\, x+
2\,r^2-12\, x^2 r^2-4\,r^3\sqrt{3} y-12\,r^3 x+3\,r^4 x^2+
r^4\,y^2+2\,r^4\sqrt{3}y\,x+12\,r^2 x-6-4\,\sqrt{3} y+
4\,r^2\,y^2 \bigr) \Bigr] \Big/ \Bigl[24 \, r^2\Bigr]
$.

\noindent \textbf{Case 10:}
\begin{multline*}
P(X_2 \in \NPE^r(X_1) \cap \G_1^r(X_1), X_1 \in T_s)=\\
\left(\int_{s_{10}}^{s_{14}}\int_{r_2(x)}^{r_{10}(x)} +
\int_{s_{14}}^{s_{13}}\int_{r_2(x)}^{r_{12}(x)} +
\int_{s_{13}}^{1/2}\int_{r_3(x)}^{r_{12}(x)} \right)
\frac{A(\msP(G_1,M_1,L_2,Q_1,N_3,L_4,L_5,M_3,G_6))}{A(\TY)^2}dydx=\\
-{\frac{6144+195456\,r^6+324\,r^{11}-
76720\,r^7-801792\,r^2+217856\,r+946432\,r^3-
239904\,r^5-275328\,r^4+39408\,r^8-11849\,r^9}{31104\,r^3}}
\end{multline*}
where
$A(\msP(G_1,M_1,L_2,Q_1,N_3,L_4,L_5,M_3,G_6))=
-\Bigl[\sqrt{3} \bigl( 4\,\sqrt{3}r^2 y+8\,\sqrt{3} y\, x+
4\,r^2\,y^2-16\,\sqrt{3} r\, y-4\,r^3\sqrt{3} y-24\,y^2+
12\, x+24\,r -6\,r^2-12\, x^2r^2-12\,r^3 x+3\,r^4 x^2+12\,r^2 x+20\,\sqrt{3
} y+2\,r^4\sqrt{3} y\, x+r^4\,y^2-24 \bigr)\Bigr] \Big/ \Bigl[24\,r^2\Bigr]
$.

\noindent \textbf{Case 11:}
\begin{multline*}
P(X_2 \in \NPE^r(X_1) \cap \G_1^r(X_1), X_1 \in T_s)=\\
\left(\int_{s_7}^{s_{11}}\int_{r_9(x)}^{\ell_{am}(x)} +
\int_{s_{11}}^{s_{10}}\int_{r_9(x)}^{r_{12}(x)} +
\int_{s_{10}}^{s_{14}}\int_{r_{10}(x)}^{r_{12}(x)} \right)
\frac{A(\msP(G_1,M_1,L_2,Q_1,Q_2,L_5,M_3,G_6))}{A(\TY)^2}dydx=\\
\Bigl[ \left( r -1 \right)  \bigl( 1080\,r^{16}+1080\,r^{15}-17820\,r^{14}-
540\,r^{13}+65394\,r^{12}-46926\,r^{11}+105435\,r^{10}-261765\,r^9+229286\,r^8-180586\,r^7+\\
101638\,r^6+40774\,r^5-46112\,r^4+24448\,r^3-20224\,r^2+10496\,r -6144 \bigr) \Bigr]
\Big/ \Bigl[10368 \, r^3 \left( 2\,r^2+1 \right)^3\Bigr]
\end{multline*}
where
$A(\msP(G_1,M_1,L_2,Q_1,Q_2,L_5,M_3,G_6))=
-\Bigl[\sqrt{3} \bigl( 6\, x+3\,r^2-4\,r^2 x\,\sqrt{3} y-4\,y^2-
6\, x^2r^2+2\,r^4\sqrt{3} y\, x+4\,\sqrt{3} y\,x-2\,r^2\,y^2-
4\,r^3\sqrt{3} y+r^4\,y^2-12\,r^3 x+3\,r^4 x^2+12\,r^2 x-6+
4\,\sqrt{3}r^2 y+2\,\sqrt{3} y \bigr) \Bigr]\Big/\Bigl[12\,r^2\Bigr]
$.

\noindent \textbf{Case 12:}
\begin{multline*}
P(X_2 \in \NPE^r(X_1) \cap \G_1^r(X_1), X_1 \in T_s)=
\int_{s_{13}}^{1/2}\int_{r_2(x)}^{r_3(x)} \frac{A(\msP(G_1,G_2,Q_1,N_3,L_4,L_5,M_3,G_6))}{A(\TY)^2}dydx=\\
-{\frac{ \left( 49\,r^6-204\,r^5+476\,r^4-768\,r^3-8\,r^2+768\,r -288 \right)  \left( -6+
5\,r \right)^2}{7776\,r^2}}
\end{multline*}
where
$A(\msP(G_1,G_2,Q_1,N_3,L_4,L_5,M_3,G_6))=
-\Bigl[\sqrt{3} \bigl( -12\, x+12\,r -3\,r^2+12\,r\, x-20\,\sqrt{3}r\, y-
12\, x^2r^2+4\,\sqrt{3}r^2 y-12\,r^3 x+3\,r^4 x^2+28\,\sqrt{3} y+
12\, x^2+12\,r^2 x-12-20\,y^2+4\,r^2\,y^2-4\,r^3 \sqrt{3} y+
r^4\,y^2+2\,r^4\sqrt{3} y\, x \bigr) \Bigr] \Big/ \Bigl[24\,r^2\Bigr]
$.

\noindent \textbf{Case 13:}
\begin{multline*}
P(X_2 \in \NPE^r(X_1) \cap \G_1^r(X_1), X_1 \in T_s)=
\int_{s_{14}}^{1/2}\int_{r_{12}(x)}^{r_{10}(x)} \frac{A(\msP(L_1,L_2,Q_1,N_3,L_4,L_5,L_6))}{A(\TY)^2}dydx=\\
{\frac{ \left( 4\,r^7+8\,r^6-37\,r^5-58\,r^4-84\,r^3+168\,r^2+336\,r -352 \right)
\left( -2+r \right)  \left( r^2+2\,r -4 \right)^2}{384\, \left( r+2 \right)^2r^2}}
\end{multline*}
where
$A(\msP(L_1,L_2,Q_1,N_3,L_4,L_5,L_6))=
-\Bigl[\sqrt{3} \bigl( -4\,r^3\sqrt{3} y-8\,\sqrt{3}r\, y+12\, x+24\,r -
8\,\sqrt{3} y\,x-12\,r^2+24\,r\, x-24-12\, x^2r^2+4\,\sqrt{3}r^2 y-
32\,y^2-12\,r^3 x+3\,r^4 x^2+20\,\sqrt{3} y-24\, x^2+
12\,r^2 x+2\,r^4\sqrt{3} y\, x+r^4\,y^2+4\,r^2\,y^2 \bigr) \Big/ \Bigl[24\,r^2\Bigr]
$.

\noindent \textbf{Case 14:}
\begin{multline*}
P(X_2 \in \NPE^r(X_1) \cap \G_1^r(X_1), X_1 \in T_s)=
\left( \int_{s_{11}}^{s_{14}}\int_{r_{12}(x)}^{\ell_{am}(x)} +
\int_{s_{14}}^{1/2}\int_{r_{10}(x)}^{\ell_{am}(x)} \right)
\frac{A(\msP(L_1,L_2,Q_1,Q_2,L_5,L_6))}{A(\TY)^2}dydx=\\
-\Bigl[ \bigl( 135\,r^{11}+675\,r^{10}-1350\,r^9-9450\,r^8+702\,r^7+39150\,r^6+24272
\,r^5-47432\,r^4-135040\,r^3+57088\,r^2+204800\,r -\\
134144 \bigr)
\left( r -1 \right) \Bigr]\Big/ \Bigl[10368 \left( r+2 \right)^2 r^2\Bigr]
\end{multline*}
where
$A(\msP(L_1,L_2,Q_1,Q_2,L_5,L_6))=
-\Bigl[\sqrt{3} \bigl( -4\,r^3\sqrt{3} y+4\,\sqrt{3}r\, y+r^4\,y^2+6\, x-
4\,\sqrt{3} y\, x+2\,r^4\sqrt{3} y\, x+12\,r\,x-
4\,r^2 x\,\sqrt{3} y-6\, x^2 r^2+4\,\sqrt{3}r^2 y-12\,r^3 x+3\,r^4 x^2+
2\,\sqrt{3} y-12\, x^2+12\,r^2 x-6-8\,y^2-2\,r^2\,y^2 \bigr) \Bigr]\Big/ \Bigl[12\,r^2\Bigr]$.
}

Adding up the $P(X_2 \in \NPE^r(X_1) \cap \G_1^r(X_1), X_1 \in T_s)$
values in the 14 possible cases above, and multiplying by 6
we get for $r \in [1,4/3)$,
$$\mu_\la(r)=
-{\frac{ \left( r -1 \right)  \left( 5\,r^5-
148\,r^4+245\,r^3-178\,r^2-232\,r+128 \right) }
{54\,r^2 \left( r+2 \right)  \left( r+1 \right) }}.$$
The $\mu_\la(r)$ values for the other intervals can be calculated similarly.
For $r=\infty $, $\mu_\la(r)=1$ follows trivially.

\subsection*{Derivation of $\nu_\la(r)$ in Theorem \ref{thm:asy-norm-under}}

By symmetry,
$P(\{X_2,X_3\} \subset \NPE^r(X_1)\cap \G_1^r(X_1))=
6\,P(\{X_2,X_3\} \subset \NPE^r(X_1)\cap \G_1^r(X_1),\; X_1 \in T_s)$.

For $r \in \bigl[6/5,\sqrt{5}-1\bigr)$,
there are 14 cases to consider for calculation of $\nu_{\la}(r)$ in the AND-underlying version:
{\small
\noindent \textbf{Case 1:}
\begin{multline*}
P(\{X_2,X_3\} \subset \NPE^r(X_1)\cap \G_1^r(X_1),\; X_1 \in T_s)=
\left(\int_{0}^{s_2}\int_{0}^{\ell_{am}(x)}+
\int_{s_2}^{s_6}\int_{0}^{r_5(x)} \right)
\frac{A(\msP(G_1,N_1,N_2,G_6))^2}{A(\TY)^3}dydx=\\
{\frac{ \left( r^2+1 \right)^2 \left( r+1 \right)^2 \left( r -1 \right)^2}{384\,r^{10}}}
\end{multline*}
where
$A(\msP(G_1,N_1,N_2,G_6))=
\sqrt{3} \left( \sqrt{3} y+3\, x \right)^2 r^2/36-
{\frac{ \left( \sqrt{3} y+3\, x \right)^2\sqrt{3}}{36\,r^2}}
$.

\noindent \textbf{Case 2:}
\begin{multline*}
P(\{X_2,X_3\} \subset \NPE^r(X_1)\cap \G_1^r(X_1),\; X_1 \in T_s)=
\left(\int_{s_5}^{s_6}\int_{r_5(x)}^{r_7(x)} +
\int_{s_6}^{s_9}\int_{0}^{r_7(x)} \right)
\frac{A(\msP(G_1,N_1,P_2,M_3,G_6))^2}{A(\TY)^3}dydx=\\
{\frac { \left( 5+38\,r+137\,r^2+320\,r^3+
552\,r^4+736\,r^5+792\,r^6+640\,r^7+407\,r^8+178\,r^9+
35\,r^{10} \right)  \left( -1+r \right)^5}{960\,r^{10} \left( r+1 \right)^5}}
\end{multline*}
where
$A(\msP(G_1,N_1,P_2,M_3,G_6))=
-{\frac {\sqrt{3} \left( -4\,r^3\sqrt{3}y-12\,r^3 x+
2\,r^4 y^2+4\,r^4\sqrt{3}y\,x+6\,r^4x^2+3\,r^2+
2\,y^2+ 4\,\sqrt{3}y\,x+6\,x^2 \right) }{24\,r^2}}
$.

\noindent \textbf{Case 3:}
\begin{multline*}
P(\{X_2,X_3\} \subset \NPE^r(X_1)\cap \G_1^r(X_1),\; X_1 \in T_s)=\\
\left(\int_{s_5}^{s_9}\int_{r_7(x)}^{r_3(x)} +
\int_{s_9}^{s_{12}}\int_{0}^{r_3(x)} +
\int_{s_{12}}^{1/2}\int_{0}^{r_6(x)}\right)
\frac{A(\msP(G_1,G_2,Q_1,P_2,M_3,G_6))^2}{A(\TY)^3}dydx=\\
-\Bigl[17496\,r^{19}-122472\,r^{18}+139968\,r^{17}+524880\,r^{16}-553095\,r^{15}-
595971\,r^{14}+368826\,r^{13}-724758\,r^{12}-543876\,r^{11}+\\
1416996\,r^{10}+1646470\,r^9+92870\,r^8+523048\,r^7-768368\,r^6-1729902\,r^5-1434990\,r^4+122185\,r^3+ 941941\,r^2+\\
573440\,r+114688\Bigr] \Big/ \Bigl[2099520\, \left( r+1 \right)^5r^{10}\Bigr]
\end{multline*}
where
$A(\msP(G_1,G_2,Q_1,P_2,M_3,G_6))=
-\Bigl[\sqrt{3} \bigl( 4\,\sqrt{3}r^2y-8\,r^3\sqrt{3}y+4\,r^2 y^2+4\,r^4 y^2+
4\,y^2+8\,r^4\sqrt{3}y\,x+6-12\,x^2 r^2-12\,x-12\,r-24\,r^3 x+
12\,r^4 x^2+9\,r^2+12\,r\,x-4\,\sqrt{3}r\,y+12\,x^2+
4\,\sqrt{3}y+12\,r^2 x \bigr)\Bigr]\Big/\Bigl[24\,r^2\Bigr]
$.

\noindent \textbf{Case 4:}
\begin{multline*}
P(\{X_2,X_3\} \subset \NPE^r(X_1)\cap \G_1^r(X_1),\; X_1 \in T_s)=\\
\left( \int_{s_8}^{s_5}\int_{r_8(x)}^{r_2(x)} +
\int_{s_5}^{s_{10}}\int_{r_3(x)}^{r_2(x)} +
\int_{s_{10}}^{s_{12}}\int_{r_3(x)}^{r_6(x)} \right)
\frac{A(\msP(G_1,M_1,L_2,Q_1,P_2,M_3,G_6))^2}{A(\TY)^3}dydx=\\
-\Bigl[32768-409264128\,r^7+1455989508\,r^{12}+680709729\,r^8-4423680\,r^3+155509\,r^2+
22889801\,r^4+202936917\,r^6+\\
6011901\,r^{20}+1060982949\,r^{16}-614739456\,r^{17}+240330993\,r^{18}-
56097792\,r^{19}-77783040\,r^5-999857664\,r^9+\\
1299257316\,r^{10}-1461851136\,r^{11}-1407624192\,r^{13}+
1414729905\,r^{14}-1352392704\,r^{15}\Bigr]\Big/\Bigl[2099520\, \left( r^2+1 \right)^5r^{10}\Bigr]
\end{multline*}
where
$A(\msP(G_1,M_1,L_2,Q_1,P_2,M_3,G_6))=
-\Bigl[\sqrt{3} \bigl( -6\,x^2r^2-3+6\,x-12\,r^3 x+6\,r^4 x^2-
4\,r^3\sqrt{3} y+4\,\sqrt{3}y\,x+4\,r^4\sqrt{3} y\,x+2\,r^4y^2+
3\,r^2+2\,\sqrt{3}r^2 y-2\,\sqrt{3} y+2\,r^2 y^2+6\,r^2 x \bigr)\Bigr]\Big/\Bigl[12 \, r^2\Bigr]
$.

\noindent \textbf{Case 5:}
\begin{multline*}
P(\{X_2,X_3\} \subset \NPE^r(X_1)\cap \G_1^r(X_1),\; X_1 \in T_s)=
\left(\int_{s_3}^{s_8}\int_{r_5(x)}^{r_2(x)} +
\int_{s_8}^{s_5}\int_{r_5(x)}^{r_8(x)}\right)
\frac{A(\msP(G_1,M_1,P_1,P_2,M_3,G_6))^2}{A(\TY)^3}dydx=\\
\Bigl[\bigl( 35361\,r^{16}-229392\,r^{15}+602820\,r^{14}-858384\,r^{13}+778848\,r^{12}-
460368\,r^{11}+277740\,r^{10}-258768\,r^9+160594\,r^8-62256\,r^7-\\
5892\,r^6-17712\,r^5+19224\,r^4+11664\,r^3+5076\,r^2+1296\,r+405 \bigr)
\left( -12\,r+7\,r^2+3 \right)^2\Bigr]\Big/\Bigl[699840\,r^{10} \left( r^2+1 \right)^5\Bigr]
\end{multline*}
where
$A(\msP(G_1,M_1,P_1,P_2,M_3,G_6))=
-{\frac {\sqrt{3} \left( -4\,r^3\sqrt{3}y-12\,r^3x+3\,r^2+6\,r^4\sqrt{3}y\,x+
9\,r^4x^2+3\,r^4y^2+y^2+2\,\sqrt{3}y\,x+3\,x^2 \right) }{12\,r^2}}
$.

\noindent \textbf{Case 6:}
\begin{multline*}
P(\{X_2,X_3\} \subset \NPE^r(X_1)\cap \G_1^r(X_1),\; X_1 \in T_s)=\\
\left(\int_{s_2}^{s_3}\int_{r_5(x)}^{\ell_{am}(x)} +
\int_{s_3}^{s_7}\int_{r_2(x)}^{\ell_{am}(x)} +
\int_{s_7}^{s_8}\int_{r_2(x)}^{r_8(x)} \right)
\frac{A(\msP(G_1,M_1,P_1,P_2,M_3,G_6))^2}{A(\TY)^3}dydx=\\
-\Bigl[3645-17496\,r+5003898912\,r^{28}+31646646384\,r^{26}+110098944\,r^{30}-1090803456\,r^{29}-
14630751360\,r^{27}+66339\,r^2-\\
99072645696\,r^{23}+79269457632\,r^{24}+66073158\,r^8-
4870743552\,r^{13}-168073488\,r^9+535086\,r^4-262440\,r^3-1737936\,r^5-\\
18592416\,r^7-107383563504\,r^{21}-41219053272\,r^{17}+58981892347\,r^{18}-78265758888\,r^{19}+95887286866\,r^{20}+\\
109053166552\,r^{22}+5500548\,r^6+466565130\,r^{10}-1070573040\,r^{11}+
2380992104\,r^{12}+9191633420\,r^{14}-16312513248\,r^{15}+\\
26801184917\,r^{16}-54759787776\,r^{25}\Bigr]\Big/\Bigl[1399680\, \left( r^2+1 \right)^5 \left( 2\,r^2+1 \right)
^5r^{10}\Bigr]
\end{multline*}
where
$A(\msP(G_1,M_1,P_1,P_2,M_3,G_6))=
-{\frac {\sqrt{3} \left( -4\,r^3\sqrt{3}y-12\,r^3 x+
3\,r^2+6\,r^4\sqrt{3}y\,x+9\,r^4 x^2+3\,r^4 y^2+y^2+
2\,\sqrt{3}y\,x+3\,x^2 \right) }{12\,r^2}}
$.

\noindent \textbf{Case 7:}
\begin{multline*}
P(\{X_2,X_3\} \subset \NPE^r(X_1)\cap \G_1^r(X_1),\; X_1 \in T_s)=\\
\left(\int_{s_7}^{s_8}\int_{r_8(x)}^{r_9(x)}+
\int_{s_8}^{s_{10}}\int_{r_2(x)}^{r_9(x)} \right)
\frac{A(\msP(G_1,M_1,L_2,Q_1,P_2,M_3,G_6))^2}{A(\TY)^3}dydx=\\
\Bigl[4\,\bigl( 162576\,r^{22}-1083456\,r^
{21}+3368016\,r^{20}-6969888\,r^{19}+11578088\,r^{18}-
15664080\,r^{17}+18796852\,r^{16}-19984824\,r^{15}+\\
19534445\,r^{14}-18170472\,r^{13}+15507752\,r^{12}-
13150464\,r^{11}+9987958\,r^{10}-7448736\,r^9+5016464\,r^8-2991768\,r^7+1857485\,r^6-\\
749160\,r^5+481804\,r^4-96720\,r^3+76160\,r^2-4032\,r+4320 \bigr)
\left( 2\,r -1 \right)^2 \left( r -1 \right)^2\Bigr]\Big/\Bigl[32805\, \left( r^2+1 \right)^5r^6 \left( 2\,r^2+1 \right)^5\Bigr]
\end{multline*}
where
$A(\msP(G_1,M_1,L_2,Q_1,P_2,M_3,G_6))=
-\Bigl[\sqrt{3} \bigl( -6\,x^2r^2-3+6\,x-12\,r^3x+6\,r^4x^2-
4\,r^3\sqrt{3}y+4\,\sqrt{3}y\,x+4\,r^4\sqrt{3}y\,x+2\,r^4y^2+
3\,r^2+2\,\sqrt{3}r^2 y-2\,\sqrt{3}y+2\,r^2 y^2+6\,r^2 x \bigr)\Bigr]\Big/\Bigl[12\,r^2\Bigr]
$.

\noindent \textbf{Case 8:}
\begin{multline*}
P(\{X_2,X_3\} \subset \NPE^r(X_1)\cap \G_1^r(X_1),\; X_1 \in T_s)=\\
\left(\int_{s_{12}}^{s_{13}}\int_{r_6(x)}^{r_3(x)} +
\int_{s_{13}}^{1/2}\int_{r_6(x)}^{r_2(x)} \right)
\frac{A(\msP(G_1,G_2,Q_1,N_3,M_C,M_3,G_6))^2}{A(\TY)^3}dydx=\\
-\Bigl[-458752+811008\,r^2+329205504\,r^8-582626304\,r^{13}-489563136\,r^9-65536\,r^4-
168708096\,r^7-57883680\,r^{17}+\\
18009258\,r^{18}-3623400\,r^{19}+352563\,r^{20}+41502720\,r^6+659111904\,r^{10}-
761846400\,r^{11}+725173376\,r^{12}+409477188\,r^{14}-\\
254829600\,r^{15}+135968852\,r^{16}\Bigr]\Big/\Bigl[8398080\,r^{10}\Bigr]
\end{multline*}
where
$A(\msP(G_1,G_2,Q_1,N_3,M_C,M_3,G_6))=
-\Bigl[\sqrt{3} \bigl( -12\,x^2r^2-12\,x-12\,r-12\,r^3x+3\,r^4x^2+
4\,\sqrt{3}r^2 y+5\,r^2+12\,r\,x+12\,x^2+2\,r^4\sqrt{3}y\,x+4\,r^2y^2-
4\,r^3\sqrt{3}y+6+4\,y^2+r^4 y^2+4\,\sqrt{3}y+12\,r^2 x-4\,\sqrt{3}r\,y \bigr)\Bigr]\Big/\Bigl[24\,r^2\Bigr]
$.

\noindent \textbf{Case 9:}
\begin{multline*}
P(\{X_2,X_3\} \subset \NPE^r(X_1)\cap \G_1^r(X_1),\; X_1 \in T_s)=\\
\left(\int_{s_{10}}^{s_{12}}\int_{r_6(x)}^{r_2(x)} +
\int_{s_{12}}^{s_{13}}\int_{r_3(x)}^{r_2(x)} \right)
\frac{A(\msP(G_1,M_1,L_2,Q_1,N_3,M_C,M_3,G_6))^2}{A(\TY)^3}dydx=
\Bigl[\bigl( 7203\,r^{16}-49392\,r^{15}+\\
170226\,r^{14}-392112\,r^{13}+680784\,r^{12}-1040256\,r^{11}+
1385628\,r^{10}-1337760\,r^9+816224\,r^8-253824\,r^7+\\
469088\,r^6-1029888\,r^5+820992\,r^4-488448\,r^3+190976\,r^2+49152\,r+8192 \bigr)
 \left( -12\,r+7\,r^2+4 \right)^2\Bigr]\Big/\Bigl[8398080\,r^{10}\Bigr]
\end{multline*}
where
$A(\msP(G_1,M_1,L_2,Q_1,N_3,M_C,M_3,G_6))=
-\Bigl[\sqrt{3} \bigl( -12\,x^2r^2-6+12\,x-12\,r^3x+3\,r^4x^2+2\,r^2+
2\,r^4\sqrt{3}y\,x+r^4y^2+8\,\sqrt{3}y\,x+4\,r^2 y^2-
4\,\sqrt{3} y+4\,\sqrt{3}r^2y+12\,r^2x-4\,r^3\sqrt{3}y \bigr)\Bigr]\Big/\Bigl[24\,r^2\Bigr]
$.

\noindent \textbf{Case 10:}
\begin{multline*}
P(\{X_2,X_3\} \subset \NPE^r(X_1)\cap \G_1^r(X_1),\; X_1 \in T_s)=\\
\left(\int_{s_{10}}^{s_{14}}\int_{r_2(x)}^{r_{10}(x)} +
\int_{s_{14}}^{s_{13}}\int_{r_2(x)}^{r_{12}(x)}+
\int_{s_{13}}^{1/2}\int_{r_3(x)}^{r_{12}(x)} \right)
\frac{A(\msP(G_1,M_1,L_2,Q_1,N_3,L_4,L_5,M_3,G_6))^2}{A(\TY)^3}dydx=\\
\Bigl[4423680-4627454976\,r^6+511684992\,r^{11}+2163142656\,r^7-660127744\,r^2-31555584\,r+
3534520320\,r^3+7647989760\,r^5+\\
7785504\,r^{15}-1313880\,r^{16}+19683\,r^{18}-7240624128\,r^4-1511047552\,r^8+
1204122240\,r^9-796453824\,r^{10}-282583320\,r^{12}+\\
107804736\,r^{13}-30362052\,r^{14}\Bigr]\Big/\Bigl[16796160\,r^6\Bigr]
\end{multline*}
where
$A(\msP(G_1,M_1,L_2,Q_1,N_3,L_4,L_5,M_3,G_6))=
-\Bigl[\sqrt{3} \bigl( -16\,\sqrt{3}r\,y+20\,\sqrt{3}y-24\,y^2-12\,x^2 r^2+
12\,x+24\,r-12\,r^3 x+3\,r^4 x^2-6\,r^2-24+4\,\sqrt{3}r^2 y+
8\,\sqrt{3}y\,x-4\,r^3\sqrt{3}y+4\,r^2y^2+r^4 y^2+2\,r^4\sqrt{3}y\,x+
12\,r^2 x \bigr) \Bigr]\Big/\Bigl[24\,r^2\Bigr]
$.

\noindent \textbf{Case 11:}
\begin{multline*}
P(\{X_2,X_3\} \subset \NPE^r(X_1)\cap \G_1^r(X_1),\; X_1 \in T_s)=\\
\left(\int_{s_7}^{s_{11}}\int_{r_9(x)}^{\ell_{am}(x)} +
\int_{s_{11}}^{s_{10}}\int_{r_9(x)}^{r_{12}(x)} +
\int_{s_{10}}^{s_{14}}\int_{r_{10}(x)}^{r_{12}(x)} \right
)\frac{A(\msP(G_1,M_1,L_2,Q_1,Q_2,L_5,M_3,G_6))^2}{A(\TY)^3}dydx=\\
-\Bigl[ \left( r -1 \right)  \bigl( -1474560+
8847360\,r+111456\,r^{26}+111456\,r^{27}-27738112\,r^2
+23311152\,r^{23}-167184\,r^{24}-808889416\,r^8-\\
2228253688\,r^{13}+366739256\,r^9-207619072\,r^4+
98557952\,r^3+397199360\,r^5+802401664\,r^7-34733448\,r^{21}-624736557\,r^{17}+\\
400615470\,r^{18}-134938386\,r^{19}+39014136\,r^{20}-18026064\,r^{22}-640058432\,r^6+
407655352\,r^{10}-1227078728\,r^{11}+\\
1996721576\,r^{12}+2033409092\,r^{14}-1681870468\,r^{15}+1064030499\,r^{16}-
2842128\,r^{25} \bigr) \Bigr]\Big/\Bigl[1866240\, \left( 2\,r^2+1 \right)^5r^6\Bigr]
\end{multline*}
where
$A(\msP(G_1,M_1,L_2,Q_1,Q_2,L_5,M_3,G_6))=
-\Bigl[\sqrt{3} \bigl( 4\,\sqrt{3}r^2y+4\,\sqrt{3}y\,x-2\,r^2 y^2-
4\,r^3\sqrt{3}y-4\,y^2-4\,\sqrt{3}r^2y\,x-6\,x^2 r^2+6\,x-12\,r^3x+
3\,r^4 x^2+3\,r^2+2\,r^4\sqrt{3}y\,x+r^4 y^2+2\,\sqrt{3}y+12\,r^2 x-6 \bigr) \Bigr]\Big/\Bigl[12\,r^2\Bigr]
$.

\noindent \textbf{Case 12:}
\begin{multline*}
P(\{X_2,X_3\} \subset \NPE^r(X_1)\cap \G_1^r(X_1),\; X_1 \in T_s)=
\int_{s_{13}}^{1/2}\int_{r_2(x)}^{r_3(x)} \frac{A(\msP(G_1,G_2,Q_1,N_3,L_4,L_5,M_3,G_6))^2}{A(\TY)^3}dydx=\\
\Bigl[ \bigl( 2322432-7554816\,r+9510912\,{
r}^2+1046068\,r^8-558720\,r^9+2444224\,r^4-5799360\,r^3-2134656\,r^5-1608672\,r^7+\\
2169696\,r^6+216300\,r^{10}-55440\,r^{11}+7095\,r^{12} \bigr)
 \left( -6+5\,r \right)^2\Bigr]\Big/\Bigl[4199040\,r^4\Bigr]
\end{multline*}
where
$A(\msP(G_1,G_2,Q_1,N_3,L_4,L_5,M_3,G_6))=
-\Bigl[\sqrt{3} \bigl( -12\,x^2r^2-12\,x+12\,r-12\,r^3x+3\,r^4x^2-3\,r^2+
12\,r\,x+28\,\sqrt{3}y+12\,x^2-20\, y^2+12\,r^2 x+r^4y^2+4\,r^2 y^2-
4\,r^3\sqrt{3}y+2\,r^4\sqrt{3} y\,x+4\,\sqrt{3}r^2 y-20\,\sqrt{3}r\, y-12 \bigr) \Bigr]\Big/\Bigl[24\,r^2\Bigr]
$.

\noindent \textbf{Case 13:}
\begin{multline*}
P(\{X_2,X_3\} \subset \NPE^r(X_1)\cap \G_1^r(X_1),\; X_1 \in T_s)=
\int_{s_{14}}^{1/2}\int_{r_{12}(x)}^{r_{10}(x)} \frac{A(\msP(L_1,L_2,Q_1,N_3,L_4,L_5,L_6))^2}{A(\TY)^3}dydx=\\
-\Bigl[\bigl( 9\,r^{14}+36\,r^{13}-132\,r^{12}-576\,r^{11}+164\,r^{10}+2512\,r^9+
4976\,r^8-1536\,r^7-13888\,r^6-17536\,r^5-3072\,r^4+79360\,r^3+\\
9216\,r^2-120832\,r+61440 \bigr)  \left( -2+r \right)
\left( r^2+2\,r -4 \right)^2\Bigr]\Big/\Bigl[7680\, \left( r+2 \right)^3r^4\Bigr]
\end{multline*}
where
$A(\msP(L_1,L_2,Q_1,N_3,L_4,L_5,L_6))=
-\Bigl[\sqrt{3} \bigl( r^4 y^2-8\,\sqrt{3} r\,y-8\,\sqrt{3}y\,x+4\,r^2y^2-
4\,r^3\sqrt{3}y-32\,y^2+2\,r^4\sqrt{3}y\,x-12\,x^2 r^2+12\,x+
24\,r-12\,r^3 x+3\,r^4 x^2-12\,r^2+4\,\sqrt{3}r^2 y+24\,r\,x-24\,x^2-24+
20\,\sqrt{3}y+12\,r^2 x \bigr) \Bigr]\Big/\Bigl[24\,r^2\Bigr]
$.

\noindent \textbf{Case 14:}
\begin{multline*}
P(\{X_2,X_3\} \subset \NPE^r(X_1)\cap \G_1^r(X_1),\; X_1 \in T_s)=\\
\left(\int_{s_{11}}^{s_{14}}\int_{r_{12}(x)}^{\ell_{am}(x)} +
\int_{s_{14}}^{1/2}\int_{r_{10}(x)}^{\ell_{am}(x)}\right)
\frac{A(\msP(L_1,L_2,Q_1,Q_2,L_5,L_6))^2}{A(\TY)^3}dydx=\\
\Bigl[\left( r -1 \right)  \bigl( 3483\,r^{18}+24381\,r^{17}-34830\,r^{16}-529416\,r^{15}-
265680\,r^{14}+4274208\,r^{13}+4999320\,r^{12}-15227352\,r^{11}-\\
25751336\,r^{10}+19466488\,r^9+62834064\,r^8+17452256\,r^7-53339200\,r^6-117114624\,r^5-
51206656\,r^4+270430208\,r^3+\\
58073088\,r^2-296222720\,r+122159104 \bigr) \Bigr]\Big/\Bigl[ 1866240\, \left( r+2 \right)^3r^4\Bigr]
\end{multline*}
where
$A(\msP(L_1,L_2,Q_1,Q_2,L_5,L_6))=
-\Bigl[\sqrt{3} \bigl( -4\,\sqrt{3}y\,x-2\,r^2 y^2+4\,\sqrt{3}r\,y-4\,r^3\sqrt{3}y-
8\,y^2-4\,\sqrt{3}r^2 y\,x-6\,x^2 r^2+6\,x-12\,r^3 x+3\,r^4 x^2+
4\,\sqrt{3}r^2y+12\,r\,x-12\, x^2+2\,r^4\sqrt{3}y\,x+r^4 y^2+
2\,\sqrt{3}y+12\,r^2x-6 \bigr) \Bigr]\Big/\Bigl[12\,r^2\Bigr]
$.
}

Adding up the $P(\{X_2,X_3\} \subset \NPE^r(X_1)\cap \G_1^r(X_1),\; X_1 \in T_s)$
values in the 14 possible cases above, and multiplying by 6
we get for $r \in \bigl[6/5,\sqrt{5}-1\bigr)$,
\begin{multline*}
\nu_\la(r)=
-\Bigl[219936\,r -3041936\,r^2-30889822
\,r^8+18084672\,r^{13}+27137438\,r^9+2364868\,r^4+2305864\,r^3-4168820\,r^5-\\
2832544\,r^7+486\,r^{21}-118850\,r^{17}-45155\,r^{18}-269\,r^{19}+3402\,r^{20}+
11101160\,r^6+24604048\,r^{10}-43009544\,r^{11}+8770788\,r^{12}-\\
13736295\,r^{14}+2751855\,r^{15}+443518\,r^{16}+49152\Bigr]
\Big/\Bigl[116640\,r^6 \left( r+2 \right)^2 \left( 2\,r^{
2}+1 \right)  \left( r+1 \right)^3\Bigr].
\end{multline*}
The $\nu_\la(r)$ values for the other intervals can be calculated similarly.

\section*{Appendix 4: Derivation of $\mu_\lo(r)$ and $\nu_\lo(r)$ under the Null Case}

\subsection*{Derivation of $\mu_\lo(r)$ in Theorem \ref{thm:asy-norm-under}}
First we find $\mu_\lo(r)$ for $r \in \Bigl[1,\infty)$.
Observe that, by symmetry,
$$\mu_\lo(r)=P\bigl( X_2 \in \NPE^r(X_1) \cup \G_1^r(X_1) \bigr)=
6\,P\bigl( X_2 \in \NY^r(X_1) \cup \G_1^r(X_1), X_1 \in T_s \bigr).$$

For $r \in [1,4/3)$,
there are 17 cases to consider for calculation of $\nu_{\lo}(r)$ in the OR-underlying version.
Each Case $j$ correspond to $R_i$ for $i=1,2,\ldots,17$ in Figure \ref{fig:cases-AND-OR}.
{\small
\noindent \textbf{Case 1:}
\begin{multline*}
P(X_2 \in \NPE^r(X_1) \cup \G_1^r(X_1), X_1 \in T_s)=
\left(\int_{0}^{s_0}\int_{0}^{\ell_{am}(x)} +
\int_{s_0}^{s_1}\int_{r_1(x)}^{\ell_{am}(x)} \right)
\frac{A(\msP(A,M_1,M_C,M_3))}{A(\TY)^2}dydx=\\
{\frac{4}{27}}\,r^2-4\,r/9+1/3
\end{multline*}
where
$A(\msP(A,M_1,M_C,M_3))=\sqrt{3}/12
$.

\noindent \textbf{Case 2:}
\begin{multline*}
P(X_2 \in \NPE^r(X_1) \cup \G_1^r(X_1), X_1 \in T_s)=\\
\left(\int_{s_0}^{s_1}\int_{0}^{r_1(x)}+
\int_{s_1}^{s_3}\int_{0}^{r_2(x)}+
\int_{s_3}^{s_4}\int_{0}^{r_5(x)} +
\int_{s_4}^{s_5}\int_{r_3(x)}^{r_5(x)} \right)
\frac{A(\msP(A,M_1,L_2,L_3,M_C,M_3))}{A(\TY)^2}dydx=\\
-{\frac{ \left( r -1 \right)  \left( 1817\,r^7-7807\,r^6+14157\,r^5-14067\,r^4+7893\,r^3-
2475\,r^2+405\,r -27 \right) }{864\,r^6}}
\end{multline*}
where
$A(\msP(A,M_1,L_2,L_3,M_C,M_3))=
{\frac{\sqrt{3} \left( -4\,\sqrt{3}r\, y-12\,r+12\,r\, x+
5\,r^2+3\,y^2+6\,\sqrt{3} y-6\,\sqrt{3} y\, x+9-18\, x+9\, x^2 \right) }
{12\,r^2}}
$.

\noindent \textbf{Case 3:}
\begin{multline*}
P(X_2 \in \NPE^r(X_1) \cup \G_1^r(X_1), X_1 \in T_s)=
\left(\int_{s_4}^{s_5}\int_{0}^{r_3(x)} +
\int_{s_5}^{s_6}\int_{0}^{r_5(x)} \right)
\frac{A(\msP(A,G_2,G_3,M_2,M_C,M_3))}{A(\TY)^2}dydx=\\
{\frac{ \left( 13\,r^4-4\,r^3+4\,r -1-2\,r^2 \right)  \left( r -1 \right)^4}{96\,r^6}}
\end{multline*}
where
$A(\msP(A,G_2,G_3,M_2,M_C,M_3))=
-{\frac{\sqrt{3} \left( \,y^2+2\,\sqrt{3} y-2\,\sqrt{3} y\, x+
3-6\, x+3\, x^2-2\,r^2 \right) }{12\,r^2}}
$.

\noindent \textbf{Case 4:}
\begin{multline*}
P(X_2 \in \NPE^r(X_1) \cup \G_1^r(X_1), X_1 \in T_s)=
\left(\int_{s_1}^{s_2}\int_{r_2(x)}^{\ell_{am}(x)} +
\int_{s_2}^{s_3}\int_{r_2(x)}^{r_5(x)} \right)
\frac{A(\msP(A,M_1,L_2,L_3,L_4,L_5,M_3))}{A(\TY)^2}dydx=\\
{\frac{ \left( 9-72\,r+192\,r^2-192\,r^3+76\,r^4 \right)
\left( 4\,r -3+\sqrt{3} \right)^2 \left( 4\,r -3-\sqrt{3} \right)^2}{10368\,r^6}}
\end{multline*}
where
$A(\msP(A,M_1,L_2,L_3,L_4,L_5,M_3))=
{\frac{\sqrt{3} \left( 4\,\sqrt{3}r\, y+9\,r^2-24\,r+12\,r\, x+15\,y^2-6\,\sqrt{3} y-
6\,\sqrt{3} y\, x+18-18\, x+9\, x^2 \right)}{12\,r^2}}
$.

\noindent \textbf{Case 5:}
\begin{multline*}
P(X_2 \in \NPE^r(X_1) \cup \G_1^r(X_1), X_1 \in T_s)=
\left(\int_{s_5}^{s_6}\int_{r_5(x)}^{r_7(x)} +
\int_{s_6}^{s_9}\int_{0}^{r_7(x)} \right)
\frac{A(\msP(A,G_2,G_3,M_2,M_C,P_2,N_2))}{A(\TY)^2}dydx=\\
{\frac{ \left( -1+2\,r+6\,r^2-6\,r^3+
22\,r^5+17\,r^6 \right)  \left( r -1 \right)^3}{96\,r^6 \left( r+1 \right)^3}}
\end{multline*}
where
$A(\msP(A,G_2,G_3,M_2,M_C,P_2,N_2))=
\Bigl[\sqrt{3} \bigl( -2\,y^2-4\,\sqrt{3} y+4\,\sqrt{3} y\, x-6+12\, x-6\, x^2+
7\,r^2-4\,r^3\sqrt{3} y-12\,r^3 x+8\,r^4\sqrt{3} y\, x+12\,r^4 x^2+
4\,r^4\,y^2 \bigr) \Bigr]\Big/\Bigl[24\,r^2\Bigr]
$.

\noindent \textbf{Case 6:}
\begin{multline*}
P(X_2 \in \NPE^r(X_1) \cup \G_1^r(X_1), X_1 \in T_s)=\\
\left(\int_{s_5}^{s_9}\int_{r_7(x)}^{r_3(x)} +
\int_{s_9}^{s_{12}}\int_{0}^{r_3(x)} +
\int_{s_{12}}^{1/2}\int_{0}^{r_6(x)} \right)
\frac{A(\msP(A,N_1,Q_1,G_3,M_2,M_C,P_2,N_2))}{A(\TY)^2}dydx=\\
-{\frac{81\,r^9-189\,r^8+561\,r^7-45\,r^6-1894\,r^5-18\,r^4+1912\,r^3+
224\,r^2-384\,r -128}{1296\, \left( r+1 \right)^3 r^4}}
\end{multline*}
where
$A(\msP(A,N_1,Q_1,G_3,M_2,M_C,P_2,N_2))=
\Bigl[\sqrt{3} \bigl( 4\,r\,y^2-4\,\sqrt{3} y+12\, x+13\,r -12+18\,r^3 x^2+
12\,r\, x-12\,r\, x^2-8\,\sqrt{3}r^2 y+4\,\sqrt{3}r\, y-24\,r^2 x+
12\,\sqrt{3}r^3 y\, x+6\,r^3\,y^2 \bigr) \Bigr]\Big/\Bigl[24\,r\Bigr]
$.

\noindent \textbf{Case 7:}
\begin{multline*}
P(X_2 \in \NPE^r(X_1) \cup \G_1^r(X_1), X_1 \in T_s)=\\
\left(\int_{s_8}^{s_5}\int_{r_8(x)}^{r_2(x)} +
\int_{s_5}^{s_{10}}\int_{r_3(x)}^{r_2(x)} +
\int_{s_{10}}^{s_{12}}\int_{r_3(x)}^{r_6(x)} \right)
\frac{A(\msP(A,N_1,Q_1,L_3,M_C,P_2,N_2))}{A(\TY)^2}dydx=\\
-\Bigl[128-1536\,r -302592\,r^7+11753\,r^{12}+346171\,r^8-28416\,r^3+
8384\,r^2+69760\,r^4+220201\,r^6-135936\,r^5-305664\,r^9+\\
186683\,r^{10}-69120\,r^{11}\Bigr]\Big/\Bigl[1944\, \left( r^2+1 \right)^3r^6\Bigr]
\end{multline*}
where
$A(\msP(A,N_1,Q_1,L_3,M_C,P_2,N_2))=
\Bigl[\sqrt{3} \bigl( -4\,\sqrt{3}r\, y+2\,\sqrt{3}
r^2 y-12\, x-12\,r+8\,r^2+12\,r\, x-6\, x^2r^2+2\,r^2\,y^2-
4\,\sqrt{3} y\, x+3\,r^4\,y^2-4\,r^3\sqrt{3} y-12\,r^3 x+9\,r^4 x^2+
4\,\sqrt{3} y+6\,r^4\sqrt{3} y\, x+6\, x^2+6\,r^2 x+6+2\,y^2 \bigr) \Bigr]\Big/\Bigl[12\,r^2\Bigr]
$.

\noindent \textbf{Case 8:}
\begin{multline*}
P(X_2 \in \NPE^r(X_1) \cup \G_1^r(X_1), X_1 \in T_s)=
\left(\int_{s_3}^{s_8}\int_{r_5(x)}^{r_2(x)} +
\int_{s_8}^{s_5}\int_{r_5(x)}^{r_8(x)} \right)
\frac{A(\msP(A,N_1,P_1,L_2,L_3,M_C,P_2,N_2))}{A(\TY)^2}dydx=\\
{\frac{ \left( 895\,r^8-2472\,r^7+3363\,r^6-2880\,r^5+2220\,r^4-1296\,r^3+675\,r^2-
216\,r+27 \right)  \left( -12\,r+7\,r^2+3 \right)^2}{7776\, \left( r^2+1 \right)^3r^6}}
\end{multline*}
where
$A(\msP(A,N_1,P_1,L_2,L_3,M_C,P_2,N_2))=
\Bigl[\sqrt{3} \bigl( 4\,r^4\,y^2+8\,r^4 \sqrt{3} y\, x+12\,r^4 x^2-
4\,r^3 \sqrt{3} y-12\,r^3 x-4\,\sqrt{3}r\, y-12\,r+12\,r\, x+8\,r^2+
3\,y^2+6\,\sqrt{3}y-6\,\sqrt{3} y\, x+9-18\, x+9\, x^2 \bigr) \Bigr]\Big/\Bigl[ 12\, r^2\Bigr]
$.

\noindent \textbf{Case 9:}
\begin{multline*}
P(X_2 \in \NPE^r(X_1) \cup \G_1^r(X_1), X_1 \in T_s)=\\
\left(\int_{s_2}^{s_3}\int_{r_5(x)}^{\ell_{am}(x)} +
\int_{s_3}^{s_7}\int_{r_2(x)}^{\ell_{am}(x)} +
\int_{s_7}^{s_8}\int_{r_2(x)}^{r_8(x)} \right)
\frac{A(\msP(A,N_1,P_1,L_2,L_3,L_4,L_5,P_2,N_2))}{A(\TY)^2}dydx=\\
-\Bigl[355328\,r^{18}-2204160\,r^{17}+
6591792\,r^{16}-13254912\,r^{15}+20639832\,r^{14}-26417664
\,r^{13}+28578916\,r^{12}-26760576\,r^{11}+\\
21960774\,r^{10}-15877152\,r^9+10180620\,r^8-5753232\,r^7+
2856483\,r^6-1222128\,r^5+438777\,r^4-128304\,r^3+28107\,r^2-\\
3888\,r+243\Bigr]\Big/\Bigl[7776\, \left( r^2+1 \right)^3 \left( 2\,r^2+1 \right)^3r^6\Bigr]
\end{multline*}
where
$A(\msP(A,N_1,P_1,L_2,L_3,L_4,L_5,P_2,N_2))=
\Bigl[\sqrt{3} \bigl( 18+4\,\sqrt{3}r\, y-18\,x-24\,r+12\,r^2+12\,r\, x-6\,\sqrt{3} y+
8\,r^4\sqrt{3} y\, x-12\,r^3 x+12\,r^4 x^2+9\, x^2+15\,y^2+4\,r^4 y^2-
4\,r^3\sqrt{3} y-6\,\sqrt{3} y\,x \bigr) \Bigr]\Big/\Bigl[12\,r^2\Bigr]
$.

\noindent \textbf{Case 10:}
\begin{multline*}
P(X_2 \in \NPE^r(X_1) \cup \G_1^r(X_1), X_1 \in T_s)=
\left(\int_{s_7}^{s_8}\int_{r_8(x)}^{r_9(x)}+
\int_{s_8}^{s_{10}}\int_{r_2(x)}^{r_9(x)} \right)
\frac{A(\msP(A,N_1,Q_1,L_3,L_4,L_5,P_2,N_2))}{A(\TY)^2}dydx=\\
\Bigl[ 8\bigl( 288\,r^{12}-864\,r^{11}+1486\,r^{10}-1896\,r^9+2056\,r^8-1608\,r^7+1189\,r^6-
654\,r^5+317\,r^4-132\,r^3+44\,r^2-12\,r+2 \bigr)  \left( 2\,r -1 \right)^2 \\
\left( r -1 \right)^2\Bigr]\Big/\Bigl[243\, \left( r^2+1 \right)^3 \left( 2\,r^2+1 \right)^3 r^4\Bigr]
\end{multline*}
where
$A(\msP(A,N_1,Q_1,L_3,L_4,L_5,P_2,N_2))=
\Bigl[\sqrt{3} \bigl( 4\,\sqrt{3}r\, y+2\,\sqrt{3}{
r}^2 y-8\,\sqrt{3} y-12\, x-24\,r+12\,r^2+12\,r\, x-6\, x^2r^2+15-
12\,r^3 x+9\,r^4 x^2+6\, x^2+6\,r^2 x+6\,r^4\sqrt{3} y\, x+2\,r^2\,y^2-
4\,\sqrt{3} y\, x+3\,r^4\,y^2-4\,r^3\sqrt{3} y+14\,y^2 \bigr) \Bigr]\Big/\Bigl[12\,r^2\Bigr]
$.

\noindent \textbf{Case 11:}
\begin{multline*}
P(X_2 \in \NPE^r(X_1) \cup \G_1^r(X_1), X_1 \in T_s)=
\left(\int_{s_{12}}^{s_{13}}\int_{r_6(x)}^{r_3(x)} +
\int_{s_{13}}^{1/2}\int_{r_6(x)}^{r_2(x)} \right)
\frac{A(\msP(A,N_1,Q_1,G_3,M_2,N_3,N_2))}{A(\TY)^2}dydx=\\
-{\frac{1536-6528\,r^2+133834\,r^8-48240\,r^9+95616\,r^4-20736\,r^3-158976\,r^5-
200064\,r^7+196680\,r^6+7107\,r^{10}}{15552\,r^4}}
\end{multline*}
where
$A(\msP(A,N_1,Q_1,G_3,M_2,N_3,N_2))=
\Bigl[\sqrt{3} \bigl( 4\,r\,y^2+12\, x+9\,r -12+9\,r^3 x^2+12\,r\, x-
12\,r\,x^2-4\,\sqrt{3}r^2 y+4\,\sqrt{3}r\, y+6\,\sqrt{3}r^3 y\, x+
3\,r^3\,y^2-12\,r^2 x-4\,\sqrt{3} y \bigr) \Bigr]\Big/\Bigl[24\,r\Bigr]
$.

\noindent \textbf{Case 12:}
\begin{multline*}
P(X_2 \in \NPE^r(X_1) \cup \G_1^r(X_1), X_1 \in T_s)=
\left(\int_{s_{10}}^{s_{13}}\int_{r_6(x)}^{r_2(x)} +
\int_{s_{12}}^{s_{13}}\int_{r_3(x)}^{r_2(x)} \right)
\frac{A(\msP(A,N_1,Q_1,L_3,N_3,N_2))}{A(\TY)^2}dydx=\\
{\frac{ \left( 147\,r^8-504\,r^7+530
\,r^6-336\,r^5+876\,r^4-1056\,r^3+896\,r^2-384\,r+64 \right)  \left( -12\,r+7\,r^2+4 \right)^2}{15552\,r^6}}
\end{multline*}
where
$A(\msP(A,N_1,Q_1,L_3,N_3,N_2))=
\Bigl[\sqrt{3} \bigl( 4\,y^2-8\,\sqrt{3} y\, x-24\, x-24\,r+8\,\sqrt{3} y+12\,r^2+
4\,\sqrt{3}r^2 y+6\,r^4\sqrt{3} y\, x+24\,r\, x-4\,r^3\sqrt{3} y+3\,r^4\,y^2-
8\,\sqrt{3}r\, y-12\, x^2r^2-12\,r^3 x+9\,r^4 x^2+12\, x^2+12\,r^2 x+
4\,r^2\,y^2+12 \bigr) \Bigr]\Big/\Bigl[24\,r^2\Bigr]
$.

\noindent \textbf{Case 13:}
\begin{multline*}
P(X_2 \in \NPE^r(X_1) \cup \G_1^r(X_1), X_1 \in T_s)=\\
\left(\int_{s_{10}}^{s_{14}}\int_{r_2(x)}^{r_{10}(x)} +
\int_{s_{14}}^{s_{13}}\int_{r_2(x)}^{r_{12}(x)} +
\int_{s_{13}}^{1/2}\int_{r_3(x)}^{r_{12}(x)} \right)
\frac{A(\msP(A,N_1,Q_1,L_3,N_3,N_2))}{A(\TY)^2}dydx=\\
\Bigl[1024-12288\,r+295680\,r^7+1053\,r^{12}-197140\,r^8+626688\,r^3-100864\,r^2-
1294848\,r^4-686528\,r^6+1282560\,r^5+\\
114336\,r^9-30930\,r^{10} \Bigr]\Big/\Bigl[31104\,r^4\Bigr]
\end{multline*}
where
$A(\msP(A,N_1,Q_1,L_3,N_3,N_2))=
\Bigl[\sqrt{3} \bigl( 4\,y^2-8\,\sqrt{3} y\, x-24\, x-24\,r+8\,\sqrt{3} y+12\,r^2+
4\,\sqrt{3}r^2 y+6\,r^4\sqrt{3} y\, x+24\,r\, x-4\,r^3\sqrt{3} y+3\,r^4\,y^2-
8\,\sqrt{3}r\, y-12\, x^2r^2-12\,r^3 x+9\,r^4 x^2+12\, x^2+12\,r^2 x+
4\,r^2\,y^2+12 \bigr)  \Bigr]\Big/\Bigl[24\,r^2\Bigr]
$.

\noindent \textbf{Case 14:}
\begin{multline*}
P(X_2 \in \NPE^r(X_1) \cup \G_1^r(X_1), X_1 \in T_s)=\\
\left(\int_{s_7}^{s_{11}}\int_{r_9(x)}^{\ell_{am}(x)} +
\int_{s_{11}}^{s_{10}}\int_{r_9(x)}^{r_{12}(x)} +
\int_{s_{10}}^{s_{14}}\int_{r_{10}(x)}^{r_{12}(x)} \right)
\frac{A(\msP(A,N_1,Q_1,L_3,L_4,Q_2,N_2))}{A(\TY)^2}dydx=\\
-\Bigl[ \left( r -1 \right)  \bigl( 1512\,r^{17}+1512\,r^{16}-16740\,r^{15}+540\,r^{14}+
84078\,r^{13}-83538\,r^{12}-164835\,r^{11}+401085\,r^{10}-487872\,r^9+\\
535728\,r^8-463124\,r^7+335596\,r^6-197440\,r^5+64640\,r^4-7936\,r^3-1792\,r^2+5632\,r-
512 \bigr)\Bigr]\Big/\Bigl[5184\, \left( 2\,r^2+1 \right)^3r^4\Bigr]
\end{multline*}
where
$A(\msP(A,N_1,Q_1,L_3,L_4,Q_2,N_2))=
\Bigl[\sqrt{3} \bigl( -6\, x-12\,r+6\,r^2+6\,r\, x+2\,\sqrt{3}r^2 y-r^2\,y^2-
2\,\sqrt{3} y\, x+r^4\,y^2+5\,y^2-2\,r^2 x\,\sqrt{3} y+2\,r^4\sqrt{3} y\, x+
2\,\sqrt{3}r\, y-2\,r^3\sqrt{3} y-3\,x^2r^2-6\,r^3 x+3\,r^4 x^2-2\,\sqrt{3} y+
3\, x^2+6\,r^2 x+6 \bigr) \Bigr]\Big/\Bigl[ 6\,r^2\Bigr]
$.

\noindent \textbf{Case 15:}
\begin{multline*}
P(X_2 \in \NPE^r(X_1) \cup \G_1^r(X_1), X_1 \in T_s)=
\int_{s_{13}}^{1/2}\int_{r_2(x)}^{r_3(x)} \frac{A(\msP(A,N_1,Q_1,G_3,M_2,N_3,N_2))}{A(\TY)^2}dydx=\\
{\frac{ \left( 147\,r^5-612\,r^4+980\,r^3-768\,r^2+744\,r -288 \right)  \left( -6+5\,r \right)^2}{7776\,r}}
\end{multline*}
where
$A(\msP(A,N_1,Q_1,L_3,L_4,Q_2,N_2))=
\Bigl[\sqrt{3} \bigl( 4\,r\,y^2+12\, x+9\,r -12+9\,r^3 x^2+12\,r\, x-12\,r\,x^2-
4\,\sqrt{3}r^2 y+4\,\sqrt{3}r\, y+6\,\sqrt{3}r^3 y\, x+3\,r^3\,y^2-12\,r^2 x-
4\,\sqrt{3} y \bigr) \Bigr]\Big/\Bigl[ 24\,r\Bigr]
$.

\noindent \textbf{Case 16:}
\begin{multline*}
P(X_2 \in \NPE^r(X_1) \cup \G_1^r(X_1), X_1 \in T_s)=
\int_{s_{14}}^{1/2}\int_{r_{12}(x)}^{r_{10}(x)} \frac{A(\msP(A,N_1,Q_1,L_3,N_3,N_2))}{A(\TY)^2}dydx=\\
-{\frac{ \left( 13\,r^8+52\,r^7+10\,r^6-184\,r^5+60\,r^4+624\,r^3-48\,r^2-832\,r+448 \right)
\left( -2+r \right)  \left( r^2+2\,r -4 \right)^2}{384\, \left( r+2 \right)^3r^2}}
\end{multline*}
where
$A(\msP(A,N_1,Q_1,L_3,N_3,N_2))=
\Bigl[\sqrt{3} \bigl( 4\,y^2-8\,\sqrt{3} y\, x-24\, x-24\,r+8\,\sqrt{3} y+12\,r^2+
4\,\sqrt{3}r^2 y+6\,r^4\sqrt{3} y\, x+24\,r\, x-4\,r^3\sqrt{3} y+3\,r^4\,y^2-
8\,\sqrt{3}r\, y-12\, x^2r^2-12\,r^3 x+9\,r^4 x^2+12\, x^2+12\,r^2 x+
4\,r^2\,y^2+12 \bigr)\Bigr]\Big/\Bigl[24\,r^2\Bigr]
$.

\noindent \textbf{Case 17:}
\begin{multline*}
P(X_2 \in \NPE^r(X_1) \cup \G_1^r(X_1), X_1 \in T_s)=
\left(\int_{s_{11}}^{s_{14}}\int_{r_{12}(x)}^{\ell_{am}(x)} +
\int_{s_{14}}^{1/2}\int_{r_{10}(x)}^{\ell_{am}(x)}\right)
\frac{A(\msP(A,N_1,Q_1,L_3,L_4,Q_2,N_2))}{A(\TY)^2}dydx=\\
\Bigl[\bigl( 189\,r^{12}+1323\,r^{11}+1026\,r^{10}-10692\,r^9-14364\,r^8+51732\,r^7+
64664\,r^6-183952\,r^5-153504\,r^4+398080\,r^3+124928\,r^2-\\
470528\,r+197632 \bigr)  \left( r -1 \right) \Bigr]\Big/\Bigl[5184\,r^2 \left( r+2 \right)^3\Bigr]
\end{multline*}
where
$A(\msP(A,N_1,Q_1,L_3,N_3,N_2))=
\Bigl[\sqrt{3} \bigl( -6\, x-12\,r+6\,r^2+6\,r\, x+2\,\sqrt{3}r^2 y-r^2\,y^2-
2\,\sqrt{3} y\, x+r^4\,y^2+5\,y^2-2\,r^2 x\,\sqrt{3} y+2\,r^4\sqrt{3} y
\, x+2\,\sqrt{3}r\, y-2\,r^3\sqrt{3} y-3\, x^2r^2-6\,r^3 x+3\,r^4 x^2-
2\,\sqrt{3} y+3\, x^2+6\,r^2 x+6 \bigr)\Bigr]\Big/\Bigl[6\,r^2\Bigr]
$.
}

Adding up the $P(X_2 \in \NPE^r(X_1) \cup \G_1^r(X_1), X_1 \in T_s)$
values in the 17 possible cases above, and multiplying by 6
we get for $r \in [1,4/3)$,
$$\nu_\lo(r)=
{\frac{860\,r^4-195\,r^5-256+720\,r -
846\,r^3-108\,r^2+47\,r^6}{108\,r^2 \left( r+2 \right)  \left( r+1 \right) }}.$$
The $\nu_\lo(r)$ values for the other intervals can be calculated similarly.

\subsection*{Derivation of $\nu_\lo(r)$ in Theorem \ref{thm:asy-norm-under}}

By symmetry,
$P(\{X_2,X_3\} \subset \NPE^r(X_1)\cup \G_1^r(X_1))=
6\,P(\{X_2,X_3\} \subset \NPE^r(X_1)\cup \G_1^r(X_1),\; X_1 \in T_s)$.
For $r \in \bigl[6/5,\sqrt{5}-1\bigr)$,
there are 17 cases to consider for calculation of $\nu_{\lo}(r)$ in the OR-underlying version
(see also Figure \ref{fig:cases-AND-OR}):
{\small
\noindent \textbf{Case 1:}
\begin{multline*}
P(\{X_2,X_3\} \subset \NPE^r(X_1)\cup \G_1^r(X_1),\; X_1 \in T_s)=
\left(\int_{0}^{s_0}\int_{0}^{\ell_{am}(x)} +
\int_{s_0}^{s_1}\int_{r_1(x)}^{\ell_{am}(x)} \right)
\frac{A(\msP(A,M_1,M_C,M_3))^2}{A(\TY)^3}dydx=\\
{\frac{4}{81}}\,r^2-{\frac{4}{27}}\,r+1/9
\end{multline*}
where
$A(\msP(A,M_1,M_C,M_3))=
1/12\,\sqrt{3}
$.

\noindent \textbf{Case 2:}
\begin{multline*}
P(\{X_2,X_3\} \subset \NPE^r(X_1)\cup \G_1^r(X_1),\; X_1 \in T_s)=\\
\left(\int_{s_0}^{s_1}\int_{0}^{r_1(x)} +
\int_{s_1}^{s_3}\int_{0}^{r_2(x)} +
\int_{s_3}^{s_4}\int_{0}^{r_5(x)} +
\int_{s_4}^{s_5}\int_{r_3(x)}^{r_5(x)}\right)
\frac{A(\msP(A,M_1,L_2,L_3,M_C,M_3))^2}{A(\TY)^3}dydx=\\
-\Bigl[ \left( r -1 \right)  \bigl( 119155\,r^{11}-845345\,r^{10}+2724777\,r^9-5206743\,r^8+
6475257\,r^7-5454855\,r^6+3155193\,r^5-1249479\,r^4+\\
332181\,r^3-56619\,r^2+5589\,r -243 \bigr) \Bigr]\Big/\Bigl[25920\,r^{10}\Bigr]
\end{multline*}
where
$A(\msP(A,M_1,L_2,L_3,M_C,M_3))=
{\frac {\sqrt{3} \left( -4\,\sqrt{3}r\,y-12\,r+12\,r\,x+
5\,r^2+3\,y^2+6\,\sqrt{3}y-6\,\sqrt{3}y\,x+9-18\,x+9\,x^2 \right)}
{12\,r^2}}
$.

\noindent \textbf{Case 3:}
\begin{multline*}
P(\{X_2,X_3\} \subset \NPE^r(X_1)\cup \G_1^r(X_1),\; X_1 \in T_s)=
\left(\int_{s_4}^{s_5}\int_{0}^{r_3(x)} +
\int_{s_5}^{s_6}\int_{0}^{r_5(x)}\right)
\frac{A(\msP(A,G_2,G_3,M_2,M_C,M_3))^2}{A(\TY)^3}dydx=\\
{\frac{ \left( 215\,r^8-136\,r^7-56\,r^6+172\,r^5-55\,r^4-60\,r^3+66\,r^2-24\,r+
3 \right)  \left( r -1 \right)^4}{2880\,r^{10}}}
\end{multline*}
where
$A(\msP(A,G_2,G_3,M_2,M_C,M_3))=
-{\frac {\sqrt{3} \left( y^2+2\,\sqrt{3}y-2\,\sqrt{3}y\,x+3-
6\,x+3\,x^2-2\,r^2 \right) }{12\,r^2}}
$.

\noindent \textbf{Case 4:}
\begin{multline*}
P(\{X_2,X_3\} \subset \NPE^r(X_1)\cup \G_1^r(X_1),\; X_1 \in T_s)=\\
\left(\int_{s_1}^{s_2}\int_{r_2(x)}^{\ell_{am}(x)} +
\int_{s_2}^{s_3}\int_{r_2(x)}^{r_5(x)} \right)
\frac{A(\msP(A,M_1,L_2,L_3,L_4,L_5,M_3))^2}{A(\TY)^3}dydx=\\
\Bigl[ \bigl( 37072\,r^8-195072\,r^7+453120\,r^6-589248\,r^5+460728\,r^4-
217728\,r^3+60480\,r^2-9072\,r+567 \bigr) \\
\left( 4\,r -3+\sqrt{3} \right)^2 \left( 4\,r -3-\sqrt{3} \right)^2\Bigr]\Big/\Bigl[1866240\,r^{10}\Bigr]
\end{multline*}
where
$A(\msP(A,M_1,L_2,L_3,L_4,L_5,M_3))=
{\frac {\sqrt{3} \left( 4\,\sqrt{3}r\,y+9\,r^2-
24\,\nu+12\,r\,x+15\,y^2-6\,\sqrt{3}y-6\,
\sqrt{3}y\,x+18-18\,x+9\,x^2 \right)}{12\,r^2}}
$.

\noindent \textbf{Case 5:}
\begin{multline*}
P(\{X_2,X_3\} \subset \NPE^r(X_1)\cup \G_1^r(X_1),\; X_1 \in T_s)=\\
\left(\int_{s_5}^{s_6}\int_{r_5(x)}^{r_7(x)} +
\int_{s_6}^{s_9}\int_{0}^{r_7(x)} \right)
\frac{A(\msP(A,G_2,G_3,M_2,M_C,P_2,N_2))^2}{A(\TY)^3}dydx=\\
{\frac{ \left( 3-12\,r -15\,r^2+84\,r^3+18\,r^4-232\,r^5+130\,r^6+504\,r^7-108\,r^8-
288\,r^9+623\,r^{10}+920\,r^{11}+373\,r^{12} \right)
\left( r -1 \right)^3}{2880\,r^{10} \left( r+1 \right)^5}}
\end{multline*}
where
$A(\msP(A,G_2,G_3,M_2,M_C,P_2,N_2))=
\Bigl[\sqrt{3} \bigl( -2\,y^2-4\,\sqrt{3}y+4\,\sqrt{3}y\,x-6+12\,x-6\,x^2+
7\,r^2-4\,r^3\sqrt{3}y-12\,r^3x+8\,r^4\sqrt{3}y\,x+12\,r^4x^2+
4\,r^4 y^2 \bigr) \Bigr]\Big/\Bigl[24\,r^2\Bigr]
$.

\noindent \textbf{Case 6:}
\begin{multline*}
P(\{X_2,X_3\} \subset \NPE^r(X_1)\cup \G_1^r(X_1),\; X_1 \in T_s)=\\
\left(\int_{s_5}^{s_9}\int_{r_7(x)}^{r_3(x)} +
\int_{s_9}^{s_{12}}\int_{0}^{r_3(x)} +
\int_{s_{12}}^{1/2}\int_{0}^{r_6(x)} \right)
\frac{A(\msP(A,N_1,Q_1,G_3,M_2,M_C,P_2,N_2))^2}{A(\TY)^3}dydx=\\
-\Bigl[19683\,r^{15}-59049\,r^{14}+83106\,r^{13}+167670\,r^{12}-211626\,r^{11}+
344466\,r^{10}-142614\,r^9-2573586\,r^8-128853\,r^7+\\
3465675\,r^6+1103824\,r^5-1473304\,r^4-730880\,r^3+107776\,r^2+
158720\,r+31744\Bigr]\Big/\Bigl[1049760\, \left( r+1 \right)^5r^6\Bigr]
\end{multline*}
where
$A(\msP(A,N_1,Q_1,G_3,M_2,M_C,P_2,N_2))=
\Bigl[\sqrt{3} \bigl( 4\,r\,y^2+12\,x+13\,r+12\,r\,x-4\,\sqrt{3}y-12+
4\,\sqrt{3}r\,y-8\,\sqrt{3}r^2y+18\,x^2r^3-12\,r\,x^2+6\,r^3y^2-
24\,r^2 x+12\,\sqrt{3}r^3y\,x \bigr)\Bigr]\Big/\Bigl[24\,r\Bigr]
$.

\noindent \textbf{Case 7:}
\begin{multline*}
P(\{X_2,X_3\} \subset \NPE^r(X_1)\cup \G_1^r(X_1),\; X_1 \in T_s)=\\
\left(\int_{s_8}^{s_5}\int_{r_8(x)}^{r_2(x)} +
\int_{s_5}^{s_{10}}\int_{r_3(x)}^{r_2(x)} +
\int_{s_{10}}^{s_{12}}\int_{r_3(x)}^{r_6(x)} \right)
\frac{A(\msP(A,N_1,Q_1,L_3,M_C,P_2,N_2))^2}{A(\TY)^3}dydx=\\
-\Bigl[6144-110592\,r -310846464\,r^7+2127553557\,r^{12}+570050560\,r^8-
5031936\,r^3+936960\,r^2+19526656\,r^4+147203072\,r^6+\\
7627473\,r^{20}+1419072042\,r^{16}-762467328\,r^{17}+288811029\,r^{18}
-68327424\,r^{19}-59166720\,r^5-923627520\,r^9+\\
1340817105\,r^{10}-1765251072\,r^{11}-2350015488\,r^{13}+
2339575338\,r^{14}-2016377856\,r^{15}\Bigr]\Big/\Bigl[ 262440\,\left( r^2+1 \right)^5r^{10}\Bigr]
\end{multline*}
where
$A(\msP(A,N_1,Q_1,L_3,M_C,P_2,N_2))=
\Bigl[\sqrt{3} \bigl( -4\,\sqrt{3}r\,y+2\,\sqrt{3}r^2 y-6\,x^2r^2-12\,x-12\,r-
12\,r^3 x+9\,r^4x^2+8\,r^2+12\,r\,x+6\,x^2+6\,r^4\sqrt{3}y\,x+2\,r^2 y^2-
4\,\sqrt{3}y\,x+3\,r^4 y^2-4\,r^3\sqrt{3}y+4\,\sqrt{3}y+2\,y^2+
6\,r^2x+6 \bigr)\Bigr]\Big/\Bigl[12\,r^2\Bigr]
$.

\noindent \textbf{Case 8:}
\begin{multline*}
P(\{X_2,X_3\} \subset \NPE^r(X_1)\cup \G_1^r(X_1),\; X_1 \in T_s)=\\
\left(\int_{s_3}^{s_8}\int_{r_5(x)}^{r_2(x)} +
\int_{s_8}^{s_5}\int_{r_5(x)}^{r_8(x)} \right)
\frac{A(\msP(A,N_1,P_1,L_2,L_3,M_C,P_2,N_2))^2}{A(\TY)^3}dydx=\\
\Bigl[\bigl( 426497\,r^{16}-2443992\,r^{15}+6726107\,r^{14}-11753232\,r^{13}+
15220771\,r^{12}-16367448\,r^{11}+15754449\,r^{10}-13773024\,r^9+\\
10839672\,r^8-7552440\,r^7+4592889\,r^6-2374272\,r^5+
1018899\,r^4-344088\,r^3+81891\,r^2-11664\,r+729 \bigr) \\
\left( -12\,r+7\,r^2+3 \right)^2\Bigr]\Big/\Bigl[699840\, \left( r^2+1 \right)^5r^{10}\Bigr]
\end{multline*}
where
$A(\msP(A,N_1,P_1,L_2,L_3,M_C,P_2,N_2))=
\Bigl[\sqrt{3} \bigl( -4\,r^3\sqrt{3}y-12\,r^3 x+8\,r^4\sqrt{3}y\,x+12\,r^4 x^2+
4\,r^4y^2-4\,\sqrt{3}r\,y-12\,r+12\,r\,x+3\,y^2+6\,\sqrt{3}y-6\,\sqrt{3}y\,x+
8\,r^2+9-18\,x+9\,x^2 \bigr)\Bigr]\Big/\Bigl[12\,r^2\Bigr]
$.

\noindent \textbf{Case 9:}
\begin{multline*}
P(\{X_2,X_3\} \subset \NPE^r(X_1)\cup \G_1^r(X_1),\; X_1 \in T_s)=\\
\left(\int_{s_2}^{s_3}\int_{r_5(x)}^{\ell_{am}(x)} +
\int_{s_3}^{s_7}\int_{r_2(x)}^{\ell_{am}(x)} +
\int_{s_7}^{s_8}\int_{r_2(x)}^{r_8(x)} \right)
\frac{A(\msP(A,N_1,P_1,L_2,L_3,L_4,L_5,P_2,N_2))^2}{A(\TY)^3}dydx=\\
-\Bigl[15309-367416\,r+60475010560\,r^{28}+437704472832\,r^{26}+1444872192\,r^{30}-
13250101248\,r^{29}-185909870592\,r^{27}+\\
4148739\,r^2-2027754648576\,r^{23}+1397612375040\,r^{24}+20429177589\,r^8-
677278256112\,r^{13}-49656902904\,r^9+\\
159963012\,r^4-30005640\,r^3-681714144\,r^5-7515142416\,r^7-3097406755584\,r^{21}-2609245249920\,r^{17}+\\
3051035360256\,r^{18}-3315184235136\,r^{19}+3337272236928\,r^{20}+2631941507968\,r^{22}+
2435971806\,r^6+\\
109069315047\,r^{10}-218273842152\,r^{11}+400534503738\,r^{12}+1059615993384\,r^{14}-
1538314485120\,r^{15}+\\
2076627064432\,r^{16}-845838600192\,r^{25}\Bigr]\Big/\Bigl[1399680\, \left(
r^2+1 \right)^5 \left( 2\,r^2+1 \right)^5r^{10}\Bigr]
\end{multline*}
where
$A(\msP(A,N_1,P_1,L_2,L_3,L_4,L_5,P_2,N_2))=
\Bigl[\sqrt{3} \bigl( 18-18\,x-24\,r-12\,r^3 x+12\,r^4x^2+12\,r^2+12\,r\,x+
4\,\sqrt{3}r\,y-4\,r^3\sqrt{3}y+4\,r^4 y^2-6\,\sqrt{3}y\,x+8\,r^4\sqrt{3} y\,x+
9\,x^2+15\,y^2-6\,\sqrt{3}y \bigr) \Bigr]\Big/\Bigl[12\,r^2\Bigr]
$.

\noindent \textbf{Case 10:}
\begin{multline*}
P(\{X_2,X_3\} \subset \NPE^r(X_1)\cup \G_1^r(X_1),\; X_1 \in T_s)=\\
\left(\int_{s_7}^{s_8}\int_{r_8(x)}^{r_9(x)} +
\int_{s_8}^{s_{10}}\int_{r_2(x)}^{r_9(x)} \right)
\frac{A(\msP(A,N_1,Q_1,L_3,L_4,L_5,P_2,N_2))^2}{A(\TY)^3}dydx=\\
\Bigl[64\,\bigl( 12-144\,r+924\,r^2-683328\,r^{23}+112976\,r^{24}+757211\,r^8-
10554918\,r^{13}-1513230\,r^9+16242\,r^4-4320\,r^3-51372\,r^5-\\
344988\,r^7-4867848\,r^{21}-18583080\,r^{17}+16493828\,r^{18}-
12883116\,r^{19}+8668124\,r^{20}+2177536\,r^{22}+141366\,r^6+2774371\,r^{10}-\\
4692510\,r^{11}+7331714\,r^{12}+14002613\,r^{14}-16948218\,r^{15}+18708475\,r^{
16} \bigr)  \left( r -1 \right)^2 \left( 2\,r -1 \right)^2\Bigr]\Big/\Bigl[32805\,\left( r^2+1 \right)^5 \left( 2\,r^2+1 \right)^5r^8\Bigr]
\end{multline*}
where
$A(\msP(A,N_1,Q_1,L_3,L_4,L_5,P_2,N_2))=
\Bigl[\sqrt{3} \bigl( 2\,\sqrt{3}r^2y+15-6\,x^2r^2-12\,x-24\,r-12\,r^3x+
9\,r^4 x^2+12\,r^2+12\,r\,x-8\,\sqrt{3} y+6\,x^2+6\,r^4\sqrt{3}y\,x+
14\,y^2-4\,\sqrt{3}y\,x+2\,r^2y^2-4\,r^3\sqrt{3}y+3\,r^4y^2+6\,r^2 x+
4\,\sqrt{3}r\,y \bigr)\Bigr]\Big/\Bigl[12\,r^2\Bigr]
$.

\noindent \textbf{Case 11:}
\begin{multline*}
P(\{X_2,X_3\} \subset \NPE^r(X_1)\cup \G_1^r(X_1),\; X_1 \in T_s)=\\
\left(\int_{s_{12}}^{s_{13}}\int_{r_6(x)}^{r_3(x)} +
\int_{s_{13}}^{1/2}\int_{r_6(x)}^{r_2(x)} \right)
\frac{A(\msP(A,N_1,Q_1,G_3,M_2,N_3,N_2))^2}{A(\TY)^3}dydx=\\
-\Bigl[-253952+1529856\,r^2+601574256\,r^8-385780320\,r^{13}-776518272\,r^9+
7803648\,r^4-70917120\,r^5-396524160\,r^7+\\
209710080\,r^6+869661288\,r^{10}-845940960\,r^{11}+668092108\,r^{12}+
147067614\,r^{14}-32610600\,r^{15}+3173067\,r^{16}\Bigr]\Big/\Bigl[8398080\,r^6\Bigr]
\end{multline*}
where
$A(\msP(A,N_1,Q_1,G_3,M_2,N_3,N_2))=
\Bigl[\sqrt{3} \bigl( 4\,r\,y^2+12\,x+4\,\sqrt{3}r\,y+9\,r-4\,\sqrt{3}y+12\,r\,x-
12+9\,x^2r^3+6\,\sqrt{3}r^3y\,x-12\,r\,x^2-4\,\sqrt{3}r^2y-12\,r^2 x+3\,r^3y^2 \bigr)\Bigr]\Big/\Bigl[24\,r\Bigr]
$.

\noindent \textbf{Case 12:}
\begin{multline*}
P(\{X_2,X_3\} \subset \NPE^r(X_1)\cup \G_1^r(X_1),\; X_1 \in T_s)=\\
\left(\int_{s_{10}}^{s_{12}}\int_{r_6(x)}^{r_2(x)} +
\int_{s_{12}}^{s_{13}}\int_{r_3(x)}^{r_2(x)} \right)
\frac{A(\msP(A,N_1,Q_1,L_3,N_3,N_2))^2}{A(\TY)^3}dydx=\\
\Bigl[\bigl( 64827\,r^{16}-444528\,r^{15}+1223334\,r^{14}-1793232\,r^{13}+
1839416\,r^{12}-2003712\,r^{11}+2286224\,r^{10}-2421504\,r^9+3095088\,r^8-\\
4428288\,r^7+5889152\,r^6-6093312\,r^5+4557056\,r^4-2340864\,r^3+774144\,r^2-147456\,r+12288 \bigr)\\
 \left( -12\,r+7\,r^2+4 \right)^2\Bigr]\Big/\Bigl[8398080\,r^{10}\Bigr]
\end{multline*}
where
$A(\msP(A,N_1,Q_1,L_3,N_3,N_2))=
\Bigl[\sqrt{3} \bigl( -12\,x^2r^2-24\,x-24\,r-12\,r^3x+9\,r^4x^2+4\,y^2-
8\,\sqrt{3}r\,y+6\,r^4\sqrt{3}y\,x+8\,\sqrt{3}y+12\,r^2+24\,r\,x+
12\,x^2-8\,\sqrt{3}y\,x+4\,r^2y^2-4\,r^3\sqrt{3}y+3\,r^4 y^2+
4\,\sqrt{3}r^2 y+12\,r^2 x+12 \bigr)\Bigr]\Big/\Bigl[24\,r^2\Bigr]
$.

\noindent \textbf{Case 13:}
\begin{multline*}
P(\{X_2,X_3\} \subset \NPE^r(X_1)\cup \G_1^r(X_1),\; X_1 \in T_s)=\\
\left(\int_{s_{10}}^{s_{14}}\int_{r_2(x)}^{r_{10}(x)} +
\int_{s_{14}}^{s_{13}}\int_{r_2(x)}^{r_{12}(x)} +
\int_{s_{13}}^{1/2}\int_{r_3(x)}^{r_{12}(x)} \right)
\frac{A(\msP(A,N_1,Q_1,L_3,N_3,N_2))^2}{A(\TY)^3}dydx=\\
\Bigl[196608-3538944\,r+8927944704\,r^7-1883996112\,r^{12}-9492593152\,r^8-
146866176\,r^3+29196288\,r^2+220250112\,r^4-\\
4486594560\,r^6+213597\,r^{20}-259250904\,r^{16}+69124752\,r^{17}-10683306\,r^{18}+864387072\,r^5+5220357120\,r^9-\\
1081136256\,r^{10}+602097408\,r^{11}+2223664128\,r^{13}-1509638512\,r^{14}+
716568768\,r^{15}\Bigr]\Big/\Bigl[16796160\,r^8\Bigr]
\end{multline*}
where
$A(\msP(A,N_1,Q_1,L_3,N_3,N_2))=
\Bigl[\sqrt{3} \bigl( -12\,x^2r^2-24\,x-24\,r-12\,r^3x+9\,r^4x^2+4\,y^2-
8\,\sqrt{3}r\,y+6\,r^4\sqrt{3}y\,x+8\,\sqrt{3}y+12\,r^2+24\,r\,x+
12\,x^2-8\,\sqrt{3}y\,x+4\,r^2 y^2-4\,r^3\sqrt{3}y+3\,r^4 y^2+
4\,\sqrt{3}r^2y+12\,r^2x+12 \bigr)\Bigr]\Big/\Bigl[24\,r^2\Bigr]
$.

\noindent \textbf{Case 14:}
\begin{multline*}
P(\{X_2,X_3\} \subset \NPE^r(X_1)\cup \G_1^r(X_1),\; X_1 \in T_s)=\\
\left(\int_{s_7}^{s_{11}}\int_{r_9(x)}^{\ell_{am}(x)} +
\int_{s_{11}}^{s_{10}}\int_{r_9(x)}^{r_{12}(x)} +
\int_{s_{10}}^{s_{14}}\int_{r_{10}(x)}^{r_{12}(x)} \right)
\frac{A(\msP(A,N_1,Q_1,L_3,L_4,Q_2,N_2))^2}{A(\TY)^3}dydx=\\
-\Bigl[\left( r -1 \right)  \bigl( -16384+278528\,r+215136\,r^{28}+40176\,r^{26}+
215136\,r^{29}-3381264\,r^{27}-2301952\,r^2-99212040\,r^{23}-\\
25050384\,r^{24}-312101312\,r^8-7215869272\,r^{13}-147586784\,r^9-
42770432\,r^4+12591104\,r^3+114049024\,r^5+345810944\,r^7+\\
55914462\,r^{21}-2082969096\,r^{17}+43443459\,r^{18}+826941555\,r^{19}-641846754\,r^{20}+
209930616\,r^{22}-232963072\,r^6+\\
1311322268\,r^{10}-3191747236\,r^{11}+
5434516904\,r^{12}+7756861008\,r^{14}-6865898928\,r^{15}+4727296416\,r^{16}+26115696\,r^{25}
\bigr)\Bigr]\Big/\\
\Bigl[466560\, \left( 2\,r^2+1 \right)^5r^8\Bigr]
\end{multline*}
where
$A(\msP(A,N_1,Q_1,L_3,L_4,Q_2,N_2))=
\Bigl[\sqrt{3} \bigl( -3\,x^2r^2-6\,x-12\,r-6\,r^3 x+3\,r^4 x^2+2\,\sqrt{3}r\,y+
6\,r^2+6\,r\,x+3\,x^2-2\,\sqrt{3}y-2\,\sqrt{3}r^2y\,x+2\,r^4\sqrt{3}y\,x+
2\,\sqrt{3}r^2y-r^2 y^2+5\,y^2-2\,r^3\sqrt{3}y+r^4 y^2-2\,\sqrt{3}y\,x+
6+6\,r^2x \bigr)\Bigr]\Big/\Bigl[6\,r^2\Bigr]
$.

\noindent \textbf{Case 15:}
\begin{multline*}
P(\{X_2,X_3\} \subset \NPE^r(X_1)\cup \G_1^r(X_1),\; X_1 \in T_s)=
\int_{s_{13}}^{1/2}\int_{r_2(x)}^{r_3(x)} \frac{A(\msP(A,N_1,Q_1,G_3,M_2,N_3,N_2))^2}{A(\TY)^3}dydx=\\
\Bigl[\bigl( 63855\,r^{10}-498960\,r^9+1650060\,r^8-3036960\,r^7+3703292\,r^6-
3657696\,r^5+3268368\,r^4-2419200\,r^3+1550448\,r^2-\\
725760\,r+155520 \bigr)  \left( -6+5\,r \right)^2\Bigr]\Big/\Bigl[4199040\,r^2\Bigr]
\end{multline*}
where
$A(\msP(A,N_1,Q_1,G_3,M_2,N_3,N_2))=
\Bigl[\sqrt{3} \bigl( 4\,r\,y^2+12\,x+4\,\sqrt{3}r\,y+9\,r-4\,\sqrt{3}y+12\,r\,x-
12+9\,x^2r^3+6\,\sqrt{3}r^3y\,x-12\,r\,x^2-4\,\sqrt{3}r^2y-12\,r^2 x+
3\,r^3y^2 \bigr)\Bigr]\Big/\Bigl[24\,r\Bigr]
$.

\noindent \textbf{Case 16:}
\begin{multline*}
P(\{X_2,X_3\} \subset \NPE^r(X_1)\cup \G_1^r(X_1),\; X_1 \in T_s)=
\int_{s_{14}}^{1/2}\int_{r_{12}(x)}^{r_{10}(x)} \frac{A(\msP(A,N_1,Q_1,L_3,N_3,N_2))^2}{A(\TY)^3}dydx=\\
-\Bigl[\bigl( 293\,r^{16}+2344\,r^{15}+4662\,r^{14}-9088\,r^{13}-32320\,r^{12}+42976\,r^{11}+
175408\,r^{10}-119680\,r^9-544144\,r^8+372352\,r^7+\\
1216512\,r^6-882688\,r^5-1564672\,r^4+1373184\,r^3+924672\,r^2-1314816\,r+380928 \bigr)\\
\left( -2+r \right)  \left( r^2+2\,r -4 \right)^2\Bigr]\Big/\Bigl[23040\, \left( r+2 \right)^5r^4\Bigr]
\end{multline*}
where
$A(\msP(A,N_1,Q_1,L_3,N_3,N_2))=
\Bigl[\sqrt{3} \bigl( -12\,x^2r^2-24\,x-24\,r-12\,r^3x+9\,r^4x^2+
4\,y^2-8\,\sqrt{3}r\,y+6\,r^4\sqrt{3}y\,x+8\,\sqrt{3}y+12\,r^2+24\,r\,x+
12\,x^2-8\,\sqrt{3}y\,x+4\,r^2y^2-4\,r^3\sqrt{3}y+3\,r^4y^2+4\,\sqrt{3}
r^2 y+12\,r^2x+12 \bigr)\Bigr]\Big/\Bigl[24\,r^2\Bigr]
$.

\noindent \textbf{Case 17:}
\begin{multline*}
P(\{X_2,X_3\} \subset \NPE^r(X_1)\cup \G_1^r(X_1),\; X_1 \in T_s)=\\
\left(\int_{s_{11}}^{s_{14}}\int_{r_{12}(x)}^{\ell_{am}(x)}+
\int_{s_{14}}^{1/2}\int_{r_{10}(x)}^{\ell_{am}(x)} \right)
\frac{A(\msP(A,N_1,Q_1,L_3,L_4,Q_2,N_2))^2}{A(\TY)^3}dydx=\\
\Bigl[\bigl( 6723\,r^{20}+73953\,r^{19}+213678\,r^{18}-433512\,r^{17}-2873232\,r^{16}+
627264\,r^{15}+20218896\,r^{14}+5675184\,r^{13}-97577924\,r^{12}-\\
39916108\,r^{11}+343932568\,r^{10}+108508576\,r^9-906967296\,r^8-96480192\,r^7+1702951296\,r^6-
293251072\,r^5-1994987520\,r^4+\\
981590016\,r^3+1118830592\,r^2-1135919104\,r+287604736 \bigr)  \left( r -1
 \right) \Bigr]\Big/\Bigl[466560\,r^4 \left( r+2 \right)^5\Bigr]
\end{multline*}
where
$A(\msP(A,N_1,Q_1,L_3,N_3,N_2))=
\Bigl[\sqrt{3} \bigl( -3\,x^2r^2-6\,x-12\,r-6\,r^3x+3\,r^4x^2+2\,\sqrt{3}r\,y+
6\,r^2+6\,r\,x+3\,x^2-2\,\sqrt{3}y-2\,\sqrt{3}r^2 y\,x+2\,r^4\sqrt{3}y\,x+
2\,\sqrt{3}r^2y-r^2 y^2+5\,y^2-2\,r^3\sqrt{3}y+r^4 y^2-2\,\sqrt{3}y\,x+
6+6\,r^2x \bigr)\Bigr]\Big/\Bigl[6\,r^2\Bigr]
$.
}

Adding up the $P(\{X_2,X_3\} \subset \NPE^r(X_1)\cup \G_1^r(X_1),\; X_1 \in T_s)$
values in the 17 possible cases above, and multiplying by 6
we get, for $r \in \bigl[6/5,\sqrt{5}-1 \bigr)$,
\begin{multline*}
\nu_\lo(r)=
-\Bigl[-413208\,r+3070468\,r^2-74801558\,r^8+75243552\,r^{13}-4883958\,r^9+
14541630\,r^4+28880-11254002\,r^3-\\
3667716\,r^5+64360782\,r^7+13122\,r^{21}-3300900\,r^{17}+156014\,r^{18}-
175011\,r^{19}+62825\,r^{20}+1458\,r^{22}-19812000\,r^6+\\
99831906\,r^{10}-120628524\,r^{11}+33155180\,r^{12}-67685050\,r^{14}+
5055135\,r^{15}+11053023\,r^{16}\Bigr]\Big/\Bigl[116640\,r^6 \left( r^2+1 \right)\\
 \left( 2\,r^2+1 \right)  \left( r+2 \right)^3 \left( r+1 \right)^3\Bigr].
\end{multline*}
The $\nu_\lo(r)$ values for the other intervals can be calculated similarly.

\section*{Appendix 5: The Asymptotic Means of Relative Edge Density
Under Segregation and Association Alternatives}
Let $\mu^S_{\la}(r,\ve)$
and
$\mu^A_{\la}(r,\ve)$
be the means of relative edge density
for the AND-underlying graph
under the segregation and association alternatives.
Define $\mu^S_{\lo}(r,\ve)$
and
$\mu^A_{\lo}(r,\ve)$ similarly.
Derivation of $\mu^S_{\la}(r,\ve)$ involves detailed geometric calculations and
partitioning of the space of $(r,\ve,x)$ for $r \in [1,\infty)$,
$\ve \in \bigl[ 0,\sqrt{3}/3 \bigr)$, and $x \in T_e$.
See Appendix 6 for the derivation of $\mu(r,\ve)$ at a demonstrative interval.

\subsection*{$\mu^S_\la(r,\ve)$ Under Segregation Alternatives}
Under segregation, we compute $\mu^S_\la(r,\ve)$ and $\mu^S_\lo(r,\ve)$ explicitly.
For $\ve \in \bigl[ 0,\sqrt{3}/8 \bigr)$,
$\mu^S_\la(r,\ve)=\sum_{i=1}^4 \varpi^\la_{i}(r,\ve)\,\I(r \in \mI_i)$
where
{\small
\begin{multline*}
\varpi^\la_{1}(r,\ve)=
-{\frac{ \left( r-1 \right)  \left( 5\,r^5+288\,r^5\ve^4+1152\,r^4\ve^4-148\,r^4+
1440\,r^3\ve^4+245\,r^3-178\,r^2+576\,r^2\ve^4-232\,r+128 \right) }
{54\,r^2 \left( 2\,\ve-1 \right)^2 \left( 2\,\ve+1 \right)^2 \left( r+2 \right)
\left( r+1 \right) }}
\end{multline*}

\begin{multline*}
\varpi^\la_{2}(r,\ve)=
-\Bigl[1152\,r^5\ve^4+101\,r^5+3456\,r^4\ve^4-801\,r^4+1302\,r^3+1152\,r^3\ve^4-
732\,r^2-3456\,r^2\ve^4-536\,r-2304\,r\ve^4+672\Bigr]\Big/\\
\Bigl[216\, \left( r+2 \right) r \left( 16\,\ve^4-8\,\ve^2+1 \right)  \left( r+1 \right) \Bigr]
\end{multline*}

\begin{multline*}
\varpi^\la_{3}(r,\ve)=
-\Bigl[-3\,r^8+128\,r^8\ve^4+384\,r^7\ve^4+39\,r^7+128\,r^6\ve^4-90\,r^6-444\,r^5-
384\,r^5\ve^4+1344\,r^4-256\,r^4\ve^4-792\,r^3-\\
864\,r^2+1104\,r-288\Bigr]\Big/\Bigl[24\,r^4 \left( 16\,{
\ve}^4-8\,\ve^2+1 \right)  \left( r+1 \right)  \left(
r+2 \right)\Bigr]
\end{multline*}

$$
\varpi^\la_{4}(r,\ve)=
-{\frac{16\,r^7\ve^4+16\,r^6\ve^4-3\,{
r}^5-16\,r^5\ve^4-3\,r^4-16\,r^4\ve^4+
9\,r^3+9\,r^2-18\,r+6}{3\, \left( r+1 \right) r^4 \left( 4\,{
\ve}^2-1 \right)^2}}
$$
}
with the corresponding intervals
$\mI_1=\Bigl[ 1,4/3 \Bigr)$,
$\mI_2=\Bigl[ 4/3,3/2\Bigr)$,
$\mI_3=\Bigl[ 3/2,2\Bigr)$, and
$\mI_4=\Bigl[ 2,\infty \Bigr)$.

For $\ve \in \Bigl[ 0,\sqrt{3}/8 \bigr)$,
$\mu^S_\lo(r,\ve)=\sum_{i=1}^4 \varpi^\lo_{i}(r,\ve)\,\I(r \in \mI_i)$
where
{\small
\begin{multline*}
\varpi^\lo_{1}(r,\ve)=
\Bigl[47\,r^6-195\,r^5+576\,r^4\ve^4-288\,r^4\ve^2+860\,r^4-846\,r^3+
1728\,r^3\ve^4-864\,r^3\ve^2-108\,r^2-576\,r^2\ve^2+1152\,r^2\ve^4+\\
720\,r-256\Bigr]\Big/\Bigl[108\,r^2 \left( 16\,r\ve^4-8\,r\ve^2+r-16\,\ve^2+
2+32\,\ve^4 \right)  \left( r+1 \right)\Bigr]
\end{multline*}

\begin{multline*}
\varpi^\lo_{2}(r,\ve)=
\Bigl[175\,r^5-579\,r^4+1450\,r^3+1152\,r^3\ve^4-576\,r^3\ve^2+3456\,r^2\ve^4-
1728\,r^2\ve^2-732\,r^2+2304\,r\ve^4-536\,r-1152\,r\ve^2+\\
672\Bigr]\Big/\Bigl[ 216\,\left( r+2 \right) r \left( 2\,\ve-1 \right)^2 \left( 2\,\ve+1 \right)^2
 \left( r+1 \right)\Bigr]
\end{multline*}

\begin{multline*}
\varpi^\la_{3}(r,\ve)=
-\Bigl[27\,r^8-63\,r^7-270\,r^6+1728\,r^6\ve^2-384\,r^6\ve^4+1024\,\ve^3\sqrt{3}r^5-
1152\,r^5\ve^4+576\,r^5\ve^2+756\,r^5+1536\,r^4\ve^3\sqrt{3}-\\
2376\,r^4-6912\,r^4\ve^2-2560\,\sqrt{3}\ve^3r^3+2304\,r^3\ve^4+2736\,r^3+
1152\,r^3\ve^2+1296\,r^2-3072\,r^2\ve^3\sqrt{3}+1536\,r^2\ve^4+
6912\,r^2\ve^2-3312\,r+\\
864\Bigr]\Big/\Bigl[72\,r^4 \left( r+1 \right) \left( 16\,r\ve^4-8\,r\ve^2+r-16\,\ve^2+2+
32\,\ve^4 \right)\Bigr]
\end{multline*}

\begin{multline*}
\varpi^\la_{4}(r,\ve)=
-\Bigl[-18-48\,r^5\ve^4-48\,r^4\ve^4+72\,r^4\ve^2-144\,r^2\ve^2-9\,r^4-32\,r^3\ve^4-
144\,r^3\ve^2+72\,r^5\ve^2-9\,r^5-32\,r^2\ve^4+54\,r+\\
64\,r^2\ve^3 \sqrt{3}+64\,\sqrt{3}\ve^3 r^3\Bigr]\Big/\Bigl[9\,r^4 \left( 4\,\ve^2-1 \right)^2 \left( r+1 \right)\Bigr]
\end{multline*}
}
with the corresponding intervals $\mI_i$ are same as before.

\subsection*{$\mu^A_\la(r,\ve)$ Under Association Alternatives}
Under association,
we compute $\mu^A_\la(r,\ve)$ and $\mu^A_\lo(r,\ve)$  explicitly.
For $\ve \in \bigl[ 0,\left( 7\,\sqrt{3}-3\,\sqrt{15} \right)/12 \approx .042 \bigr)$,
$\mu^A_\la(r,\ve)=\sum_{i=1}^4 \varsigma^\la_{i}(r,\ve)\,\I(r \in \mI_i)$
where
{\small
\begin{multline*}
\varsigma^\la_{1}(r,\ve)=
-\Bigl[-128+768\,r^6\sqrt{3}\ve^3+360\,r+8640\,\ve^4+5760\,\ve^2+393\,r^4-54\,r^2+
6912\,r^4\ve^2+5\,r^6-153\,r^5-423\,r^3-4608\,r^4\sqrt{3}\ve^3+\\
6912\,\sqrt{3}r^2\ve^3+1728\,\ve^2r-3072\,\sqrt{3}\ve^3-7776\,r^2\ve^4-
864\,r^6\ve^4-2592\,r^5\ve^4-18144\,\ve^4r^3+12960\,\ve^4 r-576\,r^6\ve^2
-3456\,r^3\ve^2+1728\,r^5\ve^2-\\
7776\,r^4\ve^4-12096\,r^2\ve^2\Bigr]\Big/\Bigl[6\,\left( \sqrt{3}+6\,
\ve \right)^2 \left( -6\,\ve+\sqrt{3} \right)^2
 \left( r+2 \right) r^2 \left( r+1 \right)\Bigr]
\end{multline*}

\begin{multline*}
\varsigma^\la_{2}(r,\ve)=
\Bigl[-672\,r+20736\,\ve^4+13824\,\ve^2-1302\,r^4+536\,r^2-101\,r^6+801\,r^5+732\,r^3-
3072\,r^6\sqrt{3}\ve^3+18432\,r^4\sqrt{3}\ve^3-9216\,\sqrt{3}r\ve^3-\\
19968\,\sqrt{3}r^2\ve^3+4608\,\sqrt{3}r^3\ve^3+31104\,r^4\ve^4+4608\,r^2\ve^2-
17280\,r^4\ve^2+58752\,\ve^4r^2-6912\,\ve^2r+3456\,r^6 \ve^4+10368\,r^5\ve^4+72576\,\ve^4r^3+\\
31104\,\ve^4r+2304\,r^6\ve^2+17280\,r^3 \ve^2-6912\,r^5\ve^2\Bigr]\Big/\Bigl[216\,\left( r+2 \right) r^2
 \left( r+1 \right)  \left( -1+12\,\ve^2 \right)^2\Bigr]
\end{multline*}

\begin{multline*}
\varsigma^\la_{3}(r,\ve)=
\Bigl[9(r^8-13\,r^7+30\,r^6-192\,r^6 \ve^2+1152\,r^6\ve^4+148\,r^5+3456\,r^5 \ve^4-
576\,r^5\ve^2-448\,r^4+2688\,r^4\ve^4-128\,r^4\ve^2+1152\,\ve^4r^3+\\
264\,r^3+768\,r^3\ve^2+512\,r^2\ve^2+768\,\ve^4r^2+288\,r^2-368\,r+96)\Bigr]\Big/
\Bigl[8\,r^4 \left( -6\,\ve+\sqrt{3} \right)^2 \left( \sqrt{3}+6\,\ve \right)^2
\left( r+1 \right)  \left( r+2 \right)\Bigr]
\end{multline*}

\begin{multline*}
\varsigma^\la_{4}(r,\ve)=
{\frac{9(r^5+6\,r+r^4-3\,r^3-3\,r^2-2+144\,r^5\ve^4+144\,r^4\ve^4+48\,\ve^4r^3+
48\,\ve^4r^2-24\,r^5\ve^2-24\,r^4\ve^2+32\,r^3\ve^2+32\,r^2\ve^2)}
{r^4 \left( r+1 \right)  \left( -\sqrt{3}+6\,\ve \right)^2
 \left( \sqrt{3}+6\,\ve \right)^2}}
\end{multline*}
}
with the corresponding intervals $\mI_i$ are same as before.

For $\ve \in \bigl[ 0,\left( 7\,\sqrt{3}-3\,\sqrt{15} \right)/12 \approx .042 \bigr)$,
$\mu^A_\lo(r,\ve)=\sum_{i=1}^4 \varsigma^\lo_{i}(r,\ve)\,\I(r \in \mI_i)$
where
{\small
\begin{multline*}
\varsigma^\lo_{1}(r,\ve)=
\Bigl[-256+720\,r-13824\,\ve^4-9216\,\ve^2+860\,r^4-108\,r^2+47\,r^6-195\,r^5-
846\,r^3+12096\,r^4\ve^4+6912\,r^2\ve^2+1152\,r^4\ve^2+\\
31104\,\ve^4r^2-6144\,\sqrt{3}\ve^3+3072\,r^6\sqrt{3}\ve^3-6144\,r^4\sqrt{3}\ve^3+
13824\,\sqrt{3}r^2\ve^3+4608\,\sqrt{3}r^5\ve^3+13824\,\ve^2r-10368\,r^5\ve^4+57024\,\ve^4r^3-\\
20736\,\ve^4r-2304\,r^6\ve^2-17280\,r^3\ve^2-3456\,r^6\ve^4\Bigr]\Big/
\Bigl[12\, \left( r+2 \right)  \left( -6\,\ve+\sqrt{3} \right)^2
\left( \sqrt{3}+6\,\ve \right)^2r^2 \left( r+1 \right)\Bigr]
\end{multline*}

\begin{multline*}
\varsigma^\lo_{2}(r,\ve)=
-\Bigl[-672+579\,r^4-1450\,r^3+536\,r+20736\,r^4\ve^4+32832\,r^2\ve^2-114048\,\ve^4 r^2-
7488\,r^3\ve^2+8064\,\ve^2r-175\,r^5+6912\,r^5\ve^4+\\
4608\,r^5\ve^2-24192\,\ve^4r^3-76032\,\ve^4r+12288\,\sqrt{3}r^3\ve^3-
9216\,r^4\sqrt{3}\ve^3+4608\,\sqrt{3}r^2\ve^3+732\,r^2-
6144\,\sqrt{3}r^5\ve^3-9216\,\sqrt{3}\ve^3-\\
19968\,\sqrt{3}r\ve^3-27648\,\ve^2\Bigr]\Big/\Bigl[216\,r \left( r+2 \right)  \left( r+1 \right)
\left( -1+12\,\ve^2 \right)^2\Bigr]
\end{multline*}

\begin{multline*}
\varsigma^\lo_{3}(r,\ve)=
-\Bigl[9(96+384\,r^4\ve^2+192\,r^6 \ve^2-2304\,r^4\ve^4-30\,r^6-1152\,r^6 \ve^4+
84\,r^5+576\,r^5\ve^2+3\,r^8-7\,r^7-368\,r+304\,r^3+144\,r^2-\\
3456\,r^5\ve^4-264\,r^4)\Bigr]\Big/\Bigl[8\,r^4 \left( r+2 \right)  \left( r+1 \right)
\left(\sqrt{3}+6\,\ve \right)^2 \left( -6\,\ve+\sqrt{3}\right)^2\Bigr]
\end{multline*}

$$
\varsigma^\lo_{4}(r,\ve)=
{\frac{9(-6\,r+r^4+r^5+2+144\,r^5\ve^4+144\,{r}
^4\ve^4-24\,r^5\ve^2-24\,r^4\ve^2)}{r^4 \left( r+1 \right)  \left( -6\,\ve+\sqrt{3} \right)^{
2} \left( \sqrt{3}+6\,\ve \right)^2}}
$$
}
with the corresponding intervals $\mI_i$ are same as before.

\section*{Appendix 6: Derivation of $\mu^S_\la(r,\ve)$ and $\mu^S_\lo(r,\ve)$}
We demonstrate the derivation of $\mu_S(r,\ve)$ for segregation with
$\ve \in \bigl[0,\sqrt{3}/8\bigr)$ and
among the intervals of $r$ that do not vanish as $\ve \rightarrow 0$.
So the resultant expressions can be used in PAE analysis.

\subsection*{Derivation of $\mu^S_\la(r,\ve)$}
\begin{figure}
\centering
\scalebox{.6}{\input{arcprobseg1und.pstex_t}}
\caption{
\label{fig:arc-prob-seg1}
The vertices for $\NPE^r(x_1,\ve) \cap \G_1^r(x_1,\ve)$ regions
for $x_1 \in T_s$
in addition to the ones given in Figure \ref{fig:vertices-AND-OR}
because of the restrictive nature of the alternatives.
}
\end{figure}
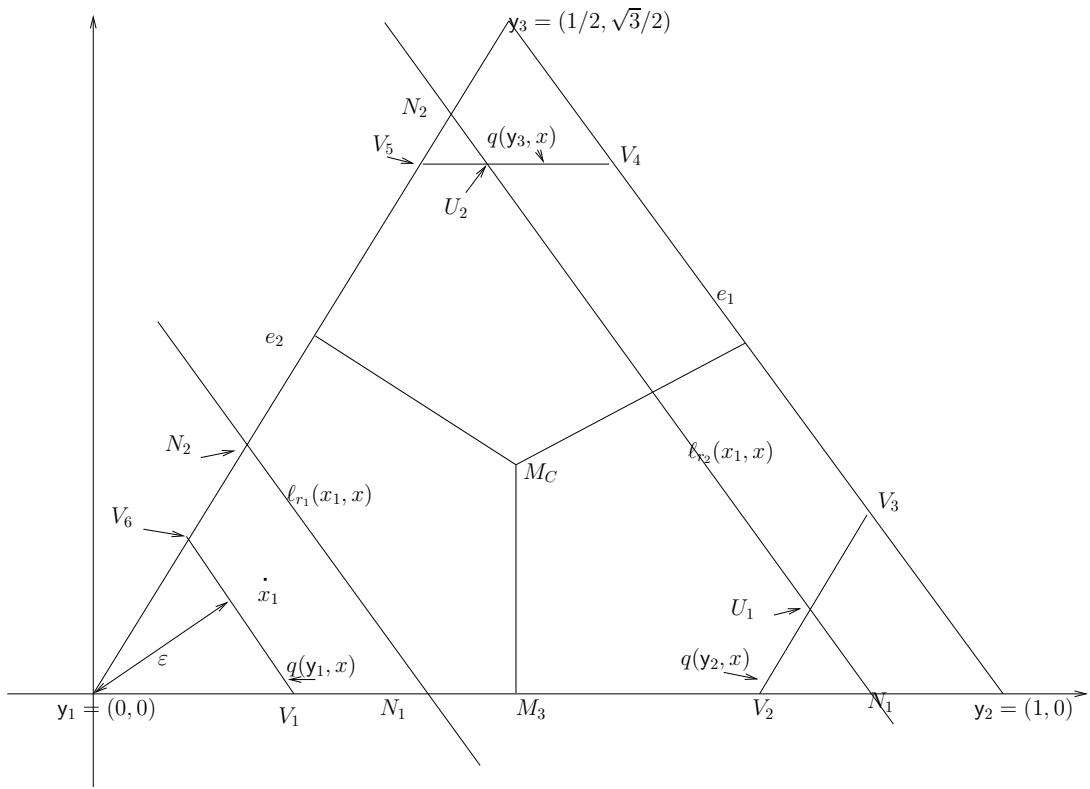

By symmetry,
$$\mu^S_\la(r,\ve)=
P\bigl(X_2 \in \NPE^r(X_1,\ve) \cap \G_1^r(X_1,\ve)\bigr)=
6\,P\bigl(X_2 \in \NY^r(X_1,\ve) \cap \G_1^r(X_1,\ve), X_1 \in T_s \setminus T(\y,\ve)\bigr).$$
Let $q(\y_i,x)$ be the line parallel to $e_i$ and crossing $\TY$ such that
$d(\y_i,q(\y_i,x))=\ve$ for $i=1,2,3$.
Furthermore, let $T_{\ve}:=\TY \setminus \cup_{j=1}^3 T(\y_i,\ve)$.
Then
$q(\y,x)=2\,\ve-\sqrt{3}\,x $,
$q(\y_2,x)=\sqrt{3}\,x-\sqrt{3}+2\,\ve$,
and $q(\y_3,x)=\sqrt{3}/2-\ve$.
Now, let
\begin{align*}
V_1&=q(\y,x)\cap \overline{\y\y_2}=\left(2\,\ve/\sqrt{3},0 \right), &
V_2&=q(\y_2,x)\cap \overline{\y\y_2}=\left(1-2\,\ve/\sqrt{3},0\right),\\
V_3&=q(\y_2,x)\cap \overline{\y_2\y_3}=\left(1-\ve/\sqrt{3},\ve\right),&
V_4&=q(\y_3,x)\cap \overline{\y_2\y_3}=\left(1/2+\ve/\sqrt{3},\sqrt{3}/2-\ve\right), \\
V_5&=q(\y_3,x)\cap \overline{\y\y_3}=\left(1/2-\ve/\sqrt{3},\sqrt{3}/2-\ve\right),&
V_6&=q(\y,x)\cap \overline{\y\y_3}=\left(\ve/\sqrt{3},\ve\right).
\end{align*}
See Figure \ref{fig:arc-prob-seg1}.

The points
$G_i$, for $i=1,2,\ldots,6$,
$P_i$, for $i=1,2$,
$L_i$, for $i=1,2,\ldots,6$,
$N_i$, for $i=1,2,3$,
$Q_i$, for $i=1,2$
and the lines
$r_i(x)$, for $i=1,2,\ldots,11$
are as in Appendix 3.

$s_0=-2\,r/3+1$,
$s_1=-r+3/2$,
$s_2=3/(8\,r)$,
$s_3=1-r/2$,
$s_4=\frac{3}{2\, \left( 2\,r^2+1 \right)}$,
$s_5={\frac{3-3\,r+2\,r^2}{6\,r}}$,
$s_6=1/(2\,r)$,
$s_7=1/(2\,r)$,
$s_8=-{\frac{-2\,r^2-6+r^3+2\,r}{4\,(r^2+1)}}$,
$s_9=-{\frac{-4-6\,r+3\,r^2}{12\,r}}$,
$s_{10}=1/\left( r+1 \right)$,
$s_{11}=-{\frac{-2\,r+r^2-1}{4\,r}}$,
$s_{12}={\frac{-3\,r+2\,r^2+4}{6\,r}}$,
$s_{13}={\frac{9-3\,r^2+2\,r^3-2\,r}{6\,(r^2+1)}}$,
$s_{14}=3\,r/8$,
$s_{15}=r -r^3/8-1/2$



$\ell_1(x)=1/3\,\sqrt{3} \left( -3\, x+2\,\ve\,\sqrt{3} \right)$,
$\ell_2(x)=-1/3\,{\frac{\sqrt{3} \left( 3\, x\,r -3+2\,\ve\,\sqrt{3}\right) }r}$,
$\ell_3(x)=-{\frac{\sqrt{3} \left(  x\,r -1 \right) }r}$,
$\ell_4(x)=1/3\,\sqrt{3} \left( -3\, x+2\,\ve\,\sqrt{3}r \right)$

$q_1=1/2\,\ve\,\sqrt{3}$,
$q_2=2/3\,\ve\,\sqrt{3}$,
$q_3=-1/4\,{\frac{-3+2\,\ve\,\sqrt{3}}r}$,
$q_4=3/4\,r^{-1}$,
$q_7=1/2\,\ve\,\sqrt{3}r$, and
$q_8=2/3\,\ve\,\sqrt{3}r$

Then $T(\y,\ve)=T(\y,Q_1,Q_6)$,
$T(\y_2,\ve)=T(Q_2,\y_2,Q_3)$, and $T(\y_3,\ve)=T(Q_4,Q_5,\y_3)$,
and for $\ve \in \bigl[0,\sqrt{3}/4 \bigr)$,
$T_{\ve}$ is the hexagon with vertices, $Q_i,\;i=1,\ldots,6$.
So we have $A(T_{\ve})=-\ve^2\sqrt{3}+\sqrt{3}/4$.

For $r \in \bigl[1,4/3\bigr)$,
since $\ve$ small enough that $q_2(x)\cap T_e=\emptyset$,
then $N(x,\ve)\subsetneq T_{\ve}$ for all $x \in T_e \setminus T(\y,\ve)$.
There are 14 cases to consider for the AND-underlying version:
{\small
\noindent \textbf{Case 1:}
\begin{multline*}
P(X_2 \in \NPE^r(X_1,\ve) \cap \G_1^r(X_1,\ve), X_1 \in T_s \setminus T(\y,\ve))=\\
\left(\int_{q_1}^{q_7}\int_{\ell_1(x)}^{\ell_{am}(x)}+
\int_{q_7}^{q_2}\int_{\ell_1(x)}^{\ell_4(x)}+
\int_{q_2}^{q_8}\int_{0}^{\ell_4(x)}\right)
\frac{A(\msP(V_1,N_1,N_2,V_6))}{A(T_{\ve})^2}dydx=
{\frac{4\,\ve^4 \left( -3\,r^2+2+r^6 \right) }
{ 9\,\left( 4\,\ve^2-1 \right)^2}}
\end{multline*}
where
$A(\msP(V_1,N_1,N_2,V_6))=
-{\frac{4\,\left( -\ve^2\sqrt{3}+1/4\,\sqrt{3} \right)^2\sqrt{3}
\left( -r^2\,y^2-2\,r^2 y\,\sqrt{3} x-3\,r^2 x^2+4\,\ve^2 \right) }
{9\,\left( 4\,\ve^2-1 \right)^2}}$.

\noindent \textbf{Case 2:}
\begin{multline*}
P(X_2 \in \NPE^r(X_1,\ve) \cap \G_1^r(X_1,\ve), X_1 \in T_s \setminus T(\y,\ve))=\\
\left(\int_{q_7}^{q_8}\int_{\ell_4(x)}^{\ell_{am}(x)}+
\int_{q_8}^{s_2}\int_{0}^{\ell_{am}(x)}+
\int_{s_2}^{s_6}\int_{0}^{r_5(x)}\right)
\frac{A(\msP(G_1,N_1,N_2,G_6))}{A(T_{\ve})^2}dydx=
-{\frac{256\,\ve^4r^{12}-256\,\ve^4r^8-9\,r^4+9}
{576\,r^6 \left( 4\,\ve^2-1 \right)^2}}
\end{multline*}
where
$A(\msP(G_1,N_1,N_2,G_6))=
{\frac{4\,\left( -\ve^2\sqrt{3}+1/4\,\sqrt{3} \right)^2
\left(  y+\sqrt{3} x \right)^2\sqrt{3} \left( r^4-1 \right) }
{9\,r^2 \left( 4\,\ve^2-1 \right)^2}}$.

\noindent \textbf{Case 3:}
\begin{multline*}
P(X_2 \in \NPE^r(X_1,\ve) \cap \G_1^r(X_1,\ve), X_1 \in T_s \setminus T(\y,\ve))=
\left(\int_{s_{11}}^{s_6}\int_{r_5(x)}^{r_7(x)}+
\int_{s_6}^{s_{10}}\int_{0}^{r_7(x)}\right)
\frac{A(\msP(G_1,N_1,P_2,M_3,G_6))}{A(T_{\ve})^2}dydx=\\
{\frac{9\,r^9-13\,r^8-14\,r^7+30\,r^6-22\,r^5+22\,r^4-6\,r^3-10\,r^2+r+3}
{96\, \left( 4\,\ve^2-1 \right)^2r^6 \left( r+1
 \right)^3}}
\end{multline*}
where
$A(\msP(G_1,N_1,P_2,M_3,G_6))=
-{\frac{ 2\,\left( -\ve^2\sqrt{3}+1/4\,\sqrt{3} \right)^2
\left( -12\,r^3 y+2\,r^4\sqrt{3}\,y^2+12\,r^4 y\, x-12\,r^3\sqrt{3} x+
6\,r^4\sqrt{3} x^2+3\,\sqrt{3}r^2+2\,\sqrt{3} y^2+12\, y\, x+6\,\sqrt{3} x^2\right) }
{9\,r^2 \left( 4\,\ve^2-1 \right)^2}}$.

\noindent \textbf{Case 4:}
\begin{multline*}
P(X_2 \in \NPE^r(X_1,\ve) \cap \G_1^r(X_1,\ve), X_1 \in T_s \setminus T(\y,\ve))=\\
\left(\int_{s_2}^{s_5}\int_{r_5(x)}^{\ell_{am}(x)}+
\int_{s_5}^{s_4}\int_{r_2(x)}^{\ell_{am}(x)} +
\int_{s_4}^{s_{13}}\int_{r_2(x)}^{r_8(x)}\right)
\frac{A(\msP(G_1,M_1,P_1,P_2,M_3,G_6))}{A(T_{\ve})^2}dydx=\\
\Bigl[243+7022682\,r^{12}-1296\,r+36612\,r^4-952704\,r^{17}+137472\,r^{18}-578976\,r^7+
7057828\,r^{14}-5116608\,r^{15}+2792712\,r^{16}-7725792\,r^{13}-\\
5484816\,r^{11}+3631995\,r^{10}-2213712\,r^9+1213271\,r^8+3051\,r^2-11664\,r^3-101952\,r^5+
292518\,r^6\Bigr]\Big/\Bigl[15552\, \left( r^2+1 \right)^3\\
\left( 2\,r^2+1 \right)^3r^6 \left( 4\,\ve^2-1 \right)^2\Bigr]
\end{multline*}
where
$A(\msP(G_1,M_1,P_1,P_2,M_3,G_6))=
-{\frac{ 4\,\left( -\ve^2\sqrt{3}+1/4\,\sqrt{3} \right)
^2 \left( -12\,r^3 y-12\,r^3\sqrt{3} x+3\,
\sqrt{3}r^2+3\,r^4\sqrt{3}\,y^2+18\,r^4{
\it y}\, x+9\,r^4\sqrt{3} x^2+\sqrt{3}{{\it
y}}^2+6\, y\, x+3\,\sqrt{3} x^2 \right) }{
9\,\left( 4\,\ve^2-1 \right)^2r^2}}$.

\noindent \textbf{Case 5:}
\begin{multline*}
P(X_2 \in \NPE^r(X_1,\ve) \cap \G_1^r(X_1,\ve), X_1 \in T_s \setminus T(\y,\ve))=
\left(\int_{s_4}^{s_{13}}\int_{r_8(x)}^{r_9(x)}+
\int_{s_{13}}^{s_{12}}\int_{r_2(x)}^{r_9(x)} \right)
\frac{A(\msP(G_1,M_1,L_2,Q_1,P_2,M_3,G_6))}{A(T_{\ve})^2}dydx=\\
-\Bigl[4(400\,r^{15}-2832\,r^{14}+8012\,r^{13}-13608\,r^{12}+16350\,r^{11}-14292\,r^{10}+
8677\,r^9-2442\,r^8-1963\,r^7+3288\,r^6-2751\,r^5+1710\,r^4-743\,r^3+\\
288\,r^2-118\,r+24)\Bigr]\Big/\Bigl[243\,r^3 \left( 2\,r^2+1 \right)^3 \left( r^2+1 \right)^3
 \left( 16\,\ve^4-8\,\ve^2+1 \right)\Bigr]
\end{multline*}
where
$A(\msP(G_1,M_1,L_2,Q_1,P_2,M_3,G_6))=
-\Bigl[4\,\left( -\ve^2\sqrt{3}+1/4\,\sqrt{3} \right)^2
\bigl( -9+42\,\sqrt{3} y\, x-45\, x^2+36\, x-15\,y^2+21\,r^2\,y^2+
2\,r^4 y^4-12\,r^4 x^2\,y^2+12\,r^4 y^2 x+18\, x^3-12\,\sqrt{3} y+
42\,y^2 x-24\,r^3\,y^2-6\,r^2 y\,\sqrt{3} x+4\,\sqrt{3}\,y^3 x+
12\, y\,x^3\sqrt{3}+54\,r^2 x^3+4\,r^4\sqrt{3}\,y^3+12\,r^4 y\,\sqrt{3} x-
12\,r^4 \sqrt{3} x^2 y+18\,r^4 x^2+6\,r^4\,y^2-36\,r^4 x^3+18\,r^4 x^4-
18\,r^2\sqrt{3} x^2 y+12\,r^3\sqrt{3} x^2 y+12\,r^2 x^3\sqrt{3} y-
4\,r^2\sqrt{3}\,y^3 x+12\,r^2 y\,\sqrt{3}-45\,r^2 x^2+9\,r^2-12\,r^3 y\,\sqrt{3}-
4\,r^3\sqrt{3}\,y^3-18\,r^2 y^2 x+12\,r^3\,y^2 x-42\, y\,x^2\sqrt{3}+
6\,r^2\sqrt{3}\,y^3+2\,r^2\,y^4-24\,y^2 x^2-18\,r^2 x^4-36\,r^3 x^3-
36\,r^3 x+72\,r^3 x^2-2\,\sqrt{3}\,y^3 \bigr)\Bigr]\Big/
\Bigl[3\,r^2 \left( -\sqrt{3} y-3+3\, x \right)  \left( - y-\sqrt{3}+
\sqrt{3} x \right)  \left( 4\,\ve^2-1 \right)^2\Bigr]
$.

\noindent \textbf{Case 6:}
\begin{multline*}
P(X_2 \in \NPE^r(X_1,\ve) \cap \G_1^r(X_1,\ve), X_1 \in T_s \setminus T(\y,\ve))=\\
\left(\int_{s_{11}}^{s_{10}}\int_{r_7(x)}^{r_3(x)} +
\int_{s_{10}}^{s_9}\int_{0}^{r_3(x)} +
\int_{s_9}^{1/2}\int_{0}^{r_6(x)} \right)
\frac{A(\msP(G_1,G_2,Q_1,P_2,M_3,G_6))}{A(T_{\ve})^2}dydx=\\
{\frac{324\,r^{11}-1620\,r^{10}-618\,r
^9+4626\,r^8+990\,r^7-2454\,r^6+2703\,r^5-
5571\,r^4-3827\,r^3+1455\,r^2+3072\,r+1024}{7776\,r^
{6} \left( r+1 \right)^3 \left( 16\,\ve^4-8\,\ve^
{2}+1 \right) }}
\end{multline*}
where
$A(\msP(G_1,G_2,Q_1,P_2,M_3,G_6))=
-\Bigl[2\,\left( -\ve^2\sqrt{3}+1/4\,\sqrt{3} \right)^2
\bigl( -9\,\sqrt{3}r^2-24\,\sqrt{3}r\, x-21\,r^2 y-8\,r^2\sqrt{3}\,y^2+
24\,r^2\sqrt{3} x^2-3\,r^2\sqrt{3} x+24\,r\, y+24\,y\, x-24\,\sqrt{3} x^2-
8\,\sqrt{3}\,y^2-6\,\sqrt{3}+18\,\sqrt{3} x-4\,y^3+12\,\sqrt{3}r+12\,\sqrt{3}r\, x^2+
4\,\sqrt{3}\,y^2r -18\, y+12\,r^4 x^2 y-24\,\sqrt{3}r^3x^2+8\,\sqrt{3}r^3\,y^2+
12\,r^4 x^3\sqrt{3}-24\, y\,r\, x-4\,r^2\,y^3+24\,r^3 y-4\,r^4\,y^2\sqrt{3} x-
12\, x^2 y+12\,r^2 x^2 y-12\,r^2 x^3\sqrt{3}+4\,y^2\sqrt{3} x-4\,r^4\sqrt{3}\,y^2-
24\,r^4 y\, x+24\,r^3\sqrt{3} x-12\,r^4\sqrt{3} x^2-4\,r^4\,y^3+4\,r^2\,y^2\sqrt{3} x+
12\,x^3\sqrt{3} \bigr)\Bigr]\Big/\Bigl[3\,r^2 \left( -\sqrt{3} y-3+3\, x \right)  \left( 4\,\ve^2-1 \right)^2\Bigr]
$.

\noindent \textbf{Case 7:}
\begin{multline*}
P(X_2 \in \NPE^r(X_1,\ve) \cap \G_1^r(X_1,\ve), X_1 \in T_s \setminus T(\y,\ve))=\\
\left(\int_{s_4}^{s_{14}}\int_{r_9(x)}^{\ell_{am}(x)} +
\int_{s_{14}}^{s_{12}}\int_{r_9(x)}^{r_{12}(x)} +
\int_{s_{12}}^{s_{15}}\int_{r_{10}(x)}^{r_{12}(x)} \right)
\frac{A(\msP(G_1,M_1,L_2,Q_1,Q_2,L_5,M_3,G_6))}{A(T_{\ve})^2}dydx=\\
\Bigl[1080\,r^{17}-18900\,r^{15}+17280\,r^{14}+65934\,r^{13}-112320\,r^{12}+152361\,r^{11}-
367200\,r^{10}+491051\,r^9-409872\,r^8+282224\,r^7-60864\,r^6-\\
86886\,r^5+70560\,r^4-44672\,r^3+30720\,r^2-16640\,r+6144\Bigr]\Big/
\Bigl[10368\,r^3 \left( 2\,r^2+1 \right)^3 \left( 16\,\ve^4-8\,\ve^2+1 \right) \Bigr]
\end{multline*}
where
$A(\msP(G_1,M_1,L_2,Q_1,Q_2,L_5,M_3,G_6))=
\Bigl[4\,\left( -\ve^2\sqrt{3}+1/4\,\sqrt{3} \right)^2
\bigl( -18+24\,\sqrt{3} y\, x-54\, x^2+54\, x-6\,y^2+21\,r^2\,y^2+
r^4 y^4-6\,r^4 x^2\,y^2+6\,r^4 y^2 x-4\,y^4+18\, x^3-6\,\sqrt{3} y+
42\,y^2 x-24\,r^3\,y^2-18\,r^2 y\,\sqrt{3} x+12\,\sqrt{3}\,y^3 x+
12\,y\, x^3\sqrt{3}+72\,r^2 x^3+2\,r^4\sqrt{3}\,y^3+6\,r^4 y\,\sqrt{3} x-
6\,r^4\sqrt{3} x^2 y+9\,r^4 x^2+3\,r^4\,y^2-18\,r^4 x^3+9\,r^4 x^4+
12\,r^3\sqrt{3} x^2 y+12\,r^2 x^2\,y^2+18\,r^2 y\,\sqrt{3}+18\,r^2 x-
81\,r^2 x^2+9\,r^2-12\,r^3 y\,\sqrt{3}-4\,r^3\sqrt{3}\,y^3-24\,r^2\,y^2 x+
12\,r^3\,y^2 x-30\,y\, x^2\sqrt{3}-2\,r^2\,y^4-36\,y^2 x^2-18\,r^2 x^4-
36\,r^3 x^3-36\,r^3 x+72\,r^3 x^2-6\,\sqrt{3}y^3 \bigr)\Bigr]\Big/
\Bigl[3\,r^2 \left( \sqrt{3} y+3-3\,x \right)  \left( - y-\sqrt{3}+\sqrt{3} x \right)
 \left( 4\,\ve^2-1 \right)^2\Bigr]
$.

\noindent \textbf{Case 8:}
\begin{multline*}
P(X_2 \in \NPE^r(X_1,\ve) \cap \G_1^r(X_1,\ve), X_1 \in T_s \setminus T(\y,\ve))=
\int_{s_9}^{1/2}\int_{r_6(x)}^{r_3(x)} \frac{A(\msP(G_1,G_2,Q_1,N_3,M_C,M_3,G_6))}{A(T_{\ve})^2}dydx=\\
-{\frac{81\,r^{12}+2048+384\,r^4-810\,r^{10}+1296\,r^8-3072\,r^2+96\,r^6}
{15552\,r^6 \left( 16\,\ve^4-8\,\ve^2+1 \right) }}
 \end{multline*}
where
$A(\msP(G_1,G_2,Q_1,N_3,M_C,M_3,G_6))=
-\Bigl[2\,\left( -\ve^2\sqrt{3}+1/4\,\sqrt{3} \right)^2
\bigl( -5\,\sqrt{3}r^2-24\,\sqrt{3}r\, x-17\,r^2 y-8\,r^2\sqrt{3}\,y^2+
24\,r^2\sqrt{3} x^2-7\,r^2\sqrt{3} x+24\,r\, y+24\,y\,x-24\,\sqrt{3} x^2-
8\,\sqrt{3}\,y^2-6\,\sqrt{3}+18\,\sqrt{3} x-4\,y^3+12\,\sqrt{3}r+
12\,\sqrt{3}r\, x^2+4\,\sqrt{3}\,y^2r -18\, y+3\,r^4 x^2 y-
12\,\sqrt{3}r^3 x^2+4\,\sqrt{3}r^3\,y^2+3\,r^4 x^3\sqrt{3}-24\, y\,r\, x-
4\,r^2\,y^3+12\,r^3 y-r^4\,y^2\sqrt{3} x-12\,x^2 y+12\,r^2 x^2 y-
12\,r^2 x^3\sqrt{3}+4\,y^2\sqrt{3} x-r^4\sqrt{3}\,y^2-6\,r^4 y\, x+
12\,r^3\sqrt{3} x-3\,r^4\sqrt{3} x^2-r^4 y^3+4\,r^2\,y^2\sqrt{3} x+
12\, x^3\sqrt{3} \bigr)\Bigr]\Big/
\Bigl[3\,r^2 \left( -\sqrt{3} y-3+3\, x \right)  \left( 4\,\ve^2-1 \right)^2\Bigr]
$.

\noindent \textbf{Case 9:}
\begin{multline*}
P(X_2 \in \NPE^r(X_1,\ve) \cap \G_1^r(X_1,\ve), X_1 \in T_s \setminus T(\y,\ve))=
\left(\int_{s_5}^{s_{13}}\int_{r_5(x)}^{r_2(x)} +
\int_{s_{13}}^{s_{11}}\int_{r_5(x)}^{r_8(x)} \right)
\frac{A(\msP(G_1,M_1,P_1,P_2,M_3,G_6))}{A(T_{\ve})^2}dydx=\\
-\Bigl[243+8673\,r^{12}-1296\,r+23571\,r^4-119712\,r^7-61488\,r^{11}+169716\,r^{10}-
246672\,r^9+216121\,r^8+1404\,r^2-3888\,r^3-35424\,r^5+\\
48816\,r^6\Bigr]\Big/\Bigl[7776\,r^6 \left( 4\,\ve^2-1 \right)^2 \left( r^2+1 \right)^3\Bigr]
 \end{multline*}
where
$A(\msP(G_1,M_1,P_1,P_2,M_3,G_6))=
-{\frac{ 4\,\left( -\ve^2\sqrt{3}+1/4\,\sqrt{3} \right)^2
\left( -12\,r^3 y-12\,r^3\sqrt{3} x+3\,\sqrt{3}r^2+3\,r^4\sqrt{3}\,y^2+
18\,r^4 y\, x+9\,r^4\sqrt{3} x^2+\sqrt{3} y^2+6\, y\, x+3\,\sqrt{3} x^2 \right) }
{9\,r^2 \left( 4\,\ve^2-1 \right)^2}}
$.

\noindent \textbf{Case 10:}
\begin{multline*}
P(X_2 \in \NPE^r(X_1,\ve) \cap \G_1^r(X_1,\ve), X_1 \in T_s \setminus T(\y,\ve))=\\
\left(\int_{s_{12}}^{s_{15}}\int_{r_2(x)}^{r_{10}(x)}+
\int_{s_{15}}^{1/2}\int_{r_2(x)}^{r_{12}(x)}\right)
\frac{A(\msP(G_1,M_1,L_2,Q_1,N_3,L_4,L_5,M_3,G_6))}{A(T_{\ve})^2}dydx=\\
-{\frac{324\,r^{11}-6949\,r^9+7248\,r^8+26896\,r^7-24960\,r^6+2160\,r^5-
259200\,r^4+645760\,r^3-552960\,r^2+155648\,r+6144}
{31104\,r^3 \left( 16\,\ve^4-8\,\ve^2+1 \right) }}
 \end{multline*}
where
$A(\msP(G_1,M_1,L_2,Q_1,N_3,L_4,L_5,M_3,G_6))=
\Bigl[2\,\left( -\ve^2\sqrt{3}+1/4\,\sqrt{3} \right)^2
\bigl( -72-24\,\sqrt{3} y\, x-144\, x^2-144\, x\,r+180\, x+24\,y^2+72\,r+
30\,r^2\,y^2+r^4\,y^4-6\,r^4 x^2 y^2+6\,r^4\,y^2 x-24\,y^4+36\, x^3+
12\,\sqrt{3} y+84\,y^2 x-24\,r^3\,y^2+12\,r^2 y\,\sqrt{3} x+56\,\sqrt{3}\,y^3 x+
24\, y\, x^3\sqrt{3}+108\,r^2 x^3+2\,r^4\sqrt{3}\,y^3+6\,r^4 y\,\sqrt{3} x-
6\,r^4\sqrt{3} x^2 y+9\,r^4 x^2+3\,r^4\,y^2-18\,r^4 x^3+9\,r^4 x^4-
36\,r^2\sqrt{3} x^2 y+12\,r^3\sqrt{3} x^2 y+24\,r^2 x^3\sqrt{3} y-
8\,r^2\sqrt{3} y^3 x-72\,r\,y^2+96\,r\,y^2 x+72\,r^2 x-126\,r^2 x^2-18\,r^2+
72\,r\, x^2-12\,r^3 y\,\sqrt{3}-4\,r^3 \sqrt{3}\,y^3+48\,r\, y\,\sqrt{3} x-
48\,r\, y\, x^2\sqrt{3}-36\,r^2\,y^2 x+12\,r^3\,y^2 x-12\, y\, x^2 \sqrt{3}+
12\,r^2\sqrt{3}\,y^3+4\,r^2\,y^4-120\,y^2 x^2-36\,r^2 x^4-36\,r^3 x^3-36\,r^3 x+
72\,r^3 x^2-28\,\sqrt{3}\,y^3-16\,r\,\sqrt{3}\,y^3 \bigr) \Bigr]\Big/
\Bigl[ 3\,r^2 \left( \sqrt{3} y+3-3\, x \right)
\left( - y-\sqrt{3}+\sqrt{3} x \right)  \left( 4\,\ve^2-1 \right)^2\Bigr]
$.

\noindent \textbf{Case 11:}
\begin{multline*}
P(X_2 \in \NPE^r(X_1,\ve) \cap \G_1^r(X_1,\ve), X_1 \in T_s \setminus T(\y,\ve))=
\int_{s_{15}}^{1/2}\int_{r_{12}(x)}^{r_{10}(x)} \frac{A(\msP(L_1,L_2,Q_1,N_3,L_4,L_5,L_6))}{A(T_{\ve})^2}dydx=\\
{\frac{4\,r^{12}+16\,r^{11}-69\,r^{10}-
260\,r^9+372\,r^8+1248\,r^7+112\,r^6-2624\,{
r}^5-8256\,r^4+12288\,r^3+13568\,r^2-27648\,r+
11264}{ 384\,\left( 16\,r^2\ve^4-8\,r^2\ve^2+
r^2+64\,r\,\ve^4-32\,r\,\ve^2+4\,r+64\,{
\ve}^4-32\,\ve^2+4 \right) r^2}}
 \end{multline*}
where
$A(\msP(L_1,L_2,Q_1,N_3,L_4,L_5,L_6))=
\Bigl[2\,\left( -\ve^2\sqrt{3}+1/4\,\sqrt{3} \right)^2
\bigl( -72+24\,\sqrt{3} y\,r -72\,\sqrt{3} y\,x-216\, x^2-72\, x\,r+180\, x+72\,r+
24\,r^2\,y^2+r^4\,y^4-6\,r^4 x^2\,y^2+6\,r^4\,y^2 x-32\,y^4-72\, x^4+180\, x^3+
12\,\sqrt{3} y+36\,y^2 x-24\,r^3\,y^2+24\,r^2 y\,\sqrt{3} x+56\,\sqrt{3}\,y^3 x+
24\, y\,x^3\sqrt{3}+108\,r^2 x^3+72\,r\,x^3+2\,r^4\sqrt{3}\,y^3+6\,r^4 y\,\sqrt{3} x-
6\,r^4\sqrt{3} x^2 y+9\,r^4 x^2+3\,r^4\,y^2-18\,r^4 x^3+9\,r^4 x^4-36\,r^2\sqrt{3} x^2 y+
12\,r^3\sqrt{3} x^2 y+24\,r^2 x^3\sqrt{3} y-8\,r^2\sqrt{3}\,y^3 x-24\,r\,y^2+72\,r\,y^2 x-
12\,r^2 y\,\sqrt{3}+108\,r^2 x-144\,r^2 x^2-36\,r^2-72\,r\, x^2-12\,r^3 y\,\sqrt{3}-
4\,r^3\sqrt{3}\,y^3+48\,r\, y\,\sqrt{3} x-72\,r\, y\, x^2\sqrt{3}-36\,r^2\,y^2 x+
12\,r^3\,y^2 x+36\,y\, x^2\sqrt{3}+12\,r^2\sqrt{3}\,y^3+4\,r^2\,y^4-72\,y^2 x^2-
36\,r^2 x^4-36\,r^3 x^3-36\,r^3 x+72\,r^3 x^2-44\,\sqrt{3}\,y^3-8\,r\,\sqrt{3}\,y^3 \bigr)\Bigr]\Big/
\Bigl[3\,r^2 \left( \sqrt{3} y+3-3\,x \right)  \left( - y-\sqrt{3}+\sqrt{3} x \right)
\left( 4\,\ve^2-1 \right)^2\Bigr]
$.

\noindent \textbf{Case 12:}
\begin{multline*}
P(X_2 \in \NPE^r(X_1,\ve) \cap \G_1^r(X_1,\ve), X_1 \in T_s \setminus T(\y,\ve))=
\left(\int_{s_{14}}^{s_{15}}\int_{r_{12}(x)}^{\ell_{am}(x)}+
\int_{s_{15}}^{1/2}\int_{r_{10}(x)}^{\ell_{am}(x)}\right)
\frac{A(\msP(L_1,L_2,Q_1,Q_2,L_5,L_6))}{A(T_{\ve})^2}dydx=\\
-\Bigl[135\,r^{12}+540\,r^{11}-2025\,r^{10}-8100\,r^9+10152\,r^8+38448\,r^7-
14878\,r^6-71704\,r^5-87608\,r^4+192128\,r^3+147712\,r^2-338944\,r+\\
134144\Bigr]\Big/\Bigl[10368\,r^2 \left( r^2+4\,r+4 \right)  \left( 16\,\ve^4-8\,\ve^2+1 \right)\Bigr]
\end{multline*}
where
$A(\msP(L_1,L_2,Q_1,Q_2,L_5,L_6))=
-\Bigl[4\,\left( -\ve^2\sqrt{3}+1/4\,\sqrt{3} \right)^2
\bigl( -18+12\,\sqrt{3} y\,r -90\, x^2+36\,x\,r+54\, x-18\,y^2+18\,r^2\,y^2+r^4\,y^4-
6\,r^4 x^2\,y^2+6\,r^4\,y^2 x-8\,y^4-36\, x^4+90\, x^3-6\,\sqrt{3} y+18\,y^2 x-
24\,r^3\,y^2-12\,r^2 y\,\sqrt{3} x+12\,\sqrt{3}\,y^3 x+12\, y\, x^3\sqrt{3}+
72\,r^2 x^3+36\,r\, x^3+2\,r^4\sqrt{3}\,y^3+6\,r^4 y\,\sqrt{3} x-
6\,r^4\sqrt{3} x^2 y+9\,r^4 x^2+3\,r^4\,y^2-18\,r^4 x^3+9\,r^4 x^4+
12\,r^3\sqrt{3} x^2 y+12\,r^2 x^2\,y^2+24\,r\,y^2-12\,r\,y^2 x+12\,r^2 y\,\sqrt{3}+
36\,r^2 x-90\,r^2 x^2-72\,r\, x^2-12\,r^3 y\,\sqrt{3}-4\,r^3\sqrt{3}\,y^3-
12\,r\, y\,x^2\sqrt{3}-24\,r^2\,y^2 x+12\,r^3\,y^2 x-6\, y\, x^2\sqrt{3}-
2\,r^2\,y^4-12\,y^2 x^2-18\,r^2 x^4-36\,r^3 x^3-36\,r^3 x+72\,r^3 x^2-
14\,\sqrt{3}\,y^3+4\,r\,\sqrt{3}y^3 \bigr)\Bigr]\Big/
\Bigl[3\,r^2 \left( -\sqrt{3} y-3+3\, x \right)  \left( - y-\sqrt{3}+\sqrt{3} x \right)
\left( 4\,\ve^2-1 \right)^2\Bigr]
$.

\noindent \textbf{Case 13:}
\begin{multline*}
P(X_2 \in \NPE^r(X_1,\ve) \cap \G_1^r(X_1,\ve), X_1 \in T_s \setminus T(\y,\ve))=\\
\left(\int_{s_{13}}^{s_{11}}\int_{r_8(x)}^{r_2(x)}+
\int_{s_{11}}^{s_{12}}\int_{r_3(x)}^{r_2(x)} +
\int_{s_{12}}^{s_9}\int_{r_3(x)}^{r_6(x)}\right)
\frac{A(\msP(G_1,M_1,L_2,Q_1,P_2,M_3,G_6))}{A(T_{\ve})^2}dydx=\\
{\frac{3654\,r^{12}-35328\,r^{11}+94802\,r^{10}-100608\,r^9-255\,r^8+138240\,r^7-
193581\,r^6+148224\,r^5-86387\,r^4+43008\,r^3-12369\,r^2+512}{7776\,r^6
\left( r^2+1 \right)^3 \left( 16\,\ve^4-8\,\ve^2+1 \right) }}
\end{multline*}
where
$A(\msP(G_1,M_1,L_2,Q_1,P_2,M_3,G_6))=
-\Bigl[4\,\left( -\ve^2\sqrt{3}+1/4\,\sqrt{3} \right)^2
\bigl( -9+42\,\sqrt{3} y\, x-45\, x^2+36\, x-15\,y^2+21\,r^2\,y^2+2\,r^4 y^4-
12\,r^4 x^2\,y^2+12\,r^4y^2 x+18\, x^3-12\,\sqrt{3} y+42\,y^2 x-24\,r^3\,y^2-
6\,r^2 y\,\sqrt{3} x+4\,\sqrt{3}\,y^3 x+12\, y\,x^3\sqrt{3}+54\,r^2 x^3+
4\,r^4\sqrt{3}\,y^3+12\,r^4 y\,\sqrt{3} x-12\,r^4\sqrt{3} x^2 y+18\,r^4 x^2+
6\,r^4\,y^2-36\,r^4 x^3+18\,r^4 x^4-18\,r^2\sqrt{3} x^2 y+12\,r^3\sqrt{3} x^2 y+
12\,r^2 x^3\sqrt{3} y-4\,r^2\sqrt{3}\,y^3 x+12\,r^2 y\,\sqrt{3}-45\,r^2 x^2+9\,r^2-
12\,r^3 y\,\sqrt{3}-4\,r^3\sqrt{3}\,y^3-18\,r^2 y^2 x+12\,r^3\,y^2 x-42\, y\,x^2\sqrt{3}+
6\,r^2\sqrt{3}\,y^3+2\,r^2\,y^4-24\,y^2 x^2-18\,r^2 x^4-36\,r^3 x^3-36\,r^3 x+72\,r^3 x^2-
2\,\sqrt{3}\,y^3 \bigr)\Bigr]\Big/\Bigl[3\,r^2 \left( -\sqrt{3} y-3+3\, x \right)
\left( - y-\sqrt{3}+\sqrt{3} x \right)  \left( 4\,\ve^2-1 \right)^2\Bigr]
$.

\noindent \textbf{Case 14:}
\begin{multline*}
P(X_2 \in \NPE^r(X_1,\ve) \cap \G_1^r(X_1,\ve), X_1 \in T_s \setminus T(\y,\ve))=
\left(\int_{s_{12}}^{s_9}\int_{r_6(x)}^{r_2(x)} +
\int_{s_9}^{1/2}\int_{r_3(x)}^{r_2(x)} \right)
\frac{A(\msP(G_1,M_1,L_2,Q_1,N_3,M_C,M_3,G_6))}{A(T_{\ve})^2}dydx=\\
{\frac{49\,r^{12}+124288\,r^4+50688\,{
r}^7+384\,r^{11}-3562\,r^{10}+13440\,r^9-36948\,{\nu
}^8+27648\,r^2-86016\,r^3-1024-89088\,r^5+160\,{
r}^6}{15552\,r^6 \left( 16\,\ve^4-8\,\ve^2+1
 \right) }}
\end{multline*}
where
$A(\msP(G_1,M_1,L_2,Q_1,N_3,M_C,M_3,G_6))=
-\Bigl[2\,\left( -\ve^2\sqrt{3}+1/4\,\sqrt{3} \right)^2
\bigl( -18+84\,\sqrt{3} y\, x-90\, x^2+72\, x-30\,y^2+38\,r^2\,y^2+r^4 y^4-
6\,r^4 x^2\,y^2+6\,r^4 y^2 x+36\, x^3-24\,\sqrt{3} y+84\,y^2 x-24\,r^3\,y^2-
4\,r^2 y\,\sqrt{3} x+8\,\sqrt{3}\,y^3 x+24\, y\,x^3\sqrt{3}+108\,r^2 x^3+
2\,r^4\sqrt{3}\,y^3+6\,r^4 y\,\sqrt{3} x-6\,r^4\sqrt{3} x^2 y+9\,r^4 x^2+
3\,r^4\,y^2-18\,r^4 x^3+9\,r^4 x^4-36\,r^2\sqrt{3} x^2 y+12\,r^3\sqrt{3} x^2 y+
24\,r^2 x^3\sqrt{3} y-8\,r^2\sqrt{3}\,y^3 x+16\,r^2 y\,\sqrt{3}+24\,r^2 x-
102\,r^2 x^2+6\,r^2-12\,r^3 y\,\sqrt{3}-4\,r^3\sqrt{3}\,y^3-36\,r^2\,y^2 x+
12\,r^3\,y^2 x-84\, y\, x^2\sqrt{3}+12\,r^2\sqrt{3}{y}^3+4\,r^2\,y^4-48\,y^2 x^2-
36\,r^2 x^4-36\,r^3 x^3-36\,r^3 x+72\,r^3 x^2-4\,\sqrt{3}\,y^3 \bigr)\Bigr]\Big/
\Bigl[3\,r^2 \left( -\sqrt{3} y-3+3\, x \right)  \left( -y-\sqrt{3}+\sqrt{3} x \right)
\left( 4\,\ve^2-1 \right)^2\Bigr]
$.
}

Adding up the $P(X_2 \in \NPE^r(X_1,\ve) \cap \G_1^r(X_1,\ve), X_1 \in T_s \setminus T(\y,\ve))$
values in the 14 possible cases above, and multiplying by 6
we get for $r \in [1,4/3)$,
$$\mu^S_\la(r,\ve)=
-{\frac{ \left( r -1 \right)  \left( 5\,r^5+288\,r^5\ve^4+1152\,r^4\ve^4-
148\,r^4+1440\,r^3\ve^4+245\,r^3+576\,r^2 \ve^4-178\,r^2-232\,r+128 \right) }
{54\,r^2 \left( 2+r \right)  \left( 2\,\ve-1 \right)^2
\left( 2\,\ve+1 \right)^2 \left( r+1 \right) }}.
$$
The $\mu^S_\la(r,\ve)$ values for the other intervals can be calculated similarly.

\subsection*{Derivation of $\mu^S_\lo(r,\ve)$}
For $r \in [1,4/3)$,
there are 16 cases to consider for the OR-underlying version:
{\small
\noindent \textbf{Case 1:}
\begin{multline*}
P(X_2 \in \NPE^r(X_1,\ve) \cup \G_1^r(X_1,\ve), X_1 \in T_s \setminus T(\y,\ve))=\\
\left(\int_{q_1}^{q_2}\int_{\ell_1(x)}^{\ell_{am}(x)}+
\int_{q_2}^{s_0}\int_{0}^{\ell_{am}(x)} +
\int_{s_0}^{s_1}\int_{r_1(x)}^{\ell_{am}(x)}\right)
\frac{A(\msP(V_1,M_1,M_C,M_3,V_6))}{A(T_{\ve})^2}dydx=
{\frac{6\,\ve^2-4\,r^2+12\,r -9}{27\,(4\,\ve^2-1)}}
\end{multline*}
where
$A(\msP(V_1,M_1,M_C,M_3,V_6))=
-{\frac{ 4\,\left( -\ve^2\sqrt{3}+1/4\,\sqrt{3} \right)^2\sqrt{3}}
{9\,(4\,\ve^2-1)}}
$.

\noindent \textbf{Case 2:}
\begin{multline*}
P(X_2 \in \NPE^r(X_1,\ve) \cup \G_1^r(X_1,\ve), X_1 \in T_s \setminus T(\y,\ve))=\\
\left(\int_{s_0}^{s_1}\int_{0}^{r_1(x)} +
\int_{s_1}^{s_5}\int_{0}^{r_2(x)} +
\int_{s_5}^{s_3}\int_{0}^{r_5(x)} +
\int_{s_3}^{s_{11}}\int_{r_3(x)}^{r_5(x)}\right)
\frac{A(\msP(V_1,M_1,L_2,L_3,M_C,M_3,V_6))}{A(T_{\ve})^2}dydx=\\
\Bigl[-2304\,r^5\ve^2+432\,r -21960\,r^4-27+9624\,r^7+5952\,r^6\ve^2+288\,r^4\ve^2+
1824\,r^8\ve^2-1817\,r^8-2880\,r^2+10368\,r^3+28224\,r^5-5760\,r^7 \ve^2-\\
21964\,r^6\Bigr]\Big/\Bigl[864\,r^6 \left( 16\,\ve^4-8\,\ve^2+1 \right)\Bigr]
\end{multline*}
where
$A(\msP(V_1,M_1,L_2,L_3,M_C,M_3,V_6))=
-\Bigl[-27+12\,\ve^2r^2 x^2+36\,\sqrt{3} y\,r+108\,\sqrt{3} y\, x-162\, x^2-
108\, x\,r -8\,\ve^2\sqrt{3}r^2 y\,x+108\, x-54\,y^2+36\,r -5\,r^2\,y^2-
3\,y^4-27\, x^4+108\, x^3-36\,\sqrt{3} y+108\,y^2 x+10\,r^2 y\,\sqrt{3} x+
12\,\sqrt{3}\,y^3 x+36\, y\,{x}^3\sqrt{3}-36\,r\, x^3+36\,r\,y^2-36\,r\,y^2 x-
10\,r^2 y\,\sqrt{3}+30\,r^2 x-15\,r^2 x^2-15\,r^2+108\,r\,x^2-
72\,r\, y\,\sqrt{3} x+36\,r\, y\,x^2\sqrt{3}+12\,r^2\ve^2-108\, y\,x^2\sqrt{3}-
54\,y^2 x^2-12\,\sqrt{3}y^3+4\,r\,\sqrt{3}\,y^3+4\,\ve^2r^2 y^2-24\,\ve^2r^2 x+
8\,\ve^2\sqrt{3}r^2 y\Bigr]\Big/\Bigl[4\,r^2 \left( -\sqrt{3} y-3+3\, x \right)
\left( - y-\sqrt{3}+\sqrt{3} x \right)\Bigr]
$.

\noindent \textbf{Case 3:}
\begin{multline*}
P(X_2 \in \NPE^r(X_1,\ve) \cup \G_1^r(X_1,\ve), X_1 \in T_s \setminus T(\y,\ve))=
\left(\int_{s_1}^{s_2}\int_{r_2(x)}^{\ell_{am}(x)} +
\int_{s_2}^{s_5}\int_{r_2(x)}^{r_5(x)}\right)
\frac{A(\msP(V_1,M_1,L_2,L_3,L_4,L_5,M_3,V_6))}{A(T_{\ve})^2}dydx=\\
-\Bigl[-3456\,r^5\ve^2+1296\,r -
65772\,r^4+26880\,r^7+9216\,r^6\ve^2+432\,r^4\ve^2+3072\,r^8\ve^2-4864\,r^8-
8640\,r^2+31104\,r^3+83808\,r^5-\\
9216\,r^7\ve^2-63744\,r^6-81\Bigr]\Big/\Bigl[2592\,r^6 \left( 16\,\ve^4-8\,\ve^2+1 \right)\Bigr]
\end{multline*}
where
$A(\msP(V_1,M_1,L_2,L_3,L_4,L_5,M_3,V_6))=
-\Bigl[-54+12\,\ve^2r^2 x^2+36\,\sqrt{3} y\,r+54\,\sqrt{3} y\, x-189\, x^2-
180\, x\,r -8\,\ve^2\sqrt{3}r^2 y\,x+162\, x-27\,y^2+72\,r -9\,r^2\,y^2-
15\,y^4-27\, x^4+108\, x^3-18\,\sqrt{3} y+108\,y^2 x+18\,r^2 y\,\sqrt{3} x+
36\,\sqrt{3}\,y^3 x+36\, y\,x^3\sqrt{3}-36\,r\, x^3+12\,r\,y^2 x-
18\,r^2 y\,\sqrt{3}+54\,r^2 x-27\,r^2 x^2-27\,r^2+144\,r\, x^2-
48\,r\, y\,\sqrt{3} x+12\,r\, y\, x^2\sqrt{3}+12\,r^2\ve^2-72\, y\, x^2\sqrt{3}-
90\,y^2 x^2-24\,\sqrt{3}\,y^3-4\,r\,\sqrt{3}\,y^3+4\,\ve^2r^2\,y^2-
24\,\ve^2 r^2 x+8\,\ve^2\sqrt{3}r^2 y\Bigr]\Big/
\Bigl[4\,r^2 \left( -\sqrt{3} y-3+3\, x \right)
\left( - y-\sqrt{3}+\sqrt{3} x \right)\Bigr]
$.

\noindent \textbf{Case 4:}
\begin{multline*}
P(X_2 \in \NPE^r(X_1,\ve) \cup \G_1^r(X_1,\ve), X_1 \in T_s \setminus T(\y,\ve))=
\left(\int_{s_3}^{s_{11}}\int_{0}^{r_3(x)} +
\int_{s_{11}}^{s_6}\int_{0}^{r_5(x)} \right)
\frac{A(\msP(V_1,G_2,G_3,M_2,M_C,M_3,V_6))}{A(T_{\ve})^2}dydx=\\
-{\frac{-8\,r+24\,r^2+56\,r^7-32\,r^3-13\,r^8+32\,r^4\ve^2+32\,r^8\ve^2-
128\,r^5\ve^2+192\,r^6\ve^2-128\,r^7\ve^2+64\,r^5-92\,r^6+1}{96\,r^6 \left( 4\,\ve^2-1 \right)^2}}
\end{multline*}
where
$A(\msP(V_1,G_2,G_3,M_2,M_C,M_3,V_6))=
-{\frac{\sqrt{3}\,y^2+6\, y-6\, y\,x+3\,\sqrt{3}-6\,\sqrt{3} x+3\,\sqrt{3} x^2-
2\,\sqrt{3}r^2+4\,\sqrt{3}r^2\ve^2}{12\,r^2}}
$.

\noindent \textbf{Case 5:}
\begin{multline*}
P(X_2 \in \NPE^r(X_1,\ve) \cup \G_1^r(X_1,\ve), X_1 \in T_s \setminus T(\y,\ve))=
\left(\int_{s_{11}}^{s_6}\int_{r_5(x)}^{r_7(x)} +
\int_{s_6}^{s_{10}}\int_{0}^{r_7(x)} \right)
\frac{A(\msP(V_1,G_2,G_3,M_2,M_C,P_2,N_2,V_6))}{A(T_{\ve})^2}dydx=\\
-{\frac{-1+32\,r^5\ve^2+5\,r+34\,r^4+15\,r^7+64\,r^6\ve^2-32\,r^4\ve^2-
32\,r^8\ve^2-17\,r^9+29\,r^8-3\,r^2-17\,r^3-2\,r^5-64\,r^7\ve^2-
43\,r^6+32\,r^9\ve^2}{ 96\,\left( r+1 \right)^3 \left( 4\,\ve^2-1 \right)^2r^6}}
\end{multline*}
where
$A(\msP(V_1,G_2,G_3,M_2,M_C,P_2,N_2,V_6))=
-\Bigl[2\,\left( -\ve^2\sqrt{3}+1/4\,\sqrt{3} \right)^2
\bigl( 2\,\sqrt{3}\,y^2+12\, y-12\, y\, x+6\,\sqrt{3}-12\,\sqrt{3} x+6\,\sqrt{3} x^2-
7\,\sqrt{3}r^2+12\,r^3 y+12\,r^3\sqrt{3}x-4\,r^4\sqrt{3}\,y^2-24\,r^4 y\,x-
12\,r^4\sqrt{3} x^2+8\,\sqrt{3}r^2\ve^2 \bigr)\Bigr]\Big/\Bigl[9\,r^2 \left( 4\,\ve^2-1 \right)^{2}\Bigr]
$.

\noindent \textbf{Case 6:}
\begin{multline*}
P(X_2 \in \NPE^r(X_1,\ve) \cup \G_1^r(X_1,\ve), X_1 \in T_s \setminus T(\y,\ve))=\\
\left(\int_{s_2}^{s_5}\int_{r_5(x)}^{\ell_{am}(x)} +
\int_{s_5}^{s_4}\int_{r_2(x)}^{\ell_{am}(x)} +
\int_{s_4}^{s_{13}}\int_{r_2(x)}^{r_8(x)}\right)
\frac{A(\msP(V_1,N_1,P_1,L_2,L_3,L_4,L_5,P_2,N_2,V_6))}{A(T_{\ve})^2}dydx=\\
\Bigl[-243-28578916\,r^{12}-1147392\,r^{15}\ve^2-1344384\,r^{11}\ve^2-
1734912\,r^{13}\ve^2-304128\,r^{17}\ve^2+989424\,r^{10}\ve^2-10368\,r^5\ve^2+3888\,r -438777\,r^4+\\
2204160\,r^{17}-355328\,r^{18}+5753232\,r^7+39312\,r^6\ve^2+1296\,r^4\ve^2+
296208\,r^8\ve^2-20639832\,r^{14}+13254912\,r^{15}-6591792\,r^{16}+\\
1693728\,r^{12}\ve^2+1507392\,r^{14}{\ve}^2+637056\,r^{16}\ve^2+26417664\,r^{13}+26760576\,r^{11}-21960774\,r^{10}+
15877152\,r^9-10180620\,r^8-28107\,r^2+\\
128304\,r^3+1222128\,r^5-120960\,r^7\ve^2-
2856483\,r^6+92160\,r^{18}\ve^2-563328\,r^9\ve^2\Bigr]\Big/\Bigl[7776\,r^6 \left( r^2+1 \right)^3
\left( 16\,\ve^4-8\,\ve^2+1 \right)  \left( 2\,r^2+1 \right)^3\Bigr]
\end{multline*}
where
$A(\msP(V_1,N_1,P_1,L_2,L_3,L_4,L_5,P_2,N_2,V_6))=
-\Bigl[4\,\left( -\ve^2\sqrt{3}+1/4\,\sqrt{3} \right)^2
\bigl( -54+12\,\ve^2r^2 x^2+36\,\sqrt{3} y\,r+54\,\sqrt{3} y\, x-189\, x^2-
180\, x\,r -8\,\ve^2\sqrt{3}r^2 y\, x+162\, x-27\,y^2+72\,r -12\,r^2\,y^2-
4\,r^4\,y^4+24\,r^4 x^2\,y^2-24\,r^4\,y^2 x-15\,y^4-27\, x^4+108\, x^3-18\,\sqrt{3} y+
108\,y^2 x+24\,r^3\,y^2+24\,r^2 y\,\sqrt{3} x +36\,\sqrt{3}\,y^3 x+36\, y\, x^3 \sqrt{3}-
36\,r\, x^3-8\,r^4\sqrt{3}\,y^3-24\,r^4 y\,\sqrt{3} x+24\,r^4\sqrt{3} x^2 y-36\,r^4 x^2-
12\,r^4\,y^2+72\,r^4 x^3-36\,r^4 x^4-12\,r^3\sqrt{3} x^2 y+12\,r\,y^2 x-24\,r^2 y\,\sqrt{3}+
72\,r^2 x-36\,r^2 x^2-36\,r^2+144\,r\, x^2+12\,r^3 y\,\sqrt{3}+4\,r^3\sqrt{3}\,y^3-
48\,r\, y\,\sqrt{3} x+12\,r\, y\, x^2\sqrt{3}-12\,r^3\,y^2 x+12\,r^2\ve^2-72\, y\,x^2\sqrt{3}-
90\,y^2 x^2+36\,r^3 x^3+36\,r^3 x-72\,r^3 x^2-24\,\sqrt{3}\,y^3-4\,r\,\sqrt{3}\,y^3+
4\,\ve^2 r^2\,y^2-24\,\ve^2r^2 x+8\,\ve^2\sqrt{3}r^2 y \bigr)\Bigr]\Big/
\Bigl[3\,r^2 \left( -\sqrt{3} y-3+3\, x \right)  \left( - y-\sqrt{3}+
\sqrt{3} x \right)  \left( 4\,\ve^2-1 \right)^2\Bigr]
$.

\noindent \textbf{Case 7:}
\begin{multline*}
P(X_2 \in \NPE^r(X_1,\ve) \cup \G_1^r(X_1,\ve), X_1 \in T_s \setminus T(\y,\ve))=
\left(\int_{s_4}^{s_{13}}\int_{r_8(x)}^{r_9(x)} +
\int_{s_{13}}^{s_{12}}\int_{r_2(x)}^{r_9(x)}\right)
\frac{A(\msP(V_1,N_1,Q_1,L_3,L_4,L_5,P_2,N_2,V_6))}{A(T_{\ve})^2}dydx=\\
-\Bigl[8(-2-55766\,r^{12}-864\,r^{15}\ve^2-4104\,r^{11}\ve^2-3024\,r^{13}\ve^2+
3690\,r^{10}\ve^2-108\,r^5\ve^2+24\,r -1833\,r^4+21576\,r^7+342\,r^6\ve^2+
18\,r^4\ve^2+\\
1710\,r^8\ve^2-20056\,r^{14}+6912\,r^{15}-1152\,r^{16}+
3816\,r^{12}\ve^2+1800\,r^{14}\ve^2+288\,r^{16}\ve^2+38376\,r^{13}+65532\,r^{11}-
63642\,r^{10}+52020\,r^9-36277\,r^8-\\
142\,r^2+576\,r^3+4848\,r^5-864\,r^7\ve^2-
10994\,r^6-2700\,r^9\ve^2)\Bigr]\Big/\Bigl[243\,r^4 \left( 2\,r^2+1 \right)^3 \left( r^2+1 \right)^3
 \left( 16\,\ve^4-8\,\ve^2+1 \right)\Bigr]
\end{multline*}
where
$A(\msP(V_1,N_1,Q_1,L_3,L_4,L_5,P_2,N_2,V_6))=
-\Bigl[4\,\left( -\ve^2\sqrt{3}+1/4\,\sqrt{3} \right)^2
\bigl( -36\,\sqrt{3}r^2-180\,\sqrt{3}r\, x+12\,\sqrt{3}r^2\ve^2-90\,r^2 y-
30\,r^2\sqrt{3}\,y^2+18\,r^2\sqrt{3} x^2+54\,r^2\sqrt{3} x+108\,r\, y+54\, y\, x-
135\,\sqrt{3} x^2-9\,\sqrt{3}\,y^2-45\,\sqrt{3}+72\,r^2 y\, x+126\,\sqrt{3} x-60\,y^3+
72\,\sqrt{3}r+144\,\sqrt{3}r\, x^2-18\, y+18\,x^4r^2\sqrt{3}-36\, x^3\sqrt{3}r+
36\,x^3r^3\sqrt{3}-3\,r^4\sqrt{3}\,y^4+54\,r^4 x^2 y-72\,\sqrt{3}r^3 x^2+
24\,\sqrt{3}r^3\,y^2+54\,r^4 x^3\sqrt{3}-144\, y\,r\, x-36\,r^3 x^2 y-
12\,y^3r+96\,y^3 x-18\,r^2\,y^3-18\, x^4\sqrt{3}+12\,r^3\,y^3+36\,r^3 y-18\,r^4\,y^2\sqrt{3} x-
108\, x^2 y+12\,\sqrt{3}\,y^2r\, x-12\,r^3\sqrt{3} x\,y^2+54\,r^2 x^2 y-54\,r^2 x^3\sqrt{3}+
72\,y^2\sqrt{3} x-9\,r^4\sqrt{3}\,y^2-54\,r^4 y\, x+36\,r^3\sqrt{3} x-27\,r^4\sqrt{3} x^2-
2\,y^4r^2\sqrt{3}-18\,r^4 y^3+18\,r^2\,y^2\sqrt{3} x+72\,x^3\sqrt{3}-72\,y^2\sqrt{3} x^2+
24\,\ve^2r^2 y-27\,r^4\sqrt{3} x^4+12\,r^2\,y^3 x-36\,r^2 x^3 y-14\,y^4\sqrt{3}+72\, x^3 y+
36\, x^2 r\, y+18\,r^4\sqrt{3}\,y^2 x^2+4\,\ve^2\sqrt{3}r^2\,y^2+12\,\ve^2\sqrt{3}r^2 x^2-
24\,\ve^2\sqrt{3}r^2 x-24\,\ve^2r^2 y\, x \bigr) \Bigr]\Big/
\Bigl[3\,r^2\left( -\sqrt{3} y-3+3\, x \right)^2 \left( 4\,\ve^2-1 \right)^2\Bigr]
$.

\noindent \textbf{Case 8:}
\begin{multline*}
P(X_2 \in \NPE^r(X_1,\ve) \cup \G_1^r(X_1,\ve), X_1 \in T_s \setminus T(\y,\ve))=\\
\left(\int_{s_{11}}^{s_{10}}\int_{r_7(x)}^{r_3(x)} +
\int_{s_{10}}^{s_9}\int_{0}^{r_3(x)} +
\int_{s_9}^{1/2}\int_{0}^{r_6(x)}\right)
\frac{A(\msP(V_1,N_1,Q_1,G_3,M_2,M_C,P_2,N_2,V_6))}{A(T_{\ve})^2}dydx=\\
\Bigl[-81\,r^9+189\,r^8-561\,r^7+1008\,r^7\ve^2+45\,r^6-432\,r^6\ve^2+1894\,r^5-
3120\,r^5\ve^2+18\,r^4-144\,r^4\ve^2-1912\,r^3+2304\,r^3\ve^2-224\,r^2+\\
768\,r^2\ve^2+384\,r+128\Bigr]\Big/\Bigl[1296\,r^4 \left( 16\,\ve^4-8\,\ve^2+1 \right)
\left( r+1 \right)^3\Bigr]
\end{multline*}
where
$A(\msP(V_1,N_1,Q_1,G_3,M_2,M_C,P_2,N_2,V_6))=
\Bigl[2\,\left( -\ve^2\sqrt{3}+1/4\,\sqrt{3} \right)^2
\bigl( -\sqrt{3}r\, x-24\,r^2 y-8\,r^2\sqrt{3}\,y^2+24\,r^2\sqrt{3} x^2-
24\,r^2\sqrt{3} x+36\,r^3 y\, x+8\,\ve^2\sqrt{3}r\, x-8\,\ve^2\sqrt{3}r+
25\,r\,y+24\, y\, x-12\,\sqrt{3} x^2-4\,\sqrt{3} y^2-8\,\ve^2r\, y-12\,\sqrt{3}+
24\,\sqrt{3} x+13\,\sqrt{3}r -24\,\sqrt{3}r\, x^2+8\,\sqrt{3}\,y^2r -24\, y+
12\, x^3\sqrt{3}r -18\,x^3r^3\sqrt{3}+18\,\sqrt{3}r^3 x^2+6\,\sqrt{3}r^3\,y^2-
18\,r^3 x^2 y+4\,y^3r+6\,r^3\,y^3-4\,\sqrt{3} y^2 r\, x+6\,r^3\sqrt{3} x\,y^2-
12\,x^2 r\, y \bigr)\Bigr]\Big/
\Bigl[3\,r\, \left( \sqrt{3} y+3-3\, x \right)  \left( 4\,\ve^2-1 \right)^2\Bigr]
$.

\noindent \textbf{Case 9:}
\begin{multline*}
P(X_2 \in \NPE^r(X_1,\ve) \cup \G_1^r(X_1,\ve), X_1 \in T_s \setminus T(\y,\ve))=\\
\left(\int_{s_4}^{s_{14}}\int_{r_9(x)}^{\ell_{am}(x)} +
\int_{s_{14}}^{s_{12}}\int_{r_9(x)}^{r_{12}(x)} +
\int_{s_{12}}^{s_{15}}\int_{r_{10}(x)}^{r_{12}(x)}\right)
\frac{A(\msP(V_1,N_1,Q_1,L_3,L_4,Q_2,N_2,V_6))}{A(T_{\ve})^2}dydx=\\
\Bigl[-512+81297\,r^{12}+55296\,r^{11}\ve^2-51264\,r^{10}\ve^2+6912\,r^5\ve^2+
6144\,r+72576\,r^4-1512\,r^{18}-798720\,r^7-35424\,r^6\ve^2-9216\,r^4\ve^2-\\
45792\,r^8\ve^2-83538\,r^{14}-17280\,r^{15}+18252\,r^{16}-51840\,r^{12}\ve^2+
6912\,r^{14}\ve^2+167616\,r^{13}-565920\,r^{11}+888957\,r^{10}-1023600\,r^9+\\
998852\,r^8-7424\,r^2-6144\,r^3-262080\,r^5+41472\,r^7\ve^2+533036\,r^6+82944\,r^9\ve^2\Bigr]
\Big/\Bigl[5184\,r^4 \left( 2\,r^2+1 \right)^3 \left( 16\,\ve^4-8\,\ve^2+1 \right)\Bigr]
\end{multline*}
where
$A(\msP(V_1,N_1,Q_1,L_3,L_4,Q_2,N_2,V_6))=
-\Bigl[8\,\left( -\ve^2\sqrt{3}+1/4\,\sqrt{3} \right)^2
\bigl( -18\,\sqrt{3}r^2-90\,\sqrt{3}r\, x+6\,\sqrt{3}r^2\ve^2-54\,r^2 y-
15\,r^2\sqrt{3}\,y^2+27\,r^2\sqrt{3} x^2+18\,r^2\sqrt{3} x+54\,r\, y+54\, y\, x-
63\,\sqrt{3} x^2-9\,\sqrt{3}\,y^2-18\,\sqrt{3}+54\,r^2 y\, x+54\,\sqrt{3} x-24\,y^3+
36\,\sqrt{3}r+72\,\sqrt{3}r\, x^2-18\, y+9\,x^4r^2\sqrt{3}-18\, x^3\sqrt{3}r+
18\,x^3r^3\sqrt{3}-r^4\sqrt{3}\,y^4+18\,r^4 x^2 y-36\,\sqrt{3}r^3 x^2+
12\,\sqrt{3}r^3\,y^2+18\,r^4 x^3\sqrt{3}-72\, y\,r\, x-18\,r^3 x^2 y-6\,y^3 r+
36\,y^3 x-9\, x^4\sqrt{3}+6\,r^3\,y^3+18\,r^3 y-6\,r^4 y^2\sqrt{3} x-72\, x^2 y-
6\,\sqrt{3}r^2\,y^2 x^2+6\,\sqrt{3}\,y^2r\, x-6\,r^3\sqrt{3} x\,y^2-
36\,r^2 x^3\sqrt{3}+36\,y^2\sqrt{3} x-3\,r^4\sqrt{3}\,y^2-18\,r^4 y\, x+
18\,r^3\sqrt{3} x-9\,r^4\sqrt{3} x^2+\,y^4 r^2\sqrt{3}-6\,r^4\,y^3+12\,r^2\,y^2\sqrt{3} x+
36\, x^3\sqrt{3}-30\,y^2\sqrt{3} x^2+12\,\ve^2r^2 y-9\,r^4\sqrt{3} x^4-5\,y^4\sqrt{3}+
36\, x^3 y+18\, x^2r\, y+6\,r^4\sqrt{3}\,y^2 x^2+2\,\ve^2\sqrt{3}r^2\,y^2+6\,\ve^2\sqrt{3}r^2 x^2-
12\,\ve^2\sqrt{3}r^2 x-12\,\ve^2r^2 y\, x \bigr)\Bigr]\Big/\Bigl[3\,r^2 \left( \sqrt{3} y+3-3\, x \right)^2
\left( 4\,\ve^2-1 \right)^2\Bigr]
$.

\noindent \textbf{Case 10:}
\begin{multline*}
P(X_2 \in \NPE^r(X_1,\ve) \cup \G_1^r(X_1,\ve), X_1 \in T_s \setminus T(\y,\ve))=
\int_{s_9}^{1/2}\int_{r_6(x)}^{r_3(x)} \frac{A(\msP(V_1,N_1,Q_1,G_3,M_2,N_3,N_2,V_6))}{A(T_{\ve})^2}dydx=\\
-{\frac{2496\,r^4+1728\,r^6\ve^2-4608\,r^4\ve^2+512-81\,r^{10}+270\,r^8-
2176\,r^2+3072\,r^2\ve^2-1080\,r^6}{5184\,r^4 \left( 16\,\ve^4-8\,\ve^2+1 \right) }}
\end{multline*}
where
$A(\msP(V_1,N_1,Q_1,G_3,M_2,N_3,N_2,V_6))=
-\Bigl[ 2\, \left( -\ve^2\sqrt{3}+1/4\,\sqrt{3} \right)^2
\bigl( 3\,\sqrt{3}r\, x-12\,r^2 y-4\,r^2 \sqrt{3}\,y^2+12\,r^2\sqrt{3} x^2-
12\,r^2\sqrt{3} x+18\,r^3 y\, x+8\,\ve^2\sqrt{3}r\, x-8\,\ve^2\sqrt{3}r+21\,r\,y+
24\, y\, x-12\,\sqrt{3} x^2-4\,\sqrt{3}\,y^2-8\,\ve^2r\, y-12\,\sqrt{3}+
24\,\sqrt{3} x+9\,\sqrt{3}r -24\,\sqrt{3}r\, x^2+8\,\sqrt{3}\,y^2r -24\, y+12\, x^3\sqrt{3}r-
9\, x^3r^3\sqrt{3}+9\,\sqrt{3}r^3 x^2+3\,\sqrt{3}r^3\,y^2-9\,r^3 x^2 y+4\,y^3r+
3\,r^3\,y^3-4\,\sqrt{3} y^2 r\, x+3\,r^3\sqrt{3} x\,y^2-12\, x^2r\, y \bigr)\Bigr]
\Big/\Bigl[3\,r\, \left( -\sqrt{3} y-3+3\, x \right)  \left( 4\,\ve^2-1 \right)^2\Bigr]
$.

\noindent \textbf{Case 11:}
\begin{multline*}
P(X_2 \in \NPE^r(X_1,\ve) \cup \G_1^r(X_1,\ve), X_1 \in T_s \setminus T(\y,\ve))=\\
\left(\int_{s_5}^{s_{13}}\int_{r_5(x)}^{r_2(x)} +
\int_{s_{13}}^{s_{11}}\int_{r_5(x)}^{r_8(x)} \right)
\frac{A(\msP(V_1,N_1,P_1,L_2,L_3,M_C,P_2,N_2,V_6))}{A(T_{\ve})^2}dydx=\\
-\Bigl[-43855\,r^{12}+14112\,r^{12}\ve^2+271488\,r^{11}-48384\,r^{11}\ve^2-
746553\,r^{10}+81792\,r^{10}\ve^2-117504\,r^9\ve^2+1230336\,r^9-1404177\,r^8+\\
123840\,r^8\ve^2+1236528\,r^7-89856\,r^7\ve^2-901350\,r^6+58752\,r^6\ve^2-
20736\,r^5\ve^2+550800\,r^5-276453\,r^4+2592\,r^4\ve^2+104976\,r^3-\\
26649\,r^2+3888\,r -243\Bigr]\Big/\Bigl[7776\,r^6 \left( 16\,\ve^4-8\,\ve^2+1 \right)  \left( r^2+1 \right)^3\Bigr]
\end{multline*}
where
$A(\msP(V_1,N_1,P_1,L_2,L_3,M_C,P_2,N_2,V_6))=
-\Bigl[4\,\left( -\ve^2\sqrt{3}+1/4\,\sqrt{3} \right)^2
\bigl( -27+12\,\ve^2r^2 x^2+36\,\sqrt{3} y\,r+108\,\sqrt{3} y\, x-162\, x^2-
108\, x\,r -8\,\ve^2\sqrt{3}r^2 y\,x+108\, x-54\,y^2+36\,r -8\,r^2\,y^2-
4\,r^4\,y^4+24\,r^4 x^2\,y^2-24\,r^4\,y^2 x-3\,y^4-27\, x^4+108\, x^3-36\,\sqrt{3} y+
108\,y^2 x+24\,r^3\,y^2+16\,r^2 y\,\sqrt{3} x+12\,\sqrt{3}\,y^3 x+36\, y\, x^3\sqrt{3}-
36\,r\, x^3-8\,r^4\sqrt{3}\,y^3-24\,r^4 y\,\sqrt{3} x+24\,r^4\sqrt{3} x^2 y-36\,r^4 x^2-
12\,r^4\,y^2+72\,r^4 x^3-36\,r^4 x^4-12\,r^3\sqrt{3} x^2 y+36\,r\,y^2-36\,r\, y^2 x-
16\,r^2 y\,\sqrt{3}+48\,r^2 x-24\,r^2 x^2-24\,r^2+108\,r\,x^2+12\,r^3 y\,\sqrt{3}+
4\,r^3\sqrt{3}\,y^3-72\,r\, y\,\sqrt{3} x+36\,r\, y\,x^2\sqrt{3}-12\,r^3\,y^2 x+12\,r^2\ve^2-
108\, y\, x^2\sqrt{3}-54\,y^2 x^2+36\,r^3 x^3+36\,r^3 x-72\,r^3 x^2-
12\,\sqrt{3}\,y^3+4\,r\,\sqrt{3} y^3+4\,\ve^2r^2\,y^2-24\,\ve^2 r^2 x+
8\,\ve^2\sqrt{3}r^2 y \bigr)\Bigr]\Big/\Bigl[3\,r^2 \left( -\sqrt{3} y-3+3\, x \right)
\left( - y-\sqrt{3}+\sqrt{3} x \right)  \left( 4\,\ve^2-1 \right)^2\Bigr]
$.

\noindent \textbf{Case 12:}
\begin{multline*}
P(X_2 \in \NPE^r(X_1,\ve) \cup \G_1^r(X_1,\ve), X_1 \in T_s \setminus T(\y,\ve))=
\left(\int_{s_{12}}^{s_{15}}\int_{r_2(x)}^{r_{10}(x)} +
\int_{s_{15}}^{1/2}\int_{r_2(x)}^{r_{12}(x)} \right)
\frac{A(\msP(V_1,N_1,Q_1,L_3,N_3,N_2,V_6))}{A(T_{\ve})^2}dydx=\\
-\Bigl[5184\,r^8\ve^2-71424\,r^6\ve^2+138240\,r^5\ve^2-73728\,r^4\ve^2-1053\,r^{12}+
16230\,r^{10}-17856\,r^9-68908\,r^8+104448\,r^7+276688\,r^6-916608\,r^5+\\
1032192\,r^4-516096\,r^3+80128\,r^2+12288\,r -1024\Bigr]\Big/\Bigl[31104\,r^4 \left( 16\,\ve^4-8\,\ve^2+1 \right)\Bigr]
\end{multline*}
where
$A(\msP(V_1,N_1,Q_1,L_3,N_3,N_2,V_6))=
-\Bigl[2\,\left( -\ve^2\sqrt{3}+1/4\,\sqrt{3} \right)^2
\bigl( -36\,\sqrt{3}r^2-216\,\sqrt{3}r\, x+24\,\sqrt{3}r^2\ve^2-108\,r^2 y-
48\,r^2\sqrt{3}\,y^2+72\,r^2\sqrt{3} x^2+36\,r^2\sqrt{3} x+216\,r\, y+432\, y\, x-
216\,\sqrt{3} x^2-72\,\sqrt{3}\,y^2-36\,\sqrt{3}+72\,r^2 y\, x+144\,\sqrt{3} x-
48\,y^3+72\,\sqrt{3}r+216\,\sqrt{3}r\, x^2+72\,\sqrt{3} y^2 r -144\, y+
36\, x^4r^2\sqrt{3}-72\,x^3\sqrt{3}r+36\, x^3r^3\sqrt{3}-3\,r^4\sqrt{3}\,y^4+
54\,r^4 x^2 y-72\,\sqrt{3}r^3 x^2+24\,\sqrt{3}r^3\,y^2+54\,r^4 x^3\sqrt{3}-
432\, y\,r\, x-36\,r^3 x^2 y+24\,y^3r+48\,y^3 x-36\,r^2\,y^3-36\, x^4\sqrt{3}+
12\,r^3\,y^3+36\,r^3 y-18\,r^4 y^2\sqrt{3} x-432\, x^2 y-72\,\sqrt{3} y^2r\, x-
12\,r^3\sqrt{3} x\,y^2+108\,r^2 x^2 y-108\,r^2 x^3\sqrt{3}+144\,y^2\sqrt{3} x-
9\,r^4\sqrt{3} y^2-54\,r^4 y\, x+36\,r^3\sqrt{3} x-27\,r^4\sqrt{3} x^2-
4\,y^4r^2 \sqrt{3}-18\,r^4\,y^3+36\,r^2\,y^2\sqrt{3} x+144\, x^3\sqrt{3}-
72\,y^2\sqrt{3} x^2+48\,\ve^2r^2 y-27\,r^4\sqrt{3} x^4+24\,r^2\,y^3 x-
72\,r^2 x^3 y-4\,y^4\sqrt{3}+144\,x^3 y+216\, x^2r\, y+18\,r^4\sqrt{3}\,y^2 x^2+
8\,\ve^2\sqrt{3}r^2 y^2+24\,\ve^2\sqrt{3}r^2 x^2-48\,\ve^2\sqrt{3}r^2 x-
48\,\ve^2r^2 y\, x \bigr)\Bigr]\Big/\Bigl[3\,r^2 \left( -\sqrt{3} y-3+3\,x \right)^2 \left( 4\,\ve^2-1 \right)^2\Bigr]
$.

\noindent \textbf{Case 13:}
\begin{multline*}
P(X_2 \in \NPE^r(X_1,\ve) \cup \G_1^r(X_1,\ve), X_1 \in T_s \setminus T(\y,\ve))=
\int_{s_{15}}^{1/2}\int_{r_{12}(x)}^{r_{10}(x)} \frac{A(\msP(V_1,N_1,Q_1,L_3,N_3,N_2,V_6))}{A(T_{\ve})^2}dydx=\\
\Bigl[-13\,r^{13}-78\,r^{12}+42\,r^{11
}+892\,r^{10}+64\,r^9\ve^2+220\,r^9-4952\,r^8+384\,r^8\ve^2-768\,r^7-3072\,r^6 \ve^2+
18048\,r^6-3136\,r^5-2048\,r^5\ve^2+\\
8192\,r^4\ve^2-39296\,r^4+20992\,r^3+4096\,r^3\ve^2+41984\,r^2-8192\,r^2 \ve^2-
48128\,r+14336\Bigr]\Big/\Bigl[ 384\,\bigl( 16\,r^3\ve^4-8\,r^3\ve^2+r^3+96\,r^2\ve^4-\\
48\,r^2\ve^2+6\,r^2+192\,r\,\ve^4-96\,r\,\ve^2+12\,r+128\,\ve^4-64\,\ve^2+8 \bigr) r^2\Bigr]
\end{multline*}
where
$A(\msP(V_1,N_1,Q_1,L_3,N_3,N_2,V_6))=
-\Bigl[2\,\left( -\ve^2\sqrt{3}+1/4\,\sqrt{3} \right)^2
\bigl( -36\,\sqrt{3}r^2-216\,\sqrt{3}r\, x+24\,\sqrt{3}r^2\ve^2-108\,r^2 y-
48\,r^2\sqrt{3}\,y^2+72\,r^2\sqrt{3} x^2+36\,r^2\sqrt{3} x+216\,r\, y+
432\, y\, x-216\,\sqrt{3} x^2-72\,\sqrt{3}\,y^2-36\,\sqrt{3}+72\,r^2 y\, x+
144\,\sqrt{3} x-48\,y^3+72\,\sqrt{3}r+216\,\sqrt{3}r\, x^2+72\,\sqrt{3} y^2 r -
144\, y+36\, x^4r^2\sqrt{3}-72\,x^3\sqrt{3}r+36\, x^3r^3\sqrt{3}-3\,r^4\sqrt{3}\,y^4+
54\,r^4 x^2 y-72\,\sqrt{3}r^3 x^2+24\,\sqrt{3}r^3\,y^2+54\,r^4 x^3\sqrt{3}-432\, y\,r\, x-
36\,r^3 x^2 y+24\,y^3r+48\,y^3 x-36\,r^2\,y^3-36\, x^4\sqrt{3}+12\,r^3\,y^3+36\,r^3 y-
18\,r^4 y^2\sqrt{3} x-432\, x^2 y-72\,\sqrt{3} y^2 r\, x-12\,r^3\sqrt{3} x\,y^2+
108\,r^2 x^2 y-108\,r^2 x^3 \sqrt{3}+144\,y^2\sqrt{3} x-9\,r^4\sqrt{3} y^2-
54\,r^4 y\, x+36\,r^3\sqrt{3} x-27\,r^4\sqrt{3} x^2-4\,y^4r^2\sqrt{3}-18\,r^4\,y^3+
36\,r^2\,y^2\sqrt{3} x+144\, x^3\sqrt{3}-72\,y^2\sqrt{3} x^2+48\,\ve^2r^2 y-
27\,r^4\sqrt{3} x^4+24\,r^2\,y^3 x-72\,r^2 x^3 y-4\,y^4\sqrt{3}+144\,x^3 y+
216\, x^2r\, y+18\,r^4\sqrt{3}\,y^2 x^2+8\,\ve^2\sqrt{3}r^2 y^2+24\,\ve^2\sqrt{3}r^2 x^2-
48\,\ve^2\sqrt{3}r^2 x-48\,\ve^2r^2 y\,x \bigr)\Bigr]
\Big/\Bigl[3\,r^2 \left( -\sqrt{3} y-3+3\,x \right)^2 \left( 4\,\ve^2-1 \right)^2\Bigr]
$.

\noindent \textbf{Case 14:}
\begin{multline*}
P(X_2 \in \NPE^r(X_1,\ve) \cup \G_1^r(X_1,\ve), X_1 \in T_s \setminus T(\y,\ve))=\\
\left(\int_{s_{14}}^{s_{15}}\int_{r_{12}(x)}^{\ell_{am}(x)} +
\int_{s_{15}}^{1/2}\int_{r_{10}(x)}^{\ell_{am}(x)}\right)
\frac{A(\msP(V_1,N_1,Q_1,L_3,L_4,Q_2,N_2,V_6))}{A(T_{\ve})^2}dydx=\\
-\Bigl[-189\,r^{13}-1134\,r^{12}+297\,r^{11}+11718\,r^{10}+864\,r^9\ve^2+3672\,r^9-
66096\,r^8+5184\,r^8\ve^2+2592\,r^7\ve^2-12932\,r^7-32832\,r^6\ve^2+248616\,r^6-\\
30448\,r^5-33408\,r^5\ve^2+76032\,r^4\ve^2-551584\,r^4+273152\,r^3+55296\,r^{3
}\ve^2+595456\,r^2-73728\,r^2\ve^2-668160\,r+197632\Bigr]\Big/\Bigl[5184\,r^2\\
\left( r^3+6\,r^2+12\,r+8 \right)  \left( 16\,\ve^4-8\,\ve^2+1 \right)\Bigr]
\end{multline*}
where
$A(\msP(V_1,N_1,Q_1,L_3,L_4,Q_2,N_2,V_6))=
-\Bigl[8\,\left( -\ve^2\sqrt{3}+1/4\,\sqrt{3} \right)^2
\bigl( -18\,\sqrt{3}r^2-90\,\sqrt{3}r\, x+6\,\sqrt{3}r^2\ve^2-54\,r^2 y-
15\,r^2\sqrt{3}\,y^2+27\,r^2\sqrt{3} x^2+18\,r^2\sqrt{3} x+54\,r\, y+
54\, y\, x-63\,\sqrt{3} x^2-9\,\sqrt{3}\,y^2-18\,\sqrt{3}+54\,r^2 y\, x+
54\,\sqrt{3} x-24\,y^3+36\,\sqrt{3}r+72\,\sqrt{3}r\, x^2-18\, y+9\,x^4r^2\sqrt{3}-
18\, x^3\sqrt{3}r+18\,x^3 r^3\sqrt{3}-r^4\sqrt{3}\,y^4+18\,r^4 x^2 y-36\,\sqrt{3}r^3 x^2+
12\,\sqrt{3}r^3\,y^2+18\,r^4 x^3\sqrt{3}-72\, y\,r\, x-18\,r^3 x^2 y-6\,y^3 r+36\,y^3 x-
9\, x^4\sqrt{3}+6\,r^3\,y^3+18\,r^3 y-6\,r^4 y^2\sqrt{3} x-72\, x^2 y-6\,\sqrt{3}r^2\,y^2 x^2+
6\,\sqrt{3}\,y^2r\, x-6\,r^3\sqrt{3} x\,y^2-36\,r^2 x^3\sqrt{3}+36\,y^2\sqrt{3} x-
3\,r^4\sqrt{3}\,y^2-18\,r^4 y\, x+18\,r^3\sqrt{3} x-9\,r^4\sqrt{3} x^2+\,y^4r^2\sqrt{3}-
6\,r^4\,y^3+12\,r^2\,y^2\sqrt{3} x+36\, x^3\sqrt{3}-30\,y^2\sqrt{3} x^2+12\,\ve^2r^2 y-
9\,r^4\sqrt{3} x^4-5\,y^4\sqrt{3}+36\, x^3 y+18\, x^2r\, y+6\,r^4\sqrt{3}\,y^2 x^2+
2\,\ve^2\sqrt{3}r^2\,y^2+6\,\ve^2\sqrt{3}r^2 x^2-12\,\ve^2\sqrt{3}r^2 x-
12\,\ve^2r^2 y\, x \bigr) \Bigr]\Big/\Bigl[3\,r^2 \left( -\sqrt{3} y-3+3\, x \right)^2
\left( 4\,\ve^2-1 \right)^2\Bigr]
$.

\noindent \textbf{Case 15:}
\begin{multline*}
P(X_2 \in \NPE^r(X_1,\ve) \cup \G_1^r(X_1,\ve), X_1 \in T_s \setminus T(\y,\ve))=\\
\left(\int_{s_{13}}^{s_{11}}\int_{r_8(x)}^{r_2(x)} +
\int_{s_{11}}^{s_{12}}\int_{r_3(x)}^{r_2(x)} +
\int_{s_{12}}^{s_9}\int_{r_3(x)}^{r_6(x)} \right)
\frac{A(\msP(V_1,N_1,Q_1,L_3,M_C,P_2,N_2,V_6))}{A(T_{\ve})^2}dydx=\\
\Bigl[4536\,r^{12}\ve^2-11753\,r^{12}-13824\,r^{11}\ve^2+69120\,r^{11}+23976\,r^{10}\ve^2-
186683\,r^{10}-34560\,r^9\ve^2+305664\,r^9+35496\,r^8\ve^2-346171\,r^8-\\
27648\,r^7\ve^2+302592\,r^7+17208\,r^6\ve^2-220201\,r^6-6912\,r^5\ve^2+
135936\,r^5+1152\,r^4\ve^2-69760\,r^4+28416\,r^3-8384\,r^2+1536\,r -128\Bigr]\Big/\\
\Bigl[1944\,r^6 \left( r^2+1 \right)^3 \left( 16\,\ve^4-8\,\ve^2+1 \right)\Bigr]
\end{multline*}
where
$A(\msP(V_1,N_1,Q_1,L_3,M_C,P_2,N_2,V_6))=
-\Bigl[4\,\left( -\ve^2\sqrt{3}+1/4\,\sqrt{3} \right)^2
\bigl( -24\,\sqrt{3}r^2-108\,\sqrt{3}r\, x+12\,
\sqrt{3}r^2\ve^2-66\,r^2 y-26\,r^2\sqrt{3}\,y^2+30\,r^2\sqrt{3} x^2+30\,r^2\sqrt{3} x+
108\,r\, y+216\, y\, x-108\,\sqrt{3} x^2-36\,\sqrt{3}\,y^2-18\,\sqrt{3}+
48\,r^2 y\, x+72\,\sqrt{3} x-24\,y^3+36\,\sqrt{3}r+108\,\sqrt{3}r\, x^2+
36\,\sqrt{3} y^2r -72\, y+18\, x^4r^2\sqrt{3}-36\,x^3\sqrt{3}r+36\, x^3r^3\sqrt{3}-
3\,r^4\sqrt{3}\,y^4+54\,r^4 x^2 y-72\,\sqrt{3}r^3 x^2+24\,\sqrt{3}r^3\,y^2+
54\,r^4 x^3\sqrt{3}-216\, y\,r\, x-36\,r^3 x^2 y+12\,y^3r+24\,y^3 x-18\,r^2\,y^3-
18\, x^4\sqrt{3}+12\,r^3\,y^3+36\,r^3 y-18\,r^4 y^2\sqrt{3} x-216\, x^2 y-
36\,\sqrt{3} y^2r\, x-12\,r^3\sqrt{3} x\,y^2+54\,r^2 x^2 y-54\,r^2 x^3\sqrt{3}+
72\,y^2\sqrt{3} x-9\,r^4\sqrt{3}{y}^2-54\,r^4 y\, x+36\,r^3\sqrt{3} x-
27\,r^4\sqrt{3} x^2-2\,y^4r^2\sqrt{3}-18\,r^4\,y^3+18\,r^2\,y^2\sqrt{3} x+
72\, x^3\sqrt{3}-36\,y^2\sqrt{3} x^2+24\,\ve^2r^2 y-27\,r^4\sqrt{3} x^4+12\,r^2\,y^3 x-
36\,r^2 x^3 y-2\,y^4\sqrt{3}+72\, x^3 y+108\,x^2r\, y+18\,r^4\sqrt{3}\,y^2 x^2+
4\,\ve^2\sqrt{3}r^2\,y^2+12\,\ve^2\sqrt{3}r^2 x^2-24\,\ve^2\sqrt{3}r^2 x-
24\,\ve^2r^2 y\, x \bigr)\Bigr]\Big/\Bigl[3\,r^2 \left( \sqrt{3} y+3-3\, x \right)^2
\left( 4\,\ve^2-1 \right)^2\Bigr]
$.

\noindent \textbf{Case 16:}
\begin{multline*}
P(X_2 \in \NPE^r(X_1,\ve) \cup \G_1^r(X_1,\ve), X_1 \in T_s \setminus T(\y,\ve))=
\left(\int_{s_{12}}^{s_9}\int_{r_6(x)}^{r_2(x)} +
\int_{s_9}^{1/2}\int_{r_3(x)}^{r_2(x)} \right)
\frac{A(\msP(V_1,N_1,Q_1,L_3,N_3,N_2,V_6))}{A(T_{\ve})^2}dydx=\\
\Bigl[-147\,r^{12}+55296\,r^5\ve^2-12288\,r+351872\,r^4-142080\,r^7+1024-73728\,r^6\ve^2-
9216\,r^4\ve^2+576\,r^8 \ve^2-1152\,r^{11}+7018\,r^{10}-20352\,r^9+\\
51188\,r^8+64000\,r^2-190464\,r^3-414720\,r^5+27648\,r^7\ve^2+
305920\,r^6\Bigr]\Big/\Bigl[15552\,r^6 \left( 16\,\ve^4-8\,\ve^2+1 \right)\Bigr]
\end{multline*}
where
$A(\msP(V_1,N_1,Q_1,L_3,N_3,N_2,V_6))=
-\Bigl[2\,\left( -\ve^2\sqrt{3}+1/4\,\sqrt{3} \right)^2
\bigl( -36\,\sqrt{3}r^2-216\,\sqrt{3}r\, x+24\,\sqrt{3}r^2\ve^2-108\,r^2 y-
48\,r^2\sqrt{3}\,y^2+72\,r^2\sqrt{3} x^2+36\,r^2\sqrt{3} x+216\,r\, y+432\, y\, x-
216\,\sqrt{3} x^2-72\,\sqrt{3}\,y^2-36\,\sqrt{3}+72\,r^2 y\, x+144\,\sqrt{3} x-
48\,y^3+72\,\sqrt{3}r+216\,\sqrt{3}r\, x^2+72\,\sqrt{3} y^2r -144\, y+36\, x^4r^2\sqrt{3}-
72\,x^3\sqrt{3}r+36\, x^3r^3\sqrt{3}-3\,r^4\sqrt{3}\,y^4+54\,r^4 x^2 y-72\,\sqrt{3}r^3 x^2+
24\,\sqrt{3}r^3\,y^2+54\,r^4 x^3\sqrt{3}-432\, y\,r\, x-36\,r^3 x^2 y+24\,y^3r+48\,y^3 x-
36\,r^2\,y^3-36\, x^4\sqrt{3}+12\,r^3\,y^3+36\,r^3 y-18\,r^4 y^2\sqrt{3} x-432\, x^2 y-
72\,\sqrt{3} y^2r\, x-12\,r^3\sqrt{3} x\,y^2+108\,r^2 x^2 y-108\,r^2 x^3\sqrt{3}+
144\,y^2\sqrt{3} x-9\,r^4\sqrt{3} y^2-54\,r^4 y\, x+36\,r^3\sqrt{3}x-27\,r^4\sqrt{3} x^2-
4\,y^4r^2\sqrt{3}-18\,r^4\,y^3+36\,r^2\,y^2\sqrt{3} x+144\, x^3\sqrt{3}-72\,y^2\sqrt{3} x^2+
48\,\ve^2r^2 y-27\,r^4\sqrt{3} x^4+24\,r^2\,y^3 x-72\,r^2 x^3 y-4\,y^4\sqrt{3}+144\,x^3 y+
216\, x^2r\, y+18\,r^4\sqrt{3}\,y^2 x^2+8\,\ve^2\sqrt{3}r^2 y^2+24\,\ve^2\sqrt{3}r^2 x^2-
48\,\ve^2\sqrt{3}r^2 x-48\,\ve^2r^2 y\, x \bigr)\Bigr]\Big/
\Bigl[3\,r^2 \left( \sqrt{3} y+3-3\,x \right)^2 \left( 4\,\ve^2-1 \right)^2\Bigr]
$.
}

Adding up the $P(X_2 \in \NPE^r(X_1,\ve) \cup \G_1^r(X_1,\ve), X_1 \in T_s \setminus T(\y,\ve))$
values in the 16 possible cases above, and multiplying by 6
we get for $r \in [1,4/3)$,

\begin{multline*}
\mu^S_\lo(r,\ve)=
\Bigl[47\,r^6-195\,r^5+576\,r^4 \ve^4+860\,r^4-288\,r^4\ve^2-864\,r^3\ve^2-
846\,r^3+1728\,r^3\ve^4-108\,r^2+1152\,r^2\ve^4-\\
576\,r^2\ve^2+720\,r -256\Bigr]\Big/
\Bigl[108\,r^2 \left( 2+r \right)  \left( 16\,\ve^4-8\,\ve^2+1 \right)  \left( r+1 \right) \Bigr].
\end{multline*}
The $\mu^S_\lo(r,\ve)$ values for the other intervals can be calculated similarly.

For $r=\infty $, it is trivial to see that $\mu(r)=1$.
In fact, for fixed $\ve>0$, $\mu(r)=1$ for $r \ge \sqrt{3}/(2\,\ve)$.

\begin{remark}
Derivation of $\mu^A_\la(r,\ve)$ and $\mu^A_\lo(r,\ve)$ is similar to
the segregation case.
\end{remark}

\subsection*{Appendix 7: Proof of Corollary \ref{cor:MT-asy-norm-NYr}:}
Recall that
$S^{\la}_n(r)=\rho^{\la}_{I,n}(r)$ is
the relative edge density of the AND-underlying graph for the multiple triangle case.
Then the expectation of $S^{\la}_n(r)$ is
$$ \E\left[S^{\la}_n(r)\right]=\frac{2}{n\,(n-1)}
\sum\hspace*{-0.1 in}\sum_{i < j \hspace*{0.25 in}}
\hspace*{-0.1 in}\,\E\left[h^{\la}_{ij}(r)\right]=\E\left[ h^{\la}_{12}(r) \right]=
P(X_2\in \NPE^r(X_1) \cap \G_1^r(X_1))=\widetilde \mu_{\la}(r).$$
But, by definition of $\NPE^r(\cdot)$ and $\G_1^r(\cdot)$, if $X_1$
and $X_2$ are in different triangles, then $P(X_2 \in \NPE^r(X_1)\cap
\G_1^r(X_1))=0$. So by the law of total probability
\begin{eqnarray*}
\widetilde \mu_{\la}(r)&:=&P(X_2 \in \NPE^r(X_1)\cap \G_1^r(X_1))\\
&=& \sum_{i=1}^{J_m}P(X_2 \in \NPE^r(X_1)\cap \G_1^r(X_1)\,|\,\{X_1,X_2\} \subset T_i)\,P(\{X_1,X_2\} \subset T_i)\\
&=& \sum_{i=1}^{J_m}\mu_{\la}(r)\,P(\{X_1,X_2\} \subset T_i)
\text{ (since $P(X_2 \in \NPE^r(X_1)\cap \G_1^r(X_1)\,|\,\{X_1,X_2\} \subset T_i)=\mu_{\la}(r)$)}\\
&=& \mu_{\la}(r) \, \sum_{i=1}^{J_m}\left(\frac{A(T_i)}{\sum_{i=1}^{J_m}A(T_i)}\right)^2
\text{ (since $P(\{X_1,X_2\} \subset T_i)=\left(\frac{A(T_i)}{\sum_{i=1}^{J_m}A(T_i)}\right)^2$).}\\
&=& \mu_{\la}(r) \, \left(\sum_{i=1}^{J_m}w_i^2\right).
\end{eqnarray*}
where $\mu_{\la}(r)$ is given by Equation \eqref{eqn:Asymean_and}.

Likewise, we get $\widetilde \mu_{\lo}(r)=\mu_{\lo}(r)\,\left(\sum_{i=1}^{J_m}w_i^2\right)$
where $\mu_{\lo}(r)$ is given by Equation \eqref{eqn:Asymean_or}.

Furthermore, the asymptotic variance is
$$\widetilde \nu_{\la}(r)=\E\left[ h^{\la}_{12}(r)h^{\la}_{13}(r) \right]-\E\left[ h^{\la}_{12}(r) \right]\E\left[ h^{\la}_{13}(r) \right]=
P\left(\{X_2,X_3\} \subset \NPE^r(X_1)\cap \G_1^r(X_1)\right)-\left(\widetilde \mu_{\la}(r)\right)^2.$$
Then for $J_m>1$, we have
\begin{multline*}
P(\{X_2,X_3\} \subset \NPE^r(X_1)\cap \G_1^r(X_1))=
\sum_{i=1}^{J_m}P(\{X_2,X_3\} \subset \NPE^r(X_1)\cap \G_1^r(X_1)\,|\, \{X_1,X_2,X_3\} \subset T_i)\, P(\{X_1,X_2,X_3\} \subset T_i)\\
 = P(\{X_2,X_3\} \subset \NPE^r(X_1)\cap \G_1^r(X_1)\,|\,\{X_1,X_2,X_3\} \subset T_e)\, \left(\sum_{i=1}^{J_m}w_i^3 \right).
\end{multline*}

Hence,
\begin{eqnarray*}
\widetilde \nu_{\la}(r)&=&P(\{X_2,X_3\} \subset \NPE^r(X_1)\cap \G_1^r(X_1)\,|\, \{X_1,X_2,X_3\} \subset T_e)\,
\left(\sum_{i=1}^{J_m}w_i^3 \right)-\left(\widetilde \mu_{\la}(r)\right)^2\\
&=&\nu_{\la}(r)\,\left(\sum_{i=1}^{J_m}w_i^3\right)+
\mu_{\la}(r)^2\,\left(\sum_{i=1}^{J_m}w_i^3-\left(\sum_{i=1}^{J_m}w_i^2\right)^2\right).
\end{eqnarray*}

Likewise, we get
$\widetilde \nu_{\lo}(r)=\nu_{\lo}(r)\,\left(\sum_{i=1}^{J_m}w_i^3\right)+
\mu_{\lo}(r)^2\,\left(\sum_{i=1}^{J_m}w_i^3-\left(\sum_{i=1}^{J_m}w_i^2\right)^2\right).$

So conditional on $\Y_m$, if $\widetilde \nu_{\la}(r)>0$ then
$\sqrt{n}\,\left(S^{\la}_n(r)-\widetilde \mu_{\la}(r)\right)
\stackrel {\mathcal L}{\longrightarrow} \mathcal N\left(0,\widetilde \nu_{\la}(r)\right)$.
A similar result holds for the OR-underlying version.

\subsection*{Appendix 8: Proof of Theorem \ref{thm:MT-asy-norm-II}:}
Recall that
$\rho^{\la}_{II,n}(r)$ is the version II of the
relative edge density of the AND-underlying graph for the multiple triangle case.
Then the expectation of $\rho^{\la}_{II,n}(r)$ is
$$ \E\left[\rho^{\la}_{II,n}(r)\right]=
\sum_{i=1}^{J_m}\frac{n_i\,(n_i-1)}{2\,n_t}\,\E\left[\rho^{\la}_{{}_{[i]}}(r)\right]=\mu_{\la}(r)$$
since by \eqref{eqn:E[rho-and-D]} we have
$$\E[\rho^{\la}_{{}_{[i]}}(r)]=\frac{2}{n_i(n_i-1)}\sum\hspace*{-0.1 in}\sum_{k < l \hspace*{0.25 in}}
\hspace*{-0.1 in}\,\E\left[h^{\la}_{kl}(r)\right]=\E\left[ h^{\la}_{12}(r) \right]=\mu_{\la}(r)$$
where
$\mu_{\la}(r)$ is given by Equation \eqref{eqn:Asymean_and}.
Likewise, we get $\widetilde \mu_{\lo}(r)=\mu_{\lo}(r)$
where $\mu_{\lo}(r)$ is given by Equation \eqref{eqn:Asymean_or}.

Next,
$$ \Var\left[\rho^{\la}_{II,n}(r)\right]=
\sum_{i=1}^{J_m}\frac{n_i^2\,(n_i-1)^2}{4\,n_t^2}\,
\Var\left[\rho^{\la}_{{}_{[i]}}(r)\right]$$
since $\rho^{\la}_{{}_{[k]}}(r)$ and $\rho^{\la}_{{}_{[l]}}(r)$ are independent for $k\not=l$.
Then by \eqref{eqn:Var[rho-and-D]} we have
$$\Var\left[\rho^{\la}_{{}_{[i]}}(r)\right]=
\frac{2}{n_i\,(n_i-1)}\Var\left[ h^{\la}_{12}(r) \right]+
\frac{4\,(n_i-2)}{n_i\,(n_i-1)} \, \Cov\left[ h^{\la}_{12}(r),h^{\la}_{13}(r) \right].$$
So,
$$ \Var\left[\rho^{\la}_{II,n}(r)\right]=
\sum_{i=1}^{J_m}\frac{n_i\,(n_i-1)}{2\,n_t^2}\, \Var\left[ h^{\la}_{12}(r) \right]+
\sum_{i=1}^{J_m}\frac{n_i\,(n_i-1)\,(n_i-2)}{n_t^2}\, \Cov\left[ h^{\la}_{12}(r),h^{\la}_{13}(r) \right].$$
Here
$\sum_{i=1}^{J_m}\frac{n_i\,(n_i-1)}{2\,n_t^2}\, \Var\left[ h^{\la}_{12}(r) \right]=
\frac{1}{n_t}\, \Var\left[ h^{\la}_{12}(r) \right]$.
Then
for large $n_i$ and $n$,
$$\frac{1}{n_t}\, \Var\left[ h^{\la}_{12}(r) \right] \approx
\frac{2}{n^2 \sum_{i=1}^{J_m} w_i^2}\, \Var\left[ h^{\la}_{12}(r) \right]$$
since $\frac{n_t^2}{n^2}=\sum_{i=1}^{J_m}\frac{n_i\,(n_i-1)}{2\,n^2}$
and $n_i/n \rightarrow w_i$ as $n_i,n \rightarrow \infty$.
Similarly,
for large $n_i$ and $n$,
$$\sum_{i=1}^{J_m}\frac{n_i\,(n_i-1)\,(n_i-2)}{n_t^2}\, \Cov\left[ h^{\la}_{12}(r),h^{\la}_{13}(r) \right]
\approx
\left[\frac{4}{n} \left(\sum_{j=1}^{J_m} w_i^3\right)\Bigg/\left(\sum_{i=1}^{J_m} w_i^2\right)^2 \right] \Cov\left[ h^{\la}_{12}(r),h^{\la}_{13}(r) \right].$$

Hence, conditional on $\Y_m$,
$\sqrt{n}\,\left(\rho^{\la}_{II,n}(r)-\widetilde \mu_{\la}(r)\right)
\stackrel {\mathcal L}{\longrightarrow} \mathcal N\left(0,4\,\breve \nu_{\la}(r)\right)$
provided that $\breve \nu_{\la}(r)>0$
where $\breve \mu_{\la}(r)=\mu_{\la}(r)$ and
$\breve \nu_{\la}(r)=\nu_{\la}(r) \left(\sum_{i=1}^{J_m} w_i^3\right)\Big/\left(\sum_{i=1}^{J_m} w_i^2\right)^2 $.
A similar result holds for the OR-underlying version.


\newpage
\clearpage

\newpage
\clearpage

\begin{figure}[ht]
\centering
\rotatebox{-90}{ \resizebox{8.5 in}{!}{\includegraphics{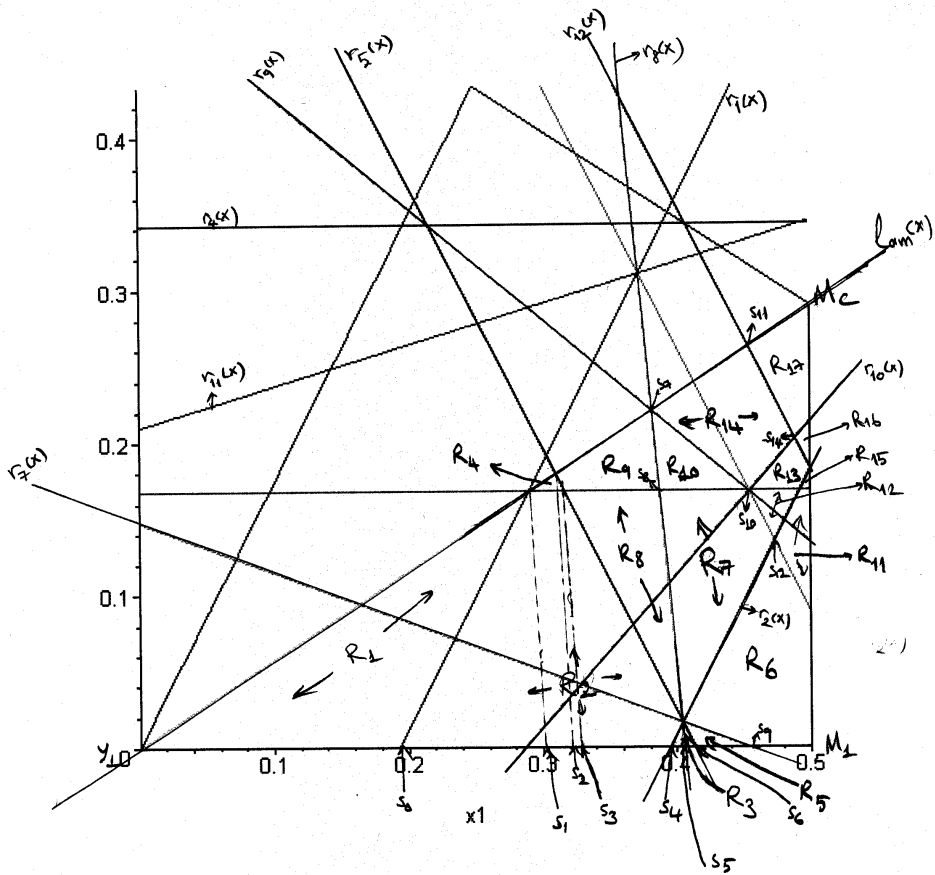}}}
\caption{
\label{fig:cases-AND-OR}
Prototype regions $R_i$ for various types of $\NPE^r(x_1) \cap \G_1^r(x_1)$
and the corresponding points whose $x$-coordinates are $s_k$ values.
}
\end{figure}


\newpage
\clearpage

\begin{figure}
\centering
\rotatebox{-90}{ \resizebox{8.5 in}{!}{\includegraphics{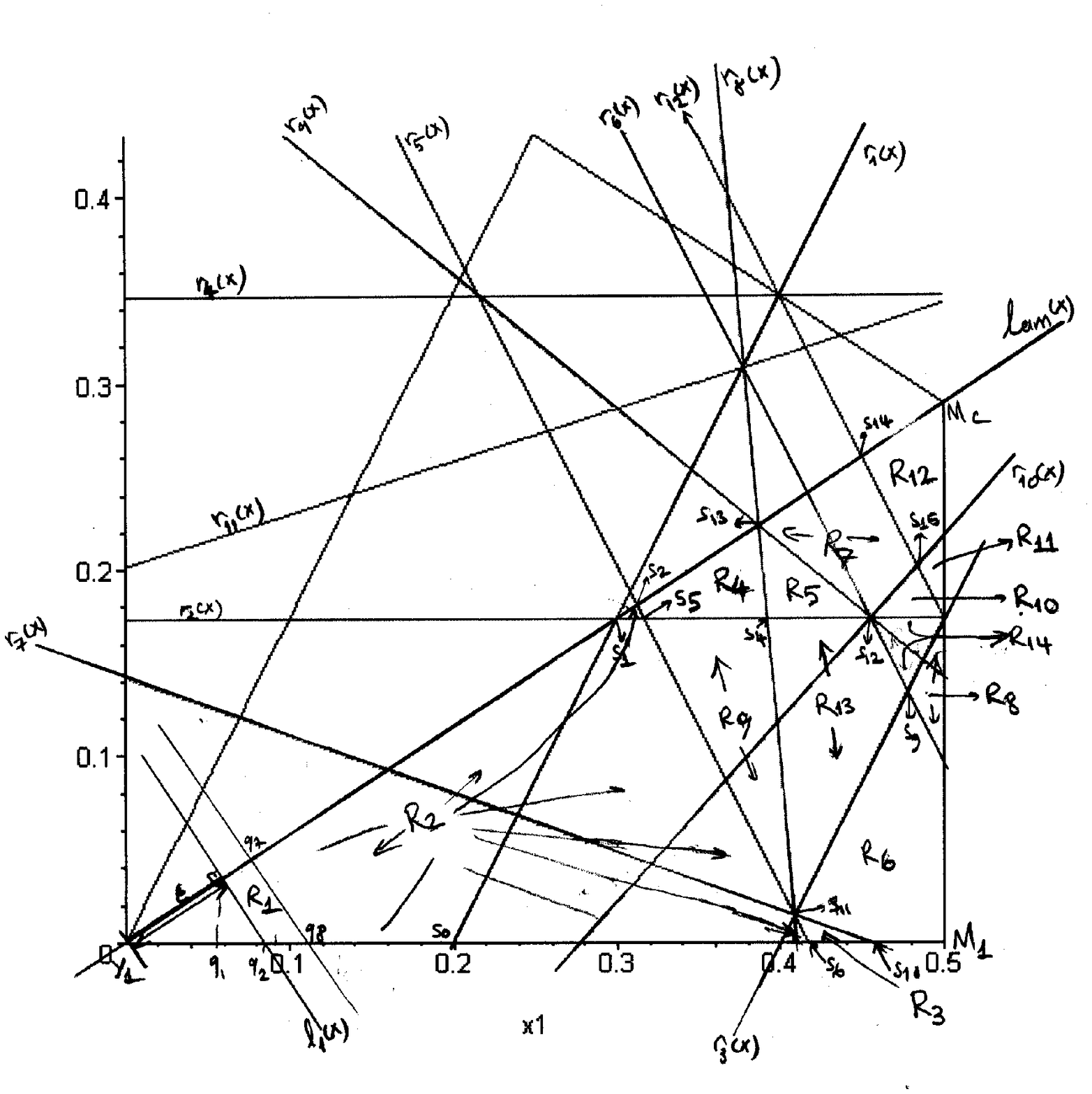}}}
\caption{
\label{fig:cases-AND-OR}
Prototype regions $R_i$ for various types of $\NPE^r(x_1) \cap \G_1^r(x_1)$
and the corresponding points whose $x$-coordinates are $s_k$ values.
}
\end{figure}

\end{document}

%% file: Nofnu2.pstex_t
\begin{picture}(0,0)%
\includegraphics{Nofnu2.pstex}%
\end{picture}%
\setlength{\unitlength}{3947sp}%
\begingroup\makeatletter\ifx\SetFigFont\undefined%
\gdef\SetFigFont#1#2#3#4#5{%
  \reset@font\fontsize{#1}{#2pt}%
  \fontfamily{#3}\fontseries{#4}\fontshape{#5}%
  \selectfont}%
\fi\endgroup%
\begin{picture}(10684,8079)(589,-7348)
\put(3376,-4936){\makebox(0,0)[lb]{\smash{{\SetFigFont{14}{16.8}{\rmdefault}{\mddefault}{\updefault}{\color[rgb]{0,0,0}$x$}%
}}}}
\put(2992,143){\rotatebox{315.0}{\makebox(0,0)[lb]{\smash{{\SetFigFont{14}{16.8}{\rmdefault}{\mddefault}{\updefault}{\color[rgb]{0,0,0}$\ell_2(v(x),x)$}%
}}}}}
\put(1405,-2508){\rotatebox{320.0}{\makebox(0,0)[lb]{\smash{{\SetFigFont{14}{16.8}{\rmdefault}{\mddefault}{\updefault}{\color[rgb]{0,0,0}$\ell(v(x),x)$}%
}}}}}
\put(6975,-1389){\rotatebox{310.0}{\makebox(0,0)[lb]{\smash{{\SetFigFont{14}{16.8}{\rmdefault}{\mddefault}{\updefault}{\color[rgb]{0,0,0}$e(x)$}%
}}}}}
\put(1051,-6511){\makebox(0,0)[lb]{\smash{{\SetFigFont{14}{16.8}{\rmdefault}{\mddefault}{\updefault}{\color[rgb]{0,0,0}$\y_1=v(x)$}%
}}}}
\put(11251,-6136){\makebox(0,0)[lb]{\smash{{\SetFigFont{14}{16.8}{\rmdefault}{\mddefault}{\updefault}{\color[rgb]{0,0,0}$\y_2$}%
}}}}
\put(4876,539){\makebox(0,0)[lb]{\smash{{\SetFigFont{14}{16.8}{\rmdefault}{\mddefault}{\updefault}{\color[rgb]{0,0,0}$\y_3$}%
}}}}
\put(3376,-6286){\rotatebox{45.0}{\makebox(0,0)[lb]{\smash{{\SetFigFont{14}{16.8}{\rmdefault}{\mddefault}{\updefault}{\color[rgb]{0,0,0}$d(v(x),\ell_2(v(x),x))=2\,d(v(x),\ell(v(x),x))$}%
}}}}}
\put(960,-5299){\rotatebox{45.0}{\makebox(0,0)[lb]{\smash{{\SetFigFont{12}{14.4}{\rmdefault}{\mddefault}{\updefault}{\color[rgb]{0,0,0} $d(v(x),\ell(v(x),x))$}%
}}}}}
\end{picture}%

%% file: Gammaofnu2.pstex_t
\begin{picture}(0,0)%
\includegraphics{Gammaofnu2.pstex}%
\end{picture}%
\setlength{\unitlength}{3947sp}%
\begingroup\makeatletter\ifx\SetFigFont\undefined%
\gdef\SetFigFont#1#2#3#4#5{%
  \reset@font\fontsize{#1}{#2pt}%
  \fontfamily{#3}\fontseries{#4}\fontshape{#5}%
  \selectfont}%
\fi\endgroup%
\begin{picture}(10684,8043)(589,-7348)
\put(3376,-4936){\makebox(0,0)[lb]{\smash{{\SetFigFont{12}{14.4}{\rmdefault}{\bfdefault}{\updefault}{\color[rgb]{0,0,0}$x$}%
}}}}
\put(2626,-2011){\makebox(0,0)[lb]{\smash{{\SetFigFont{12}{14.4}{\rmdefault}{\mddefault}{\updefault}{\color[rgb]{0,0,0}$\xi_3(x)$}%
}}}}
\put(3867,-6839){\rotatebox{320.0}{\makebox(0,0)[lb]{\smash{{\SetFigFont{12}{14.4}{\rmdefault}{\mddefault}{\updefault}{\color[rgb]{0,0,0}$\xi_1(x)$}%
}}}}}
\put(7533,-5107){\rotatebox{65.0}{\makebox(0,0)[lb]{\smash{{\SetFigFont{12}{14.4}{\rmdefault}{\mddefault}{\updefault}{\color[rgb]{0,0,0}$\xi_2(x)$}%
}}}}}
\put(2101,-7036){\rotatebox{45.0}{\makebox(0,0)[lb]{\smash{{\SetFigFont{12}{14.4}{\rmdefault}{\mddefault}{\updefault}{\color[rgb]{0,0,0}$d(y_1,\xi_1(x))$}%
}}}}}
\put(1405,-2508){\rotatebox{320.0}{\makebox(0,0)[lb]{\smash{{\SetFigFont{12}{14.4}{\rmdefault}{\mddefault}{\updefault}{\color[rgb]{0,0,0}$\ell(y_1,x)$}%
}}}}}
\put(960,-5299){\rotatebox{45.0}{\makebox(0,0)[lb]{\smash{{\SetFigFont{12}{14.4}{\rmdefault}{\mddefault}{\updefault}{\color[rgb]{0,0,0}$d(y_1,\ell(y_1,x))=r\,d(y_1,\xi_1(x))$}%
}}}}}
\put(1051,-6511){\makebox(0,0)[lb]{\smash{{\SetFigFont{12}{14.4}{\rmdefault}{\mddefault}{\updefault}{\color[rgb]{0,0,0}$\y_1$}%
}}}}
\put(4876,539){\makebox(0,0)[lb]{\smash{{\SetFigFont{12}{14.4}{\rmdefault}{\mddefault}{\updefault}{\color[rgb]{0,0,0}$\y_3$}%
}}}}
\put(11251,-6136){\makebox(0,0)[lb]{\smash{{\SetFigFont{12}{14.4}{\rmdefault}{\mddefault}{\updefault}{\color[rgb]{0,0,0}$\y_2$}%
}}}}
\end{picture}%

%% file: means3.pstex_t
\begin{picture}(0,0)%
\includegraphics{means3.pstex}%
\end{picture}%
\setlength{\unitlength}{3947sp}%
\begingroup\makeatletter\ifx\SetFigFont\undefined%
\gdef\SetFigFont#1#2#3#4#5{%
  \reset@font\fontsize{#1}{#2pt}%
  \fontfamily{#3}\fontseries{#4}\fontshape{#5}%
  \selectfont}%
\fi\endgroup%
\begin{picture}(10211,8344)(-5,-9933)
\put(5655,-4426){\rotatebox{20.0}{\makebox(0,0)[lb]{\smash{{\SetFigFont{25}{14.4}{\rmdefault}{\mddefault}{\updefault}{\color[rgb]{0,0,0}$\mu_{\la}(r)$}%
}}}}}
\put(5490,-3367){\rotatebox{10.0}{\makebox(0,0)[lb]{\smash{{\SetFigFont{25}{14.4}{\rmdefault}{\mddefault}{\updefault}{\color[rgb]{0,0,0}$\mu(r)$}%
}}}}}
\put(5483,-2686){\rotatebox{5.0}{\makebox(0,0)[lb]{\smash{{\SetFigFont{25}{14.4}{\rmdefault}{\mddefault}{\updefault}{\color[rgb]{0,0,0}$\mu_{\lo}(r)$}%
}}}}}
\end{picture}%

%% file: var3.pstex_t
\begin{picture}(0,0)%
\includegraphics{var3.pstex}%
\end{picture}%
\setlength{\unitlength}{3947sp}%
\begingroup\makeatletter\ifx\SetFigFont\undefined%
\gdef\SetFigFont#1#2#3#4#5{%
  \reset@font\fontsize{#1}{#2pt}%
  \fontfamily{#3}\fontseries{#4}\fontshape{#5}%
  \selectfont}%
\fi\endgroup%
\begin{picture}(10211,8344)(-5,-9933)
\put(8551,-4786){\makebox(0,0)[lb]{\smash{{\SetFigFont{25}{14.4}{\rmdefault}{\mddefault}{\updefault}{\color[rgb]{0,0,0}$\nu_{\la}(r)$}%
}}}}
\put(6451,-6886){\makebox(0,0)[lb]{\smash{{\SetFigFont{25}{14.4}{\rmdefault}{\mddefault}{\updefault}{\color[rgb]{0,0,0}$\nu(r)$}%
}}}}
\put(5551,-8311){\makebox(0,0)[lb]{\smash{{\SetFigFont{25}{14.4}{\rmdefault}{\mddefault}{\updefault}{\color[rgb]{0,0,0}$\nu_{\lo}(r)$}%
}}}}
\end{picture}%

%% file: pae_seg_und.pstex_t
\begin{picture}(0,0)%
\includegraphics{pae_seg_und.pstex}%
\end{picture}%
\setlength{\unitlength}{3947sp}%
\begingroup\makeatletter\ifx\SetFigFont\undefined%
\gdef\SetFigFont#1#2#3#4#5{%
  \reset@font\fontsize{#1}{#2pt}%
  \fontfamily{#3}\fontseries{#4}\fontshape{#5}%
  \selectfont}%
\fi\endgroup%
\begin{picture}(10204,10004)(-2,-10763)
\put(6673,-4737){\rotatebox{45.0}{\makebox(0,0)[lb]{\smash{{\SetFigFont{25}{14.4}{\rmdefault}{\mddefault}{\updefault}{\color[rgb]{0,0,0}$\PAE_{\la}^S(r)$}%
}}}}}
\put(7131,-5348){\rotatebox{45.0}{\makebox(0,0)[lb]{\smash{{\SetFigFont{25}{14.4}{\rmdefault}{\mddefault}{\updefault}{\color[rgb]{0,0,0}$\PAE^S(r)$}%
}}}}}
\put(7388,-7205){\rotatebox{25.0}{\makebox(0,0)[lb]{\smash{{\SetFigFont{25}{14.4}{\rmdefault}{\mddefault}{\updefault}{\color[rgb]{0,0,0}$\PAE_{\lo}^S(r)$}%
}}}}}
\end{picture}%

%% file: pae_agg_und.pstex_t
\begin{picture}(0,0)%
\includegraphics{pae_agg_und.pstex}%
\end{picture}%
\setlength{\unitlength}{3947sp}%
\begingroup\makeatletter\ifx\SetFigFont\undefined%
\gdef\SetFigFont#1#2#3#4#5{%
  \reset@font\fontsize{#1}{#2pt}%
  \fontfamily{#3}\fontseries{#4}\fontshape{#5}%
  \selectfont}%
\fi\endgroup%
\begin{picture}(10204,10004)(-2,-10763)
\put(3976,-6811){\makebox(0,0)[lb]{\smash{{\SetFigFont{25}{14.4}{\rmdefault}{\mddefault}{\updefault}{\color[rgb]{0,0,0}$\PAE_{\lo}^A(r)$}%
}}}}
\put(3976,-7486){\makebox(0,0)[lb]{\smash{{\SetFigFont{25}{14.4}{\rmdefault}{\mddefault}{\updefault}{\color[rgb]{0,0,0}$\PAE_{\la}^A(r)$}%
}}}}
\put(4126,-8161){\makebox(0,0)[lb]{\smash{{\SetFigFont{25}{14.4}{\rmdefault}{\mddefault}{\updefault}{\color[rgb]{0,0,0}$\PAE^A(r)$}%
}}}}
\end{picture}%

%% file: ls_lam_cases.pstex_t
\begin{picture}(0,0)%
\includegraphics{ls_lam_cases.pstex}%
\end{picture}%
\setlength{\unitlength}{3947sp}%
\begingroup\makeatletter\ifx\SetFigFont\undefined%
\gdef\SetFigFont#1#2#3#4#5{%
  \reset@font\fontsize{#1}{#2pt}%
  \fontfamily{#3}\fontseries{#4}\fontshape{#5}%
  \selectfont}%
\fi\endgroup%
\begin{picture}(11349,8280)(289,-7498)
\put(10426,-6586){\makebox(0,0)[lb]{\smash{{\SetFigFont{22}{16.8}{\rmdefault}{\mddefault}{\updefault}{\color[rgb]{0,0,0}$\y_2=(1,0)$}%
}}}}
\put(7201,-6586){\makebox(0,0)[lb]{\smash{{\SetFigFont{22}{16.8}{\rmdefault}{\mddefault}{\updefault}{\color[rgb]{0,0,0}$e_3$}%
}}}}
\put(5101,-6586){\makebox(0,0)[lb]{\smash{{\SetFigFont{22}{16.8}{\rmdefault}{\mddefault}{\updefault}{\color[rgb]{0,0,0}$M_3$}%
}}}}
\put(3301,-6586){\makebox(0,0)[lb]{\smash{{\SetFigFont{22}{14.4}{\rmdefault}{\mddefault}{\updefault}{\color[rgb]{0,0,0}$s_1$}%
}}}}
\put(1201,-2611){\makebox(0,0)[lb]{\smash{{\SetFigFont{22}{16.8}{\rmdefault}{\mddefault}{\updefault}{\color[rgb]{0,0,0}$\ell_s(r=4,x)$}%
}}}}
\put(1801,-736){\makebox(0,0)[lb]{\smash{{\SetFigFont{22}{16.8}{\rmdefault}{\mddefault}{\updefault}{\color[rgb]{0,0,0}$\ell_s(r=1.75,x)$}%
}}}}
\put(3226,314){\makebox(0,0)[lb]{\smash{{\SetFigFont{22}{16.8}{\rmdefault}{\mddefault}{\updefault}{\color[rgb]{0,0,0}$\ell_s\bigl(r=\sqrt{2},x\bigr)$}%
}}}}
\put(5926,614){\makebox(0,0)[lb]{\smash{{\SetFigFont{22}{16.8}{\rmdefault}{\mddefault}{\updefault}{\color[rgb]{0,0,0}$\y_3=\bigl(1/2,\sqrt{3}/2\bigr)$}%
}}}}
\put(7876,-2236){\makebox(0,0)[lb]{\smash{{\SetFigFont{22}{16.8}{\rmdefault}{\mddefault}{\updefault}{\color[rgb]{0,0,0}$e_1$}%
}}}}
\put(3976,-6586){\makebox(0,0)[lb]{\smash{{\SetFigFont{22}{14.4}{\rmdefault}{\mddefault}{\updefault}{\color[rgb]{0,0,0}$s_2$}%
}}}}
\put(826,-6736){\makebox(0,0)[lb]{\smash{{\SetFigFont{22}{16.8}{\rmdefault}{\mddefault}{\updefault}{\color[rgb]{0,0,0}$\y_1=(0,0)$}%
}}}}
\put(5851,-4111){\makebox(0,0)[lb]{\smash{{\SetFigFont{22}{16.8}{\rmdefault}{\mddefault}{\updefault}{\color[rgb]{0,0,0}$M_{C}$}%
}}}}
\put(3001,-2761){\makebox(0,0)[lb]{\smash{{\SetFigFont{22}{16.8}{\rmdefault}{\mddefault}{\updefault}{\color[rgb]{0,0,0}$e_2$}%
}}}}
\end{picture}%

%% file: G1ofxCase1.pstex_t
\begin{picture}(0,0)%
\includegraphics{G1ofxCase1.pstex}%
\end{picture}%
\setlength{\unitlength}{3947sp}%
\begingroup\makeatletter\ifx\SetFigFont\undefined%
\gdef\SetFigFont#1#2#3#4#5{%
  \reset@font\fontsize{#1}{#2pt}%
  \fontfamily{#3}\fontseries{#4}\fontshape{#5}%
  \selectfont}%
\fi\endgroup%
\begin{picture}(11349,8130)(289,-7348)
\put(10426,-6586){\makebox(0,0)[lb]{\smash{{\SetFigFont{25}{16.8}{\rmdefault}{\mddefault}{\updefault}{\color[rgb]{0,0,0}$\y_2=(1,0)$}%
}}}}
\put(826,-6586){\makebox(0,0)[lb]{\smash{{\SetFigFont{25}{16.8}{\rmdefault}{\mddefault}{\updefault}{\color[rgb]{0,0,0}$\y_1=(0,0)$}%
}}}}
\put(5326,614){\makebox(0,0)[lb]{\smash{{\SetFigFont{25}{16.8}{\rmdefault}{\mddefault}{\updefault}{\color[rgb]{0,0,0}$\y_3=(1/2,\sqrt{3}/2)$}%
}}}}
\put(7726,-2236){\makebox(0,0)[lb]{\smash{{\SetFigFont{25}{16.8}{\rmdefault}{\mddefault}{\updefault}{\color[rgb]{0,0,0}$e_1$}%
}}}}
\put(3001,-2686){\makebox(0,0)[lb]{\smash{{\SetFigFont{25}{16.8}{\rmdefault}{\mddefault}{\updefault}{\color[rgb]{0,0,0}$e_2$}%
}}}}
\put(7201,-6586){\makebox(0,0)[lb]{\smash{{\SetFigFont{25}{16.8}{\rmdefault}{\mddefault}{\updefault}{\color[rgb]{0,0,0}$e_3$}%
}}}}
\put(5626,-6586){\makebox(0,0)[lb]{\smash{{\SetFigFont{25}{16.8}{\rmdefault}{\mddefault}{\updefault}{\color[rgb]{0,0,0}$M_3$}%
}}}}
\put(2026,-5986){\makebox(0,0)[lb]{\smash{{\SetFigFont{25}{16.8}{\rmdefault}{\mddefault}{\updefault}{\color[rgb]{0,0,0}$\xi_1(r,x)$}%
}}}}
\put(5701,-4111){\makebox(0,0)[lb]{\smash{{\SetFigFont{25}{16.8}{\rmdefault}{\mddefault}{\updefault}{\color[rgb]{0,0,0}$M_{C}$}%
}}}}
\put(1426,-6211){\makebox(0,0)[lb]{\smash{{\SetFigFont{25}{16.8}{\rmdefault}{\mddefault}{\updefault}{\color[rgb]{0,0,0}$x_1$}%
}}}}
\end{picture}%

%% file: G1ofxCase2.pstex_t
\begin{picture}(0,0)%
\includegraphics{G1ofxCase2.pstex}%
\end{picture}%
\setlength{\unitlength}{3947sp}%
\begingroup\makeatletter\ifx\SetFigFont\undefined%
\gdef\SetFigFont#1#2#3#4#5{%
  \reset@font\fontsize{#1}{#2pt}%
  \fontfamily{#3}\fontseries{#4}\fontshape{#5}%
  \selectfont}%
\fi\endgroup%
\begin{picture}(11349,8130)(289,-7348)
\put(10426,-6586){\makebox(0,0)[lb]{\smash{{\SetFigFont{25}{16.8}{\rmdefault}{\mddefault}{\updefault}{\color[rgb]{0,0,0}$\y_2=(1,0)$}%
}}}}
\put(826,-6586){\makebox(0,0)[lb]{\smash{{\SetFigFont{25}{16.8}{\rmdefault}{\mddefault}{\updefault}{\color[rgb]{0,0,0}$\y_1=(0,0)$}%
}}}}
\put(5326,614){\makebox(0,0)[lb]{\smash{{\SetFigFont{25}{16.8}{\rmdefault}{\mddefault}{\updefault}{\color[rgb]{0,0,0}$\y_3=(1/2,\sqrt{3}/2)$}%
}}}}
\put(7726,-2236){\makebox(0,0)[lb]{\smash{{\SetFigFont{25}{16.8}{\rmdefault}{\mddefault}{\updefault}{\color[rgb]{0,0,0}$e_1$}%
}}}}
\put(3001,-2686){\makebox(0,0)[lb]{\smash{{\SetFigFont{25}{16.8}{\rmdefault}{\mddefault}{\updefault}{\color[rgb]{0,0,0}$e_2$}%
}}}}
\put(7201,-6586){\makebox(0,0)[lb]{\smash{{\SetFigFont{25}{16.8}{\rmdefault}{\mddefault}{\updefault}{\color[rgb]{0,0,0}$e_3$}%
}}}}
\put(5626,-6586){\makebox(0,0)[lb]{\smash{{\SetFigFont{25}{16.8}{\rmdefault}{\mddefault}{\updefault}{\color[rgb]{0,0,0}$M_3$}%
}}}}
\put(2026,-5986){\makebox(0,0)[lb]{\smash{{\SetFigFont{25}{16.8}{\rmdefault}{\mddefault}{\updefault}{\color[rgb]{0,0,0}$\xi_1(r,x)$}%
}}}}
\put(5701,-4111){\makebox(0,0)[lb]{\smash{{\SetFigFont{25}{16.8}{\rmdefault}{\mddefault}{\updefault}{\color[rgb]{0,0,0}$M_{C}$}%
}}}}
\put(3601,-6136){\makebox(0,0)[lb]{\smash{{\SetFigFont{25}{16.8}{\rmdefault}{\mddefault}{\updefault}{\color[rgb]{0,0,0}$x_1$}%
}}}}
\put(6226,-4636){\makebox(0,0)[lb]{\smash{{\SetFigFont{25}{16.8}{\rmdefault}{\mddefault}{\updefault}{\color[rgb]{0,0,0}$\xi_2(r,x)$}%
}}}}
\end{picture}%

%% file: G1ofxCase3.pstex_t
\begin{picture}(0,0)%
\includegraphics{G1ofxCase3.pstex}%
\end{picture}%
\setlength{\unitlength}{3947sp}%
\begingroup\makeatletter\ifx\SetFigFont\undefined%
\gdef\SetFigFont#1#2#3#4#5{%
  \reset@font\fontsize{#1}{#2pt}%
  \fontfamily{#3}\fontseries{#4}\fontshape{#5}%
  \selectfont}%
\fi\endgroup%
\begin{picture}(11349,8130)(289,-7348)
\put(10426,-6586){\makebox(0,0)[lb]{\smash{{\SetFigFont{25}{16.8}{\rmdefault}{\mddefault}{\updefault}{\color[rgb]{0,0,0}$\y_2=(1,0)$}%
}}}}
\put(826,-6586){\makebox(0,0)[lb]{\smash{{\SetFigFont{25}{16.8}{\rmdefault}{\mddefault}{\updefault}{\color[rgb]{0,0,0}$\y_1=(0,0)$}%
}}}}
\put(5326,614){\makebox(0,0)[lb]{\smash{{\SetFigFont{25}{16.8}{\rmdefault}{\mddefault}{\updefault}{\color[rgb]{0,0,0}$\y_3=(1/2,\sqrt{3}/2)$}%
}}}}
\put(7726,-2236){\makebox(0,0)[lb]{\smash{{\SetFigFont{25}{16.8}{\rmdefault}{\mddefault}{\updefault}{\color[rgb]{0,0,0}$e_1$}%
}}}}
\put(3001,-2686){\makebox(0,0)[lb]{\smash{{\SetFigFont{25}{16.8}{\rmdefault}{\mddefault}{\updefault}{\color[rgb]{0,0,0}$e_2$}%
}}}}
\put(7201,-6586){\makebox(0,0)[lb]{\smash{{\SetFigFont{25}{16.8}{\rmdefault}{\mddefault}{\updefault}{\color[rgb]{0,0,0}$e_3$}%
}}}}
\put(5626,-6586){\makebox(0,0)[lb]{\smash{{\SetFigFont{25}{16.8}{\rmdefault}{\mddefault}{\updefault}{\color[rgb]{0,0,0}$M_3$}%
}}}}
\put(5701,-4111){\makebox(0,0)[lb]{\smash{{\SetFigFont{25}{16.8}{\rmdefault}{\mddefault}{\updefault}{\color[rgb]{0,0,0}$M_{C}$}%
}}}}
\put(7726,-5236){\makebox(0,0)[lb]{\smash{{\SetFigFont{25}{16.8}{\rmdefault}{\mddefault}{\updefault}{\color[rgb]{0,0,0}$\xi_2(r,x)$}%
}}}}
\put(4651,-6136){\makebox(0,0)[lb]{\smash{{\SetFigFont{25}{16.8}{\rmdefault}{\mddefault}{\updefault}{\color[rgb]{0,0,0}$x_1$}%
}}}}
\put(2851,-5536){\makebox(0,0)[lb]{\smash{{\SetFigFont{25}{16.8}{\rmdefault}{\mddefault}{\updefault}{\color[rgb]{0,0,0}$\xi_1(r,x)$}%
}}}}
\end{picture}%

%% file: G1ofxCase4.pstex_t
\begin{picture}(0,0)%
\includegraphics{G1ofxCase4.pstex}%
\end{picture}%
\setlength{\unitlength}{3947sp}%
\begingroup\makeatletter\ifx\SetFigFont\undefined%
\gdef\SetFigFont#1#2#3#4#5{%
  \reset@font\fontsize{#1}{#2pt}%
  \fontfamily{#3}\fontseries{#4}\fontshape{#5}%
  \selectfont}%
\fi\endgroup%
\begin{picture}(11349,8130)(289,-7348)
\put(10426,-6586){\makebox(0,0)[lb]{\smash{{\SetFigFont{25}{16.8}{\rmdefault}{\mddefault}{\updefault}{\color[rgb]{0,0,0}$\y_2=(1,0)$}%
}}}}
\put(826,-6586){\makebox(0,0)[lb]{\smash{{\SetFigFont{25}{16.8}{\rmdefault}{\mddefault}{\updefault}{\color[rgb]{0,0,0}$\y_1=(0,0)$}%
}}}}
\put(5326,614){\makebox(0,0)[lb]{\smash{{\SetFigFont{25}{16.8}{\rmdefault}{\mddefault}{\updefault}{\color[rgb]{0,0,0}$\y_3=(1/2,\sqrt{3}/2)$}%
}}}}
\put(7201,-6586){\makebox(0,0)[lb]{\smash{{\SetFigFont{25}{16.8}{\rmdefault}{\mddefault}{\updefault}{\color[rgb]{0,0,0}$e_3$}%
}}}}
\put(2026,-5986){\makebox(0,0)[lb]{\smash{{\SetFigFont{25}{16.8}{\rmdefault}{\mddefault}{\updefault}{\color[rgb]{0,0,0}$\xi_1(r,x)$}%
}}}}
\put(5701,-4111){\makebox(0,0)[lb]{\smash{{\SetFigFont{25}{16.8}{\rmdefault}{\mddefault}{\updefault}{\color[rgb]{0,0,0}$M_{C}$}%
}}}}
\put(3826,-5536){\makebox(0,0)[lb]{\smash{{\SetFigFont{25}{16.8}{\rmdefault}{\mddefault}{\updefault}{\color[rgb]{0,0,0}$x_1$}%
}}}}
\put(3901,-1111){\makebox(0,0)[lb]{\smash{{\SetFigFont{25}{16.8}{\rmdefault}{\mddefault}{\updefault}{\color[rgb]{0,0,0}$e_2$}%
}}}}
\put(6901,-1111){\makebox(0,0)[lb]{\smash{{\SetFigFont{25}{16.8}{\rmdefault}{\mddefault}{\updefault}{\color[rgb]{0,0,0}$e_1$}%
}}}}
\put(5326,-6586){\makebox(0,0)[lb]{\smash{{\SetFigFont{25}{16.8}{\rmdefault}{\mddefault}{\updefault}{\color[rgb]{0,0,0}$M_3$}%
}}}}
\put(3151,-6661){\makebox(0,0)[lb]{\smash{{\SetFigFont{25}{16.8}{\rmdefault}{\mddefault}{\updefault}{\color[rgb]{0,0,0}$G_1$}%
}}}}
\put(1501,-4936){\makebox(0,0)[lb]{\smash{{\SetFigFont{25}{16.8}{\rmdefault}{\mddefault}{\updefault}{\color[rgb]{0,0,0}$G_6$}%
}}}}
\put(2851,-2686){\makebox(0,0)[lb]{\smash{{\SetFigFont{25}{16.8}{\rmdefault}{\mddefault}{\updefault}{\color[rgb]{0,0,0}$M_2$}%
}}}}
\put(3976,-2986){\makebox(0,0)[lb]{\smash{{\SetFigFont{25}{16.8}{\rmdefault}{\mddefault}{\updefault}{\color[rgb]{0,0,0}$L_5$}%
}}}}
\put(5251,-2986){\makebox(0,0)[lb]{\smash{{\SetFigFont{25}{16.8}{\rmdefault}{\mddefault}{\updefault}{\color[rgb]{0,0,0}$\xi_3(r,x)$}%
}}}}
\put(6901,-2986){\makebox(0,0)[lb]{\smash{{\SetFigFont{25}{16.8}{\rmdefault}{\mddefault}{\updefault}{\color[rgb]{0,0,0}$L_4$}%
}}}}
\put(7651,-3061){\makebox(0,0)[lb]{\smash{{\SetFigFont{25}{16.8}{\rmdefault}{\mddefault}{\updefault}{\color[rgb]{0,0,0}$L_3$}%
}}}}
\put(5626,-6211){\makebox(0,0)[lb]{\smash{{\SetFigFont{25}{16.8}{\rmdefault}{\mddefault}{\updefault}{\color[rgb]{0,0,0}$L_2$}%
}}}}
\put(6526,-4786){\makebox(0,0)[lb]{\smash{{\SetFigFont{25}{16.8}{\rmdefault}{\mddefault}{\updefault}{\color[rgb]{0,0,0}$\xi_2(r,x)$}%
}}}}
\put(7951,-2686){\makebox(0,0)[lb]{\smash{{\SetFigFont{25}{16.8}{\rmdefault}{\mddefault}{\updefault}{\color[rgb]{0,0,0}$M_1$}%
}}}}
\end{picture}%

%% file: G1ofxCase5.pstex_t
\begin{picture}(0,0)%
\includegraphics{G1ofxCase5.pstex}%
\end{picture}%
\setlength{\unitlength}{3947sp}%
\begingroup\makeatletter\ifx\SetFigFont\undefined%
\gdef\SetFigFont#1#2#3#4#5{%
  \reset@font\fontsize{#1}{#2pt}%
  \fontfamily{#3}\fontseries{#4}\fontshape{#5}%
  \selectfont}%
\fi\endgroup%
\begin{picture}(11349,8130)(289,-7348)
\put(10426,-6586){\makebox(0,0)[lb]{\smash{{\SetFigFont{25}{16.8}{\rmdefault}{\mddefault}{\updefault}{\color[rgb]{0,0,0}$\y_2=(1,0)$}%
}}}}
\put(826,-6586){\makebox(0,0)[lb]{\smash{{\SetFigFont{25}{16.8}{\rmdefault}{\mddefault}{\updefault}{\color[rgb]{0,0,0}$\y_1=(0,0)$}%
}}}}
\put(5326,614){\makebox(0,0)[lb]{\smash{{\SetFigFont{25}{16.8}{\rmdefault}{\mddefault}{\updefault}{\color[rgb]{0,0,0}$\y_3=(1/2,\sqrt{3}/2)$}%
}}}}
\put(7726,-2236){\makebox(0,0)[lb]{\smash{{\SetFigFont{25}{16.8}{\rmdefault}{\mddefault}{\updefault}{\color[rgb]{0,0,0}$e_1$}%
}}}}
\put(3001,-2686){\makebox(0,0)[lb]{\smash{{\SetFigFont{25}{16.8}{\rmdefault}{\mddefault}{\updefault}{\color[rgb]{0,0,0}$e_2$}%
}}}}
\put(7201,-6586){\makebox(0,0)[lb]{\smash{{\SetFigFont{25}{16.8}{\rmdefault}{\mddefault}{\updefault}{\color[rgb]{0,0,0}$e_3$}%
}}}}
\put(5626,-6586){\makebox(0,0)[lb]{\smash{{\SetFigFont{25}{16.8}{\rmdefault}{\mddefault}{\updefault}{\color[rgb]{0,0,0}$M_3$}%
}}}}
\put(5701,-4111){\makebox(0,0)[lb]{\smash{{\SetFigFont{25}{16.8}{\rmdefault}{\mddefault}{\updefault}{\color[rgb]{0,0,0}$M_{C}$}%
}}}}
\put(4801,-2986){\makebox(0,0)[lb]{\smash{{\SetFigFont{25}{16.8}{\rmdefault}{\mddefault}{\updefault}{\color[rgb]{0,0,0}$\xi_3(r,x)$}%
}}}}
\put(4951,-5386){\makebox(0,0)[lb]{\smash{{\SetFigFont{25}{16.8}{\rmdefault}{\mddefault}{\updefault}{\color[rgb]{0,0,0}$x_1$}%
}}}}
\put(7576,-5236){\makebox(0,0)[lb]{\smash{{\SetFigFont{25}{16.8}{\rmdefault}{\mddefault}{\updefault}{\color[rgb]{0,0,0}$\xi_2(r,x)$}%
}}}}
\put(3076,-5536){\makebox(0,0)[lb]{\smash{{\SetFigFont{25}{16.8}{\rmdefault}{\mddefault}{\updefault}{\color[rgb]{0,0,0}$\xi_1(r,x)$}%
}}}}
\end{picture}%

%% file: G1ofxCase6.pstex_t
\begin{picture}(0,0)%
\includegraphics{G1ofxCase6.pstex}%
\end{picture}%
\setlength{\unitlength}{3947sp}%
\begingroup\makeatletter\ifx\SetFigFont\undefined%
\gdef\SetFigFont#1#2#3#4#5{%
  \reset@font\fontsize{#1}{#2pt}%
  \fontfamily{#3}\fontseries{#4}\fontshape{#5}%
  \selectfont}%
\fi\endgroup%
\begin{picture}(11349,8130)(289,-7348)
\put(10426,-6586){\makebox(0,0)[lb]{\smash{{\SetFigFont{25}{16.8}{\rmdefault}{\mddefault}{\updefault}{$\y_2=(1,0)$}%
}}}}
\put(826,-6586){\makebox(0,0)[lb]{\smash{{\SetFigFont{25}{16.8}{\rmdefault}{\mddefault}{\updefault}{$\y_1=(0,0)$}%
}}}}
\put(5326,614){\makebox(0,0)[lb]{\smash{{\SetFigFont{25}{16.8}{\rmdefault}{\mddefault}{\updefault}{$\y_3=(1/2,\sqrt{3}/2)$}%
}}}}
\put(7726,-2236){\makebox(0,0)[lb]{\smash{{\SetFigFont{25}{16.8}{\rmdefault}{\mddefault}{\updefault}{$e_1$}%
}}}}
\put(3001,-2686){\makebox(0,0)[lb]{\smash{{\SetFigFont{25}{16.8}{\rmdefault}{\mddefault}{\updefault}{$e_2$}%
}}}}
\put(7201,-6586){\makebox(0,0)[lb]{\smash{{\SetFigFont{25}{16.8}{\rmdefault}{\mddefault}{\updefault}{$e_3$}%
}}}}
\put(5626,-6586){\makebox(0,0)[lb]{\smash{{\SetFigFont{25}{16.8}{\rmdefault}{\mddefault}{\updefault}{$M_3$}%
}}}}
\put(5701,-4111){\makebox(0,0)[lb]{\smash{{\SetFigFont{25}{16.8}{\rmdefault}{\mddefault}{\updefault}{$M_{C}$}%
}}}}
\put(5251,-4411){\makebox(0,0)[lb]{\smash{{\SetFigFont{25}{16.8}{\rmdefault}{\mddefault}{\updefault}{$x_1$}%
}}}}
\put(5851,-4636){\makebox(0,0)[lb]{\smash{{\SetFigFont{25}{16.8}{\rmdefault}{\mddefault}{\updefault}{$\xi_2(r,x)$}%
}}}}
\put(5176,-3436){\makebox(0,0)[lb]{\smash{{\SetFigFont{25}{16.8}{\rmdefault}{\mddefault}{\updefault}{$\xi_3(r,x)$}%
}}}}
\put(4201,-4711){\makebox(0,0)[lb]{\smash{{\SetFigFont{25}{16.8}{\rmdefault}{\mddefault}{\updefault}{$\xi_1(r,x)$}%
}}}}
\put(5476,-6961){\makebox(0,0)[lb]{\smash{{\SetFigFont{25}{14.4}{\rmdefault}{\mddefault}{\updefault}{case-6}%
}}}}
\end{picture}%

%% file: arcprobseg1und.pstex_t
\begin{picture}(0,0)%
\includegraphics{arcprobseg1und.pstex}%
\end{picture}%
\setlength{\unitlength}{3947sp}%
\begingroup\makeatletter\ifx\SetFigFont\undefined%
\gdef\SetFigFont#1#2#3#4#5{%
  \reset@font\fontsize{#1}{#2pt}%
  \fontfamily{#3}\fontseries{#4}\fontshape{#5}%
  \selectfont}%
\fi\endgroup%
\begin{picture}(11349,8168)(289,-7348)
\put(10426,-6586){\makebox(0,0)[lb]{\smash{{\SetFigFont{14}{16.8}{\rmdefault}{\mddefault}{\updefault}{\color[rgb]{0,0,0}$\y_2=(1,0)$}%
}}}}
\put(826,-6586){\makebox(0,0)[lb]{\smash{{\SetFigFont{14}{16.8}{\rmdefault}{\mddefault}{\updefault}{\color[rgb]{0,0,0}$\y_1=(0,0)$}%
}}}}
\put(7726,-2236){\makebox(0,0)[lb]{\smash{{\SetFigFont{14}{16.8}{\rmdefault}{\mddefault}{\updefault}{\color[rgb]{0,0,0}$e_1$}%
}}}}
\put(3001,-2686){\makebox(0,0)[lb]{\smash{{\SetFigFont{14}{16.8}{\rmdefault}{\mddefault}{\updefault}{\color[rgb]{0,0,0}$e_2$}%
}}}}
\put(5626,-6586){\makebox(0,0)[lb]{\smash{{\SetFigFont{14}{16.8}{\rmdefault}{\mddefault}{\updefault}{\color[rgb]{0,0,0}$M_3$}%
}}}}
\put(5701,-4111){\makebox(0,0)[lb]{\smash{{\SetFigFont{14}{16.8}{\rmdefault}{\mddefault}{\updefault}{\color[rgb]{0,0,0}$M_{C}$}%
}}}}
\put(3226,-4336){\makebox(0,0)[lb]{\smash{{\SetFigFont{14}{16.8}{\rmdefault}{\mddefault}{\updefault}{\color[rgb]{0,0,0}$\ell_{r_1}(x_1,x)$}%
}}}}
\put(7426,-3886){\makebox(0,0)[lb]{\smash{{\SetFigFont{14}{16.8}{\rmdefault}{\mddefault}{\updefault}{\color[rgb]{0,0,0}$\ell_{r_2}(x_1,x)$}%
}}}}
\put(2926,-5386){\makebox(0,0)[lb]{\smash{{\SetFigFont{14}{16.8}{\rmdefault}{\mddefault}{\updefault}{\color[rgb]{0,0,0}$x_1$}%
}}}}
\put(1876,-6061){\makebox(0,0)[lb]{\smash{{\SetFigFont{14}{16.8}{\rmdefault}{\mddefault}{\updefault}{\color[rgb]{0,0,0}$\varepsilon$}%
}}}}
\put(3226,-6136){\makebox(0,0)[lb]{\smash{{\SetFigFont{14}{16.8}{\rmdefault}{\mddefault}{\updefault}{\color[rgb]{0,0,0}$q(\y_1,x)$}%
}}}}
\put(7351,-6061){\makebox(0,0)[lb]{\smash{{\SetFigFont{14}{16.8}{\rmdefault}{\mddefault}{\updefault}{\color[rgb]{0,0,0}$q(\y_2,x)$}%
}}}}
\put(5326,-586){\makebox(0,0)[lb]{\smash{{\SetFigFont{14}{16.8}{\rmdefault}{\mddefault}{\updefault}{\color[rgb]{0,0,0}$q(\y_3,x)$}%
}}}}
\put(4426,-286){\makebox(0,0)[lb]{\smash{{\SetFigFont{14}{16.8}{\rmdefault}{\mddefault}{\updefault}{\color[rgb]{0,0,0}$N_2$}%
}}}}
\put(4201,-6586){\makebox(0,0)[lb]{\smash{{\SetFigFont{14}{16.8}{\rmdefault}{\mddefault}{\updefault}{\color[rgb]{0,0,0}$N_1$}%
}}}}
\put(9301,-6511){\makebox(0,0)[lb]{\smash{{\SetFigFont{14}{16.8}{\rmdefault}{\mddefault}{\updefault}{\color[rgb]{0,0,0}$N_1$}%
}}}}
\put(5551,614){\makebox(0,0)[lb]{\smash{{\SetFigFont{14}{16.8}{\rmdefault}{\mddefault}{\updefault}{\color[rgb]{0,0,0}$\y_3=(1/2,\sqrt{3}/2)$}%
}}}}
\put(4876,-1336){\makebox(0,0)[lb]{\smash{{\SetFigFont{14}{16.8}{\rmdefault}{\mddefault}{\updefault}{\color[rgb]{0,0,0}$U_2$}%
}}}}
\put(7876,-5536){\makebox(0,0)[lb]{\smash{{\SetFigFont{14}{16.8}{\rmdefault}{\mddefault}{\updefault}{\color[rgb]{0,0,0}$U_1$}%
}}}}
\put(1951,-3811){\makebox(0,0)[lb]{\smash{{\SetFigFont{14}{16.8}{\rmdefault}{\mddefault}{\updefault}{\color[rgb]{0,0,0}$N_2$}%
}}}}
\put(8101,-6586){\makebox(0,0)[lb]{\smash{{\SetFigFont{14}{16.8}{\rmdefault}{\mddefault}{\updefault}{\color[rgb]{0,0,0}$V_2$}%
}}}}
\put(9414,-4411){\makebox(0,0)[lb]{\smash{{\SetFigFont{14}{16.8}{\rmdefault}{\mddefault}{\updefault}{\color[rgb]{0,0,0}$V_3$}%
}}}}
\put(6714,-786){\makebox(0,0)[lb]{\smash{{\SetFigFont{14}{16.8}{\rmdefault}{\mddefault}{\updefault}{\color[rgb]{0,0,0}$V_4$}%
}}}}
\put(4126,-674){\makebox(0,0)[lb]{\smash{{\SetFigFont{14}{16.8}{\rmdefault}{\mddefault}{\updefault}{\color[rgb]{0,0,0}$V_5$}%
}}}}
\put(3139,-6661){\makebox(0,0)[lb]{\smash{{\SetFigFont{14}{16.8}{\rmdefault}{\mddefault}{\updefault}{\color[rgb]{0,0,0}$V_1$}%
}}}}
\put(1389,-4574){\makebox(0,0)[lb]{\smash{{\SetFigFont{14}{16.8}{\rmdefault}{\mddefault}{\updefault}{\color[rgb]{0,0,0}$V_6$}%
}}}}
\end{picture}%